\documentclass[10pt]{amsart}
\usepackage{amsmath}
\usepackage{amsxtra}
\usepackage{amscd}
\usepackage{amsthm}
\usepackage{amsfonts}
\usepackage{amssymb}
\usepackage{eucal}

\usepackage{latexsym}

\newtheorem{cor}[subsubsection]{Corollary}
\newtheorem{lem}[subsubsection]{Lemma}
\newtheorem{prop}[subsubsection]{Proposition}

\newtheorem{conj}[subsubsection]{Conjecture}
\newtheorem{mainconj}[subsubsection]{Main Conjecture}
\newtheorem{thm}[subsubsection]{Theorem}
\newtheorem{mainthm}[subsubsection]{Main Theorem}

\newtheorem{lemconstr}[subsubsection]{Lemma-Construction}
\theoremstyle{remark}
\newtheorem{remark}[subsubsection]{Remark}

\theoremstyle{definition}

\numberwithin{equation}{section}

\newcommand{\thmref}[1]{Theorem~\ref{#1}}
\newcommand{\mainthmref}[1]{Main Theorem~\ref{#1}}
\newcommand{\secref}[1]{Sect.~\ref{#1}}
\newcommand{\lemref}[1]{Lemma~\ref{#1}}
\newcommand{\propref}[1]{Proposition~\ref{#1}}
\newcommand{\corref}[1]{Corollary~\ref{#1}}
\newcommand{\conjref}[1]{Conjecture~\ref{#1}}
\newcommand{\mainconjref}[1]{Main Conjecture~\ref{#1}}

\newcommand{\nc}{\newcommand}
\nc{\renc}{\renewcommand}
\nc{\ssec}{\subsection}
\nc{\sssec}{\subsubsection}
\nc{\on}{\operatorname}

\nc\ol{\overline}
\nc\ul{\underline}
\nc\wt{\widetilde}
\nc\tboxtimes{\wt{\boxtimes}}
\nc{\wh}{\widehat}
\nc{\mc}{\mathcal}

\nc{\CM}{{\mathcal M}}
\nc{\CN}{{\mathcal N}}
\nc{\CF}{{\mathcal F}}
\nc{\D}{{\mathcal D}}
\nc{\CQ}{{\mathcal Q}}
\nc{\CY}{{\mathcal Y}}
\nc{\CX}{{\mathcal X}}
\nc{\CG}{{\mathcal G}}
\nc{\CE}{{\mathcal E}}
\nc{\CC}{{\mathcal C}}
\nc{\CO}{{\mathcal O}}
\renc{\CC}{{\mathcal C}}
\nc{\CT}{{\mathcal T}}
\nc{\CK}{{\mathcal K}}
\nc{\CS}{{\mathcal S}}
\nc{\CH}{{\mathcal H}}
\nc{\CU}{{\mathcal U}}
\nc{\CV}{{\mathcal V}}
\nc{\CA}{{\mathcal A}}
\nc{\CB}{{\mathcal B}}
\nc{\CW}{{\mathcal W}}
\nc{\CL}{{\mathcal L}}
\nc{\CP}{{\mathcal P}}
\nc{\CI}{{\mathcal I}}
\nc{\CJ}{{\mathcal J}}
\nc{\CR}{{\mathcal R}}

\nc{\BA}{{\mathbb{A}}}
\nc{\BC}{{\mathbb{C}}}
\nc{\BG}{{\mathbb{G}}}
\nc{\BM}{{\mathbb{M}}}
\nc{\BN}{{\mathbb{N}}}
\nc{\BP}{{\mathbb{P}}}
\nc{\BR}{{\mathbb{R}}}
\nc{\BZ}{{\mathbb{Z}}}
\nc{\BV}{{\mathbb{V}}}
\nc{\BW}{{\mathbb{W}}}
\nc{\BS}{{\mathbb{S}}}
\nc{\BD}{{\mathbb{D}}}
\nc{\BQ}{{\mathbb{Q}}}
\nc{\BL}{{\mathbb{L}}}
\renc{\BW}{{\mathbb{W}}}

\nc{\fa}{{\mathfrak{a}}}
\nc{\fb}{{\mathfrak{b}}}
\nc{\fg}{{\mathfrak{g}}}
\nc{\fgl}{{\mathfrak{gl}}}
\nc{\fh}{{\mathfrak{h}}}
\nc{\fj}{{\mathfrak{j}}}
\nc{\fm}{{\mathfrak{m}}}
\nc{\fl}{{\mathfrak{l}}}
\nc{\fn}{{\mathfrak{n}}}
\nc{\fu}{{\mathfrak{u}}}
\nc{\fp}{{\mathfrak{p}}}
\nc{\fr}{{\mathfrak{r}}}
\nc{\fs}{{\mathfrak{s}}}
\nc{\fsl}{{\mathfrak{sl}}}

\nc{\hsl}{{\widehat{\mathfrak{sl}}}}
\nc{\hgl}{{\widehat{\mathfrak{gl}}}}
\nc{\hg}{{\widehat{\mathfrak{g}}}}
\nc{\hb}{{\widehat{\mathfrak{b}}}}
\nc{\hn}{{\widehat{\mathfrak{n}}}}

\nc{\fA}{{\mathfrak{A}}}
\nc{\fB}{{\mathfrak{B}}}
\nc{\fO}{{\mathfrak{O}}}
\nc{\fD}{{\mathfrak{D}}}
\nc{\fE}{{\mathfrak{E}}}
\nc{\fF}{{\mathfrak{F}}}
\nc{\fG}{{\mathfrak{G}}}
\nc{\fK}{{\mathfrak{K}}}
\nc{\fL}{{\mathfrak{L}}}
\nc{\fC}{{\mathfrak{C}}}
\nc{\fM}{{\mathfrak{M}}}
\nc{\fN}{{\mathfrak{N}}}
\nc{\fH}{{\mathfrak{H}}}
\nc{\fP}{{\mathfrak{P}}}
\nc{\fU}{{\mathfrak{U}}}
\nc{\fV}{{\mathfrak{V}}}
\nc{\fZ}{{\mathfrak{Z}}}
\nc{\fz}{{\mathfrak{z}}}

\nc{\bc}{{\mathbf{c}}}
\nc{\bd}{{\mathbf{d}}}
\nc{\bh}{{\mathbf{h}}}
\nc{\be}{{\mathbf{e}}}
\nc{\bj}{{\mathbf{j}}}
\nc{\bn}{{\mathbf{n}}}
\nc{\bp}{{\mathbf{p}}}
\nc{\bg}{{\mathbf{g}}}
\nc{\bq}{{\mathbf{q}}}
\nc{\bs}{{\mathbf{s}}}
\nc{\bu}{{\mathbf{u}}}
\nc{\bv}{{\mathbf{v}}}
\nc{\bx}{{\mathbf{x}}}
\nc{\by}{{\mathbf{y}}}
\nc{\bw}{{\mathbf{w}}}
\nc{\bA}{{\mathbf{A}}}
\nc{\bK}{{\mathbf{K}}}
\nc{\bB}{{\mathbf{B}}}
\nc{\bC}{{\mathbf{C}}}
\nc{\bD}{{\mathbf{D}}}
\nc{\bH}{{\mathbf{H}}}
\nc{\bM}{{\mathbf{M}}}
\nc{\bN}{{\mathbf{N}}}
\nc{\bV}{{\mathbf{V}}}
\nc{\bW}{{\mathbf{W}}}
\nc{\bL}{{\mathbf{L}}}
\nc{\bU}{{\mathbf{U}}}
\nc{\bX}{{\mathbf{X}}}
\nc{\bI}{{\mathbf{I}}}
\nc{\bZ}{{\mathbf{Z}}}
\nc{\bS}{{\mathbf{S}}}

\nc{\sA}{{\mathsf{A}}}
\nc{\sB}{{\mathsf{B}}}
\nc{\sC}{{\mathsf{C}}}
\nc{\sD}{{\mathsf{D}}}
\nc{\sF}{{\mathsf{F}}}
\nc{\sH}{{\mathsf{H}}}
\nc{\sG}{{\mathsf{G}}}
\nc{\sK}{{\mathsf{K}}}
\nc{\sM}{{\mathsf{M}}}
\nc{\sO}{{\mathsf{O}}}
\nc{\sQ}{{\mathsf{Q}}}
\nc{\sP}{{\mathsf{P}}}
\nc{\sV}{{\mathsf{V}}}
\nc{\sZ}{{\mathsf{Z}}}
\nc{\sfp}{{\mathsf{p}}}
\nc{\sr}{{\mathsf{r}}}
\nc{\sg}{{\mathsf{g}}}
\nc{\sk}{{\mathsf{k}}}
\nc{\ssf}{{\mathsf{f}}}
\nc{\ssh}{{\mathsf{h}}}
\nc{\sse}{{\mathsf{e}}}
\nc{\sfb}{{\mathsf{b}}}
\nc{\sfc}{{\mathsf{c}}}
\nc{\sd}{{\mathsf{d}}}

\nc{\Av}{\on{Av}}
\nc{\act}{\on{act}}
\nc{\Hom}{\on{Hom}}
\nc{\End}{\on{End}}
\nc{\Lie}{\on{Lie}}
\nc{\Loc}{\on{Loc}}
\nc{\IC}{\on{IC}}
\nc{\Aut}{\on{Aut}}
\nc{\rk}{\on{rk}}
\nc{\Sh}{\on{Sh}}
\nc{\Perv}{\on{Perv}}
\nc{\pos}{{\on{pos}}}
\nc{\Conv}{\on{Conv}}
\nc{\Sph}{\on{Sph}}
\nc{\Sym}{\on{Sym}}
\nc{\Rep}{{\mc R}ep(\cG)}
\nc{\RepH}{{\mc R}ep(H)}
\nc{\Fun}{\on{Fun}}
\nc{\Id}{\on{Id}}
\nc{\id}{\on{id}}
\renc{\mod}{\on{--mod}}

\nc{\oG}{\overset{\circ}{G}{}}
\nc{\oGB}{{\overset{\circ}{G/B}{}}}
\nc{\oGN}{{\overset{\circ}{G/N}{}}}
\nc{\uBC}{\underline{\BC}}

\nc{\crit}{{\on{crit}}}
\nc{\reg}{{\on{reg}}}
\nc{\nilp}{{\on{nilp}}}
\nc{\ord}{\on{ord}}
\nc{\nil}{\wt{\on{reg}}}
\nc{\mb}{\mathbf}
\nc{\ren}{\on{ren}}
\nc{\res}{\on{res}}
\nc{\RS}{{\on{RS}}}
\nc{\Dist}{\on{Dist}}
\nc{\semiinf}{{\frac{\infty}{2}}}
\nc{\semiinfi}{{\frac{\infty}{2}+i}}
\nc{\semiinfb}{{\frac{\infty}{2}+\bullet}}
\nc{\torsemiinf}{{\overset{\semiinf}\otimes}}
\nc{\Hitch}{\on{Hitch}}

\nc{\hl}{\overset{\leftarrow}h}
\nc{\hr}{\overset{\rightarrow}h}
\nc\Dh{\widehat{\D}}
\nc{\Gr}{\on{Gr}}
\nc{\Fl}{\on{Fl}}
\nc{\Flt}{\wt{\Fl}{}}
\nc{\Pic}{\on{Pic}}
\nc{\Bun}{\on{Bun}}

\nc{\bDR}{\mathbf {DR}}
\nc{\uV}{\underline{V}}
\nc{\arrowtimes}{\overset{\to}\otimes}
\nc{\hattimes}{\widehat\otimes}
\nc{\larrowtimes}{\overset{\leftarrow}\otimes}
\nc{\shriektimes}{\overset{!}\otimes}
\nc{\startimes}{\overset{*}\otimes}
\nc{\sCliff}{\mathsf {Cliff}}
\nc{\sSpin}{\mathsf {Spin}}

\nc{\one}{{\mathbf{1}}}

\nc\Spec{\on{Spec}}
\nc{\Pro}{\on{Pro}}
\nc{\QCoh}{\on{QCoh}}
\nc{\uHom}{\underline{\on{Hom}}}
\nc{\RHom}{\on{RHom}}
\nc{\uRHom}{\underline{\on{RHom}}}
\nc{\CHom}{{\mathcal Hom}}
\nc{\uCHom}{\underline{{\mathcal Hom}}}
\nc{\uCRHom}{\underline{{\mathcal R}{\mathcal Hom}}}

\nc{\cg}{\check \fg}
\nc{\Op}{\on{Op}_{\cg}}
\nc{\nOp}{\on{Op}^{\nilp}_{\cg}}
\nc{\nMOp}{\on{MOp}^{\nilp}_{\cg}}
\nc{\rOp}{\on{Op}^{\reg}_{\cg}}

\nc{\tg}{\wt{\check \fg}}
\nc{\cn}{\check \fn}
\nc{\tn}{\wt{\cn}}
\nc{\cG}{\check G}
\nc{\cB}{\check B}
\nc{\cb}{\check \fb}
\nc{\MOp}{\on{MOp}_{\check \cg}}
\nc{\cN}{\CN_{\cG}}
\nc{\tN}{\wt{\CN}_{\cG}}
\nc{\dIsom}{{\mathsf{Isom}}_{\Op}}
\nc{\disom}{{\mathsf{isom}}_{\Op}}
\nc{\Kdv}{{\mathsf{Isom}}_{\Op^\reg}}
\nc{\kdv}{{\mathsf{isom}}_{\Op^\reg}}
\nc{\Isom}{{\mathsf{Isom}}_{\on{Op}_\fg}}
\nc{\isom}{{\mathsf{isom}}_{\on{Op}_\fg}}

\nc{\wcosta}{j_{\wt{w},*}}
\nc{\wsta}{j_{\wt{w},!}}
\nc{\wcost}{j_{w,*}}
\nc{\wst}{j_{w,!}}

\nc{\epsi}{{\mathbf e}^\psi}
\nc{\epsip}{{\mathbf e}^{\psi'}}

\nc{\Ppi}{{\mathbf \Pi}}

\nc{\hCO}{{\hat{\CO}}}
\nc{\hCK}{{\hat{\CK}}}

\nc{\CPreg}{\CP_{G,\on{Op}^\reg}}
\nc{\CPBreg}{\CP_{B,\on{Op}^\reg}}
\nc{\CPnilp}{\CP_{G,\on{Op}^\nilp}}
\nc{\CPBnilp}{\CP_{B,\on{Op}^\nilp}}
\nc{\CPla}{\CP_{G,\on{Op}_{\cla}}}
\nc{\CPBla}{\CP_{B,\on{Op}_{\cla}}}

\nc{\Cat}{\hg_\crit\mod^{I,m}_\nilp}
\nc{\Catf}{{}^f\hg_\crit\mod^{I,m}_\nilp}
\nc{\DCat}{D^b(\hg_\crit\mod_\nilp)^{I^0}}
\nc{\DCatf}{{}^f D^b(\hg_\crit\mod_\nilp)^{I^0}}
\nc{\Catr}{\hg_\crit\mod^{I,m}_\reg}
\nc{\Catrf}{{}^f\hg_\crit\mod^{I,m}_\reg}
\nc{\DCatr}{D^b(\hg_\crit\mod_\reg)^{I^0}}
\nc{\DCatrf}{{}^f D^b(\hg_\crit\mod_\reg)^{I^0}}

\nc{\ch}{\mbox{ch}}
\nc{\Z}{{\mathbb Z}}
\nc{\C}{{\mathbb C}}
\nc{\pone}{{\mathbb C}{\mathbb P}^1}
\nc{\pa}{\partial}
\nc{\F}{{\mathcal F}}
\nc{\arr}{\rightarrow}
\nc{\larr}{\longrightarrow}
\nc{\al}{\alpha}
\nc{\ri}{\rangle}
\nc{\lef}{\langle}
\nc{\W}{{\mathcal W}}
\nc{\la}{\lambda}
\nc{\ep}{\epsilon}
\nc{\su}{\widehat{{\mathfrak s}{\mathfrak l}}_2}
\nc{\sw}{{\mathfrak s}{\mathfrak l}}
\nc{\g}{{\mathfrak g}}
\nc{\h}{{\mathfrak h}}
\nc{\n}{{\mathfrak n}}
\nc{\N}{\widehat{\n}}
\nc{\De}{\Delta}
\nc{\gt}{\widetilde{\g}}
\nc{\Ga}{\Gamma}
\nc{\z}{{\mathfrak Z}}
\nc{\La}{\Lambda}
\nc{\cri}{_{\kappa_c}}
\nc{\kk}{h^\vee}
\nc{\sun}{\widehat{\sw}_N}
\nc{\si}{\sigma}
\nc{\el}{\ell}
\nc{\bi}{\bibitem}
\nc{\om}{\omega}
\nc{\ds}{\displaystyle}
\nc{\dzz}{\frac{dz}{z}}
\nc{\Res}{\on{Res}}
\nc{\Cal}{\mathcal}
\nc{\bb}{{\mathfrak b}}
\nc{\ot}{\otimes}
\nc{\R}{{\mc R}}
\nc{\yy}{{\mc Y}}
\nc{\ga}{\gamma}

\nc{\us}{\underset}
\nc{\opl}{\oplus}
\nc{\beq}{\begin{equation}}
\nc{\Fq}{{\mathbb F}_q}
\nc{\Mq}{{\mathcal M}}
\nc{\lan}{\langle}
\nc{\ran}{\rangle}

\nc{\Vect}{\on{Vect}}
\nc{\ghat}{\wh\fg}
\nc{\T}{\mc T}
\nc{\Tloc}{\T^\g_{\on{loc}}}
\nc{\vac}{|0\ran}
\nc{\Wick}{{\mb :}}
\nc{\delz}{\partial_z}
\nc{\K}{{\cali K}}
\nc{\cali}{\mathcal}
\nc{\li}{\mathfrak l}
\nc{\lt}{\widetilde{\li}}
\nc{\astar}{a^*}
\nc{\cA}{{\mc A}}
\nc{\ka}{\kappa}

\nc{\OO}{{\mc O}}
\nc{\AutO}{\on{Aut}\OO}
\nc{\DerO}{\on{Der}\OO}
\nc{\DerpO}{\on{Der}_+\OO}
\nc{\Au}{{\mc A}ut}
\nc{\mf}{\mathfrak}
\nc{\V}{{\mc V}}
\nc{\hh}{\wh{\h}}

\nc{\pp}{{\mathfrak p}}
\nc{\mm}{{\mathfrak m}}
\nc{\rr}{{\mathfrak r}}
\nc{\ket}{\rangle}
\nc{\zz}{{\mathfrak z}}
\nc{\gr}{\on{gr}}
\nc{\Spe}{\on{Spec}}
\nc{\rv}{\crho}
\nc{\can}{\on{can}}
\nc{\Db}{{\mathbb D}}
\nc{\ww}{w}

\nc{\RR}{\on{R}}
\nc{\PPi}{{\mathbf \Pi}}
\nc{\M}{{\mathbb M}}
\nc{\Mv}{{\mathbb M}^\vee}
\nc{\VV}{{\mathbb V}}
\nc{\bsl}{\backslash}

\nc{\bchi}{{\mathbf {\chi}}}
\nc{\anch}{{\mathbf {anch}}}

\nc{\cla}{{\check{\la}}}
\nc{\cmu}{{\check{\mu}}}
\nc{\crho}{{\check{\rho}}}
\nc{\com}{{\check{\omega}}}
\nc{\DD}{{\mc D}}
\nc{\E}{{\mc E}}
\nc{\Ll}{{\mc L}}

\nc{\ConnX}{\on{Conn}_H(\omega_X^{\crho})}
\nc{\ConHX}{\on{Conn}_{\check{H}}(\omega_X^{\rho})}
\nc{\ConnD}{\on{Conn}_H(\omega_{\D}^{\crho})}
\nc{\ConHD}{\on{Conn}_{\check{H}}(\omega_{\D}^{\rho})}
\nc{\ConnDt}{\on{Conn}_H(\omega_{\D^\times}^{\crho})}
\nc{\ConHDt}{\on{Conn}_{\check{H}}(\omega_{\D^\times}^{\rho})}

\nc{\Ind}{\on{Ind}}

\nc{\CTop}{{\mathcal Top}}

\nc{\ppart}{(\!(t)\!)}

\nc{\qu}{/\!/}

\nc{\gen}{\on{gen}}

\begin{document}

\renewcommand{\thepart}{\Roman{part}}

\renewcommand{\partname}{\hspace*{30mm} Part}

\title[Local Langlands correspondence and affine Kac-Moody
algebras]{Local geometric Langlands correspondence and affine Kac-Moody
algebras}

\dedicatory{Dedicated to Vladimir Drinfeld on his 50th birthday}

\author{Edward Frenkel}\thanks{Both authors were supported by the
DARPA grant HR0011-04-1-0031 and E.F. was also supported by the NSF
grant DMS-0303529.}

\address{Department of Mathematics, University of California,
  Berkeley, CA 94720, USA}

\email{frenkel@math.berkeley.edu}

\author{Dennis Gaitsgory}

\address{Department of Mathematics, Harvard University,
Cambridge, MA 02138, USA}

\email{gaitsgde@math.harvard.edu}

\date{August 2005; Revised: November 2005}

\maketitle

\setcounter{tocdepth}{1}

\tableofcontents

\section*{Introduction}

Let $\fg$ be a simple Lie algebra over $\BC$ and $G$ a connected
algebraic group with Lie algebra $\fg$. The affine Kac-Moody algebra
$\wh\fg$ is the universal central extension of the formal loop agebra
$\fg\ppart$. Representations of $\wh\fg$ have a parameter, an
invariant bilinear form on $\fg$, which is called the
level. Representations corresponding to the bilinear form which is
equal to minus one half of the Killing form are called representations
of {\em critical level}. Such representations can be realized in
spaces of global sections of twisted D-modules on the quotient of the
loop group $G\ppart$ by its "open compact" subgroup $K$, such as
$G[[t]]$ or the Iwahori subgroup $I$.

This is the first in a series of papers devoted to the study of the
categories of representations of the affine Kac-Moody algebra $\hg$ of
the critical level and D-modules on $G\ppart/K$ from the point of view
of a geometric version of the local Langlands correspondence. Let us
explain what we mean by that.

\ssec{} First of all, we recall the classical setting of local
Langlands correspondence. Let $\hCK$ be a local non-archimedian field
such as ${\mathbb F}_q\ppart$ and $G$ a connected reductive algebraic
group over $\hCK$. The local Langlands correspondence sets up a
relation between two different types of data. Roughly speaking, the
first data are the equivalence classes of homomorphisms, denoted by
$\sigma$, from the Galois group of $\hCK$ (more precisely, its
version, called the Weil-Deligne group) to $\cG$, the Langlands dual
group of $G$. The second data are the isomorphism classes of
irreducible smooth representations, denoted by $\pi$, of the locally
compact group $G(\hCK)$ (we refer the reader to \cite{Vogan} for a
precise formulation of this correspondence).

A naive analogue of this correspondence in the geometric situation is
as follows. Since the geometric analogue of the Galois group is the
fundamental group, the geometric analogue of a homomorphism from the
Galois group of $\hCK$ to $\cG$ is a $\cG$-local system on $\on{Spec}
\hCK$. Now we wish to replace $\hCK = {\mathbb F}_q\ppart$ by
$\C\ppart$. Then $\Spec \BC\ppart$ is the formal punctured disc
$\D^\times$. By a $\cG$-local system on $\D^\times$ we will always
understand its de Rham version: a $\cG$-bundle on $\D^\times$ with a
meromorphic connection that may have a pole of an arbitrary order at
the origin. By analogy with the classical local Langlands
correspondence, we would like to attach to such a local system a
representation of the formal loop group $G\ppart =
G(\C\ppart)$. However, we will argue in this paper that in contrast to
the classical setting, this representation of $G\ppart$ should be
defined not on a vector space, but on a {\em category} (see
\secref{HCh action of groups} where the notion of a group acting on a
category is spelled out).

Thus, to each $\cG$-local system $\sigma$ we would like to attach an
abelian category $\CC_\sigma$ equipped with an action of the ind-group
$G\ppart$. This is what we will mean by a geometric local Langlands
correspondence for the formal loop group $G\ppart$. This
correspondence may be viewed as a "categorification" of the
classical local Langlands correspondence, in the sense that we expect
the Grothendieck groups of the categories $\CC_\sigma$ to "look
like" irreducible smooth representations of $G(\hCK)$. At the moment
we cannot characterize $\CC_\sigma$ in local terms. Instead, we 
shall now explain how this local correspondence fits
in with the pattern of the global geometric Langlands correspondence.

In the global geometric Langlands correspondence we start with a
smooth projective connected curve $X$ over $\C$ with distinct marked
points $x_1,...,x_n$. Let $\sigma^{\on{glob}}$ be a $\cG$-local
system on $\overset\circ{X} = X \backslash \{x_1,...,x_n\}$, i.e., a
$\cG$-bundle on $X \backslash \{x_1,...,x_n\}$ with a connection
which may have poles of arbitrary order at the points
$x_1,\ldots,x_n$. Let $\Bun_G^{x_1,...,x_n}$ be the moduli stack
classifying $G$-bundles on $X$ with the full level structure at
$x_1,...,x_n$ (i.e., trivializations on the formal discs $\D_{x_i}$
around $x_i$). Let $\fD(\Bun_G^{x_1,...,x_n})\mod$ be the category of
D-modules on $\Bun_G^{x_1,...,x_n}$. One defines, as in \cite{BD},
the Hecke correspondence between $\Bun_G^{x_1,...,x_n}$ and
$\overset\circ{X} \times \Bun_G^{x_1,...,x_n}$ and the
notion of a Hecke "eigensheaf" on $\Bun_G^{x_1,...,x_n}$ with the
"eigenvalue" $\sigma^{\on{glob}}$.

The {\em Hecke correspondence} is the following moduli space:
$${\mc Hecke}=\{(\CP,\CP',x,\phi)\; | \: \Bun_G^{x_1,...,x_n},x \in
\overset\circ{X}, \phi:\CP|_{\circ{X}\backslash x}
\overset\sim\to \CP'|_{\circ{X}\backslash x}\}.$$ It is
equipped with the projections
$$
\begin{array}{ccccc}
& & {\mc Hecke} & & \\
& \stackrel{\hl}\swarrow & & \stackrel{\hr}\searrow & \\
\Bun_G^{x_1,...,x_n} & & & & \overset\circ{X} \times \Bun_G^{x_1,...,x_n}
\end{array}
$$
where $\hl(\CP,\CP',x,\phi) = \CP$ and $\hr(\CP,\CP',x,\phi) =
(x,\CP')$. The fiber of ${\mc Hecke}$ over $(x,{\mc P}')$ is
isomorphic to $\Gr_x^{\CP'}$, the twist of the affine Grassmannian
$\Gr_x = G(\hCK_x)/G(\hCO_x)$ by the $G(\hCO_x)$--torsor of
trivializations of $\CP'|_{\D_x}$ (here we denote by $\hCO_x$ and
$\hCK_x$ the completed local ring of $X$ at $x$ and its field of
fractions, respectively). The stratification of $\Gr_x$ by
$G(\hCO_x)$--orbits induces a stratification of ${\mc Hecke}$. The
strata are parametrized by the set of isomorphism classes of
irreducible representations of the Langlands dual group $\cG$. To each
such isomorphism class $V$ therefore corresponds an irreducible
D-module on ${\mc Hecke}$ supported on the closure of the orbit
labeled by $V$. We denote it by $\CF^{\on{glob}}_V$.

One defines the Hecke functors $H_V, V \in \on{Irr}(\Rep)$ from the
derived category of D-modules on $\Bun_G^{x_1,...,x_n}$ to the derived
category of D-modules on $\overset\circ{X} \times
\Bun_G^{x_1,...,x_n}$ by the formula
$$
H_V({\mc F}) = \hr_!(\hl{}^*({\mc F}) \otimes \CF^{\on{glob}}_V).
$$
A D-module on $\Bun_G^{x_1,...,x_n}$ is called a {\em Hecke eigensheaf
with eigenvalue} $\sigma^{\on{glob}}$ if we are given isomorphisms
\begin{equation}    \label{compat}
H_V({\mc F}) \simeq V_{\sigma^{\on{glob}}} \boxtimes {\mc F}
\end{equation}
of D-modules on $\overset\circ{X} \times \Bun_G^{x_1,...,x_n}$ which
are compatible with the tensor product structure on the category of
representations of $\cG$ (here $V_{\sigma^{\on{glob}}}$ is the
associated vector bundle with a connection on $\overset\circ{X}$
corresponding to $\sigma^{\on{glob}}$ and $V$).

The aim of the global geometric Langlands correspondence is to
describe the category
$\fD(\Bun_G^{x_1,...,x_n})^{\on{Hecke}}_{\sigma^{\on{glob}}}\mod$ of
such eigensheaves.

For example, if there are no marked points, and so
$\sigma^{\on{glob}}$ is unramified everywhere, it is believed that
this category is equivalent to the category of vector spaces, provided
that $\sigma^{\on{glob}}$ is sufficiently generic. In particular, in
this case $\fD(\Bun_G)^{\on{Hecke}}_{\sigma^{\on{glob}}}\mod$ should
contain a unique, up to an isomorphism, irreducible object, and all
other objects should be direct sums of its copies.  The irreducible
Hecke eigensheaf may be viewed as a geometric analogue of an
unramified automorphic function from the classical global Langlands
correspondence. This Hecke eigensheaf has been constructed by
A. Beilinson and V. Drinfeld in \cite{BD} in the case when
$\sigma^{\on{glob}}$ has an additional structure of an "oper".

In order to explain what we expect from the category
$\fD(\Bun_G^{x_1,...,x_n})^{\on{Hecke}}_{\sigma^{\on{glob}}}\mod$ when
the set of marked points is non-empty, let us revisit the classical
situation. Denote by ${\mathbb A}$ the ring of ad\`eles of the field
of rational functions on $X$. Let $\pi_{\sigma^{\on{glob}}}$ be an
irreducible automorphic representation of the ad\`elic group
$G({\mathbb A})$ corresponding to $\sigma^{\on{glob}}$ by the
classical global Langlands correspondence. Denote by
$\left(\pi_{\sigma^{\on{glob}}}\right)_{x_1,...,x_n}$ the subspace of
$\pi_{\sigma^{\on{glob}}}$ spanned by vectors unramified away from
$x_1,...,x_n$. Then
$\left(\pi_{\sigma^{\on{glob}}}\right)_{x_1,...,x_n}$ is a
representation of the locally compact group $\underset{i=1,...,n}\Pi\,
G(\hat\CK_{x_i})$ (here $\hat\CK_{x_i}$ denotes the local field at
$x_i$). A basic compatibility between the local and global classical
Langlands correspondences is that this representation should be
isomorphic to the tensor product of local factors
$$\bigotimes_{i=1,\ldots,n} \pi_{\sigma_i},$$ where $\pi_{\sigma_i}$ is
the irreducible representation of $G(\hat\CK_i)$, attached via the
local Langlands correspondence to the restriction $\sigma_i$ of
$\sigma^{\on{glob}}$ to the formal punctured disc around $x_i$.

In the geometric setting we view the category
$\fD(\Bun_G^{x_1,...,x_n})^{\on{Hecke}}_{\sigma^{\on{glob}}}\mod$ as a
"categorification" of the representation
$\left(\pi_{\sigma^{\on{glob}}}\right)_{x_1,...,x_n}$. Based on this,
we expect that there should be a natural functor
\begin{equation} \label{loc to glob}
\underset{i=1,...,n}\bigotimes\, \CC_{\sigma_i}\to
\fD(\Bun_G^{x_1,...,x_n})^{\on{Hecke}}_{\sigma^{\on{glob}}}\mod,
\end{equation}
relating the local and global categories. Moreover, we expect this
functor to be an equivalence when $\sigma^{\on{glob}}$ is sufficiently
generic. This gives us a basic compatibility between the local and
global geometric Langlands correspondences.

\ssec{}

How can we construct the categories $\CC_\sigma$ and the corresponding
functors to the global categories? At the moment we see two ways to do
that. In order to explain them, we first illustrate the main idea on a
toy model.

\medskip

Let $G$ be a split reductive group over $\Z$, and $B$ its Borel
subgroup. A natural representation of the finite group $G(\Fq)$ is
realized in the space of complex (or $\ol{\mathbb Q}_\ell$-) valued
functions on the quotient $G(\Fq)/B(\Fq)$. We can ask what is the
"correct" analogue of this representation when we replace the field
$\Fq$ by the complex field and the group $G(\Fq)$ by $G(\C)$. This may
be viewed as a simplified version of our quest, since instead of
considering $G(\Fq\ppart)$ we now look at $G(\Fq)$.

The quotient $G(\Fq)/B(\Fq)$ is the set of $\Fq$-points of the
algebraic variety $G/B$ defined over $\Z$ called the flag variety of
$G$. Let us recall the Grothendieck {\it faisceaux-fonctions} dictionary: if
${\mc F}$ is an $\ell$-adic sheaf on an algebraic variety $V$ over
${\mathbb F}_q$ and $x$ is an ${\mathbb F}_q$-point of $V$, then we
have the Frobenius conjugacy class $\on{Fr}_x$ acting on the stalk
${\mc F}_x$ of ${\mc F}$ at $x$. Hence we can define a $\ol{\mathbb
Q}_\ell$-valued function $\text{\tt f}_{q}({\mc F})$ on the set of
${\mathbb F}_{q}$-points of $X$, whose value at $x$ is
$\on{Tr}(\on{Fr}_x,{\mc F}_x)$. We also obtain in the same way a
function on the set $V({\mathbb F}_{q^n})$ of ${\mathbb
F}_{q^n}$-points of $V$ for $n>1$. This passage from $\ell$-adic
sheaves to functions satisfies various natural
properties. This construction identifies the Grothendieck
group of the category of $\ell$-adic sheaves on $V$ with a subgroup of
the direct product of the spaces of functions on $V({\mathbb F}_{q^n}), n>0$ (see
\cite{Laumon}). Thus, the category of $ell$-adic sheaves (or its
derived category) may be viewed as a categorification of this space of
functions.

This suggests that in order to pass from $\Fq$ to $\C$ we first need
to replace the notion of a function on $(G/B)(\Fq)$ by the notion of
an $\ell$-adic sheaf on the variety $(G/B)_{\Fq} = G/B
\underset{\Z}\otimes \Fq$.

Next, we replace the notion of an $\ell$-adic sheaf on $G/B$
considered as an algebraic variety over $\Fq$, by a similar notion of
a constructible sheaf on $(G/B)_{\C} = G/B \underset{\Z}\otimes \C$
which is an algebraic variety over $\C$. The group $G_\BC$ naturally
acts on $(G/B)_{\C}$ and hence on this category. We shall now apply
two more metamorphoses to this category.

Recall that for a smooth complex algebraic variety $V$ we have a
Riemann-Hilbert correspondence which is an equivalence between the
derived category of constructible sheaves on $V$ and the derived
category of D-modules on $V$ that are holonomic and have regular
singularities. Thus, over $\C$ we may pass from constructible sheaves
to D-modules. Generalizing this, we consider the category of all
D-modules on the flag variety $(G/B)_\C$. This category carries a
natural action of $G_\BC$.

Let us also recall that by taking global sections we obtain a functor
from the category of D-modules on $(G/B)_\C$ to the category of
$\fg$-modules. Moreover, A. Beilinson and J. Bernstein have proved
\cite{BB} that this functor is an equivalence between the category of
D-modules on $(G/B)_\C$ and the category of $\g$-modules on which the
center of the universal enveloping algebra $U(\g)$ acts through the
augmentation character.  Observe that the latter category also carries
a natural $G_\BC$-action that comes from the adjoint action of $G_\BC$
on $\fg$.

We arrive at the following conclusion: a meaningful geometric analogue
of the notion of representation of $G(\Fq)$ is that of a {\em
category} equipped with an action of $G_\BC$. In particular, an
analogue of the space of functions on $G(\Fq)/B(\Fq)$ is the category
$\fD((G/B)_\C)\mod$, which can be also realized as the category of
$\g$-modules with a fixed central character.

\medskip

Our challenge is to find analogues of the above two categories in
the case when the reductive group $G$ is replaced by its loop group
$G\ppart$. The exact relation between them will be given by a loop
group analogue of the Beilinson-Bernstein equivalence, and will be by
itself of great interest to us.

As the previous discussion demonstrates, one possibility is to
consider representations of the complex loop group $G\ppart$ on
various categories of D-modules on the ind-schemes $G\ppart/K$,
where $K$ is an "open compact" subgroup of $G\ppart$, such as $G[[t]]$ or
the Iwahori subgroup $I$ (the preimage of a Borel subgroup $B \subset
G$ under the homomorphism $G[[t]] \to G$). The other possibility is to
consider various categories of representations of the Lie algebra
$\g\ppart$, or of its universal central extension $\ghat$, because the
group $G\ppart$ still acts on $\ghat$ via the adjoint action.

\ssec{}    \label{Satake}

To explain the main idea of this paper, we consider an important
example of a category of D-modules which may be viewed as a
"categorification" of an irreducible {\em unramified} representation
of the group $G(\hCK)$, where $\hCK = \Fq\ppart$. We recall that a
representation $\pi$ of $G(\hCK)$ is called unramified if it contains
a non-zero vector $v$ such that $G(\hCO) v = v$, where $\hCO =
\Fq[[t]]$. The spherical Hecke algebra $H(G(\hCK),G(\hCO))$ of
bi-$G(\hCO)$-invariant compactly supported functions on $G(\hCK)$ acts
on the subspace spanned by such vectors.

The Satake isomorphism identifies $H(G(\hCK),G(\hCO))$ with the
representation ring $\on{Rep}(\cG)$ of finite-dimensional
representations of the Langlands dual group $\cG$ \cite{La}. This
implies that equivalence classes of irreducible unramified
representations of $G(\hCK)$ are parameterized by semi-simple conjugacy
classes in the Langlands dual group $\cG$. This is in fact a baby
version of the local Langlands correspondence mentioned above, because
a semi-simple conjugacy class in $\cG$ may be viewed as an equivalence
class of unramified homomorphisms from the Weil group $W_{\hCK}$ to
$\cG$ (i.e., one that factors through the homomorphism $W_{\hCK} \to
W_{\Fq} \simeq \Z$).

For a semi-simple conjugacy class $\gamma$ in $\cG$ denote by
$\pi_\gamma$ the corresponding irreducible unramified representation
of $G(\hCK)$. It contains a unique, up to a scalar, vector $v_\ga$
such that $G(\hCO) v_\gamma = v_\gamma$. It also satisfies the
following property. For a finite-dimensional representation $V$ of
$\cG$ denote by $F_V$ the element of $H(G(\hCK),G(\hCO))$
corresponding to $[V] \in \on{Rep}(\cG)$ under the Satake
isomorphism. Then we have $F_V \cdot v_\gamma = \on{Tr}(\gamma,V)
v_\gamma$ (to simplify our notation, we omit a $q$-factor in this
formula).

Now we embed $\pi_\gamma$ into the space of locally constant functions
on $G(\hCK)/G(\hCO)$, by using matrix coefficients, as follows:
$$
u \in \pi_\gamma \mapsto f_u, \qquad f_u(g) = \langle u,g v_\gamma
\rangle,
$$
where $\langle,\rangle$ is an invariant bilinear form on
$\pi_\gamma$. Clearly, the functions $f_u$ are right
$G(\hCO)$-invariant and satisfy the condition
\begin{equation}    \label{hecke for functions}
f \star F_V = \on{Tr}(\gamma,V) f,
\end{equation}
where $\star$ denotes the convolution product. Let
$C(G(\hCK)/G(\hCO))_\gamma$ be the space of locally constant functions
on $G(\hCK)/G(\hCO)$ satisfying \eqref{hecke for functions}. We have
constructed an injective map $\pi_\gamma \to
C(G(\hCK)/G(\hCO))_\gamma$, and one can show that for generic $\gamma$
it is an isomorphism.

Thus, we obtain a realization of irreducible unramified
representations of $G(\hCK)$ in functions on the quotient
$G(\hCK)/G(\hCO)$. According to the discussion in the previous
subsection, a natural complex geometric analogue of the space of
functions on $G(\hCK)/G(\hCO)$ is the category of (right) D-modules on
$G\ppart/G[[t]]$. The latter has the structure of an ind-scheme over
$\C$ which is called the affine Grassmannian and is denoted by
$\on{Gr}_G$.

The classical Satake isomorphism has a categorical version due to
Lusztig, Drinfeld, Ginzburg and Mirkovi\'c-Vilonen (see \cite{MV})
which may be formulated as follows: the category of
$G[[t]]$-equivariant D-modules on $\on{Gr}_G$, equipped with the
convolution tensor product, is equivalent to the category 
${\mathcal R}ep(\cG)$ of finite-dimensional representations of $\cG$ as a
tensor category. For a representation $V$ of $\cG$ let $\CF_V$ be the
corresponding $D$-module on $\Gr_G$. A D-module $\CF$ on $\on{Gr}_G$
satisfies the geometric analogue of the property \eqref{hecke for
functions} if we are given isomorphisms
\begin{equation}    \label{hecke for sheaves}
\alpha_V: \CF \star \CF_V \overset{\sim}\longrightarrow \underline{V}
\otimes \CF, \qquad V \in \on{Ob} {\mathcal R}ep(\cG)
\end{equation}
satisfying a natural compatibility with tensor products.
In other words, observe that we now have two monoidal
actions of the tensor category ${\mathcal R}ep(\cG)$ on the category
$\fD(\Gr_G)\mod$ of right D-modules on $\Gr_G$: one is given by
tensoring D-modules with $\underline{V}$, the vector space underlying
a representation $V$ of $\cG$, and the other is given by convolution
with the D-module $\CF_V$. The collection of isomorphisms $\alpha_V$
in \eqref{hecke for sheaves} should give us an isomorphism between
these two actions applied to the object $\CF$.

Let $\fD(\Gr_G)^{\on{Hecke}}\mod$ be the category whose objects are
the data $(\CF,\{ \alpha_V\})$, where $\CF$ is a D-module on $\Gr_G$
and $\{ \alpha_V \}$ are the isomorphisms \eqref{hecke for sheaves}
satisfying the above compatibility. This category carries a natural
action of the loop group $G\ppart$ that is induced by the (left)
action of $G\ppart$ on the Grassmannian $\Gr_G$. We believe that the
category $\fD(\Gr_G)^{\on{Hecke}}\mod$, together with this action of
$G\ppart$, is the "correct" geometric analogue of the unramified
irreducible representations of $G(\Fq\ppart)$ described above. Thus,
we propose that
\begin{equation}    \label{conj 0}
\CC_{\sigma_0} \simeq \fD(\Gr_G)^{\on{Hecke}}\mod,
\end{equation}
where $\sigma_0$ is the trivial $\cG$-local system on $\D^\times$.
This is our simplest example of the conjectural categories
$\CC_{\sigma}$, and indeed its Grothendieck group "looks like" an
unramified irreducible representation of $G(\hCK)$.

\ssec{}

Next, we attempt to describe the category $\CC_{\sigma_0}$ 
in terms of representations of the affine Kac-Moody algebra $\ghat$. 
Since the affine analogue of the Beilinson-Bernstein equivalence
is {\it a priori} not known, the answer is not as obvious as in the
finite-dimensional case. However, the clue is provided by the 
Beilinson-Drinfeld construction of the Hecke eigensheaves.

The point of departure is a theorem of \cite{FF} which states that the
completed universal enveloping algebra of $\wh\fg$ at the critical
level has a large center. More precisely, according to \cite{FF}, it
is isomorphic to the algebra of functions of the affine ind-scheme
$\Op(\D^\times)$ of $\cg$-{\it opers} over the formal punctured disc
(where $\cg$ is the Langlands dual of the Lie algebra $\g$). Thus,
each point $\chi \in \Op(\D^\times)$ defines a character of the
center, and hence the category $\hg_\crit\mod_\chi$ of discrete
$\wh\fg$-modules of critical level on which the center acts according
to the character $\chi$.

We recall that a $\cg$-oper (on a curve or on a disc) is a $\cG$-local
system plus some additional data. This notion was introduced in
\cite{DS,BD} (see \secref{def of oper} for the definition). Thus, we
have a natural forgetful map $\Op(\D^\times) \to
\on{LocSys}_{\cG}(\D^\times)$, where $\on{LocSys}_{\cG}(\D^\times)$ is
the stack of $\cG$-local systems on $\D^\times$.\footnote{Note that it
is not an algebraic stack, but in this paper we will work with its
substacks which are algebraic.} In this subsection we will restrict our 
attention to those
opers, which extend regularly to the formal disc $\D$; they correspond
to points of a closed subscheme of regular opers $\Op^\reg\subset
\Op(\D^\times)$. In particular, the local systems on $\D^\times$ defined by
such opers are unramified, i.e., they extend to local systems on $\D$,
which means that they are isomorphic to the trivial local system
(non-canonically, since the group $\cG$ acts by automorphisms of the
trivial local system).

For a fixed point $x\in X$ Beilinson and Drinfeld construct a
local-to-global functor $\hg_\crit\mod\to \fD(\Bun^x_G)\mod$ as a
Beilinson-Bernstein type localization functor by realizing $\Bun^x_G$
as the quotient $G\ppart/G(X-x)$.

Given a regular oper on the formal disc $\D$ around $x$, consider the
restriction of this localization functor to the category
$\hg_\crit\mod_\chi$. It was shown in \cite{BD} that the latter
functor is non-zero if and only if $\chi$ extends to an oper on the
global curve $X$, and in that case it gives rise to a functor
$$\hg_\crit\mod_\chi\to
\fD(\Bun^x_G)^{\on{Hecke}}_{\sigma^{\on{glob}}}\mod,$$ where
$\sigma^{\on{glob}}$ is the $\cG$-local system on $X$ corresponding to
the above oper.

This construction, combined with \eqref{loc to glob}, suggests that
for every regular oper $\chi$ on $\D$ we should have an equivalence of
categories

\begin{equation}   \label{single oper}
\CC_{\sigma_0}\simeq \hg_\crit\mod_\chi.
\end{equation}

Thus, we have two conjectural descriptions of the category
$\CC_{\sigma^0}$: one is given by \eqref{conj 0}, and the other by
\eqref{single oper}. Comparing the two, we obtain a conjectural
analogue of the Beilinson-Bernstein equivalence for the affine
Grassmannian:
\begin{equation}    \label{bb}
\fD(\Gr_G)^{\on{Hecke}}\mod \simeq \hg_\crit\mod_{\chi}
\end{equation}
for any $\chi \in \Op^{\on{reg}}$. In fact, as we shall see later, we
should have an equivalence as in \eqref{bb} for every trivialization
of the local system on $\D$ corresponding to the oper $\chi$. In
particular, the group of automorphisms of such a local system, which
is non-canonically isomorphic to $\cG$, should act on the category
$\hg_\crit\mod_{\chi}$ by automorphisms. In a sense, it is this action
that replaces the Satake parameters of irreducible unramified
representations of $G(\hCK)$ in the geometric setting.

Let us note that the equivalence conjectured in \eqref{bb} does not
explicitly involve the Langlands correspondence. Thus, our attempt to 
describe the simplest of the categories $\CC_\sigma$ has already paid 
dividends: it has led us to a formulation of Beilinson-Bernstein type 
equivalence for $\Gr_G$.

It is instructive to compare it with the Beilinson-Bernstein
equivalence for a finite-dimensional flag variety $(G/B)_{\BC}$, which
says that the category of D-modules on $(G/B)_{\BC}$ is equivalent to
the category of $\fg$-modules with a fixed central character. Naively,
one might expect that the same pattern holds in the affine case as
well, and the category $\fD(\Gr_G)\mod$ is equivalent to the category
of $\ghat_\crit$-modules with a fixed central character. However, in
contrast to the finite-dimensional case, the category $\fD(\Gr_G)\mod$
carries an additional symmetry, namely, the monoidal action of the
category $\Rep$ (which can be traced back to the action of the
spherical Hecke algebra in the classical setting). The existence of
this symmetry means that, unlike the category $\hg_\crit\mod_{\chi}$,
the category $\fD(\Gr_G)\mod$ is a $\cG$-equivariant category (in
other words, $\fD(\Gr_G)\mod$ is a category over the stack
$\on{pt}/\cG$, see below). From the point of view of Langlands correspondence, 
this equivariant structure is related to the fact that $\cG$ is the group
of automorphism of the trivial local system $\cG$. In order to obtain a
Beilinson-Bernstein type equivalence, we need to de-equivariantize
this category and replace it by $\fD(\Gr_G)^{\on{Hecke}}\mod$.

\ssec{}

Our next goal is to try to understand in similar terms what the
categories $\CC_\sigma$ look like for a general local system
$\sigma$. Unfortunately, unlike the unramified case, we will not be
able to construct them directly as categories of D-modules on some
homogeneous space of $G\ppart$. The reason for this can be traced to
the classical picture. If $\sigma$ is ramified, then the corresponding
irreducible representation $\pi$ of the group $G(\Fq\ppart)$ does not
contain non-zero vectors invariant under $G(\Fq[[t]])$, but it
contains vectors invariant under a smaller compact subgroup $K \subset
G(\Fq[[t]])$. As in the ramified case, we can realize $\pi$, by taking
matrix coefficients, in the space of functions on $G(\Fq\ppart)/K$
with values in the space $\pi^K$ of $K$-invariant vectors in $\pi$
satisfying a certain Hecke property. However, unlike the case of
unramified representations, $\pi^K$ generically has dimension greater
than one. When we pass to the geometric setting, we need, roughly
speaking, to find a proper "categorification" not only for the space
of functions on $G(\Fq\ppart)/K$ (which is the category of D-modules
on the corresponding ind-scheme, as explained above), but also for
$\pi^K$ and for the Hecke property. In the case when $\sigma$ is
tamely ramified, we can take as $K$ the Iwahori subgroup $I$. Then the
desired categorification of $\pi^I$ and the Hecke property can be
constructed following R. Bezrukavnikov's work \cite{Bez}, as we will
see below. This will allow us to relate the conjectural category
$\CC_\sigma$ to the category of D-modules on $G\ppart/I$. But we do
not know how to do that for more general local systems.

Therefore, we try first to describe the categories $\CC_\sigma$ in
terms of the category of representations of $\ghat$ at the critical
level rather than categories of D-modules on homogeneous spaces of
$G\ppart$.

A hint is once again provided by the Beilinson-Drinfeld construction
of Hecke eigensheaves from representations of $\ghat$ at the critical
level described above, because it may be applied in the ramified
situation as well. Extending \eqref{single oper}, we conjecture that
for {\em any} oper $\chi$ on $\D^\times$ and the corresponding local
system $\sigma$, we have an equivalence of categories
\begin{equation} \label{meta pointwise}
\CC_\sigma \simeq \hg_\crit\mod_{\chi}
\end{equation}
equipped with an action of $G\ppart$. This statement implies, in
particular, that the category $\hg_\crit\mod_{\chi}$ depends not on
the oper $\chi$, but only on the underlying local system! This is by
itself a deep conjecture about representations of $\ghat$ at the
critical level.

At this point, in order to elaborate more on what this conjecture
implies and to describe the results of this paper, we will need to
discuss a more refined version of the local geometric Langlands
correspondence indicated above. For that we have to use the notion of
an abelian or a triangulated category over a stack. In the abelian
case this is an elementary notion, introduced, e.g., in \cite{Ga1}. It
amounts to a sheaf (in the faithfully-flat topology) of abelian
categories over a given stack $\CY$. When $\CY$ is an affine scheme
$\Spec(A)$, this amounts to the notion of $A$-linear abelian
category. In the triangulated case, some extra care is needed, and we
refer the reader to \cite{Ga2} for details. The only property of this
notion needed for the discussion that follows is that whenever $\CC$
is a category over $\CY$ and $\CY'\to \CY$ is a map of stacks, we can
form the base-changed category $\CC' = \CC \underset{\CY}\times \CY'$;
in particular, for a point $y\in \CY(\BC)$ we have the category-fiber
$\CC_y$.

A refined version of the local geometric Langlands correspondence
should attach to any $\CY$-family of $\cG$-local systems $\wt\sigma$
on $\D^\times$ a category $\CC_{\wt\sigma}$ over $\CY$, equipped with
an action of $G\ppart$, in a way compatible with the above base change
property.  Such an assignment may be viewed as a category $\CC$ over
the stack $\on{LocSys}_{\cG}(\D^\times)$ equipped with a "fiberwise"
action of $G\ppart$. Then the categories $\CC_\sigma$ discussed above
may be obtained as the fibers of $\CC$ at $\C$-points $\sigma$ of
$\on{LocSys}_{\cG}(\D^\times)$.

We shall now present a refined version of \eqref{meta pointwise}. Namely,
although at the moment we cannot construct $\CC$, the following
meta-conjecture will serve as our guiding principle:
\begin{equation} \label{meta-conj}
\hg_\crit\mod \simeq \CC\underset{\on{LocSys}_{\cG}(\D^\times)}\times
\Op(\D^\times).
\end{equation}
We will not even try to make this meta-conjecture precise in this
paper. Instead we will derive from it some more concrete conjectures,
and the goal of this paper will be to give their precise formulation
and provide evidence for their validity.

\ssec{}

Let us first revisit the unramified case discussed above. Since the
trivial local system $\sigma_0$ has $\cG$ as the group of its
automorphisms, we have a natural map from the stack $\on{pt}/\cG$ to
$\on{LocSys}_{\cG}(\D^\times)$. Let us denote by $\CC_\reg$ the base
change of (the still conjectural category) $\CC$ under the above map.
Then, by definition, we have an equivalence:
\begin{equation}    \label{by def}
\CC_{\sigma_0}\simeq \CC_\reg\underset{\on{pt}/\cG}\times \on{pt}.
\end{equation}

Now observe that the geometric Satake equivalence of \secref{Satake}
gives us an action of the tensor category ${\mathcal R}ep(\cG)$ on
$\fD(\Gr_G)\mod$, $V,\CF \mapsto \CF \star \CF_V$. This precisely
amounts to saying that $\fD(\Gr_G)\mod$ {\it is a category over
the stack} $\on{pt}/\cG$. Moreover, we then have the following base
change equivalence:
\begin{equation}    \label{b change}
\fD(\Gr_G)^{\on{Hecke}}\mod \simeq \fD(\Gr_G)\mod
\underset{\on{pt}/\cG}\times \on{pt}.
\end{equation}

Combining \eqref{by def} and \eqref{b change}, we arrive at the
following generalization of \eqref{conj 0}:
\begin{equation} \label{conj reg}
\CC_\reg\simeq \fD(\Gr_G)\mod.
\end{equation}

Let us now combine this with \eqref{meta-conj}. Let us denote by
$\hg_\crit\mod_{\on{reg}}$ the subcategory of $\wh\fg$-modules at the 
critical level on which the center acts in such a way that their
scheme-theoretic support in $\Op(\D^\times)$ belongs to
$\Op^{\on{reg}}$.  By the definition of the map $\Op(\D^\times)\to
\on{LocSys}_{\cG}(\D^\times)$, its restriction to $\Op^{\on{reg}}$
factors through a map $\Op^{\on{reg}}\to \on{pt}/\cG$, which assigns
to an oper $\chi$ the $\cG$-torsor on $\Op^{\on{reg}}$ obtained by
taking the fiber of $\chi$ at the origin in $\D$.

Thus, combining \eqref{meta-conj} with the identification
\begin{equation}    \label{ident}
\hg_\crit\mod_{\on{reg}} \simeq
\hg_\crit\mod\underset{\Op(\D^\times)}\times \Op^\reg,
\end{equation}
we obtain the following statement:
\begin{equation}  \label{conj 2}
\hg_\crit\mod_\reg \simeq \fD(\Gr_G)\mod\underset{\on{pt}/\cG}\times
\Op^\reg,
\end{equation}

By making a further base change with respect to an embedding
of the point-scheme into $\Op^\reg$ corresponding to some regular
oper $\chi$, we obtain \eqref{bb}. Thus, \eqref{conj 2} is a family
version of \eqref{bb}.

Let us now comment on one more aspect of the conjectural equivalence
proposed in \eqref{conj 2}. With any category $\CC$ acted on by
$G\ppart$ and an "open compact" subgroup $K\subset G\ppart$ we can associate
the category $\CC^K$ of $K$-equivariant objects. Applying this to
$\hg_\crit\mod_\reg$ we obtain the category, consisting of those
representations, which are $K$-integrable (i.e., those, for which the
action of $\on{Lie} K$ may be exponentiated to that of $K$). In the
case of $\fD(\Gr_G)\mod$ we obtain the category of $K$-equivariant
D-modules in the usual sense.

Let us take $K=G[[t]]$, and compare the categories obtained from the
two sides of \eqref{conj 2}:
\begin{equation}  \label{conj 2 equiv}
\hg_\crit\mod_\reg^{G[[t]]}\simeq
\fD(\Gr_G)^{G[[t]]}\mod\underset{\on{pt}/\cG}\times \Op^\reg.
\end{equation}

However, the Satake equivalence mentioned above says that
$\fD(\Gr_G)^{G[[t]]}\simeq {\mathcal R}ep(\cG)$, implying that the
RHS of \eqref{conj 2 equiv} is equivalent to the category of
quasi-coherent sheaves on $\Op^\reg$:
\begin{equation}    \label{unram case}
\hg_\crit\mod_{\on{reg}}^{G[[t]]} \simeq \QCoh(\Op^{\on{reg}}).
\end{equation}

The latter equivalence is not conjectural, but has already been
established in \cite{FG}, Theorem 6.3 (see also \cite{BD}).

Thus, we obtain a description of the category of modules at the critical
level with a specified integrability property and a condition on the
central character as a category of quasi-coherent sheaves on a scheme
related to the Langlands dual group. Such a description is a prototype
for the main conjecture of this paper, described below.

\ssec{}

The main goal of this paper is to develop a picture similar to the one
presented above, for {\em tamely ramified} local systems $\sigma$
on $\D^\times$, i.e., those with regular singularity at the origin and
unipotent monodromy. The algebraic stack classifying such local
systems is isomorphic to $\cN/\cG$, where $\cN\subset \cg$ is the
nilpotent cone. Let $\CC_\nilp$ be the corresponding hypothetically
existing category over $\cN/\cG$ equipped with a fiberwise action of
$G\ppart$.

We shall first formulate a conjectural analogue of Theorem
\eqref{unram case} in this set-up. As we will see, one essential
difference with the unramified case is the necessity to consider
derived categories.

Denote by $I\subset G[[t]]$ the Iwahori subgroup; it is the preimage
of a once and for all fixed Borel subgroup of $G$ under the
homomorphism $G[[t]] \to G$. We wish to give a description of the
$I$-monodromic part $D^b(\CC_\nilp)^{I,m}$ of the derived category
$D^b(\CC_\nilp)$ that is similar in spirit to the one obtained in the
unramified case (the notion of Iwahori-monodromic derived category
will be introduced in \secref{I-categories}).

In \secref{nilpotent opers} we will introduce a subscheme $\nOp\subset
\Op(\D^\times)$ of {\em opers with nilpotent singularities}.  Note
that $\nOp$ contains as a closed subscheme the scheme $\Op^\reg$ of
regular opers. Denote by $\hg_\crit\mod_\nilp$ the subcategory of
$\hg_\crit\mod$ whose objects are the $\hg_\crit$-modules whose
scheme-theoretic support in $\Op(\D^\times)$ is contained in $\nOp$.

Let $\tN$ be the Springer resolution of $\cN$. We will show in
\secref{nilpotent opers} that the composition $\nOp\to
\Op(\D^\times)\to \on{LocSys}_{\cG}(\D^\times)$ factors as
$$\nOp\overset{\Res^\nilp}\longrightarrow \tN/\cG\to
\cN/\cG\hookrightarrow \on{LocSys}_{\cG}(\D^\times).$$ The first map,
denoted by $\Res^\nilp$, is smooth.

Our \mainconjref{main} describes the (bounded) derived category of
$\hg_\crit\mod_\nilp$ as follows:
\begin{equation} \label{conj 3}
D^b(\hg_\crit\mod_\nilp)^{I,m}\simeq
D^b(\QCoh(\tg/\cG\underset{\cg/\cG}\times \nOp)),
\end{equation}
where $\tg\to \cg$ is Grothendieck's alteration. This is an analogue
for nilpotent opers of Theorem \eqref{unram case}.

As will be explained below, the scheme
$\tg/\cG\underset{\cg/\cG}\times \nOp$, appearing in the RHS of
\eqref{conj 3}, has a natural interpretation as the moduli space of
Miura opers with nilpotent singularities (see \secref{Miura}).  The
main motivation for the above conjecture came from the theory of
Wakimoto modules introduced in \cite{FF:si,F:wak}. Namely, to each
Miura oper with nilpotent singularity one can attach a Wakimoto module
which is an object of the category $\hg_\crit\mod_\nilp^{I,m}$. Our
\conjref{conj 3} extends this "pointwise" correspondence to an
equivalence of categories.

\ssec{}

Next, we would like to formulate conjectures concerning $\CC_\nilp$
that are analogous to \eqref{conj reg} and \eqref{conj 2}, and relate
them to Conjecture \eqref{conj 3} above.

The main difficulty is that we do not have an explicit description of 
$\CC_\nilp$ in terms of D-modules as the one for $\CC_\reg$,
given by \eqref{conj reg}. Instead, we will be able to describe
a certain base change of $\CC_\nilp$, suggested by the work of 
S.~Arkhipov and R.~Bezrukavnikov \cite{Bez,AB}.

Let $\Fl_G=G\ppart/I$ be the affine flag variety and the affine
Grassmannian, corresponding to $G$ and $\fD(\Fl_G)\mod$ the category
of right D-modules on $\Fl_G$. The group $G\ppart$ naturally acts on
$\fD(\Fl_G)\mod$. According to \cite{AB}, the triangulated category
$D^b(\fD(\Fl_G)\mod)$ is a category over the stack $\tN/\cG$.

We propose the following conjecture, describing the hypothetically
existing category $\CC_\nilp$, which generalizes \eqref{conj reg}:

\begin{equation}  \label{C nilp and reg}
D^b(\CC_\nilp)\underset{\cN/\cG}\times \tN/\cG\simeq
D^b(\fD(\Fl_G)\mod).
\end{equation}

Combining \eqref{C nilp and reg} with our meta-conjecture
\eqref{meta-conj}, we arrive at the statement
\begin{equation}   \label{conj 1}
D^b(\hg_\crit\mod_\nilp) \simeq
D^b(\fD(\Fl_G)\mod)\underset{\tN/\cG}\times \nOp
\end{equation}
(the RHS of the above equivalence uses the formalism of 
triangulated categories over stacks from \cite{Ga2}). Note
that Conjecture \eqref{conj 1} is an analogue for opers with 
nilpotent singularities of Conjecture \eqref{conj 2}
for regular opers.

We would now like to explain the relation of Conjectures 
\eqref{C nilp and reg} and \eqref{conj 1} to the description
of $D^b(\hg_\crit\mod_\nilp)$ via quasi-coherent sheaves,
given by Conjecture \eqref{conj 3} once we pass to the
$I$-monodromic category.

We propose the following description of the category
$D(\CC_\nilp)^{I,m}$:
\begin{equation} \label{C as coherent}
D^b(\CC_\nilp)^{I,m}\simeq D^b(\QCoh(\tg\underset{\cg}\times \cN/\cG)).
\end{equation}

Let us note that Conjecture \eqref{C as coherent} is compatible
with \eqref{C nilp and reg}. Namely, by combining the two
we obtain the following:
\begin{equation}    \label{I-mon}
D^b(\fD(\Fl_G)\mod)^{I,m} \simeq D^b(\QCoh(\tg\underset{\cg}\times \tN/\cG)).
\end{equation}

However, this last statement is in fact a theorem, which is one
of the main results of Bezrukavnikov's work \cite{Bez}.

Finally, combining \eqref{I-mon} and \eqref{conj 1},
we arrive to the statement of Conjecture \eqref{conj 3}, providing
another piece of motivation for it, in addition to the one via
Wakimoto modules given above.

\ssec{}

The principal objective of our project is to prove conjectures \eqref{conj 2}
and \eqref{conj 3}. In the present paper we review some background
material necessary to introduce the objects we are studying and
formulate the above conjectures precisely. We also prove two results
concerning the category of representations of affine Kac-Moody
algebras at the critical level which provide us with additional
evidence for the validity of these conjectures.

\medskip

Our first result is \mainthmref{equiv of quot}, and it deals with a
special case of \mainconjref{main}. Namely, in \secref{partially
integrable} we will explain that if $\CC$ is a category endowed with
an action of $G\ppart$, the corresponding category $\CC^{I,m}$ of
Iwahori-monodromic objects admits a Serre quotient, denoted
$^f\CC^{I,m}$, by the subcategory, consisting of the so-called
partially integrable objects. (Its classical analog is as follows:
given a representation $\pi$ of a locally compact group $G(\hCK)$, we
first take the subspace $\pi^I$ is Iwahori-invariant vectors, and then
inside $\pi^I$ we take the subspace of vectors corresponding to the
sign character of the Iwahori-Hecke algebra.)

Performing this procedure on the two sides of \eqref{conj 1}, we
should arrive at an equivalence of the corresponding triangulated
categories:
\begin{equation} \label{thm1}
^f D^b(\hg_\crit\mod_\nilp)^{I,m} \simeq {}^f
{}D^b(\fD(\Fl_G)\mod)^{I,m}\underset{\tN/\cG}\times \nOp
\end{equation}

However, using Bezrukavnikov's result (see \thmref{bezr quot}), the
RHS of the above expression can be rewritten as
$D^b(\QCoh(\Spec(h_0)\times \nOp))$, where $h_0$ is a
finite-dimensional commutative algebra isomorphic to
$H(\cG/\cB,\BC)$. The resulting description of $^f
D^b(\hg_\crit\mod_\nilp)^{I,m}$ is our \mainthmref{equiv of quot}. In
fact, we show that at the level of quotient categories by partially
integrable objects, the equivalence holds not only at the level of
triangulated categories, but also at the level of abelian ones:
$$\Catf\simeq \QCoh(\Spec(h_0)\times \nOp).$$

We note that while we use \cite{Bez} for motivational purposes, the 
proofs presented in this paper are independent of the results of
\cite{Bez}.

\medskip

Our second main result is \thmref{Hecke fully faithful}. We construct
a natural functor from the RHS of \eqref{conj 2} to the LHS and prove
that it is fully faithful at the level of derived categories.

\ssec{}

Let us now describe the structure of the paper. It is logically
divided into 5 parts.

\medskip

Part I is a review of results about opers and Miura opers.

\medskip

In Part II we discuss various categories of representations of affine
Kac-Moody algebras at the critical level. We give more precise
formulations of the conjectural equivalences that we mentioned above
and the interrelations between them. In particular, we prove one of
our main results, \thmref{Hecke fully faithful}.

\medskip

In Part III we review the Wakimoto modules. We present a definition of
Wakimoto modules by means of a kind of semi-infinite induction
functor. We also describe various important properties of these
modules.

\medskip

In Part IV we prove \mainthmref{equiv of quot} which establishes a
special case of our conjectural equivalence of categories \eqref{conj
3}.

\medskip

Part V is an appendix, most of which is devoted to the formalism of
group-action on categories.

Finally, a couple of comments on notation.

We will write $X \underset{Z}\times Y$ for the fiber product of
schemes $X$ and $Y$ equipped with morphisms to a scheme $Z$. To
distinguish this notation from the notation for associated fiber
bundles, we will write $\CY\overset{K}\times \CF_K$ for the fiber
bundle associated to a principal $K$-bundle $\CF_K$ over some base,
where $K$ is an algebraic group and $\CY$ is a $K$-space. We also
denote this associated bundle by $\CY_{\CF_K}$.

If a group $G$ acts on a variety $X$, we denote by $X/G$ the
stack-theoretic quotient. If $X$ is affine, we denote by $X\qu G$ the
GIT quotient, i.e., the spectrum of the algebra of invariant
functions. We have a natural morphism $X/G \to X\qu G$.

\bigskip

\bigskip

\noindent {\bf Acknowledgments}. We thank A. Beilinson,
R. Bezrukavnikov and B. Feigin for useful discussions.

\newpage

\vspace*{10mm}

{\Large \part{Opers and Miura opers}}

\vspace*{10mm}


We this Part we collect the definitions and results on opers and
Miura opers. As a mathematical object, opers first appeared in
\cite{DS}, and their connection to representations of affine Kac-Moody
algebras at the critical level was discovered in \cite{FF}.

In \secref{opers} we recall the definition of opers following
\cite{BD} and the explicit description of the scheme classifying them
as a certain affine space.

In \secref{opers with sing} we study opers on the formal punctured
disc with a prescribed form of singularity at the closed point. After
reviewing some material from \cite{BD}, we show that the subscheme of
opers with regular singularities and a specific value of the residue
can be interpreted as a scheme of opers with nilpotent singularities,
which we denote by $\on{Op}_\fg^\nilp$. We show that the scheme
$\on{Op}_\fg^\nilp$ admits a natural secondary residue map to the
stack $\fn/B\simeq \wt{\CN}/G$.

In \secref{Miura} we study Miura opers. The notion of Miura oper was
introduced in \cite{F:wak}, following earlier work of Feigin and
Frenkel. By definition, a Miura oper on a curve $X$ is an oper plus a
reduction of the underlying $G$-local system to a Borel subgroup $B^-$
opposite to the oper Borel subgroup $B$. The functor $\on{MOp}_\fg(X)$
of Miura opers admits a certain open subfunctor, denoted by
$\on{MOp}_{\fg,\gen}(X)$ that corresponds to generic Miura opers. The
(D-) scheme classifying the latter is affine over $X$, and as was
shown in \cite{F:wak}, it is isomorphic to the (D-) scheme of
connections on some fixed $H$-bundle over $X$, where $H$ is the Cartan
quotient of $B$. \footnote{The corresponding space of $H$-connections
on the formal punctured disc and its map to the space of opers were
introduced in \cite{DS} as the phase space of the generalized mKdV
hierarchy and the Miura transformation from this space to the phase
space of the generalized KdV hierarchy.}  The new results in this
section are \propref{MT and residue} which describes the forgetful map
from generic Miura opers to opers over the locus of opers with regular
singularities and \thmref{four versions of Miura} which describes the
behavior of Miura opers and generic Miura opers over
$\on{Op}^\nilp_\fg$.

In \secref{algebroid on opers} we introduce the isomonodromy groupoid
over the ind-scheme $\on{Op}_{\fg}(\D^\times)$ and its various subschemes. We
recall the definition of Poisson structure on the space
$\on{Op}_\fg(\D^\times)$ of opers on the formal punctured disc
introduced in \cite{DS}.  Following \cite{BD} and \cite{CHA}, we
interpret this Poisson structure as a structure of Lie algebroid on
the cotangent sheaf $\Omega^1(\on{Op}_\fg(\D^\times))$ and following
\cite{DS} we show that it is isomorphic to the Lie algebroid of the
isomonodromy groupoid on the space of opers. The new results in this
section concern the behavior of this algebroid along the subscheme
$\on{Op}^\nilp_\fg$.

\bigskip

\section{Opers}    \label{opers}

\ssec{Definition of opers}    \label{def of oper}

Throughout Part I (except in \secref{non-adjoint}), we will assume
that $G$ is a simple algebraic group of {\em adjoint type}. Let $B$ be
its Borel subgroup and $N = [B,B]$ its unipotent radical, with the
corresponding Lie algebras $\fn \subset \fb \subset \fg$. There is an
open $B$-orbit ${\bf O}\subset \fg/\fb$, consisting of vectors which
are invariant with respect to the radical $N\subset B$, and such that
all of their negative simple root components, with respect to the
adjoint action of $H = B/N$, are non-zero. This orbit may also be
described as the $B$-orbit of the sum of the projections of simple
root generators $f_\imath$ of any nilpotent subalgebra $\fn^-$, which
is in generic position with $\fb$, onto $\fg/\fb$. The torus $H = B/N$
acts simply transitively on ${\bf O}$, so ${\bf O}$ is an
$H$-torsor. Note in addition that ${\bf O}$ is invariant with respect
to the action of $\BG_m$ on $\fg$ by dilations.

\medskip

Let $X$ be a smooth curve, or the formal disc $\D = \on{Spec} (\hCO)$,
where $\hCO$ is a one-dimensional complete local ring, or a formal
punctured disc $\D^\times=\Spec (\hCK)$, where $\hCK$ is the field of
fractions of $\hCO$. We will denote by $\omega_X$ the canonical line
bundle on $X$; by a slight abuse of notation we will identify it with
the corresponding $\BG_m$-torsor on $X$.

\medskip

Suppose we are given a principal $G$-bundle $\CF_G$ on $X$, together
with a connection $\nabla$ (automatically flat) and a reduction
$\CF_B$ of $\CF_G$ to the Borel subgroup $B$ of $G$. Then we define
the relative position of $\nabla$ and $\CF_B$ (i.e., the failure of
$\nabla$ to preserve $\CF_B$) as follows. Locally, choose any
connection $\nabla'$ on $\CF$ preserving $\CF_B$, and take the
difference $\nabla - \nabla'\in \fg_{\CF_G}\simeq \fg_{\CF_B}$.  It is
clear that the projection of $\nabla - \nabla'$ to $(\fg/\fb)_{\CF_B}
\otimes \omega_X$ is independent of $\nabla'$; we will denote it by
$\nabla/\CF_B$. This $(\fg/\fb)_{\CF_B}$-valued one-form on $X$ is by
definition the relative position of $\nabla$ and $\CF_B$.

\medskip

Following Beilinson and Drinfeld, see \cite{BD}, Sect. 3.1, and
\cite{BD1}, one defines a $\fg$-{\em oper} on $X$ to be a triple
$(\CF_G,\nabla,\CF_B)$, where $\CF_G$ is a principal $G$-bundle
$\CF_G$ on $X$, $\nabla$ is a connection on $\CF_G$, and $\CF_B$ is a
$B$-reduction of $\CF_G$ such that the one-form $\nabla/\CF_B$ takes
values in
$${\bf O}_{\CF_B,\omega_X}:=
{\bf O}\overset{B\times \BG_m}\times (\CF_B\times \omega_X)
\subset(\fg/\fb)_{\CF_B}\otimes \omega_X.$$

\medskip

Consider the $H$-bundle $\omega^{\crho}_X$ on $X$, induced
from the line bundle $\omega_X$ by means of the homomorphism
$\crho:\BG_m\to H$. (The latter is well-defined, since $G$ was assumed
to be of adjoint type.) 

\begin{lem}    \label{FH}
For an oper $(\CF_G,\nabla,\CF_B)$, 
the induced $H$-bundle $\CF_H:=N\backslash \CF_B$
is canonically isomorphic to $\omega^{\crho}_X$.
\end{lem}

\begin{proof}

We have to show that for every simple root $\alpha_\imath:B\to \BG_m$,
the line bundle $\BC^{\alpha_\imath}_{\CF_B}$ is canonically
isomorphic to $\omega_X$.

Decomposing $\nabla/\CF_B$ with respect to negative simple roots, we
obtain for every $\imath$ a non-vanishing section of the line bundle
$$\BC^{-\alpha_\imath}_{\CF_B}\otimes \omega_X.$$
This provides the required identification.

\end{proof}

Here is an equivalent way to think about opers. Let us choose a
trivialization of the $B$-bundle $\CF_B$, and let $\nabla^0$ be the
tautological connection on it. Then an oper is given by a connection
$\nabla$ of the form
\begin{equation}    \label{form of nabla}
\nabla = \nabla^0 + \sum_{\imath} \phi_\imath \cdot f_\imath + \bq,
\end{equation}
where each $\phi_\imath$ is a nowhere vanishing one-form on $X$, and
$\bq$ is a $\fb$-valued one-form. If we change the trivialization of
$\CF_B$ by $\bg:X\to B$, the connection will get transformed by the
corresponding gauge transformation:
\begin{equation} \label{gauge action}
\nabla\mapsto \on{Ad}_{\bg}(\nabla):=\nabla^0 +
\on{Ad}_{\bg}\left(\sum_{\imath} \phi_\imath \cdot f_\imath + \bq\right)-
\bg^{-1}\cdot d(\bg).
\end{equation}

The following will be established in the course of the proof of
\propref{free}:

\begin{lem}    
If $\on{Ad}_{\bg}(\nabla)=\nabla$, then $\bg=1$.
\end{lem}

\medskip

In a similar way one defines the notion of an $R$-family of opers on
$X$, where $R$ is an arbitrary commutative $\BC$-algebra. We shall
denote this functor by $\on{Op}_\fg(X)$. For $X=\D$ (resp.,
$X=\D^\times$) some extra care is needed when one defines the notion
of $R$-family of bundles. To simplify the notation we will choose a
coordinate $t$ on $\D$, thereby identifying $\hCO\simeq \BC[[t]]$ and
$\hCK\simeq \BC\ppart$. Although this choice of the coordinate
trivializes $\omega_X$ by means of $dt$, we will keep track of the
distinction between functions and forms by denoting the $\hCO$-module
$\omega_\D$ (resp., the $\hCK$-vector space $\omega_{\D^\times}$) by
by $\BC[[t]]dt$ (resp., $\BC\ppart dt$).

\medskip

By definition, an $R$-family of $G$-bundles on $X=\D$ is a $G$-bundle
on $\Spec(R[[t]])$, or which is the same, a compatible family of
$G$-bundles on $\Spec(R[t]/t^i)$; such a family is always locally
trivial in the \'etale topology on $\Spec(R)$.

An $R$-family of $G$-bundles on $\D^\times$ is a $G$-bundle on
$\Spec\left(R\ppart\right)$, which {\it we require to be locally trivial 
in the \'etale topology in $\Spec(R)$}. 

Connections on the trivial $R$-family of $G$-bundles on $\D$ and
$\D^\times$ are expressions of the form $\nabla^0+\phi$, where $\phi$
is an element of $\fg\otimes R[[t]]dt$ and $\fg\otimes R\ppart dt$,
respectively. Gauge transformations are elements of $G(R[[t]])$ and
$G(R\ppart)$, respectively, and they act on connections by the formula
\eqref{gauge action}.

Thus, $\on{Op}_\fg(\D)$ and $\on{Op}_\fg(\D^\times)$ are well-defined
as functors on the category of $\BC$-algebras. Following \cite{BD},
Sect. 3.1.10, we will prove below that these functors are
representable by a scheme and ind-scheme, respectively.

\ssec{D-scheme picture}  \label{D-scheme picture}

When $X$ is a curve of finite type, a natural way to think of $\fg$-opers
on $X$ is in terms of D-schemes (we refer the reader to \cite{CHA}, Sect. 2.3
for the general discussion of D-schemes, and to \cite{CHA}, Sect. 2.6.8
for the discussion of opers in this context).

Namely,  let us notice that the notion of $R$-family of $\fg$-opers on $X$
makes sense when $R$ is a $\fD_X$-algebra, i.e., a quasi-coherent
sheaf of algebras over $X$, endowed with a connection. 

Repeating the argument of \propref{free} (see below), one obtains that
the above functor on the category of $\fD_X$-algebras is
representable; the corresponding affine $\fD_X$-scheme, denoted
$\on{Op}_\fg(X)^{\fD}$, is isomorphic to the $\fD_X$-scheme of jets
into a finite-dimensional vector space.

By definition, for a $\BC$-algebra $R$ we have:
\begin{equation} \label{usual from D}
\on{Op}_\fg(X)(R)\simeq
\Hom_{\fD_X}\left(\on{\Spec}(R\otimes \CO_X),\on{Op}_\fg(X)^{\fD}\right).
\end{equation}
If $\D$ is the formal neighborhood of a point $x\in X$ with a local
coordinate $t$, the functors $\on{Op}_\fg(\D)$ and
$\on{Op}_\fg(\D^\times)$ are reconstructed as
\begin{equation}  \label{formal from D}
R\mapsto
\Hom_{\fD_X}\left(\on{\Spec}(R[[t]]),\on{Op}_\fg(X)^{\fD}\right)
\text{ and } R\mapsto
\Hom_{\fD_X}\left(\on{\Spec}\left(R\ppart\right),\on{Op}_\fg(X)^{\fD}
\right),
\end{equation}
respectively.

In addition, one also has an isomorphism between the scheme
$\on{Op}_\fg(\D)$ and the fiber of $\on{Op}_\fg(X)^{\fD}$, regarded as
a mere scheme over $X$, at $x\in X$.

\ssec{Explicit description and canonical representatives}
\label{canon}

To analyze opers on $\D^\times$ (resp., $\D$) more explicitly we will
continue to use an identification $\CO\simeq \BC[[t]]$, and we will
think of opers as equivalence classes of connections of the form
\begin{equation} \label{exp form of nabla}
\nabla = \nabla^0 + \sum_{\imath} \phi_\imath(t)dt \cdot f_\imath +
\bq(t)dt,
\end{equation}
where now $\phi_\imath$ and $\bq$ are now elements of $R\ppart$ and
$\fb\otimes R\ppart$, respectively (resp., $R[[t]]$ and $\fb\otimes
R[[t]]$) such that each $\phi_\imath$ is invertible. Two such
connections are equivalent, if they can be conjugated one into another
by a gauge transformation by an element of
$\Hom\bigl(\Spec\left(R\ppart\right),B\bigr)$ (resp.,
$\Hom\bigl(\Spec\left(R[[t]]\right),B\bigr)$.

Let us observe that since $H\simeq B/N$ acts simply-transitively on
${\bf O}$, any connection as above can be brought to the form when all
the functions $\phi_\imath(t)$ are equal to $1$. Moreover, this can be
done uniquely, up to a gauge transformation by means of
$\Hom(\Spec(R\ppart,N)$.

\medskip

The operator $\on{ad} \crho$ defines the principal grading on $\fb$,
with respect to which we have a direct sum decomposition $\fb=
\underset{d\geq 0}\oplus\, \fb_d$. Set
$$p_{-1} = \sum_{\imath} f_\imath;$$
we shall call this element the negative principal nilpotent.

Let $p_1$ be the unique element of $\fn$ such that $\{
p_{-1},2\rv,p_1\}$ is an $\sw_2$-triple.  Let $V_{\can} =
\underset{d>0}\oplus\, V_{\can,d}$ be the space of $\on{ad}
p_1$-invariants in $\fn$. The operator $\on{ad} p_{-1}$ acts from
$\fb_{d+1}$ to $\fb_{d}$ injectively for all $d\geq 0$, and we have
$\fb_d = [p_{-1},\fb_{d+1}] \oplus V_{\can,d}$.

We will call the $\BG_m$-action on $V_{\can}$, resulting from the
above grading, "principal".  We will call the $\BG_m$-action on
$V_{\can}$, obtained by multiplying the principal one by the standard
character, "canonical". Recall that by a theorem of Kostant, the map
\begin{equation} \label{Kostant}
V_{\can}\overset{\bc\mapsto p_{-1}+\bc}\longrightarrow \fg\to
\fg\qu G\simeq \fh\qu W
\end{equation}
is an isomorphism. This map is compatible with the canonical
$\BG_m$-action on $V_{\can}$ and the action on $\fh\qu W$, induced by
the standard $\BG_m$-action on $\fh$.

\begin{prop}[\cite{DS}]    \label{free}  
The gauge action of $\Hom(\Spec(R\ppart),B)$ on the set of connections
of the form \eqref{exp form of nabla} is free. Each gauge equivalence
class contains a unique representative of the form
\begin{equation} \label{can form}
\nabla = \nabla^0 + p_{-1}dt+ \bv(t)dt, \,\, \bv(t)\in V_{\can}\otimes
R\ppart.
\end{equation}
\end{prop}

As we shall see, the same assertion with the same proof is valid if we
replace $R\ppart$ by $R[[t]]$. In what follows we will refer to
\eqref{can form} as the canonical representative of an oper.

\begin{proof}

We already know that we can bring a connection \eqref{exp form of
nabla} to the form
$$\nabla^0+p_{-1}dt+\bq(t)dt,$$ uniquely up to an element in
$\Hom(\Spec(R\ppart),N)$. We need to show now that there exists a
unique element $u(t)\in \fn\otimes R\ppart$ such that
$$\on{Ad}_{\exp(u(t))}\left(\nabla^0 + p_{-1}dt+\bq(t)dt\right) =
\nabla^0 + p_{-1}dt + \bv(t)dt,$$ $\bv(t)\in V_{\can}\otimes R\ppart$.

Let us decompose the unknown element $u(t)$ as $\underset{d}\Sigma\,
u_d(t)$, where $u_d(t)\in \fn_d\otimes R\ppart$, and we claim that we
can find the elements $u_d(t)$ by induction $d$. Indeed, let us assume
that $\bq_{d'}\in V_{\can,d'}$ for $d'< d$. Then $u_{d+1}(t)$ must
satisfy:
$$[u_{d+1}(t),p_{-1}]+\bq_d(t)\in V_{\can,d},$$
and this indeed has a unique solution.

\end{proof}

\begin{cor}    \label{can rep general}
The set of $R$-families of opers on $\D$ and $\D^\times$ is isomorphic
to $V_{\can}\otimes R[[t]]$ and $V_{\can}\otimes R\ppart$,
respectively. In particular, the functor $\on{Op}_\fg(\D)$ (resp.,
$\on{Op}_\fg(\D^\times)$) is representable by the scheme (resp.,
ind-scheme), isomorphic to $V_{\can}[[t]]$ (resp., $V_{\can}\ppart$).
\end{cor}

We should note, however, that the isomorphisms $\on{Op}_\fg(\D)\simeq
V_{\can}[[t]]$ and $\on{Op}_\fg(\D^\times)\simeq V_{\can}\ppart$ are
not canonical, since they depend on the choice of the coordinate $t$
on $\D$.

\medskip

By the very definition, on the scheme
$$
\on{Op}_\fg(\D)\widehat\times \D := \on{Spec}
(\on{Fun}\on{Op}_\fg(\D))[[t]]
$$
there exists a universal $G$-bundle $\CF_{G,\on{Op}_\fg(\D)}$ with a
reduction to a $B$-bundle $\CF_{B,\on{Op}_\fg(\D)}$ and a connection
$\nabla_{\on{Op}_\fg(\D)}$ in the $\D$-direction such that the triple
$(\CF_{G,\on{Op}_\fg(\D)},\nabla_{\on{Op}_\fg(\D)},\CF_{B,\on{Op}_\fg(\D)})$
is a $\on{Op}_\fg(\D)$-family of $\fg$-opers on $\D$.  By the above,
when we identify $\hCO\simeq \BC[[t]]$, the $G$-bundle
$\CF_{B,\on{Op}_\fg(\D)}$, and hence $\CF_{G,\on{Op}_\fg(\D)}$,
becomes trivialized. But this trivialization depends on the choice of
the coordinate.

In what follows we will denote by $\CP_{G,\on{Op}_\fg(\D)}$ (resp.,
$\CP_{B,\on{Op}_\fg(\D)}$) the restriction of
$\CF_{G,\on{Op}_\fg(\D)}$ (resp., $\CF_{B,\on{Op}_\fg(\D)}$) to the
subscheme $\on{Op}_\fg(\D)\subset \on{Op}_\fg(\D)\widehat\times \D$,
corresponding to the closed point of $\D$. Note that
$\CP_{G,\on{Op}_\fg(\D)}$ can also be defined as the torsor of
horizontal, with respect to the connection along $\D$, sections of
$\CF_{G,\on{Op}_\fg(\D)}$.

\ssec{Action of $\on{Aut}(\D)$}

Let $\on{Aut}(\D)$ (resp., $\on{Aut}(\D^\times)$) be the group-scheme 
(resp., group ind-scheme) of automorphisms of $\D$ (resp., $\D^\times$).
\footnote{Note that $\on{Aut}(\D)$ is not reduced, see \cite{BD},
Sect. 2.6.5}.  Since $\on{Op}_\fg(\D)$ (resp.,
$\on{Op}_\fg(\D^\times)$) is canonically attached to $\D$ (resp.,
$\D^\times$), it carries an action of $\on{Aut}(\D)$ (resp.,
$\on{Aut}(\D^\times)$), see \secref{topological algebroids} for the
definition of the latter notion.

By transport of structure, the action of $\on{Aut}(\D)$ on
$\on{Op}_\fg(\D)\widehat\times \D$ lifts onto
$\CF_{G,\on{Op}_\fg(\D)}$ and $\CF_{B,\on{Op}_\fg(\D)}$.  The
interpretation of $\CP_{G,\on{Op}_\fg(\D)}$ as the space of horizontal
sections of $\CF_{G,\on{Op}_\fg(\D)}$ implies that the action of
$\on{Aut}(\D)$ on $\on{Op}_\fg(\D)$ lifts also onto the $G$-torsor
$\CP_{G,\on{Op}_\fg(\D)}$.

\medskip

To a choice of a coordinate $t$ on $\D$ there corresponds a
homomorphism $\BG_m\to \on{Aut}(\D)$ that acts by the "loop rotation",
i.e., $t\mapsto c\cdot t$.  We shall now describe the resulting action
of $\BG_m$ on $\on{Op}_\fg(\D^\times)$ in terms of the isomorphism of
\corref{can rep general}:

\begin{lem}  \label{grading}  \hfill

\smallskip

\noindent{\em (1)} The trivialization $\CP_{G,\on{Op}_\fg(\D)}\simeq
G\times \on{Op}_\fg(\D)$ corresponding to the given choice of a
coordinate is compatible with the $\BG_m$-action, via the homomorphism
$\crho:\BG_m\to G$.

\smallskip

\noindent{\em (2)} The action of $c\in \BG_m$ on
$$\on{Op}_\fg(\D^\times)\simeq V_{\can}\ppart\simeq
\underset{d}\oplus\, V_{\can,d}\ppart$$ is given by
$$\bv_d(t)\in V_{\can,d}\ppart\mapsto c^{d+1}\cdot \bv_d(c\cdot t).$$
\end{lem}

\begin{proof}

By definition, the action of $c\in \BG_m$ on a connection in form
\eqref{can form} transforms it to
\begin{equation} \label{transformed connection}
\nabla^0 + p_{-1}d(c\cdot t)+\bv(c\cdot t)d(c\cdot t)=
\nabla^0 + c\cdot p_{-1}dt+c\cdot \bv(c\cdot t)dt.
\end{equation}
In order to bring it back to the form \eqref{can form} we need to apply
a gauge transformation by means of the constant $H$-valued function
$\crho(c)$. This implies point (1) of the lemma.

This gauge transformation transforms \eqref{transformed connection}
to $$\nabla^0 + p_{-1}dt+c\cdot \on{Ad}_{\crho(c)}(\bv(c\cdot t))dt.$$
This implies point (2) of the lemma, since
$\on{Ad}_{\crho(c)}(\bv_d(c\cdot t))=c^d\cdot \bv_d(c\cdot t)$.

\end{proof}

\ssec{Quasi-classics: the Hitchin space}    \label{Hitch}

Recall that the Hitchin space $\Hitch_\fg(X)$ corresponding to the Lie
algebra $\fg$ and a curve $X$ is a functor on the category of algebras
that attaches to $R$ the set of sections of the pull-back to
$\Spec(R)\times X$ of the fiber bundle $$(\fh\qu W)\overset{\BG_m}\times
\omega_X,$$ where $\fh\qu W:=\Spec\left(\Sym(\fh^*)^W\right)\simeq
\Spec\left(\Sym(\fg^*)^G\right)$ is endowed with a canonical action of
$\BG_m$.

When $X=\D$ or $X=\D^\times$, in the above definition we replace 
$\Spec(R)\times X$ by $\Spec(R[[t]])$ and $\Spec\left(R\ppart\right)$, 
respectively.

For $X=\D$ the Hitchin space is a scheme, isomorphic to
$\underset{d}\oplus\, V_{\can,d}\otimes \omega^{\otimes
{d+1}}_\D$. For $X=\D^\times$ this is an ind-scheme, isomorphic to
$\underset{d}\oplus\, V_{\can,d}\otimes \omega^{\otimes
{d+1}}_{\D^\times}$. In particular, $\Hitch_\fg(\D)$ (resp.,
$\Hitch_\fg(\D^\times)$) has a natural structure of group-scheme
(resp., group ind-scheme).

According to \cite{BD}, Sect. 2.4.1, the natural map
$$\Spec\left(\Sym(\check\fg\ppart/\check\fg[[t]])^{\check{\fg}[[t]]}
\right)\to \Hitch_\fg(\D)$$ is an isomorphism, where $\check\fg$ is
the Langlands dual Lie algebra.  This implies that the maps
$$\Fun\left(\Hitch_\fg(\D^\times)\right)\to
\biggl(\underset{k}{\underset{\longleftarrow}{\lim}}\,
\Sym\left(\check\fg\ppart/t^k\cdot
\check\fg[[t]]\right)\biggr)^{\check{\fg}\ppart}\to
\underset{k}{\underset{\longleftarrow}{\lim}}\,
\Sym\left(\check\fg\ppart/t^k\cdot
\check\fg[[t]]\right)^{\check{\fg}[[t]]}$$ are also isomorphisms.

\medskip

By \propref{free} the scheme $\on{Op}_\fg(\D)$ (resp.,
$\on{Op}_\fg(\D^\times)$) is non-canonically isomorphic to
$\Hitch_\fg(\D)$ (resp., $\Hitch_\fg(\D^\times)$).  However, one can
deduce from the proof (see \cite{BD}, Sect. 3.10.11) that
$\on{Op}_\fg(\D)$ (resp., $\on{Op}_\fg(\D^\times)$) is {\it
canonically} a torsor over $\Hitch_\fg(\D)$ (resp.,
$\Hitch_\fg(\D^\times)$).

In particular, the algebra $\Fun(\on{Op}_\fg(\D))$ acquires a
filtration, whose associated graded is $\Fun(\Hitch_\fg(\D))$. This
filtration can also be defined as follows, see \cite{BD},
Sect. 3.11.14:

We claim that there exists a flat $\BG_m$-equivariant family of
schemes over $\BA^1\simeq \Spec(\BC[\hslash])$, whose fiber over $1\in
\BA^1$ is $\on{Op}_\fg(\D)$, and whose fiber over $0\in \BA^1$ is
$\Hitch_\fg(\D)$.

Indeed, this family is obtained from $\on{Op}_\fg(\D)$ by replacing
the word "connection" by "$\hslash$-connection". The identification at
the special fiber results from Kostant's theorem that the adjoint
action of $B$ on the pre-image of ${\mb O}$ in $\fg$ is free, and the
quotient projects isomorphically onto $\fh\qu W$.

\ssec{The case of groups of non-adjoint type}    \label{non-adjoint}

In the rest of the paper we will have to consider the case when the
group $G$ is not necessarily of adjoint type. Let $Z(G)$ be the center
of $G$.

The notion of $R$-family of $G$-opers in this case is formally the
same as in the adjoint case, i.e., a triple $(\CF_G,\nabla,\CF_B)$,
where $\CF_G$ is an $R$-family of $G$-bundle on $X$, $\CF_B$ is its
reduction to $B$, and $\nabla$ is a connection on $\CF_G$ in the
$X$-direction, which satisfies the same condition on $\nabla/\CF_B$.

We will denote the functor of $R$-families of $G$-opers on $X$ by
$\on{Op}_G(X)$.  The difference with the adjoint case is that now
$\on{Op}_G(\D)$ is not representable by a scheme, but rather by a
Deligne-Mumford stack, which is non-canonically isomorphic to
$\on{Op}_\fg(\D)\times \on{pt}/Z(G)$, see \cite{BD}, Sect. 3.4.

The following statement established in \cite{BD}, Sect. 3.4, will
suffice for our purposes:

\begin{lem}
Every choice of the square root $\omega_X^{\frac{1}{2}}$ of the canonical
bundle gives a map of functors $\on{Op}_\fg(X)\to \on{Op}_G(X)$.
\end{lem}

In particular, the lemma implies that for every choice of a square 
root of $\omega_\D$, there exists a canonically defined family of
$G$-opers over $\on{Op}_\fg(\D)$. A similar statement holds
for $\D$ replaced by $\D^\times$.

\begin{proof}

One only has to show how to lift the $B/Z(G)$-bundle $\CF_{B/Z(G)}$ to
a $B$-bundle. This is equivalent to lifting the $H$-bundle $\CF_{H/Z(G)}$
to an $H$-bundle $\CF_H$. 

We set $\CF_H$ to be the bundle induced by means of the homomorphism
$2\crho:\BG_m\to H$ from the line bundle $\omega_X^{\frac{1}{2}}$. By
\lemref{FH}, it satisfies our requirement.

\end{proof}

\section{Opers with singularities}   \label{opers with sing}

\ssec{}   \label{order k}

For $X$ a curve of finite type over $\BC$ we shall fix $x\in X$ to be any
closed point. For $X=\D$, we let $x$ to be the unique closed point of
$\Spec(\BC[[t]])$. We shall now define the notion of  {\em $\fg$-oper on $X$ 
with singularity of order $k$ at $x$}. 

By definition, this a triple $(\CF_G,\nabla,\CF_B)$, where
$(\CF_G,\CF_B)$ are as in the definition of opers, but the connection
$\nabla$on $\CF_G$ is required to have a pole of order $k$ such that
\begin{equation}  \label{defn of ord k} 
(\nabla-\nabla')\on{mod} \fb_{\CF_B} \otimes \omega_X(k\cdot x)
\in {\bf O}_{\CF_B,\omega_X(k\cdot x)}\subset 
(\fg/\fb)_{\CF_B}\otimes \omega_X(k\cdot x),
\end{equation}
for any regular connection $\nabla'$ on $\CF_G$ that preserves $\CF_B$.

Again, if we trivialize $\CF_B$ and choose a coordinate $t$ near $x$,
the set of opers with singularity of order $k$ at $x$ identifies with
the set of equivalence classes of connections of the form
\begin{equation} \label{ord k form}
\nabla^0+ t^{-k} \left( \sum_{\imath} \phi_\imath(t)dt \cdot f_\imath
+ \bq(t)dt\right),
\end{equation}
where $\phi_\imath(t)$ are nowhere vanishing functions on $X$, and
$\bq(t)$ is a $\fb$-valued function. Two such connections are
equivalent if they are conjugate by means of an element of
$\Hom(X,B)$. Equivalently, opers with singularity of order $k$ at $x$
is the set of $\Hom(X,N)$-equivalence classes of connections of the
form $$\nabla^0 + t^{-k} \left( p_{-1}dt + \bq(t)dt\right)$$ for
$\bq(t)$ as above.

As in \lemref{FH} one has:

\begin{lem}
For $(\CF_G,\nabla,\CF_B)$-an oper on $X$ with singularity of order
$k$ at $x$, the $H$-bundle $\CF_H=N\backslash \CF_B$ is canonically
isomorphic to $\left(\omega_X(x)\right)^\crho$.
 \end{lem}

One defines the notion of $R$-family of opers on $X$ with singularity
of order $k$ at $x$ in a straightforward way. The corresponding
functor on the category of $R$-algebras will be denoted
$\on{Op}^{\ord_k}_\fg(X)$. We will be mainly concerned with the case
when $X=\D$; if no confusion is likely to occur, we will denote the
corresponding functor simply by $\on{Op}^{\ord_k}_\fg$. We have an
evident morphism of functors $\on{Op}^{\ord_k}_\fg\to
\on{Op}_\fg(\D^\times)$, and below we will see (see \corref{can form
ord k}) that $\on{Op}^{\ord_k}_\fg$ is a representable by a closed
subscheme of $\on{Op}_\fg(\D^\times)$.  Note that for $k=0$ we recover
$\on{Op}_\fg(\D)$, and we will often use the notation
$\on{Op}^\reg_\fg$ for it.

\medskip

As in the case of usual opers, there exists a naturally defined
functor $\on{Op}^{\ord_k}_\fg(X)^{\fD}$ on the category of
$\fD_X$-algebras.  The functors $\on{Op}^{\ord_k}_\fg(X)$ and
$\on{Op}^{\ord_k}_\fg$ are reconstructed by the analogs of
\eqref{usual from D} and \eqref{formal from D},
respectively. \corref{can form ord k} implies that
$\on{Op}^{\ord_k}_\fg(X)^{\fD}$ is representable by an affine
$\fD_X$-scheme, which over the curve $(X-x)$ is isomorphic to
$\on{Op}_\fg(X)^{\fD}$.
  
\ssec{Changing $k$}

\begin{prop}[\cite{BD}, 4.3]    \label{inductive system of opers}
For every $k$ there is a natural morphism of functors
$\on{Op}_\fg^{\ord_k}\to\on{Op}_\fg^{\ord_{k+1}}$. We have:
$$\on{Op}_\fg(\D^\times)\simeq
\underset{k}{\underset{\longrightarrow}{\lim}}\,
\on{Op}_\fg^{\ord_k}.$$
\end{prop}

\begin{proof}

Given a triple $(\CF_G,\nabla,\CF_B)\in \on{Op}_\fg^{\ord_k}(R)$ we
define the corresponding $(\CF'_G,\nabla',\CF'_B)\in
\on{Op}_\fg^{\ord_{k+1}}(R)$ as follows:

Let us choose (locally) a trivialization of $\CF_B$, and let us apply
the gauge transformation by means of $t^{\crho}\in H\ppart$, where $t$
is any uniformizer on $\D$. We thus obtain a different extension of
$\CF_B$ from $\D^\times$ to $\D$, and let it be our $\CF'_B$. It is
clear that $\CF'_B$ is independent of both the choice of the
trivialization and the coordinate.

Let $\CF'_G$ be the induced $G$-bundle, and $\nabla'$ the resulting
meromorphic connection on it. By \eqref{gauge action}, $\nabla'$ has
the form required by the \eqref{ord k form}.

\medskip

Let now $(\CF_G,\nabla,\CF_B)$ be an $R$-point of
$\on{Op}_\fg(\D^\times)$, represented as a gauge equivalence class of
some connection written in the form \eqref{exp form of nabla},
$\phi_\imath(t)\in \left(R\ppart\right)^\times$.

Consider the $R$-point of $H\ppart$ equal to $\left(\underset{\imath}\Pi\, 
(\check\omega_\imath)(\phi_\imath)\right)^{-1}\cdot t^{k\cdot \crho}$,
where each $\check\omega_\imath$ is regarded as a homomorphism
$\BG_m\to H$.
It is clear from \eqref{gauge action} that for $k$ large enough, the
resulting connection will be of the form \eqref{ord k form}.

\end{proof}

\ssec{Description in terms of canonical representatives}

By repeating the proof of  \propref{free}, we obtain:

\begin{lem}   \label{FHk}
For every $R$-point of $\on{Op}_\fg^{\ord_k}$, the canonical form of
its image in $\on{Op}_\fg(\D^\times)(R)$ is such that each homogeneous
component $\bv_d(t)$ has a pole in $t$ of order $\leq k\cdot (d+1)$.
\end{lem}

\begin{cor}    \label{can form ord k}
The morphisms of functors $\on{Op}_\fg^{\ord_k}\to
\on{Op}_\fg^{\ord_{k+1}}$ and $\on{Op}_\fg^{\ord_k}\to
\on{Op}_\fg(\D^\times)$ are closed embeddings.  The latter identifies
with the subscheme
$$\underset{d}\oplus\, t^{-k\cdot (d+1)}\cdot V_{\can,d}[[t]]\subset
V_{\can,d}\ppart.$$
\end{cor}

\begin{proof}
Evidently, given a point of $\on{Op}_\fg(\D^\times)$, written in the
canonical form \eqref{can form} such that $t^{k\cdot (d+1)}\cdot
\bv_d(t)\in V_{\can,d}[[t]]$, by applying the gauge transformation by
means of $t^{k\cdot \crho}$, we bring it to the form \eqref{ord k
form}. Thus, we obtain the maps
$$\on{Op}_\fg^{\ord_k}\rightleftarrows \underset{d}\oplus\, t^{-k\cdot
(d+1)}\cdot V_{\can,d}[[t]].$$ Finally, by induction on $d$ it is easy
to see that if some $\bg\in B\ppart(R)$ conjugates an $R$-point of
$\on{Op}_\fg^{\ord_k}$ to another point of $\on{Op}_\fg^{\ord_k}$,
then $\bg\in B[[t]](R)$.

\end{proof}


\medskip

Recall the ind-scheme $\Hitch_\fg(\D^\times)$, and let us denote by
$\Hitch_\fg^{\ord_k}$ its subscheme, corresponding to sections of
$(\fh\qu W)\overset{\BG_m}\times \omega_\D(k\cdot x)$. This is a scheme,
canonically isomorphic to $\underset{d}\oplus\, V_{\can,d}\otimes
\left(\omega_\D(k\cdot x)\right)^{\otimes d+1}$, which gives it a
structure of group-scheme. Evidently, $\Hitch_\fg(\D^\times)$ is
isomorphic to
$\underset{k}{\underset{\longrightarrow}{\lim}}\,\Hitch_\fg^{\ord_k}$. We
also have an isomorphism:
$$\Fun(\Hitch_\fg^{\ord_k})\simeq \Sym\left(\check\fg\ppart/t^k\cdot
\check\fg[[t]]\right)^{\check{\fg}[[t]]}.$$

\medskip

As in the case of $k=0$, the scheme $\Hitch_\fg^{\ord_k}$ acts simply
transitively on $\on{Op}_\fg^{\ord_k}$. Moreover,
$$\on{Op}_\fg^{\ord_k}\simeq
\on{Op}_\fg(\D)\overset{\Hitch_\fg(\D)}\times \Hitch_\fg^{\ord_k}.$$
This defines a filtration on the algebra $\Fun(\on{Op}_\fg^{\ord_k})$,
whose associated graded is $\Fun(\Hitch_\fg^{\ord_k})$. This
filtration can be alternatively described by the deformation procedure
mentioned at the end of \secref{Hitch}.

\ssec{Opers with regular singularities}    \label{reg sing} 

In the context of the previous subsection let us set $k=1$, in which
case we will replace the superscript $\ord_1$ by $\RS$, and call the
resulting scheme $\on{Op}_\fg^\RS$ "the scheme of opers with regular
singularities". The terminology is partly justified by the following
assertion, which will be proved in the next section:

\begin{prop}  \label{char RS}
If a $\BC$-point $(\CF_G,\nabla,\CF_B)$ of $\on{Op}_\fg(\D^\times)$
has regular singularities as a $G$-bundle with connection, then it
belongs to $\on{Op}_\fg^\RS$.
\end{prop}

\medskip

We claim now that there exists a canonical map
$\Res^\RS:\on{Op}_\fg^\RS\to \fh\qu W$, see \cite{BD}, Sect. 3.8.11:

\medskip

{}Recall first that if $(\CF_G,\nabla)$ is an $R$-family of
$G$-bundles on $X$ with a connection that has a pole of order $1$ at
$x$, its residue (or polar part) is well-defined as a section of
$\fg_{\CP_G}$, where $\CP_G$ is the restriction of $\CF_G$ to
$\Spec(R)\times x\subset \Spec(R)\times X$.  In other words, we obtain
an $R$-point of the stack $\fg/G$.
 
Given an $R$-point $(\CF_G,\nabla,\CF_B)\in \on{Op}_\fg^\RS$, we 
compose the above map with $\fg/\on{Ad}(G)\to \fh\qu W$. The resulting map
$\Spec(R)\to \fh\qu W$ is the map $\Res^\RS$.

\medskip

Explicitly, to a connection written as
\begin{equation}    \label{RS form}
 \nabla^0 + t^{-1} \left( \sum_{\imath} \phi_\imath(t)dt \cdot
 f_\imath + \bq(t)dt\right),
\end{equation}
we attach the projection to $\fh\qu W$ of the element 
$\underset{\imath}\Sigma \phi_\imath(0) \cdot f_\imath+\bq(0)$.

\medskip

Let $\varpi$ denote the tautological projection $\fh\to \fh\qu W$. For
$\cla\in \fh$ we will denote by $\on{Op}_\fg^{\RS,\varpi(\cla)}$ the
preimage under $\Res^\RS$ of the point $\varpi(\cla)\in \fh\qu W$.  From
the proof of \propref{inductive system of opers} for $k=0$ we obtain
that the subscheme $\on{Op}_\fg(\D)=:\on{Op}_\fg^\reg\subset
\on{Op}_\fg^\RS$ is contained in $\on{Op}_\fg^{\RS,\varpi(-\crho)}$.

Let us now describe the map $\Res^\RS$ in terms of the isomorphism of
\corref{can form ord k}:

\begin{lem}  \label{residue of canonical}
The composition
$$\underset{d}\oplus\, t^{-d-1}\cdot V_{\can,d}[[t]]\simeq
\on{Op}_\fg^\RS\overset{\Res^\RS}\to \fh\qu W$$ equals the map
$$\underset{d}\oplus\, t^{-d-1}\cdot V_{\can,d}[[t]]\to
\underset{d}\oplus\, V_{\can,d}\simeq V_{\can}\simeq \fh\qu W,$$ where
the last arrow is given by \eqref{Kostant}, and the first arrow is
defined as follows:

\begin{itemize}

\item For $d\neq 1$, this is the projection on the top polar part.

\item For $d=1$, this is the projection on the top polar part,
followed by the affine shift by $\frac{p_1}{4}$.

\end{itemize}
\end{lem}

\begin{proof}

By the proof of \propref{free}, we have to check that for any $\bv'\in
\underset{d\neq 1}\oplus\, V_{\can,d}$ and $\bv''\in V_{\can,1}$, the
elements of $\fg$ given by $p_{-1}+\bv'+\bv''-\crho$ and
$p_{-1}+\bv'+\bv''+\frac{p_1}{4}$ project to the same element of
$\fh\qu W$. However, this follows from the fact that
$\exp(\frac{p_1}{2})$ conjugates one to the other.

\end{proof}

\begin{cor}   \label{can form nilp}
Under the isomorphism of \corref{can form ord k} the subscheme
$\on{Op}_\fg^{\varpi(-\crho)}\subset \on{Op}_\fg^\RS$ identifies with
$$\underset{d}\oplus\, t^{-d}\cdot V_{\can,d}[[t]]\subset 
\underset{d}\oplus\, t^{-d-1}\cdot V_{\can,d}[[t]].$$
\end{cor}

Of course, as in the case of \corref{can form ord k}, the isomorphism
$\on{Op}_\fg^{\varpi(-\crho)}\simeq \underset{d}\oplus\, t^{-d}\cdot
V_{\can,d}[[t]]$ depends on the choice of the coordinate
$t$. Canonically, $\on{Op}_\fg^{\varpi(-\crho)}$ can be described in
terms of the Hitchin space as follows.

Let us denote $\Hitch_\fg^{\ord_1}$ by $\Hitch_\fg^\RS$, and let us
note that we have a natural homomorphism $\Hitch_\fg^\RS\to
V_{\can}$. Let $\Hitch_\fg^\nilp\subset \Hitch_\fg^\RS$ be the
preimage of $0$. 
The algebra of functions on $\Hitch_\fg^\nilp$ also admits the
following description, see \cite{F:wak}, Lemma 9.4:
$$\Fun(\Hitch_\fg^\nilp)\simeq
\Sym\left(\check\fg\ppart)/\Lie(\check{I})\right)^{\check{I}},$$ where
$\check{I}\subset \check{G}[[t]]$ is the Iwahori subgroup.

We have:
$$\on{Op}_\fg^{\varpi(-\crho)}\simeq
\on{Op}_\fg(\D)\overset{\Hitch_\fg(\D)}\times \Hitch_\fg^\nilp.$$

\medskip

Consider now the gradings on the algebras $\Fun(\on{Op}_\fg^{\ord_k})$,
$\Fun(\on{Op}_\fg^\RS)$ and $\Fun(\on{Op}_\fg^{\varpi(-\crho)})$, coming 
from $\BG_m\to \on{Aut}(\D)$, corresponding to some choice of a coordinate
$t$ on $\D$. From \lemref{grading} we obtain:

\begin{lem}  \label{descr of grading} \hfill

\smallskip

\noindent{\em (1)} The algebra $\Fun(\on{Op}_\fg^\RS)$ is
non-positively graded. The subalgebra, consisting of degree $0$
elements, identifies with $\Fun(\fh\qu W)$ under the map $\Res^\RS$.

\smallskip

\noindent{\em (2)} For every $k\geq 2$, the ideal of
$\Fun(\on{Op}_\fg^{\ord_k})\to \Fun(\on{Op}_\fg^\RS)$ is freely
generated by finitely many elements, each having a positive degree.

\smallskip

\noindent{\em (3)} The algebra $\Fun(\on{Op}_\fg^{\varpi(-\crho)})$ is
freely generated by elements of strictly negative degrees.

\end{lem}

\ssec{Opers with nilpotent singularities}  \label{nilpotent opers}

Let $X$ and $x$ be as above.  We define a $\fg$-oper on $X$ with a
nilpotent singularity at $x$ to be a triple $(\CF_G,\nabla,\CF_B)$,
where $(\CF_G,\CF_B)$ are as in the definition of opers, and the
connection $\nabla$ has a pole of order $1$ at $x$ such that for some
(or any) regular connection $\nabla'$ that preserves $B$, we have:
\begin{itemize}

\item(i) $(\nabla-\nabla')$, which a priori is an element of
$\fg_{\CF_G}\otimes \omega_X(x)$, is in fact contained in
$\fb_{\CF_B}\otimes \omega_X(x)+\fg_{\CF_G}\otimes \omega_X \subset
\fg_{\CF_G}\otimes \omega_X(x)$. Once this condition is satisfied, we
impose the following two:

\item(ii) $(\nabla-\nabla')\on{mod} \fg_{\CF_G}\otimes \omega_X$,
which is an element of $\fb_{\CP_B}\simeq \fb_{\CF_B}\otimes
\omega_X(x)/\fb_{\CF_B}\otimes \omega_X$, must be contained in
$\fn_{\CP_B}\subset \fb_{\CP_B}$, where $\CP_B$ is the fiber of
$\CF_B$ at $x$.

\item(iii) $(\nabla-\nabla')\on{mod} \fb_{\CF_B} \otimes \omega_X(x)$,
which is a section of $(\fg/\fb)_{\CF_B}\otimes \omega_X$, must be
contained in ${\bf O}_{\CF_B,\omega_X}\subset (\fg/\fb)_{\CF_B}\otimes
\omega_X$.

\end{itemize}

In other words, we are looking at gauge equivalence classes with
respect to $\Hom(X,B)$ of connections of the form
\begin{equation} \label{nilp form}
\nabla^0 + \sum_{\imath} \phi_\imath \cdot f_\imath  + \bq,
\end{equation}
where $\phi_\imath$ are as in \eqref{form of nabla}, and $\bq$ is a
$\fb$-valued one-form on $X$ with a pole of order $1$ at $x$, whose
residue belongs to $\fn$.

This definition makes sense for $R$-families, so we obtain a functor
on the category of $\BC$-algebras, which we will denote by
$\on{Op}^\nilp_\fg(X)$.  For $X=\D$ we will denote the corresponding
functor simply by $\on{Op}^\nilp_\fg$.

As in the previous cases, one can define the functor
$\on{Op}^\nilp_\fg(X)^\fD$ on the category of $\fD_X$-algebras. Once
we prove its representability (see below), this functor will be
related to $\on{Op}^\nilp_\fg(X)$, $\on{Op}^\nilp_\fg$ and
$\on{Op}_\fg(\D^\times)$ in the same way as in the case of
$\on{Op}^\RS_\fg(X)^\fD$.

\ssec{}

We have an evident morphism of functors $\on{Op}_\fg^\nilp\to
\on{Op}_\fg(\D^\times)$.

\begin{thm}  \label{descr of nilp opers}
The above map is a closed embedding of functors, and an isomorphism 
onto $\on{Op}_\fg^{\RS,\varpi(-\crho)}$.
\end{thm}

Since the assertion is local, a similar statement holds for any pair
$(X,x)$. Before proving this theorem let us make the following
observation, which implies in particular that the map in question is
injective at the level of $\BC$-points.

\medskip

Let $(\CF_G,\nabla,\CF_B)$ be a $\BC$-point of
$\on{Op}_\fg(\D^\times)$, and let us first regard it as a $G$-local
system on $\D^\times$. Recall that if a local system $(\CF_G,\nabla)$
on $\D^\times$ admits an extension to a bundle on $\D$ with a
meromorphic connection with a pole of order $1$ {\it and nilpotent
residue}, then such extension is unique; we will refer to it as
Deligne's extension.
\footnote{Such an extension exists if and only if $(\CF_G,\nabla)$ has regular 
singularities, and when regarded analytically, has a unipotent monodromy.}

Thus, a necessary condition for $(\CF_G,\nabla,\CF_B)$ to come from
$\on{Op}_\fg^\nilp$ is that it admits such an extension. Since the flag
variety $G/B$ is compact, the $B$-bundle $\CF_B$, which is a 
priori defined on $\D^\times$, admits a unique extension to $\D$, 
compatible with the above extension of $\CF_G$. 

Having fixed this extension, our point comes from $\on{Op}_\fg^\nilp$ 
if and only if conditions (i) and (iii) from the definition of opers with
nilpotent singularities hold (condition (ii) is automatic from (i) and
the nilpotency assumption on the residue).

\ssec{Proof of \thmref{descr of nilp opers}}

To an oper with nilpotent singularities, written in the form 
\begin{equation} \label{exp nilp form}
\nabla^0 + \sum_{\imath} \phi_\imath(t)dt \cdot f_\imath +
\frac{\bq(t)}{t}dt,
\end{equation}
$\phi_\imath(t)\in (R[[t]])^\times$, $\bq(t)\in \fb\otimes R[[t]]$
with $\bq(0)\in \fn\otimes R$, we associate a point of
$\on{Op}_\fg^{\RS,\varpi(-\crho)}$ by applying the gauge
transformation by means of $t^{\crho}$.

The gauge action of $B[[t]](R)$ on connections of the form \eqref{exp
nilp form} gets transformed into the gauge action on connections with
regular singularities by means of $\on{Ad}_{t^{\crho}}(B[[t]])(R)$,
which is a subgroup of $B[[t]]$. This shows that the map
$\on{Op}_\fg^\nilp\to \on{Op}_\fg(\D^\times)$ factors through
$\on{Op}_\fg^{\RS,\varpi(-\crho)}$.

\medskip

To prove the theorem we must show that any connection written as
\begin{equation} \label{form RS nilp}
\nabla=\nabla^0 + \frac{p_{-1} -\crho + \bq(t)}{t}dt,
\end{equation} 
with $\bq(t)\in \fb[[t]]$ such that the image of $p_{-1} -\crho +
\bq(0)$ in $\fh\qu W$ equals $\varpi(-\crho)$ can be conjugated by means
of $N[[t]]$ into a connection of similar form such that $t^{-d}\cdot
\bq_d(t)\in \fb_d[[t]]$, uniquely up to $\on{Ad}_{t^{\crho}}(N[[t]])$.

\medskip

Note first of all that by applying a gauge transformation by means of
a {\it constant} loop into $N$, we can assume that $\bq(0)=0$.
By induction on $d$ we will prove the following statement:

{\it Every connection as in \eqref{form RS nilp} can be conjugated
by means of $N[[t]]$ to one which satisfies:
\begin{equation} \label{ind hyp}
t^{-d'}\cdot \bq_{d'}(t)\in \fb_{d'}[[t]] \text{ for } d'\leq d,
\end{equation} 
and $t^{-d}\cdot \bq_{d''}(t)\in \fb_{d''}[[t]]$ for $d''\geq d$.}

By the above, the statement holds for $d=1$. To perform the
induction step we will use a {\it descending} inductive argument. 
We assume that $\nabla$ satisfies \eqref{ind hyp}, and
that for some $k\geq d+1$,
$$t^{-d-1}\cdot \bq_{k'}(t)\in \fb_{k'}[[t]] \text{ for } k' \text{
satisfying } k'>k,$$ and
$$t^{-d}\cdot \bq_{k''}(t)\in \fb_{k''}[[t]] \text{ for } d+1\leq k''\leq k.$$

We will show how to modify $\nabla$ so that it continues to satisfy
\eqref{ind hyp}, and in addition:
$$t^{-d-1}\cdot \bq_{k'}(t)\in \fb_{k'}[[t]] \text{ for } k' 
\text{ satisfying } k'\geq k,$$
and
$$t^{-d}\cdot \bq_{k''}(t)\in \fb_{k''}[[t]] \text{ for } d+1\leq k'' < k.$$

Namely, we will replace $\nabla$ by
$$\nabla':=\on{Ad}_{\exp(t^d\cdot u_k)}(\nabla)= \nabla^0+
\frac{p_{-1} -\crho + \bq'(t)}{t}dt$$
for a certain element $u_k\in \fb_k$.

For any such $u_k$ the conditions, involving $\bq'_{k'}(t)$ for
$k'$ with either $k'<k$ or $k'>k$, hold automatically. The condition
on $\bq'_k(t)$ reads as follows:
\begin{equation} \label{key equation}
-d\cdot u_k+[u_k,-\crho]=-t^{-d}\bq_k(t)\on{mod} t.
\end{equation}
However, $[u_k,-\crho]=k\cdot u_k$, and since $k>d$ the above
condition is indeed solvable uniquely.

\medskip

This finishes the proof of the fact that any connection as in
\eqref{form RS nilp} can be conjugated by means of $N[[t]]$ to one
satisfying $t^{-d}\cdot \bq_d(t)\in \fb_d[[t]]$. The uniqueness of the
solution of \eqref{key equation} implies that the conjugation is
unique modulo $\on{Ad}_{t^{\crho}}(N[[t]])$. Thus, the proof of
\thmref{descr of nilp opers} is complete.  Let us note that the same
argument proves the following generalization:

\begin{prop}   \label{lambda opers}
Let $\cla$ be an anti-dominant coweight. Then a data of an $R$-point of
$\on{Op}_\fg^{\RS,\varpi(\lambda-\rho)}$ is equivalent to a data of
$B[[t]]$-conjugacy class of connections of the form
$$\nabla^0+\sum_{\imath} \phi_\imath(t)dt \cdot
f_\imath+\frac{\bq(t)}{t}dt,$$ where $\phi_\imath(t)$ are as in
\eqref{exp nilp form}, and $\bq(t)\in \fb\otimes R[[t]]$ is such that
$\bq(0)\on{mod}\fn=\cla$.
\end{prop}

\ssec{The secondary residue map}

Note that by definition the scheme $\on{Op}_\fg^\nilp\wh\times\D$
carries a universal oper with nilpotent singularities. Let us denote
by $\CP_{G,\on{Op}_\fg^\nilp}$ (resp., $\CP_{B,\on{Op}_\fg^\nilp}$)
the resulting $G$-bundle (resp., $B$-bundle) on
$\on{Op}_\fg^\nilp$. In particular, we obtain a map
$\on{Op}_\fg^\nilp\to \on{pt}/B$.

By taking the residue of the connection (see \secref{reg sing}), we
obtain a map from $\on{Op}_\fg^\nilp$ to the stack $\fn/B$, where $B$
acts on $\fn$ by means of the adjoint action; we will denote this map
by $\Res^\nilp$.

\begin{lem} \label{nilp res}  \hfill

\smallskip

\noindent{\em (1)} The map $\Res^\nilp$ is smooth. Moreover, the
$B$-scheme $\fn\underset{\fn/B}\times \on{Op}_\fg^\nilp$ can be
represented as a product of an infinite-dimensional affine space by a
finite-dimensional variety with a free action of $B$.

\smallskip

\noindent{\em (2)}
We have a natural identification:
$$\on{pt}/B\underset{\fn/B}\times \on{Op}_\fg^\nilp\simeq
\on{pt}/B\underset{\on{pt}/G}\times \on{Op}^\reg_G,$$ where
$\on{pt}/B\to \fn/B$ corresponds to $0\in \fn$, and the map
$\on{Op}_\fg^\reg:= \on{Op}_\fg(\D)\to \on{pt}/G$ is given by
$\CP_{G,\on{Op}_\fg(\D)}$.

\end{lem}

\begin{proof}

The second point of the lemma results from the definitions. To prove
the first point, note that $\fn\underset{\fn/B}\times
\on{Op}_\fg^\nilp$ identifies with the quotient of the space of
connections of the form \eqref{exp nilp form} by gauge transformations
by means of $B(\BC[[t]])$. As in the proof of \propref{free}, we
obtain that any such connection can be uniquely, up to the action of
$B(t\BC[[t]])$, brought into the form
$$\nabla^0+\bigl(\sum_{\imath} a_\imath \cdot
f_\imath+\frac{\bq'}{t}+\bq''+t\cdot \bv(t)\bigr)dt,$$ where $0\neq
a_\imath\in \BC$, $\bq'\in \fn$, $\bq''\in \fb$ and $\bv(t)\in
V_{\can}[[t]]$.  This scheme projects onto the variety of expressions
of the form
$$\sum_{\imath} a_\imath \cdot f_\imath+\frac{\bq'}{t}+\bq'',$$ on
which $B$ acts freely.

\end{proof}

\ssec{Opers with an integral residue}   \label{opers with integral residue}

For completeness, we shall now give a description of the scheme
$\on{Op}_\fg^{\RS,\varpi(-\cla-\crho)}$ when $\cla$ is an integral
coweight with $\cla+\crho$ dominant, similar to the one given by
\thmref{descr of nilp opers} in the case when $\cla=0$.

Let $\CJ$ be the subset of the set $\CI$ of vertices of the Dynkin
diagram, corresponding to those simple roots, for which $\langle
\alpha_{\jmath},\cla\rangle=-1$.  Let $\fp_\CJ\subset \fg$ be the
corresponding standard parabolic subalgebra, $\fn_\CJ\subset \fn$ its
unipotent radical, and $\fm_\CJ$ the Levi factor.

We introduce the notion of oper with $\cla$-nilpotent singularity to be
a triple $(\CF_G,\nabla,\CF_B)$ as in the definition of nilpotent
opers, where conditions (i)--(iii) are replaced by the following ones:

\begin{itemize}

\item(i) 
$(\nabla-\nabla')\on{mod} \fb_{\CF_B} \otimes \omega_X(x)$,
which is a section of $(\fg/\fb)_{\CF_B}\otimes \omega_X(x)$,
must be contained in 
${\bf O}_{\CF_B}\overset{\BG_m}\times \omega^\crho_X(-\cla\cdot x)$.

\item(ii) 
$\on{Res}(\nabla):=(\nabla-\nabla')\on{mod} \fg_{\CF_G}\otimes \omega_X$,
which is {\it a priori} an element of $\fg_{\CP_B}$, is contained in 
$\bigl(\fp_\CJ\bigr)_{\CP_B}$.

\item(iii) 
The image of $\on{Res}(\nabla)$ under 
$\bigl(\fp_\CJ\bigr)_{\CP_B}\to \bigl(\fm_\CJ\bigr)_{\CP_B}$ is nilpotent.

\end{itemize}

As in the case of $\cla=0$, this definition makes sense for
$R$-families, where $R$ is a $\BC$-algebra or a $\fD_X$-algebra. We
will denote by $\on{Op}_\fg^{\cla,\nilp}$ the resulting functor for
$X=\D$. Explicitly, an $R$-point of $\on{Op}_\fg^{\cla,\nilp}$ is a
$B[[t]](R)$-equivalence classes of connections of the form
\begin{equation} \label{form for lambda} 
\nabla^0 + \sum_{\imath} t^{\langle\alpha_\imath,\cla\rangle}
\cdot \phi_\imath(t)dt \cdot f_\imath  + \frac{\bq(t)}{t}dt,
\end{equation}
where $\phi_\imath(t)$ and $\bq(t)$ are as in \eqref{exp nilp form},
subject to the condition that the element
$$\underset{\jmath\in \CJ}\Sigma\,
\phi_{\jmath}(0)+\bq(0)\on{mod}\fn_\CJ\in \fm_\CJ$$
be nilpotent.

As in the case of $\cla=0$, there exists a natural map
of functors $\on{Op}_\fg^{\cla,\nilp}\to \on{Op}_\fg(\D^\times)$.

\begin{thm}
The above map is an isomorphism onto the subscheme
$\on{Op}_\fg^{\RS,\varpi(-\cla-\crho)}$.
\end{thm}

The proof of this theorem repeats that of \thmref{descr of nilp
opers}, where instead of the principal grading on $\fn$ we use the
one, defined by the adjoint action of $\cla+\crho$.

\medskip

Consider the subvariety of $\fp_\CI$, denoted ${\mathbf O}_\CJ$,
consisting of elements of the form
$$\underset{\jmath\in \CJ}\Sigma\, c_\jmath\cdot f_\jmath+\bq,\,\,
c_\jmath\neq 0,\,\, \bq\in \fb,$$ that are nilpotent. We have a
natural action of $B$ on ${\mathbf O}_\CJ$. 

As in the case of nilpotent opers, i.e, $\cla=0$, there exists a
natural smooth map
$$\Res^{\cla,\nilp}:\on{Op}_\fg^{\cla,\nilp}\to {\mathbf O}_\CJ/B,$$
obtained by taking the polar part of a connection as in 
\eqref{form for lambda}.

\medskip

Finally, let us consider the case when $\cla$ itself is dominant. In
this case $\CJ=\emptyset$, and ${\mathbf O}_\CJ=\fn$. Let us denote by
$\on{Op}_\fg^{\cla,\reg}$ the preimage of $\on{pt}/B\subset \fn/B$
under the map $\Res^{\cla,\nilp}$, where $\on{pt}\to \fn$ corresonds
to the point $0$.

The scheme $\on{Op}_\fg^{\cla,\reg}$ is the scheme of $\cla$-opers
introduced earlier by Beilinson and Drinfeld. As in the case of 
$\on{Op}_\fg^{\cla,\nilp}$, we have the notion of (an $R$-family) of
regular $\cla$-opers over any curve. By definition, this is a triple 
$(\CF_G,\nabla,\CF_B)$, where $\CF_G$ and $\nabla$ are a principal 
$G$-bundle and a connection on it, defined on the entire $X$,
and $\CF_B$ is a reduction of $\CF_G$ to $B$, such that $\nabla/\CF_B$,
as a section of $(\fg/\fb)_{\CF_B}\otimes \omega_X$, belongs to
${\bf O}_{\CF_B}\overset{\BG_m}\times \omega^\crho_X(-\cla\cdot x)$.

\section{Miura opers}    \label{Miura}

\ssec{}

Ler $R$ be a $\fD_X$-algebra. Let us fix once and for all another
Borel subgroup $B^-$ of $G$ which is in generic relative position with
$B$. The definition of Miura opers given below uses $B^-$. However,
the resulting scheme of Miura opers is defined canonically and is
independent of this choice.

Following \cite{F:wak}, Sect. 10.3, one defines a {\em Miura oper}
over $R$ to be a quadruple $$(\CF_G,\nabla,\CF_B,\CF_{B^-}),$$ where:

\begin{itemize}

\item
$(\CF_G,\nabla,\CF_B)$ are as in the definition of opers, i.e.,
$(\CF_G,\CF_B)$ is a $G$-bundle on $\Spec(R)$ with a reduction to $B$,
and $\nabla$ is a connection on $\CF_G$ along $X$ such that
$\nabla/\CF_B\in {\bf O}_{\CF_B,\omega_X}$.

\item $\CF_{B^-}$ is a reduction of $\CF_G$ to the opposite Borel
subgroup $B^-$ which is {\it preserved by the connection} $\nabla$.

\end{itemize}

We will denote the functor of Miura opers on the category of
$\fD_X$-algebras by $\on{MOp}_\fg(X)^\fD$, and the resulting functor
on the category of $\BC$-algebras by $\on{MOp}_\fg(X)$, i.e.,
$$\on{MOp}_\fg(X)(R):=\on{MOp}_\fg(X)^\fD(R\otimes \CO_X).$$

\begin{lem}   \label{Miura repr}
The functor $\on{MOp}_\fg(X)^\fD$ is representable by a $\fD_X$-scheme.
\end{lem}

\begin{proof}

Since the functor $\on{Op}_\fg(X)^\fD$ is known to be representable,
it suffices to show that the morphism $\on{MOp}_\fg(X)^\fD\to
\on{Op}_\fg(X)^\fD$ is representable as well.

Consider another functor on the category of $\fD_X$-algebras that
associates to $R$ the set of quadruples
$(\CF_G,\nabla,\CF_B,\CF_{B^-})$, {\it but without the condition that
$\CF_{B^-}$ be compatible with the connection}. The latter functor is
clearly representable over $\on{Op}_\fg(X)^\fD$, and it contains
$\on{MOp}_\fg(X)^\fD$ as a closed subfunctor.

\end{proof}

We will denote by $\on{MOp}_\fg(\D)$ or $\on{MOp}^\reg_\fg$ the
resulting scheme of Miura opers over $\D$. Note, however, that since
the flag variety $G/B^-$ is non-affine, the $\fD_X$-scheme
$\on{MOp}_\fg(X)$ is non-affine over $X$. Hence, the functor
$\on{MOp}_\fg(\D^\times)$ on the category of $\BC$-algebras is
ill-behaved; in particular, it cannot be represented by an ind-scheme.

\medskip

We define the $\fD_X$-schemes $\on{MOp}_\fg^{\ord_k}(X)^\fD$ (resp.,
$\on{MOp}_\fg^{\RS}(X)^\fD$, $\on{MOp}_\fg^{\nilp}(X)^\fD$) to
classify quadruples $(\CF_G,\nabla,\CF_B,\CF_{B^-})$, where the first
three pieces of data are as in the definition of
$\on{Op}_\fg^{\ord_k}(X)^\fD$ (resp., $\on{Op}_\fg^{\RS}(X)^\fD$,
$\on{Op}_\fg^{\nilp}(X)^\fD$), and $\CF_{\CB^-}$ is a reduction of the
$G$-bundle $\CF_G$, which is defined on the entire $X$, to the
subgroup $\CF_{B^-}$, compatible with the connection $\nabla$. The
last condition means, in the case of $\on{Op}_\fg^{\ord_k}(X)^\fD$,
that the operator $\nabla_{t^k\pa_t}$ preserves $\CF_{B^-}$. For the
other $\fD_X$-schemes this condition is defined similarly.

Each of these $\fD_X$-schemes is isomorphic to $\on{MOp}_\fg(X)^\fD$
over the curve $(X-x)$.

We will denote the corresponding schemes for $X=\D$ simply by
$\on{MOp}_\fg^{\ord_k}$, $\on{MOp}_\fg^{\RS}$ and
$\on{MOp}_\fg^{\nilp}$, respectively.

Note that there are {\it no} natural maps from $\on{MOp}_\fg^{\ord_k}$
to $\on{MOp}_\fg^{\ord_{k+1}}$ or from $\on{MOp}_\fg^{\nilp}$ to
$\on{MOp}_\fg^{\RS}$.

\ssec{Generic Miura opers}  \label{tr Miura}

Following \cite{F:wak}, Sect. 10.3, we shall say that a Miura oper
$(\CF_G,\nabla,\CF_B,\CF_{B^-})$ is {\em generic} if the given
reductions of $\CF_G$ to $B$ and $B^-$ are in generic relative
position. More precisely, observe that given a $G$-bundle on a scheme
$X$ with two reductions to $B$ and $B^-$, we obtain a morphism $X \to
B\backslash G/B^-$. The Miura oper is called generic if this morphism
takes values in the open part $B \cdot B^-$ of $B\backslash G/B^-$.

\begin{lem}  \label{gen trans}
Let $(\CF_G,\nabla,\CF_B)$ be a $\BC$-valued oper on $\D^\times$, and 
let $\CF_{B^-}$ be any horizontal reduction of $\CF_G$ to $B^-$. Then
it is in generic relative position with $\CF_B$.
\end{lem}

\begin{proof}

The following short argument is due to Drinfeld. The $G$-bundle
$\CF_G$ can be assumed to be trivial, and we can think of $\CF_B$ and
$\CF_{B^-}$ as two families of Borel subalgebras
$$\fb_1^{-}\subset \fg\ppart\supset \fb_2.$$
The connection on $\CF_G$ has the form $\nabla^0+{\mathbf q}(t)$,
where ${\mathbf q}(t)\in \fb_1^-$. 

Let $\fh'\subset \fg\ppart$ be any Cartan subalgebra, contained in both
$\fb_1^{-}$ and $\fb_2$. Let us decompose ${\mathbf q}(t)$ with
respect to the characters of $\fh'$, acting on $\fg\ppart$, i.e., with
respect to the roots. 

Then, on the one hand, each ${\mathbf q}(t)_\alpha$ belongs to $\fb_1^{-}$. 
I.e., if ${\mathbf q}(t)_\alpha\neq 0$ then $\alpha$ is positive
with respect to $\fb_1^-$.

On the other hand, if $\alpha_{\imath}$ be a simple root of $\fh'$
with respect to $\fb_2$, then by the oper condition, ${\mathbf
q}(t)_{-\alpha_\imath}\neq 0$.  Hence, every positive simple root with
respect to $\fb_2$ is negative with respect to $\fb_1^-$. This implies
that $\fb_1^-\cap \fb_2=\fh'$, i.e., the two reductions are
in generic relative position.

\end{proof}

Evidently, generic Miura opers form an open $\fD_X$-subscheme of
$\on{MOp}_\fg(X)^\fD$; we will denote it by
$\on{MOp}_{\fg,\gen}(X)^\fD$.

\begin{lem}  \label{tr Miura is affine}
The $\fD_X$-scheme $\on{MOp}_{\fg,\gen}(X)^\fD$ is affine over $X$.
\end{lem}

\begin{proof}

We know that the $\fD_X$-scheme $\on{Op}_\fg(X)^\fD$ is affine. Hence,
it is sufficient to show that $\on{MOp}_{\fg,\gen}(X)^\fD$ is affine
over $\on{Op}_\fg(X)$.

By definition, $\on{MOp}_{\fg,\gen}(X)^\fD$ is a closed subfunctor of
the functor that associates to a $\fD_X$-algebra $R$ the set of
quadruples $(\CF_G,\nabla,\CF_B,\CF_{B^-})$, where
$(\CF_G,\nabla,\CF_B)$ are as above, and $\CF_{B^-}$ is a reduction of
$\CF_G$ to $B^-$, which is in generic relative position with $\CF_B$,
and {\it not necessarily compatible with the connection}.

Since the big cell $B^-\cdot 1\subset G/B$ is affine, the latter
functor is evidently affine over $\on{Op}_\fg(X)^\fD$, implying our
assertion.

\end{proof}

We will denote by $\on{MOp}_{\fg,\gen}^{\ord_k}(X)^\fD$ (resp.,
$\on{MOp}_{\fg,\gen}^{\RS}(X)^\fD$, $\on{MOp}_{\fg,\gen}^{\nilp}(X)^\fD$)
the corresponding open $\fD_X$-subscheme of
$\on{MOp}_\fg^{\ord_k}(X)^\fD$ (resp., $\on{MOp}_\fg^{\RS}(X)^\fD$,
$\on{MOp}_\fg^{\nilp}(X)^\fD$). We will denote by
$$\on{MOp}_{\fg,\gen}^\reg:=\on{MOp}_{\fg,\gen}(\D),
\,\,\on{MOp}_{\fg,\gen}^{\nilp},\,\, \on{MOp}_{\fg,\gen}^{\RS},\,\,
\on{MOp}_{\fg,\gen}^{\ord_k}$$ the corresponding open subschemes of
$$\on{MOp}_{\fg}^\reg:=\on{MOp}_{\fg}(\D), \,\,
\on{MOp}_{\fg}^{\nilp}, \,\, \on{MOp}_{\fg}^{\RS},\,\,
\on{MOp}_{\fg}^{\ord_k},$$ respectively. By \lemref{tr Miura is
affine}, it makes sense also to consider the ind-scheme
$\on{MOp}_{\fg,\gen}(\D^\times)$.

\ssec{Miura opers and $H$-connections}   \label{Miura and H}

We will now establish a crucial result that connects generic Miura
opers with another, very explicit, $\fD_X$-scheme.

Consider the $H$-bundle $\omega_X^\crho$, and let $\ConnX^\fD$ be
the $\fD_X$-scheme of connections on it, i.e., it associates to a
$\fD_X$-algebra $R$ the set of connections on the pull-back of
$\omega_X^\crho$ to $\Spec(R)$ along $X$. This is a principal
homogeneous space with respect to the group $\fD_X$-scheme that
associates to $R$ the set of $\fh$-valued sections of the pull-back of
$\omega_X$ to $\Spec(R)$. In particular, $\ConnX^\fD$ is affine over
$X$.

We will denote the resulting functor on $\BC$-algebras by $\ConnX$.
For $X=\D$ (resp., $\D^\times$) this functor is evidently
representable by a scheme (resp., ind-scheme), which we will denote by
$\ConnD=:\ConnD^\reg$ (resp., $\ConnDt$). This scheme
(resp., ind-scheme) is a principal homogeneous space with respect to
$\fh\otimes \omega_\D$ (resp., $\fh\otimes \omega_{\D^\times}$).

\medskip

Note we have a natural map of $\fD_X$-schemes
\begin{equation} \label{Miura and conn}
\on{MOp}_{\fg,\gen}(X)^\fD\to \ConnX^\fD.
\end{equation}

Indeed, given an $R$-point of $\on{MOp}_{\fg,\gen}(X)^\fD$, let $\CF'_H$
be the $H$-bundle with connection, induced by means of
$\CF_{B^-}$. However, the assumption that the Miura oper is generic
implies that $\CF'_H\simeq \CF'_B \cap \CF_{B^-} \simeq \CF_H$, where
$\CF_H$ is the $H$-bundle induced by $\CF_B$. Now, let us recall
that by \lemref{FH}, we have a canonical isomorphism
$\CF_H\simeq \omega_X^\crho$.

\begin{prop}[\cite{F:wak}, Prop. 10.4]    \label{Miura isom conn}
The map \eqref{Miura and conn} is an isomorphism.
\end{prop}

\begin{proof}

We construct the inverse map $\ConnX^\fD\to
\on{MOp}_{\fg,\gen}(X)^\fD$ as follows.  Recall first of all that a data
of a $G$-bundle with two reductions to $B$ and $B^-$ in generic
position is equivalent to a data of an $H$-bundle. Thus, from
$\CF_H:=\omega_X^\crho$ we obtain the data $(\CF_G,\CF_B,\CF_{B^-})$
from the Miura oper quadruple.

A connection on (the pull-back of) $\omega_X^\crho$ (to some
$\fD_X$-scheme) induces a connection, that we will call $\nabla_H$ on
$\CF_G$, compatible with both reductions. We produce the desired
connection $\nabla$ on $\CF_G$ by adding to $\nabla_H$ the
$\fn^-_{\CF_H}\otimes \omega_X$-valued 1-form equal to
$\underset{\imath}\Sigma\, \phi_\imath$, where each $\phi_\imath$ is
the tautological trivialization of
$(\fn^-_{-\alpha_\imath})_{\CF_H}\otimes \omega_X$.

The resulting connection preserves $\CF_{B^-}$ and satisfies the oper
condition with respect to $\CF_B$. Hence
$(\CF_G,\nabla,\CF_B,\CF_{B^-})$ is a generic Miura oper.  Clearly,
two maps $$\on{MOp}_{\fg,\gen}(X)^\fD\rightleftarrows \ConnX^\fD$$
are mutually inverse.

\end{proof}

This proposition immediately implies the isomorphisms of $\fD_X$-schemes
$$\ConnD^\reg\simeq \on{MOp}_{\fg,\gen}^\reg \text{ and }
\ConnDt\to \on{MOp}_{\fg,\gen}(\D^\times).$$

Let us denote by $\bigl(\ConnX^{\ord_k}\bigr)^\fD$ (resp.,
$\bigl(\ConnX^\RS\bigr)^\fD$) the $\fD_X$-scheme of meromorphic
connections on $\omega_X^\crho$ over $\D$ with pole of order $\leq k$
(resp., $\leq 1$). Each of these $\fD_X$-schemes is isomorphic to
$\ConnX^\fD$ over $(X-x)$. We will denote by $\ConnD^\RS$ and
$\ConnD^{\ord_k}$ the resulting schemes of connections on $\D$.

Using the fact that connections on $\omega_X^{\crho}$ with a pole of
order $k$, $k\geq 1$, are in a canonical bijection with those on
$\omega_X^{\crho}(\cla\cdot x)$ for any coweight $\cla$, from the
above proposition we obtain also the isomorphisms
$$\bigl(\ConnX^{\ord_k}\bigr)^\fD\simeq \on{MOp}_{\fg,\gen}^{\ord_k}(X)^\fD
\text{ and } \bigl(\ConnX^\RS\bigr)^\fD \simeq
\on{MOp}_{\fg,\gen}^{\RS}(X)^\fD,$$ implying that
\begin{equation} \label{conn and Miura RS}
\ConnD^{\ord_k}\simeq \on{MOp}_{\fg,\gen}^{\ord_k} \text{ and }
\ConnD^\RS\simeq \on{MOp}_{\fg,\gen}^\RS.
\end{equation}

\medskip

We call the composed map of $\fD_X$-schemes
\begin{equation} \label{Miura transformation}
\ConnX^\fD\to \on{MOp}_{\fg,\gen}(X)^\fD\to \on{Op}_\fg(X)^\fD
\end{equation}
the {\em Miura transformation} and denote it by $\on{MT}$. By a
slight abuse of notation, we will denote by the same symbol $\on{MT}$
the corresponding maps
$$\ConnD^{\RS}\to \on{Op}_\fg^{\RS},\,\, \ConnD^{\ord_k}\to
\on{Op}_\fg^{\ord_k},\,\, \ConnDt\to
\on{Op}_\fg(\D^\times).$$

\ssec{An application: proof of \propref{char RS}}

Let $(\CF_G,\nabla,\CF_B)$ be an oper on $\D^\times$ such that the
$G$-bundle with connection $(\CF_G,\nabla)$ has regular singularities,
i.e., $\CF_G$ can be extended to a $G$-bundle $\CF'_G$ on $\D$, so
that $\nabla$ has a pole of order $\leq 1$.

Then it is known that $(\CF_G,\nabla)$ admits at least one horizontal
reduction to $B^-$, call it $\CF_{B^-}$. By \lemref{gen trans}, the
quadruple $(\CF_G,\nabla,\CF_B,\CF_{B^-})$ is a generic Miura oper on
$\D^\times$.

By the compactness of $G/B^-$, the above reduction extends uniquely on
the entire $\D$.  The connection on the resulting $B^-$-bundle
$\CF'_{B^-}$ has a pole of order $\leq 1$.  Hence, the connection on
the $H$-bundle $\CF'_H$, induced from $\CF'_{B^-}$, also has a pole of
order $\leq 1$.

Therefore, the above point $(\CF_G,\nabla,\CF_B,\CF_{B^-})\in
\on{MOp}_{\fg,\gen}(\D^\times)$, viewed as a point of
$\ConnDt$, belongs to $\ConnD^\RS$. Hence, the triple
$(\CF_G,\nabla,\CF_B)$, being the image of the above point under the
map $\on{MT}$, belongs to $\on{Op}_\fg^\RS$.

\ssec{Miura opers with regular singularities}

Consider the map $\ConnD^\RS\to \fh$, that assigns to a connection with
a pole of order $1$ its residue; we will denote it by $\Res^\fh$. For
$\cla\in \fh$ we will denote by $\ConnD^{\RS,\cla}$ the preimage of
$\cla$ under $\Res^\fh$.

A coweight $\cla$ such that $\langle \alpha_\imath,\cla\rangle\notin
\BZ^{<0}$ (resp., $\notin \BZ^{>0}$) for $\alpha\in \Delta^+$, will be
called {\em dominant} (resp., {\em anti-dominant}).

\begin{prop}  \label{MT and residue} \hfill

\smallskip

\noindent{\em (1)}
We have a commutative diagram:
$$
\CD
\ConnD^\RS @>{\on{MT}}>> \on{Op}_\fg^\RS \\
@V{\Res^\fh}VV  @V{\Res^\RS}VV  \\
\fh @>>> \fh\qu W,
\endCD
$$
where the bottom arrow is $\cla\mapsto \varpi(\cla-\crho)$.

\smallskip

\noindent{\em (2)} If $\cla$ is dominant with respect to $B$,
then the map $\on{MT}:\ConnD^{\RS,\cla+\crho}\to
\on{Op}_\fg^{\RS,\varpi(\cla)}$ is an isomorphism.

\end{prop}

The rest of this subsection is devoted to the proof of this
proposition. Part (1) follows from the construction:

Given a generic Miura oper with regular singularities
$(\CF_G,\nabla,\CF_B,\CF_{B^-})$, the induced $H$-bundle $\CF_H$ is
$\omega_X^\crho(\crho\cdot x)$, by \lemref{FHk}. The polar part of
$\nabla$ is a section $\bq\in \fb^-_{\CP_{B^-}}$. Let $\cla$ denote
the projection of $\bq$ onto
$\fb^-_{\CP_{B^-}}/\fn^-_{\CP_{B^-}}\simeq \fh$, which equals the
polar part of the connection on $\CF_H$.

Then $\Res^\RS(\CF_G,\nabla,\CF_B)$ equals the projection of $\bq$
under $\fg/B^-\to \fg/G\to \fh\qu W$, and hence it equals
$\varpi(\cla)$. It remains to notice that the resulting connection on
$\omega_X^\crho$ has the polar part equal to $\cla+\crho$.

\medskip

To prove part (2), we will use the following general assertion.

\begin{lem} \label{reduction to B}
Let $(\CF_G,\nabla)$ be an $R$-family of $G$-connections on $\D$ with
a pole of order $1$, and let $\CP_G$ be the fiber of $\CF_G$ at the
closed point of the disc. Let $\CP_{B^-}$ be a reduction of $\CP_G$ to
$G$ with the property that the residue $\bq$ of $\nabla$, which is
{\it a priori} an element of $\fg_{\CP_G}$, belongs to
$\fb^-_{\CP_{B^-}}$.  Assume that the projection of $\bq$ to
$\fb^-_{\CF_{B^-}}/\fn^-_{\CF_{B^-}}\simeq \fh$ is constant and
anti-dominant with respect to $B^-$.

Then there exists a unique $B^-$-reduction $\CF_{B^-}$ of $\CF_G$,
which is compatible with $\nabla$ and whose fiber at $x$ equals $\CP_{B^-}$.
\end{lem}

Let us first show how this lemma implies the proposition.  Consider
the subvariety $$(p_{-1}+\fb)^{\cla}\subset (p_{-1}+\fb)\subset \fg,$$
consisting of elements whose image in $\fh\qu W$ equals $\cla$. This
is the $N$-orbit of the element $p_{-1}+\cla$. We claim that each
point of this orbit is contained in a {\em unique} Borel subalgebra of
$\fg$ that is in generic relative position with $\fb$.

More precisely, consider the {\em Grothendieck alteration} $\wt\fg \to \fg$ defined 
as the subvariety of $\fg\times G/B^-$ consisting of the pairs
\begin{equation}    \label{Groth}
\wt\fg = \{ \bq \in\fg,\fb'{}^-\in G/B^-\,|\,\bq \in \fb'{}^- \}.
\end{equation}

Let $\widetilde{(p_{-1}+\fb)}{}^{\cla}$ be the scheme-theoretic
intersection of the preimages of $(p_{-1}+\fb)^{\cla}\subset \fg$ 
and the the big cell $B\cdot 1\subset G/B^-$ in $\wt\fg$.

\begin{lem} \label{alt sec}
The projection $\widetilde{(p_{-1}+\fb)}{}^{\cla}\to (p_{-1}+\fb)^{\cla}$
is an isomorphism.
\end{lem}

\begin{proof}

The inverse map $(p_{-1}+\fb)^{\cla}\to \widetilde{(p_{-1}+\fb)}{}^{\cla}$
is obtained by conjugating the element
$\{p_{-1}+\cla,\fb^-\}$ by means of $N$.

\end{proof}

Let us denote by $\CF_{G,\on{Op}_\fg^{\RS,\cla}}$ (resp.,
$\CP_{G,\on{Op}_\fg^{\RS,\cla}}$) the universal $G$-bundle with
connection on $\on{Op}_\fg^{\RS,\cla}\widehat\times\D$ (resp., its
restriction to $\on{Op}_\fg^{\RS,\cla}\times x\subset
\on{Op}_\fg^{\RS,\cla}\widehat\times\D$). Let
$\CF_{B,\on{Op}_\fg^{\RS,\cla}}$ and $\CP_{B,\on{Op}_\fg^{\RS,\cla}}$
be their reductions to $B$ given by the oper structure.

{}From the above lemma we obtain that the $G$-bundle
$\CP_{G,\on{Op}_\fg^{\RS,\cla}}$ over $\on{Op}_\fg^{\RS,\cla}$ admits
a unique reduction to $B^-$ such that the polar part of the connection
belongs to $\fb^-_{\CF_{B^-}}$ and its image in $\fh$ equals
$\cla$. Moreover, the resulting $B^-$-bundle
$\CP_{B^-,\on{Op}_\fg^{\RS,\cla}}$ is in generic relative position
with $\CP_{B,\on{Op}_\fg^{\RS,\cla}}$.

\medskip

Note that if $\cla$ is dominant with respect to $B$, then it is
anti-dominant with respect to $B^-$. Hence by our assumption on $\cla$
and \lemref{reduction to B}, the $G$-bundle
$\CF_{G,\on{Op}_\fg^{\RS,\cla}}$ on $\on{Op}_\fg^{\RS,\cla}\wh\times
\D$ admits a unique horizontal reduction to $B^-$. This reduction is
automatically in generic position with
$\CF_{B,\on{Op}_\fg^{\RS,\cla}}$, because this is so over the closed
point $x\in \D$. Thus, we have constructed the inverse map
$$\on{Op}_\fg^{\RS,\varpi(\cla)}\to \ConnD^{\RS,\cla+\crho}.$$

This map is evidently a left inverse of the map $\on{MT}$. The
uniqueness assertion of \lemref{reduction to B}, combined with
\lemref{alt sec}, implies that it is also a right inverse.
This completes the proof of part (2) of \propref{MT and residue}.

Let us now prove \lemref{reduction to B}.

\begin{proof} (Drinfeld)

With no restriction of generality, we can assume that our $G$-bundle
$\CF_G$ is trivial, and the connection has the form
$\nabla=\nabla^0+\frac{\bq(t)}{t}$, where $\bq(t)\in \fg[[t]]$ and
$\bq(0)\in \fb^-$. We must show that there exists an element $\bg\in
\on{ker}(G[[t]]\to G)$, unique modulo $B^-$, such that
$$\on{Ad}_{\bg}(\nabla^0+\frac{\bq(t)}{t}) =
:\nabla'=\nabla^0+\frac{\bq'(t)}{t}$$ is such that $\bq'(t)\in
\fb^-[[t]]$.

Assume by induction that $\bq(t)\on{mod} t^k\in \fb^-[t]/t^k$. We must
show that there exists an element $u\in \fg$, unique modulo $\fb^-$,
so that
$$t\cdot \left(\on{Ad}_{\exp(t^k\cdot u)}(\frac{\bq(t)}{t})-k\cdot
t^{k-1}\cdot u\right) \on{mod} t^{k+1} \in \fb^-[t]/t^{k+1}.$$ This
can be rewritten as
$$k\cdot u + [\bq_0,u]=\bq_k.$$
However, this equation is indeed solvable uniquely in $\fg/\fb^-$,
since by assumption, negative integers are not among the eigenvalues
of the adjoint action of $\bq_0$ on $\fg/\fb^-$.

\end{proof}

We shall now describe the behavior of the map $\on{MT}$,
restricted to $\ConnD^{\RS,\cmu}$, for $\cmu$ anti-dominant and
integral. This is the case which is in a sense opposite to the one
considered in \propref{MT and residue}(2).

\begin{prop} \label{regular Miura}
Let $\cla$ be a dominant integral weight. Then the image of the map
$\on{MT}|_{\ConnD^{\RS,-\cla}}$ belongs to the closed subscheme
$\on{Op}_\fg^{\cla,\reg}\subset \on{Op}_\fg^{\cla,\nilp}\simeq
\on{Op}_\fg^{\RS,\varpi(-\cla-\crho)}$.  Moreover, we have a Cartesian
square:
\begin{equation} \label{diag reg Miura}
\CD
\ConnD^{\RS,-\cla}@ >>> (\overset{\circ}{B^-\backslash G})/B  \\
@V{\on{MT}}VV    @VVV  \\
\on{Op}_\fg^{\cla,\reg} @>{\Res^{\lambda,\nilp}}>> \on{pt}/B,
\endCD
\end{equation}
where $\overset{\circ}{B^-\backslash G}$ denotes the open $B$-orbit in the
flag variety $B^-\backslash G$.
\end{prop}

\begin{proof}

Choosing a coordinate on $\D$, and thus trivializing $\omega_X$, a
point of $\ConnD^{\RS,-\cla}$ can be thought of as a connection on the
trivial bundle of the form
$$\nabla^0+\frac{\bq(t)}{t}dt$$ with $\bq(t)\in \fh[[t]]$ and
$\bq(0)=-\cla$.  The oper, corresponding to the Miura transformation
of the above connection, equals
$$\nabla^0+p_{-1}dt+\frac{\bq(t)}{t}dt.$$ Conjugating this connection
by means of $t^{-\cla}$ we obtain a connection of the form \eqref{form
for lambda}. Let us denote by $(\CF_G,\nabla,\CF_B)$ the resulting
point of $\on{Op}_\fg^{\cla,\reg}$.

Note that the horizontal generic reduction to $B^-$ of $\CF_G$, which
was defined over $\D^\times$, extends to one over $\D$.  Indeed, under
the above trivialization of $\CF_G$, the reduction to $B$ corresponds
to the subgroup $B$ itself, and the reduction to $B^-$ corresponds to
$B^-$, which are manifestly in the generic position.
This defines the upper horizontal map in \eqref{diag reg Miura}. 

To show that this diagram is indeed Cartesian, it suffices to show that
given a ($R$-) point $(\CF_G,\nabla,\CF_B)$ of
$\on{Op}_\fg^{\cla,\reg}$, any reduction to $B^-$ of the fiber $\CP_G$
at $x$ of $\CF_G$, which is in the generic position with respect
$\CP_B$ (the latter being the fiber of $\CF_B$ at $x$), comes from a
unique reduction of $\CF_G$ to $B^-$. However, this immediately
follows from \lemref{reduction to B}, since $\cla$ was assumed
dominant with respect to $B$, and, hence, anti-dominant with respect
to $B^-$.

\end{proof}

\ssec{Miura opers with nilpotent singularities} \label{Miura with nilp sing}

Let us observe that we have four geometric objects that may be called
"Miura opers with nilpotent singularities"
$$\on{MOp}_\fg^\nilp,\,\,
\on{MOp}_\fg^\RS\underset{\on{Op}_\fg^\RS}\times
\on{Op}_\fg^\nilp,\,\,
\on{MOp}_{\fg,\gen}^\RS\underset{\on{Op}_\fg^\RS}\times
\on{Op}_\fg^\nilp,\text{ and }
\on{MOp}_{\fg,\gen}(\D^\times)\underset{\on{Op}_\fg(\D^\times)}\times
\on{Op}_\fg^\nilp.$$ The first three of the above objects are schemes,
and the fourth is an ind-scheme.  In this section we will study the
relationship between them.

First, we have the following:

\begin{lem} \label{naive bijection}
The sets of $\BC$-points of the four objects above are in a natural
bijection.
\end{lem}

\begin{proof}

In all the four cases the set in question classifies the data of an
oper with a nilpotent singularity on $\D$, and its horizontal
reduction to $B^-$ over $\D^\times$ (which is necessarily generic by
\lemref{gen trans}).

\end{proof}

We will establish the following:

\begin{thm}  \label{four versions of Miura}
There exist natural maps
$$
\CD \on{MOp}_{\fg,\gen}^\RS\underset{\on{Op}_\fg^\RS}\times
\on{Op}_\fg^\nilp @>{{\bf 1}}>>
\on{MOp}_{\fg,\gen}(\D^\times)\underset{\on{Op}_\fg(\D^\times)} \times
\on{Op}_\fg^\nilp @>{{\bf 2}}>> \on{MOp}_\fg^\nilp \\ @V{{\bf 3}}VV \\
\on{MOp}_\fg^\RS\underset{\on{Op}_\fg^\RS}\times \on{Op}_\fg^\nilp,
\endCD
$$
which commute with the projection to $\on{Op}_\fg^\nilp$, and which on
the level of $\BC$-points induce the bijection of \lemref{naive
bijection}.  Moreover, the map ${\bf 1}$ is a closed embedding, the
map ${\bf 2}$ is formally smooth, and the map ${\bf 3}$ is an
isomorphism.
\end{thm}

The rest of this section is devoted to the proof of this
theorem. Note, however, that the existence of the map ${\bf 1}$ and
the fact that it is a closed embedding is immediate from the fact that
$\on{MOp}_{\fg,\gen}^\RS\to \on{MOp}_{\fg,\gen}(\D^\times)$ is a
closed embedding.

Also, the map ${\bf 3}$ comes from the tautological map
$\on{MOp}_{\fg,\gen}^\RS\to \on{MOp}_{\fg}^\RS$. Since the latter is
an open embedding, the map ${\bf 3}$ is one too. Since it induces a
bijection on the set of $\BC$-points by \lemref{naive bijection}, we
obtain that it is an isomorphism.

\medskip

To construct the other maps appearing in \thmref{four versions of
Miura}, we need to describe the corresponding schemes more
explicitly. First, by \eqref{conn and Miura RS}, we have an
isomorphism:
$$\on{MOp}_{\fg,\gen}^\RS\underset{\on{Op}_\fg^\RS}\times
\on{Op}_\fg^\nilp\simeq \ConnD^\RS\underset{\on{Op}_\fg^\RS}\times
\on{Op}_\fg^\nilp,$$ and the latter identifies, 
by \thmref{descr of nilp opers}, with
$$\ConnD^\RS\underset{\fh}\times (\fh\underset{\fh\qu W}\times \on{pt}),$$
where $\on{pt}\to \fh\qu W$ corresponds to the point $\varpi(-\crho)$. Hence,
by \propref{MT and residue}, we obtain an isomorphism
$$\on{MOp}_{\fg,\gen}^\RS\underset{\on{Op}_\fg^\RS}\times
\on{Op}_\fg^\nilp\simeq \underset{w\in W}\cup\,
\ConnD^{\RS,\crho-w(\crho)}.$$

Since the map ${\bf 1}$ in \thmref{four versions of Miura} is a closed
embedding and an isomorphism at the level of $\BC$-points, the
ind-scheme
$\on{MOp}_{\fg,\gen}(\D^\times)\underset{\on{Op}_\fg(\D^\times)}
\times \on{Op}_\fg^\nilp$ also splits into connected components,
numbered by elements of $W$; we will denote by
$\on{MOp}_{\fg,\gen}(\D^\times)^w$ the component corresponding to a
given $w\in W$.

\medskip

Let now $\wt\fg$ be the Grothendieck alteration of $\fg$ defined in
\eqref{Groth}.  Let $\wt\fn$ be the scheme-theoretic preimage of
$\fn\subset \fg$ under the forgetful map $\wt\fg\to \fg$; this is a
scheme acted on by $B$.  Note that $\wt\fn$ is connected and
non-reduced.

\medskip

By \lemref{reduction to B} we have:

\begin{cor} \label{nilp Miura opers as cart}
There exists a canonical isomorphism
$$\on{MOp}_\fg^\nilp\simeq \on{Op}_\fg^\nilp\underset{\fn/B}\times
\wt\fn/B.$$
\end{cor} 

Let now $\wt\fn^w$ be the subvariety of $\wt\fn$, obtained by
requiring that the pair $(\bq\in\fg,\fb'{}^-)\in \wt\fg$ be such that
the Borel subalgebra $\fb'{}^-$ is in position $w$ with respect to
$\fb$, i.e., the corresponding point of $G/B^-$ belongs to the
$B$-orbit $B \cdot w^{-1} \cdot B^-$. This is a reduced scheme
isomorphic to the affine space of dimension $\dim(\fn)$. Let us denote
by $\wt\fn^{w,\on{th}}$ the formal neighborhood of $\wt\fn^w$ in
$\wt\fn$, regarded as an ind-scheme.  Clearly, the action of $B$ on
$\wt\fn$ preserves each $\wt\fn^w$, and
$$\wt\fn(\BC)\simeq \underset{w\in W}\cup\, \wt\fn^w(\BC).$$

Let us denote by $\on{MOp}_\fg^{\nilp,w}$ the subscheme of
$\on{MOp}_\fg^\nilp$ equal to $\on{Op}_\fg^\nilp\underset{\fn/B}\times
\wt\fn^w/B$ in terms of the isomorphism of \corref{nilp Miura opers as
cart}. Let us denote by $\on{MOp}_\fg^{\nilp,w,\on{th}}$ the
ind-scheme $\on{Op}_\fg^\nilp\underset{\fn/B}\times
\wt\fn^{w,\on{th}}/B$.

\begin{thm}   \label{two versions of Miura}
For every $w\in W$ there exists an isomorphism
$$\on{MOp}_{\fg,\gen}(\D^\times)^w\simeq \on{MOp}_\fg^{\nilp,w,\on{th}},$$
compatible with the forgetful map to $\on{Op}_\fg^\nilp$ and
the bijection of \lemref{naive bijection}.
\end{thm}

Clearly, \thmref{two versions of Miura} implies the remaining assertions
of \thmref{four versions of Miura}. In addition, by passing to reduced
schemes underlying the isomorphism of \thmref{two versions of Miura},
and using \lemref{nilp res}(1), we obtain the following:

\begin{cor}
There exists a canonical isomorphism
$\ConnD^{\RS,\crho-w(\crho)}\simeq \on{MOp}_\fg^{\nilp,w}$.
\end{cor}

\ssec{Proof of \thmref{two versions of Miura}}

We begin by constructing the map
\begin{equation} \label{from Conn to Miura}
\on{MOp}_{\fg,\gen}(\D^\times)^w\to \on{MOp}_\fg^\nilp
\end{equation}

Given an $R$-point $(\CF_G,\nabla,\CF_B,\CF_{B^-})$ of
$\on{MOp}_{\fg,\gen}(\D^\times)^w$, let $(\CF'_G,\CF'_B)$ be an
extension of the pair $(\CF_G,\CF_B)$ onto $\D$ such that the
resulting triple $(\CF'_G,\nabla,\CF'_B)$ is a point of
$\on{Op}_\fg^\nilp$.  Such an extension exists, according to
\thmref{descr of nilp opers}. We claim that the reduction to $B^-$ of
$\CF_G$, given by $\CF_{B^-}$, gives rise to a reduction of $\CF'_G$
to $B^-$:

Let us think of a reduction to $B^-$ in the Pl\"ucker picture (see
\cite{FGV}). Let $V^\la$ is the irreducible representation
representation of $\fg$ with highest weight $\la$. Then our point of
$\on{MOp}_{\fg,\gen}(\D^\times)^w$ gives rise to a system of
meromorphic maps
$$V^\lambda_{\CF'_G}\to \omega_X^{\langle \lambda,\crho\rangle},$$ for
dominant weights $\la$, compatible with the (meromorphic) connections
on the two sides. Note that the
connection on $\omega_\D^\crho$, corresponding to
$(\CF_G,\nabla,\CF_B,\CF_{B^-})$, restricted to the subscheme
$$\Spec(R)\underset{\on{MOp}_{\fg,\gen}(\D^\times)^w}\times
\ConnD^{\RS,\crho-w(\crho)}$$ has the property that its pole is of
order $1$ and the residue equals $\crho-w(\crho)$.  We apply the
following:

\begin{lem}
Let $(\CF_H,\nabla_H)$ be an $R$-family of $H$-bundles with
meromorphic connections on $\D$. Assume that there exists a quotient
$R\twoheadrightarrow R'$ by a nilpotent ideal such that connection on
the resulting $R'$-family has a pole of order $1$ and a fixed residue
integral $\cla\in \fh$. Then there exists a unique modification
$\CF'_H$ of $\CF_H$ at $x$ such that the resulting connection on
$\CF'_H$ is regular.
\end{lem}

The lemma produces an $R$-family of $H$-bundles $\CF'_H$ with a
regular connection, and a horizontal system of {\it a priori
meromorphic} maps
$$\fs^\la:V^\lambda_{\CF'_G}\to \BC^\lambda_{\CF'_H},$$
satisfying the Pl\"ucker equations. We claim that each of these maps 
$\fs^\la$ is in fact regular and surjective. This is particular case of the 
following lemma:

\begin{lem}
Let $\CV$ and $\CL$ be $R$-families of vector bundles and a line
bundles on $\D$, respectively, both equipped with connections such
that on $\CV$ it has a pole of order $1$ and nilpotent residue, and on
$\CL$ the connection is regular. Let $\CV\to \CL$ be a non-zero
meromorphic map, compatible with connections. Then this map is regular
and surjective.
\end{lem}

Thus, we obtain a horizontal reduction $\CF'_{B^-}$ of
$\CF'_G$ to $B^-$, and the desired map in \eqref{from Conn to Miura}.

\medskip

Consider the restriction of the map \eqref{from Conn to Miura} to
$\ConnD^{\RS,\crho-w(\crho)}\subset \on{MOp}_{\fg,\gen}(\D^\times)^w$.
Since the former scheme is reduced and irreducible, the image of this
map is contained in $\on{MOp}_\fg^{\nilp,w'}$ for some $w'\in W$. This
implies that the map \eqref{from Conn to Miura} itself factors through
$\on{MOp}_\fg^{\nilp,w',\on{th}}$ for the same $w'$.

We have to show that $w'=w$ and that the resulting map is an 
isomorphism. We claim that for that purpose it is sufficient to
construct a map in the opposite direction
\begin{equation} \label{from Miura to Conn}
\on{MOp}_\fg^{\nilp,w}\to \ConnD^{\RS,\crho-w(\crho)},
\end{equation}
compatible with the identification of \lemref{naive bijection}. This
follows from the next observation:

\begin{lem}
Let $(\CF_G,\nabla,\CF_B)$ be an $R$-point of $\on{Op}_\fg^\nilp$, and
let $R'$ be a quotient of $R$ by a nilpotent ideal. Let
$(\CF'_G,\nabla',\CF'_B,\CF'_{B^-})$ be a lift of the induced
$R'$-family to a point of $\on{MOp}_{\fg,\gen}(\D^\times)$.  Then the
sets of extensions of this lift to $R$-points of
$\on{MOp}_{\fg,\gen}(\D^\times)$ and $\on{MOp}_\fg^\nilp$ are in
bijection.
\end{lem}

The lemma follows from the fact a deformation over a nilpotent base
of a generic Miura oper remains generic.

\medskip

Given an $R$-point of $\on{MOp}_\fg^{\nilp,w}$ and a dominant weight
$\lambda$, consider the diagram
\begin{equation} 
\omega_X^{\langle \lambda,\crho\rangle}\overset{\fs'}\to 
V^\lambda_{\CF_G}\overset{\fs}\to \CL,
\end{equation}
where the map $\fs'$ corresponds to the reduction of $\CF_G$ to $B$,
and $\CL$ is {\it some} line bundle on $\Spec(R[[t]])$ with a regular
connection $\nabla_\CL$ in the $t$-direction, and the map $\fs$ is a
surjective bundle map, compatible with connections, corresponding to
the reduction of $\CF_G$ to $B^-$. We will denote by $\nabla(\pa_t)$
(resp., $\nabla_\CL(\pa_t)$) the action of the vector field $\pa_t$ on
$\D$ on sections of $V^\lambda_{\CF_G}$ (resp., $\CL$), given by the
connection.

To construct the map as in \eqref{from Miura to Conn}, we have to show
that the composition $\fs\circ \fs'$ has a zero of order $\langle
\lambda, \crho-w(\crho)\rangle$. This is equivalent to the following:
let $\bv$ be a non-vanishing section of $\omega_X^{\langle
\lambda,\crho\rangle}$, thought of as a section of $V^\lambda_{\CF_G}$
by means of $\fs'$. We need to show that the section
$\nabla_\CL(\pa_t)^{n'}\left(\fs(\bv)\right)$ of $\CL$ is regular and
non-vanishing for $n'=n:=\langle \lambda, \crho-w(\crho)\rangle$, and
has a zero at $x$ if $n'<n$.  Since the map $\fs$ is compatible with
connections, we have to calculate
$\fs\left(\nabla(\pa_t)^n(\bv)\right)$.

\medskip

Let $F^j(V^\lambda)$ be the increasing $B$-stable filtration on
$V^\la$, defined by the condition that a vector $v\in V^\lambda$ of
weight $\lambda'$ belongs to $F^j(V^\lambda)$ if and only if $\langle
\lambda-\lambda',\crho\rangle\leq j$. Let $F^j(V^\lambda_{\CF_B})$ be
the corresponding induced filtration on the vector bundle
$V^\lambda_{\CF_B}\simeq V^\lambda_{\CF_G}$. Each successive quotient
$F^j(V^\lambda_{\CF_B})/F^{j-1}(V^\lambda_{\CF_B})$ is isomorphic to
$$\underset{\la', \langle \lambda-\lambda',\crho\rangle= j}\oplus\, 
F^j(V^\lambda)/F^{j-1}(V^\lambda)\otimes \omega_\D^{\langle \lambda,\crho\rangle}.$$

By the condition on $\nabla$,
\begin{equation} \label{pole of nabla}
\nabla(\pa_t)(F^j(V^\lambda_{\CF_B}))\subset
F^{j-1}(V^\lambda_{\CF_B})(x)+ F^{j+1}(V^\lambda_{\CF_B}),
\end{equation}
and the induced map
\begin{equation} \label{Grif}
\nabla(\pa_t):F^j(V^\lambda_{\CF_B})/F^{j-1}(V^\lambda_{\CF_B})\to
F^{j+1}(V^\lambda_{\CF_B})/F^j(V^\lambda_{\CF_B})
\end{equation}
comes from the map $F^j(V^\lambda)/F^{j-1}(V^\lambda)\to
F^{j+1}(V^\lambda)/F^j(V^\lambda)$, given by $p_{-1}$. (The latter
makes sense, since the vector field $\pa_t$ trivializes the line
bundle $\omega_\D$.) 

Let us denote by $n''$ the maximal integer such that the composition
$\fs\circ \fs'$ vanishes to the order $n''$ along $\Spec(R)\times
x\subset \Spec(R[[t]])$.  By induction on $j$, from \eqref{pole of
nabla}, we obtain that the map
$$\fs:F^j(V^\lambda_{\CF_B})\to \CL\to \CL_x$$ vanishes to the order
$n''-j$, where $\CL_x$ is the restriction of $\CL$ to $\Spec(R)\times
x$.

Assume first that $n''<n$. Then, by the maximality assumption on
$n''$, the image of $\fs(\nabla(\pa_t)^{n''}(\bv))$ in $\CL_x$ is
non-zero. However, this is impossible since the composition
$$F^j(V^\lambda)\hookrightarrow V^\lambda\to (V^\lambda)_{\fb'{}^-}$$
vanishes for any $\fb'{}^-\in G/B^-$ in relative position $w$ with
respect to $B$ and $j<\langle \lambda-w(\lambda),\crho\rangle=n$.

\medskip

Thus, $\fs\left(\nabla(\pa_t)^n(\bv)\right)$ is regular, and it
remains to show that its image in $\CL_x$ is nowhere
vanishing. However, this follows from \eqref{Grif}, since for a
highest weight vector $v\in V^\lambda$ and $n$ and $\fb'{}^-$ as
above, the image of $p_{-1}^n(v)$ in $(V^\lambda)_{\fb'{}^-}$ is
non-zero.

\section{Groupoids and Lie algebroids associated to opers}
\label{algebroid on opers}

\ssec{The isomonodromy groupoid} \label{identification of algebroids}

Let us recall that a {\em groupoid} over a scheme $S$ is a scheme
$\CG$ equipped with morphisms $l: \CG \to S, r: \CG \to S$, $m: \CG
\underset{r,S,l}\times \CG \to \CG$, an involution $\gamma: \CG \to
\CG$, and a morphism $u: S \to \CG$ that satisfy the following
conditions:

\begin{itemize}

\item associativity: $m \circ (m \times \on{id}) = m \circ (\on{id}
  \times m)$ as morphisms $\CG \underset{r,S,l}\times \CG
  \underset{r,S,l}\times \CG \to \CG$;

\item unit: $r \circ u = l \circ u = \on{id}_S$.

\item inverse: $l \circ \gamma = r, r \circ \gamma = l$, $m\circ
(\gamma\times \on{id}_\CG)=u\circ r$, $m\circ (\on{id}_\CG\times
\gamma)=u\circ l$.

\end{itemize}

If $S_1\subset S$ is a subscheme, we will denote by $\CG|_{S_1}$ the
restriction of $\CG$ to $S_1$, i.e., the subscheme of $\CG$ equal to
$(l\times r)^{-1}(S_1\times S_1)$.  This is a groupoid over $S_1$.

\medskip

The normal sheaf to $S$ inside $\CG$ acquires a structure of {\em Lie
algebroid}; we will denote it by $\fG$, and by $\anch$ the anchor map
$\fG \to T(S)$, where $T(S)$ is the tangent algebroid of $S$. (We
refer to \cite{Ma} for more details on groupoids and Lie algebroids).

The notion of groupoid generalizes in a straightforward way to the
case when both $S$ and $\CG$ are ind-schemes. However, to speak about
a Lie algebroid attached to a Lie groupoid, we will need to assume
that $\CG$ is formally smooth over $S$ (with respect to either, or
equivalently, both projections). In this case $\fG$ will be a Tate
vector bundle over $S$; we refer the reader to \secref{topological
algebroids} for details.

\medskip

We now define the {\em isomonodromy groupoid} $\Isom$ over the
ind-scheme $\on{Op}_\fg(\D^\times)$. Points of the ind-scheme $\Isom$
over an algebra $R$ are triples $(\chi,\chi',\phi)$, where 
$\chi=(\CF_G,\nabla,\CF_B)$ and $\chi'=(\CF'_G,\nabla',\CF'_B)$ 
are both $R$-points of $\on{Op}_\fg(\D^\times)$, and 
$\phi$ is an isomorphism of $G$-bundles with connections
$(\CF_G,\nabla) \simeq (\CF'_G,\nabla')$.

Explicitly, if $\chi$ and $\chi'$ are connections $\nabla$ and $\nabla'$,
respectively, on the trivial bundle, both of the form
\begin{equation} \label{again nabla}
\nabla^0+p_{-1}dt+\phi(t)dt,\,\, \phi(t)\in \fb\otimes R\ppart,
\end{equation}
then a point of $\Isom(R)$ over $(\chi,\chi)$ is an element 
$\bg\in G\left(R\ppart\right)$ such that $\on{Ad}_{\bg}(\nabla)=\nabla'$.
Two triples $(\chi_1,\chi'_1,\bg_1)$ and $(\chi_2,\chi'_2,\bg_2)$
are equivalent of there exist elements $\bg,\bg'\in N\left(R\ppart\right)$,
such that $\nabla_1=\on{Ad}_{\bg}(\nabla_2)$,
$\nabla'_2=\on{Ad}_{\bg'}(\nabla'_1)$ and $\bg_2=\bg'\cdot \bg_1\cdot \bg$.

\medskip

The morphisms $l$ and
$r$ send $(\chi,\chi',\phi)$ to $\chi$ and $\chi'$, respectively. The
morphism $m$ sends the pair $(\chi,\chi',\phi),(\chi',\chi'',\phi')$
to $(\chi,\chi'',\phi' \circ \phi)$, the morphisms $\gamma$ sends
$(\chi,\chi',\phi)$ to $(\chi',\chi,\phi^{-1})$ and the morphism $u:
\on{Op}_\fg(\D^\times) \to \Isom$ sends $\chi$ to $(\chi,\chi,\on{id})$.

\medskip

We call $\Isom$ the isomonodromy groupoid for the following
reason. In the analytic context two connections on the trivial bundle 
on a punctured disc are called isomonodromic if they have the same
monodromy and the Stokes data (in case of irregular singularity). 
In the case of connections on the formal punctured disc the appropriate 
analogue of the notion of isomonodromy is the notion of gauge 
equivalence of connections.

\begin{prop}   \label{Isom is smooth}
The groupoid $\Isom$ is formally smooth over $\on{Op}_\fg(D^\times)$.
\end{prop}

\ssec{Description of tangent space and proof of 
\propref{Isom is smooth}}  \label{descr of tang}

Let $R'\to R$ be a homomorphism of rings such that its kernel $\bI$
satisfies $\bI^2=0$. Let $\chi'=(\CF'_G,\nabla',\CF'_B)$ be an
$R'$-point of $\on{Op}_\fg(\D^\times)$, and let
$\chi=(\CF_G,\nabla,\CF_B)$ be the corresponding $R$-point. Let $\bg$
be an automorphism of $\CF_G$ such that the quadruple
$(\CF_G,\nabla,\CF_B,\bg)$ is an $R$-point of $\Isom$ over $\chi$. We
need to show that it can be lifted to an $R'$-point
$(\CF'_G,\nabla',\CF'_B,\bg')$ of $\Isom$.

Since the ind- scheme $G\ppart$ is formally smooth, we can always find
some automorphism $\bg'_1$ of $\CF'_G$, lifting $\bg$. To show the
existence of the required lift we must find an element $\bu\in
\fg_{\CF_G}\underset{R\ppart}\otimes \bI\ppart$ such that the point
$\bg'=\bg'_1\cdot (1+\bu)$ satisfies
$$\on{Ad}_{\bg'}(\nabla')-\nabla'\in
\fb_{\CF'_B}\underset{R'\ppart}\otimes R'\ppart dt.$$ By assumption,
$\bq:=\on{Ad}_{\bg'_1}(\nabla')-\nabla'$ belongs to the subspace
$$\fb_{\CF'_B}\underset{R'\ppart}\otimes \bI\ppart dt\simeq 
\fb_{\CF_B}\underset{R\ppart}\otimes \bI\ppart dt.$$
Therefore, the desired element $\bu$ must satisfy:
$$\nabla(\bu)=\bq\in (\fg/\fb)_{\CF_B}\underset{R\ppart}\otimes
\bI\ppart dt.$$

Hence, it is sufficient to show that the map
$$\fg_{\CF_B}\underset{R\ppart}\otimes
\bI\ppart \overset{\nabla}\to \fg_{\CF_B}\underset{R\ppart}\otimes
\bI\ppart dt\to (\fg/\fb)_{\CF_B}\underset{R\ppart}\otimes \bI\ppart
dt$$ is surjective. But this follows from the oper condition on
$\nabla/\CF_B$. \footnote{The above description makes it explicit that
both $\Omega^1(\on{Op}_\fg(\D^\times))$ and the conormal
$N^*_{\on{Op}_\fg(\D^\times)/\Isom}$ are Tate vector bundles on
$\on{Op}_\fg(\D^\times)$, i.e., we do not have to use the general
\thmref{Tate families} to prove this fact.}

\medskip

Thus, \propref{Isom is smooth} is proved. In particular, the Lie
algebroid $\isom$, corresponding to the groupoid $\Isom$, is
well-defined. Let us write down an explicit expression for $\isom$ and
for the anchor map.

Since the ind-scheme $\on{Op}_\fg(\D^\times)$ is reasonable and
formally smooth, its tangent $T(\on{Op}_\fg(\D^\times))$ is a Tate
vector bundle. For a $R$-point $(\CF_G,\nabla,\CF_B)$ of
$\on{Op}_\fg(\D^\times)$ we have:
\begin{equation} \label{tangent space to opers}
T(\on{Op}_\fg(\D^\times))|_{\Spec(R)}\simeq 
\on{coker}(\nabla):\fn_{\CF_B}\to
\fb_{\CF_B}\underset{\BC\ppart}\otimes \BC\ppart dt,
\end{equation}
and 
\begin{equation} \label{exp algebroid}
\isom|_{\Spec(R)}\simeq
\on{ker}(\nabla):(\fg/\fn)_{\CF_G}\to
(\fg/\fb)_{\CF_B}\underset{\BC\ppart}\otimes \BC\ppart dt.
\end{equation}

The anchor map $\anch:\isom\to T(\on{Op}_\fg(\D^\times))$ acts as
follows:
$$\bu\in \fg_{\CF_G}\mapsto \nabla(\bu)\in
\fb_{\CF_B}\underset{\BC\ppart}\otimes \BC\ppart dt.$$

\medskip

Consider the cotangent sheaf $\Omega^1(\on{Op}_\fg(\D^\times))$; this
is also a Tate vector bundle on $\on{Op}_\fg(\D^\times)$.  From
\eqref{tangent space to opers}, we obtain that once we identify $\fg$
with its dual by means of any invariant form $\kappa:\fg\otimes\fg\to
\BC$, we obtain an isomorphism:
\begin{equation}  \label{pre-ident}
\Omega^1(\on{Op}_\fg(\D^\times))\simeq \isom.
\end{equation}

As we shall see in the next subsection, a choice of $\kappa$ defines a
Poisson structure on $\on{Op}_\fg(\D^\times)$, and in particular makes
$\Omega^1(\on{Op}_\fg(\D^\times))$ into a Lie algebroid. We will show
that the above identification of bundles is compatible with the Lie
algebroid structure.

\medskip

\begin{remark}
In the analytic context this Poisson structure
is used to define the KdV flow on $\on{Op}_\fg(\D^\times)$ as the
system of evolution equations corresponding to a certain
Poisson-commuting system functions on the space of opers. The
isomorphism with $\isom$ implies in particular that the KdV flows
preserve gauge equivalence classes.
\end{remark}

\ssec{The Drinfeld-Sokolov reduction and Poisson structure}

Consider the space of all connections on the trivial $G$-bundle on
$\D^\times$, i.e., the the space $\on{Conn}_G(\D^\times)$ of operators
of the form
\begin{equation} \label{just connections}
\nabla^0+\phi(t), \phi(t)\in \fg\otimes \omega_{\D^\times}.
\end{equation}

This is an ind-scheme, acted on by the group $G\ppart$ by gauge
transformations.  We can consider the natural isomonodromy ind-groupoid
over $\on{Conn}_G(\D^\times)$:
$${\mathsf {Isom}}_{\on{Conn}_G(\D^\times)}:=\{\bg,\nabla,\nabla'\,|\,
\on{Ad}_{\fg}(\nabla)=\nabla'\},\,\, l(\bg,\nabla,\nabla')=\nabla,
r(\bg,\nabla,\nabla')=\nabla'.$$ Since ${\mathsf
{Isom}}_{\on{Conn}_G(\D^\times)}\simeq \on{Conn}_G(\D^\times)\times
G\ppart$, it is formally smooth over $\on{Conn}_G(\D^\times)$.

\medskip

Let us choose a symmetric invariant form $\kappa:\fg\otimes \fg\to
\BC$, and let $\ghat_\kappa$ be the corresponding Kac-Moody extension
of $\fg\ppart$.  Using the form $\kappa$, we can identify the space
$\on{Conn}_G(\D^\times)$ with hyperplane in $\ghat_\kappa^*$ equal to
the preimage of $1\in \BC$ under the natural map $\ghat_\kappa\to
\BC$. It is well-known that under this identification the coadjoint
action of $G\ppart$ on $\ghat_\kappa^*$ corresponds to the gauge
action of $G\ppart$ on $\on{Conn}_G(\D^\times)$.

The space $\ghat_\kappa^*$ carries a canonical Poisson structure,
which induces a Poisson structure also on
$\on{Conn}_G(\D^\times)$.

\begin{lem}
We have a canonical isomorphism of Lie algebroids
\begin{equation} \label{two algebroids upstairs}
\Omega^1(\on{Conn}_G(\D^\times))\simeq {\mathsf
{isom}}_{\on{Conn}_G(\D^\times)},
\end{equation}
where ${\mathsf {isom}}_{\on{Conn}_G(\D^\times)}$ is the Lie algebroid
of ${\mathsf {Isom}}_{\on{Conn}_G(\D^\times)}$. 
\end{lem}

\begin{proof}

We claim that (global sections of) both the LHS and the RHS identify
with
$$\fg\ppart\shriektimes \Fun(\on{Conn}_G(\D^\times))$$
with the natural bracket (we refer to \secref{topological vector spaces}, 
where the notation $\shriektimes$ is introduced).

The assertion concerning ${\mathsf {isom}}_{\on{Conn}_G(\D^\times)}$ follows
from the fact that ${\mathsf {Isom}}_{\on{Conn}_G(\D^\times)}$ is the
product of $\on{Conn}_G(\D^\times)$ and the group $G\ppart$ acting on
it, and $\fg\ppart$ is the Lie algebra of $G\ppart$.

The assertion concerning $\Omega^1(\on{Conn}_G(\D^\times))$ 
follows from the identification of $\on{Conn}_G(\D^\times)$ with
a hyperplane in $\hg_\kappa^*$, and the description of the
Poisson structure on the dual space to a Lie algebra.

\end{proof}

For any group ind-subscheme $K\subset G\ppart$ such that
$\ghat_\kappa$ is split over $\sk\subset \fg\ppart$, the map
$\ghat_\kappa^*\to \sk^*$ is a moment map for the action of $K$ on
$\ghat_\kappa^*$, and, in particular, on $\on{Conn}_G(\D^\times)$.

We take $K=N\ppart)$, and we obtain a moment map 
$$\mu:\on{Conn}_G(\D^\times)\to \left(\fn\ppart\right)^*\simeq
\fg/\fb\otimes \omega_{\D^\times},$$ where we identify $\fn^*\simeq
\fg/\fb$ using $\kappa$.

We have an identification
\begin{equation} \label{DS}
\on{Op}_\fg(\D^\times)\simeq \left(\mu^{-1}(p_{-1}dt)\right)/N(\hCK),
\end{equation}
where the action of $N(\hCK)$ on $\mu^{-1}(p_{-1}dt)$ is free. It is
in this fashion that $\on{Op}_\fg(\D^\times)$ was originally
introduced in \cite{DS} and this is why this Hamiltonian reduction is
called the Drinfeld-Sokolov reduction.

\begin{lem}   \label{two algebroids}
There exists a canonical isomorphism of Lie algebroids over
$\on{Op}_\fg(\D^\times)$
$$\Omega^1(\on{Op}_\fg(\D^\times))\simeq \isom.$$
\end{lem}

\begin{proof}

Note that the action of $N\ppart$ on $\on{Conn}_G(\D^\times)$ lifts
naturally to an action of the group $N\ppart\times N\ppart$ on
${\mathsf {Isom}}_{\on{Conn}_G(\D^\times)}$.  We have a canonical
identification of $\Isom$ with the two-sided quotient of ${\mathsf
{Isom}}_{\on{Conn}_G(\D^\times)}$:
\begin{equation} \label{DS for groupoid}
\Isom\simeq 
\biggl(\left((\mu\times \mu)\circ (l\times r)\right)^{-1}((p_{-1}\cdot
dt)\times (p_{-1}\cdot dt))\biggr) /N\ppart\times N\ppart.
\end{equation}

Hence, $\isom$ is obtained as a reduction with respect to $N\ppart$ of
the Lie algebroid ${\mathsf {isom}}_{\on{Conn}_G(\D^\times)}$. By the
definition of the Poisson structure on $\on{Op}_\fg(\D^\times)$, the
Lie algebroid $\Omega^1(\on{Op}_\fg(\D^\times))$ is the reduction of
the Lie algebroid $\Omega^1(\on{Conn}_G(\D^\times))$ on
$\on{Conn}_G(\D^\times)$.

Hence, the assertion of the lemma follows from \eqref{two algebroids
upstairs}.

\end{proof}

\ssec{The groupoid and Lie algebroid over regular opers}
\label{identification of algebroids, reg} \hfill

Let $S$ be an ind-scheme with a Poisson structure and $S_1\subset S$
be a reasonable subscheme, which is co-isotropic, i.e., the ideal
$\bI=\on{ker}\left(\Fun(S)\to \Fun(S_1)\right)$ satisfies
$[\bI,\bI]\subset \bI$. We will assume that both $S$ and $S_1$ are
formally smooth; we will also assume that the normal bundle
$N_{S_1/S}$ (which by our assumption is discrete), is locally
projective. \footnote{We do not know whether this follows directly
from the formal smoothness assumption.}

In this case the conormal $N^*_{S_1/S}$ acquires a structure of Lie
algebroid, and the sheaf $\Omega^1(S_1)$ is a module over
it. Moreover, we have the following commutative diagram:
\begin{equation}   \label{abstract coisotrop}
\CD
0 @>>> N^*_{S_1/S} @>>> \Omega^1(S)|_{S_1} @>>> \Omega^1(S_1) @>>> 0 \\
& & @VVV @V{\anch}VV @VVV & & \\
0 @>>> T(S_1) @>>> T(S)|_{S_1} @>>> N_{S_1/S} @>>> 0
\endCD
\end{equation}
such that the right vertical arrow is a map of modules over $N^*_{S_1/S}$.

We claim:

\begin{lem}  \label{regular coisotrop}
The subscheme $\on{Op}_\fg^\reg\subset \on{Op}_\fg(\D^\times)$ is
co-isotropic.
\end{lem}

\begin{proof}

Consider the subscheme $\on{Conn}_G^\reg$ of $\on{Conn}_G(\D^\times)$
obtained by imposing the condition that $\phi(t)$ belongs to
$\fg\otimes \omega_\hCO$.  It is co-isotropic, since the corresponding
ideal in $\Fun(\on{Conn}_G(\D^\times))$ is generated by $\fg\otimes
\omega_\hCO\subset \ghat_\kappa$, which is a subalgebra.

By \secref{canon}, the scheme $\on{Op}_\fg^\reg$ can be realized as
$$\left(\mu^{-1}(p_{-1}\cdot dt)\cap \on{Conn}_G^\reg\right)/N[[t]],$$
which implies the assertion of the lemma.

\end{proof}

\medskip

Let $\Isom^\reg$ be groupoid over the scheme
$\on{Op}_\fg^\reg=\on{Op}_\fg(D)$ whose $R$-points are triples
$(\chi,\chi',\phi)$, where $\chi=(\CF_G,\nabla,\CF_B)$ and
$\chi=(\CF_G,\nabla,\CF_B)$ are $R$-points of $\on{Op}_\fg^\reg$ and
$\phi$ is an isomorphism of $R$-families of $G$-bundles on $\D$ with
connections $(\CF_G,\nabla)\to (\CF'_G,\nabla')$.

Recall now the principal $G$-bundle $\CP_{G,\on{Op}_\fg^\reg}$ over
$\on{Op}_\fg^\reg$ obtained by restriction to $\on{Op}_\fg^\reg\times
x$ from the tautological $G$-bundle $\CF_G$ on
$\on{Op}_\fg^\reg\hat\times\D$.  This $G$-bundle defines a map
$\on{Op}_\fg^\reg\to \on{pt}/G$.

\begin{lem}  \label{ident of groupoid over regular}  \hfill

\smallskip

\noindent{\em (1)} The natural map
$\Isom^\reg\to \Isom|_{\on{Op}_\fg^\reg}$ is an isomorphism.

\smallskip

\noindent{\em (2)} The groupoid $\Isom^\reg$ is naturally isomorphic
to
$$\on{Op}_\fg^\reg \underset{\on{pt}/G}\times \on{Op}_\fg^\reg.$$
\end{lem}

\begin{proof}

The assertion of the lemma amounts to the following. Let $S=\Spec(R)$
be an affine scheme and let $(\CF_G,\nabla)$, $(\CF'_G,\nabla')$ be
two $G$-bundles on $\Spec(R[[t]])$ with a regular connection along
$t$. Let $\CP_G$, $\CP'_G$ be their restrictions to $\Spec(R)$,
respectively. Then the set of connection-preserving isomorphisms
$\CF_G\to \CF'_G$ maps isomorphically to both the set of
connection-preserving isomorphisms $\CF_G|_{\D^\times}\to
\CF'_G|_{\D^\times}$ and the set of isomorphisms $\CP_G\to \CP'_G$.

\end{proof}

Let $\isom^\reg$ be the Lie algebroid of $\Isom^\reg$.  \lemref{ident
of groupoid over regular}(2) implies that $\isom^\reg$ is identifies
with the Atiyah algebroid $\on{At}(\CPreg)$ of infinitesimal
symmetries of the $G$-bundle $\CPreg$. Therefore it fits in the exact
sequence
$$
0 \to \fg_{\on{Op}^\reg_G} \to \isom^\reg \to T (\on{Op}_\fg^\reg) \to 0,
$$
where $\fg_{\on{Op}^\reg_G}:=\fg_{\CP_{\on{Op}^\reg_G}}$. In what
follows we will denote by $\fb_{\on{Op}^\reg_G}$ (resp.,
$\fn_{\on{Op}^\reg_G}$) the subbundle of $\fg_{\on{Op}^\reg_G}$,
corresponding to the reduction $\CP_{B,\on{Op}^\reg_G}$ of
$\CP_{G,\on{Op}^\reg_G}$ to $B$.

\medskip

Note that by \lemref{two algebroids} and \lemref{ident of groupoid
over regular}(1) we have a natural map of algebroids
\begin{equation} \label{map of alg on reg}
N^*_{\on{Op}_\fg^\reg/\on{Op}_\fg(\D^\times)}\to \isom^\reg.
\end{equation}

Following \cite{BD}, Sect. 3.7.16, we have:

\begin{prop}     \label{Atiyah alg}
The map of \eqref{map of alg on reg} is an isomorphism.
\end{prop}

\begin{proof}

The assertion of the proposition amounts to the fact that the map
\begin{equation}  \label{map of normals}
\Omega^1(\on{Op}_\fg^\reg)\to N_{\on{Op}_\fg^\reg/\on{Op}_\fg(\D^\times)}
\end{equation}
from \eqref{abstract coisotrop} is an injective bundle map.

\medskip

Since the scheme $\on{Op}_\fg^\reg$ is smooth, for an $R$-point
$(\CF_G,\nabla,\CF_B)$ of $\on{Op}_\fg^\reg$, the restrictions of
$T(\on{Op}_\fg^\reg)$ and $\Omega^1(\on{Op}_\fg^\reg)$ to $\Spec(R)$
can be canonically identified with
$$\on{coker}(\nabla):\fn_{\CF_B}\to
\fb_{\CF_B}\underset{\BC[[t]]}\otimes \BC[[t]]dt$$ and
$$\on{ker}(\nabla):(\fg/\fn)_{\CF_B}\underset{\BC[[t]]}\otimes
(\BC\ppart/\BC[[t]])\to (\fg/\fb)_{\CF_B}\underset{\BC[[t]]}\otimes
(\BC\ppart dt/\BC[[t]]dt),$$ respectively, where we have used the
identification $\fg^*\simeq \fg$ given by $\kappa$.

Hence, the the restriction of
$N_{\on{Op}_\fg^\reg/\on{Op}_\fg(\D^\times)}$ to $\Spec(R)$ can be
identified with
$$\on{coker}(\nabla):\fn_{\CF_B}\underset{\BC[[t]]}\otimes
(\BC\ppart/\BC[[t]])\to \fb_{\CF_B} \underset{\BC[[t]]}\otimes
(\BC\ppart dt/\BC[[t]]dt),$$ and the map of \eqref{map of normals} is
given by
$$\bu\in \fg_{\CF_B}\underset{\BC[[t]]}\otimes
(\BC\ppart/\BC[[t]])\mapsto \nabla(\bu)\in
\fb_{\CF_B}\underset{\BC[[t]]}\otimes (\BC\ppart dt/\BC[[t]]dt).$$ The
injectivity of the map in question is now evident from the oper
condition on $\nabla/\CF_B$.

\end{proof}

\begin{cor}   \label{kernel and cokernel of anchor}
The kernel and the cokernel of the anchor map 
$$\anch:\Omega^1(\on{Op}_\fg(\D^\times))|_{\on{Op}_\fg^\reg}\to
 T(\on{Op}_\fg(\D^\times))|_{\on{Op}_\fg^\reg}$$ are both isomorphic
 to $\fg_{\on{Op}^\reg_G}$ as
 $N^*_{\on{Op}_\fg^\reg/\on{Op}_\fg(\D^\times)}$-modules.
\end{cor}

\begin{proof}
The isomorphism concerning the kernel follows by combining
\propref{Atiyah alg} and \lemref{ident of groupoid over regular}. The
isomorphism concerning the cokernel follows from the first one by a
general D-scheme argument, see \cite{CHA}, Sect. 2.5.22.

Let us, however, reprove both isomorphisms directly. We have:
$$\on{ker}(\anch|_{\on{Op}_\fg^\reg})\simeq
\on{ker}(\nabla):\fg_{\CF_G}\underset{\BC[[t[[}\otimes \BC\ppart\to
\fg_{\CF_G}\underset{\BC[[t]]}\otimes \BC\ppart dt,$$ which is easily
seen to identify with $\fg_{\on{Op}^\reg_G}$.

The assertion concerning $\on{coker}(\anch|_{\on{Op}_\fg^\reg})$
follows by Serre's duality on $\D^\times$. Indeed, the dual of
$\Omega^1(\on{Op}_\fg(\D^\times))$ is canonically isomorphic to
$T(\on{Op}_\fg(\D^\times))$, and under this isomorphism, the dual of
the map $\anch$ goes to itself. Hence,
$$\left(\on{coker}(\anch)\right)^*\simeq \on{ker}(\anch)\simeq
\fg_{\on{Op}^\reg_G},$$ which we identify with
$\fg^*_{\on{Op}^\reg_G}$ using the form $\kappa$.

\end{proof}

To summarize, we obtain  the following commutative diagram:
$$
\CD & & 0 & & 0 & & & & \\ & & @VVV @VVV & & \\ 0 @>>>
\fg_{\on{Op}_\fg^\reg} @>{\on{id}}>> \fg_{\on{Op}_\fg^\reg} @>>> 0 & &
\\ & & @VVV @VVV @VVV & \\ 0 @>>>
N^*_{\on{Op}_\fg^\reg/\on{Op}_\fg(\D^\times)} @>>>
\Omega^1(\on{Op}_\fg(\D^\times)) @>>> \Omega^1(\on{Op}^\reg_\fg) @>>>
0 \\ & & @VVV @VVV @VVV & \\ 0 @>>> T (\on{Op}_\fg^\reg) @>>> T
(\on{Op}_\fg(\D^\times))|_{\on{Op}_\fg^\reg} @>>>
N_{\on{Op}_\fg^\reg/\on{Op}_\fg(\D^\times)} @>>> 0 \\ & & @VVV @VVV
@VVV & \\ & & 0 @>>> \fg_{\on{Op}_\fg^\reg} @>{\on{id}}>>
\fg_{\on{Op}_\fg^\reg}. & & \\ & & & & @VVV @VVV \\ & & & & 0 & & 0 &
& \endCD
$$

We will conclude this subsection by the following remark. Let
$(\CF_G,\nabla,\CF_B)$ be an $R$-point of $\on{Op}_\fg^\reg$, and let
$\bu$ be an element of $\fb_{\CF_B}\underset{\BC[[t]]}\otimes
\BC\ppart dt$, giving rise to a section of
$T(\on{Op}_\fg(\D^\times))|_{\Spec(R)}$ by \eqref{tangent space to
opers}.

{}From the proof of \corref{kernel and cokernel of anchor}, we obtain
the following:

\begin{lem} \label{exp normal}
The image of $\bu$ in
$N_{\on{Op}_\fg^\reg/\on{Op}_\fg(\D^\times)}/\isom^\reg\simeq
\fg_{\on{Op}_\fg^\reg}$ equals the image of $\bu$ under the
composition
$$\fb_{\CF_B}\underset{\BC[[t]]}\otimes \BC\ppart dt\to
\fg_{\CF_G}\underset{\BC[[t]]}\otimes \BC\ppart dt\to
H^0_{DR}(\D^\times,\fg_{\CF_G})\simeq \fg_{\on{Op}_\fg^\reg}.$$
\end{lem}

\ssec{The groupoid and algebroid on opers with nilpotent singularities}

Consider now the subscheme $\on{Op}_\fg^\nilp\subset
\on{Op}_\fg(\D^\times)$.  As in \lemref{regular coisotrop}, it is easy
to see that $\on{Op}_\fg^\nilp$ is co-isotropic, since
$$\on{Op}_\fg^\nilp\simeq \left(\mu^{-1}(p_{-1}dt)\cap
\on{Conn}_G^\nilp\right)/N[[t]],$$ where $\on{Conn}_G^\nilp$ is the
subscheme of $\on{Conn}_G(\D^\times)$, consisting of connections as in
\eqref{just connections}, for which $\phi(t)\in \fg[[t]]+\fn\otimes
t^{-1}\BC[[t]]$, and the latter is the orthogonal complement to the
Iwahori subalgebra in $\ghat_\kappa$.

Let us consider the groupoid
$$\Isom^\nilp:=\on{Op}_\fg^\nilp\underset{\fn/B}\times \on{Op}_\fg^\nilp$$
over $\on{Op}_\fg^\nilp$, and let $\isom^\nilp$ be the corresponding Lie
algebroid. 

\begin{lem}  \label{map of groupoids, nilp}
There exists a natural closed embedding $\Isom^\nilp\to
\Isom|_{\on{Op}_\fg^\nilp}$.
\end{lem}

\begin{proof}

The lemma is proved in the following general framework. Let
$(\CF_G,\nabla)$ and $(\CF'_G,\nabla')$ be two $R$-families of bundles
with connections on $\D$ with poles of order $1$ and nilpotent
residues. Let $\CP_G$ and $\CP'_G$ be the resulting $G$-bundles on
$\Spec(R)$, and $\Res(\nabla)$ (resp., $\Res(\nabla')$) be the
residue, which is an element in $\fg_{\CP_G}$ (resp., $\fg_{\CP'_G}$).

Then there is a bijection between the set of connection-preserving
isomorphisms $\CF_G\to \CF'_G$ of bundles on $\Spec(R[[t]])$ and
isomorphisms $\CP_G\to \CP'_G$, which map $\Res(\nabla)$ to
$\Res(\nabla')$.

\end{proof}

Note, however, that unlike the case of regular opers, the map of
\lemref{map of groupoids, nilp} is {\it not} an isomorphism. Indeed,
the restriction of $\Isom^\nilp$ to $\on{Op}_\fg^\reg$ is
$\on{Op}_\fg^\reg\underset{\on{pt}/B}\times \on{Op}_\fg^\reg$, which
is strictly contained in $\on{Op}_\fg^\reg\underset{\on{pt}/G}\times
\on{Op}_\fg^\reg\simeq \Isom^\reg$.

We shall now establish the following:

\begin{prop}  \label{Atiyah, nilp}
The map of \eqref{abstract coisotrop} induces an isomorphism of Lie
algebroids
$$N^*_{\on{Op}_\fg^\nilp/\on{Op}_\fg(\D^\times)}\simeq \isom^\nilp.$$
\end{prop}

\begin{proof}

For an $R$-point $(\CF_G,\nabla,\CF_B)$ of $\on{Op}_\fg^\nilp$ let us
describe the restrictions of
$N^*_{\on{Op}_\fg^\nilp/\on{Op}_\fg(\D^\times)}$ and $\isom^\nilp$ to
$\Spec(R)$ as subspaces of the restriction of
$\Omega^1(\on{Op}_\fg(\D^\times))\simeq \isom$. We have:
$$T(\on{Op}_\fg^\nilp)|_{\Spec(R)}=\on{coker}(\nabla): \fn_{\CF_B}\to
\left(\fb_{\CF_B}\underset{\BC[[t]]}\otimes \BC[[t]]dt+
\fn_{\CF_B}\underset{\BC[[t]]}\otimes t^{-1}\BC[[t]]dt\right).$$

Hence, $N_{\on{Op}_\fg^\nilp/\on{Op}_\fg(\D^\times)}|_{\Spec(R)}$ is
isomorphic to the cokernel of $\nabla$:
$$\fn_{\CF_B}\underset{\BC[[t]]}\otimes (\BC\ppart dt/\BC[[t]]dt)\to
\fb_{\CF_B}\underset{\BC[[t]]}\otimes \BC\ppart dt/  
(\fb_{\CF_B}\underset{\BC[[t]]}\otimes \BC[[t]]dt+
\fn_{\CF_B}\underset{\BC[[t]]}\otimes t^{-1}\BC[[t]]dt).$$

Finally,
$$N^*_{\on{Op}_\fg^\nilp/\on{Op}_\fg(\D^\times)}|_{\Spec(R)}\simeq
\on{ker}(\nabla): \left((\fg/\fn)_{\CF_B}(-x)+
(\fb/\fn)_{\CF_B}\right) \to (\fg/\fb)_{\CF_B}\underset{\hCO}\otimes
\omega_\D.$$

In other words, we can identify
$N^*_{\on{Op}_\fg^\nilp/\on{Op}_\fg(\D^\times)}|_{\Spec(R)}$ as a subset of
$\Omega^1(\on{Op}_\fg(\D^\times))|_{\Spec(R)}$ with
$$\{\bu\in \fg_{\CF_G}(-x)+ \fb_{\CF_B}\subset \fg_{\CF_G}\,|\, \nabla(\bu)\in 
\fb_{\CF_B}(x)\underset{\hCO}\otimes \omega_\D\}/\{\bu\in \fn_{\CF_B}\}.$$

The latter is easily seen to be the image of $\isom^\nilp|_{\Spec(R)}$ inside
$\isom|_{\Spec(R)}$.

\end{proof}

We shall now study the behavior of the restriction of $\isom^\nilp$
to the subscheme $\on{Op}_\fg^\reg\subset \on{Op}_\fg^\nilp$. The
above proposition combined with \lemref{nilp res}(2) implies:

\begin{cor}   \label{nilp dir}
The Lie algebroid $N^*_{\on{Op}_\fg^\nilp/\on{Op}_\fg(\D^\times)}$
preserves the subscheme $\on{Op}_\fg^\reg$.  The restriction
$N^*_{\on{Op}_\fg^\nilp/\on{Op}_\fg(\D^\times)}|_{\on{Op}_\fg^\reg}$
identifies with the Atiyah algebroid
$\on{At}(\CP_{B,\on{Op}_\fg^\reg})$ of the $B$-bundle
$\CP_{B,\on{Op}_\fg^\reg}$, and we have a commutative diagram
$$
\CD 0 @>>> \fb_{\on{Op}_\fg^\reg} @>>>
N^*_{\on{Op}_\fg^\nilp/\on{Op}_\fg(\D^\times)} |_{\on{Op}_\fg^\reg}
@>>> T(\on{Op}_\fg^\nilp)|_{\on{Op}_\fg^\reg} \\ & & @VVV @VVV @AAA &
& \\ 0 @>>> \fg_{\on{Op}_\fg^\reg} @>>>
N^*_{\on{Op}_\fg^\reg/\on{Op}_\fg(\D^\times)} @>>> T(\on{Op}_\fg^\reg)
@>>> 0 \endCD
$$
\end{cor}

\begin{cor}   \label{nilpotent directions}
The composition
$$N_{\on{Op}_\fg^\reg/\on{Op}_\fg^\nilp}\to
N_{\on{Op}_\fg^\reg/\on{Op}_\fg(\D^\times)}\to
N_{\on{Op}_\fg^\reg/\on{Op}_\fg(\D^\times)}/\Omega^1(\on{Op}_\fg^\reg)\simeq
\fg_{\on{Op}_\fg^\reg}$$ is an injective bundle map, and its image
coincides with $\fn_{\on{Op}_\fg^\reg}\subset \fg_{\on{Op}_\fg^\reg}$.
\end{cor}

\begin{proof}

We claim that it is enough to show that the natural surjection
$N_{\on{Op}_\fg^\reg/\on{Op}_\fg(\D^\times)}\to
N_{\on{Op}_\fg^\nilp/\on{Op}_\fg(\D^\times)}|_{\on{Op}_\fg^\reg}$ fits
into a commutative diagram with exact rows
\begin{equation} \label{nilp diag}
\CD & & \Omega^1(\on{Op}_\fg^\nilp)|_{\on{Op}_\fg^\reg} @>{\anch}>>
N_{\on{Op}_\fg^\nilp/\on{Op}_\fg(\D^\times)}|_{\on{Op}_\fg^\reg} @>>>
(\fg/\fn)_{\on{Op}_\fg^\reg} @>>> 0 \\ & & @VVV @AAA @AAA \\ 0 @>>>
\Omega^1(\on{Op}_\fg^\reg) @>>>
N_{\on{Op}_\fg^\reg/\on{Op}_\fg(\D^\times)} @>>>
\fg_{\on{Op}_\fg^\reg} @>>> 0.  \endCD
\end{equation}

Indeed, this would imply that the map
$N_{\on{Op}_\fg^\reg/\on{Op}_\fg^\nilp} \to \fg_{\on{Op}_\fg^\reg}$
appearing in the corollary is a surjective bundle map onto
$\fn_{\on{Op}_\fg^\reg}$; hence it must be an isomorphism because of
the equality of the ranks.

By Serre duality, the existence of the diagram \eqref{nilp diag} is
equivalent to the diagram appearing in the previous corollary.

\end{proof}

Let us consider now the sequence of embeddings of schemes:
$$\on{Op}_\fg^\reg\hookrightarrow \on{Op}_\fg^\nilp\hookrightarrow
\on{Op}_\fg^\RS.$$ By \thmref{descr of nilp opers}, the normal bundle
$N_{\on{Op}_\fg^\nilp/\on{Op}_\fg^\RS}$ is canonically trivialized and
its fiber isomorphic to the tangent space to $\fh\qu W$ at the point
$-\crho$; this tangent space is in turn canonically isomorphic to
$\fh$.

\begin{lem}   \label{h-directions}
The composition
$$\fh\to N_{\on{Op}_\fg^\nilp/\on{Op}_\fg^\RS}|_{\on{Op}_\fg^\reg}\to
N_{\on{Op}_\fg^\nilp/\on{Op}_\fg(\D^\times)}|_{\on{Op}_\fg^\reg}
\twoheadrightarrow (\fg/\fn)_{\on{Op}_\fg^\reg}$$ equals the canonical
map
$$\fh\simeq (\fb/\fn)_{\on{Op}_\fg^\reg}\hookrightarrow
(\fg/\fn)_{\on{Op}_\fg^\reg}.$$
\end{lem}

\begin{proof}

Let $(\CF_G,\nabla,\CF_B)$ be an $R$-point of $\on{Op}_\fg^\reg$,
written in the form $\nabla^0+p_{-1}dt+\phi(t)dt$, $\phi(t)\in
\fb\otimes R[[t]]$. Then by \propref{lambda opers} the map
$$\fh\to N_{\on{Op}_\fg^\nilp/\on{Op}_\fg^\RS}|_{\Spec(R)}\to 
N_{\on{Op}_\fg^\nilp/\on{Op}_\fg(\D)}|_{\Spec(R)}$$ 
can be realized by 
$$\cla\mapsto \frac{\cla}{t}\in \fb\otimes \BC\ppart dt\subset
T(\on{Op}_\fg(\D^\times))|_{\Spec(R)}.$$

To prove the lemma it would be enough to show that the image of
$\frac{\cla}{t}$ under
$$T(\on{Op}_\fg(\D^\times))|_{\Spec(R)}\to
N_{\on{Op}_\fg^\reg/\on{Op}_\fg(\D^\times)}\to
\fg_{\on{Op}_\fg^\reg}$$ equals $\lambda$. But this follows from
\lemref{exp normal}.

\end{proof}

\ssec{The case of opers with an integral residue}

For completeness, let us describe the behaviour of the groupoid $\Isom$
and the algebroid $\isom$, when restricted to the subscheme 
$$\on{Op}_\fg^{\cla,\nilp}\simeq \on{Op}_\fg^{\RS,\varpi(-\cla-\crho)}
\subset \on{Op}_\fg(\D^\times)$$ when 
$\cla+\crho$ dominant and integral.

Recall that to $\cla$ as above there corresponds a subset $\CJ$ of vertices of
the Dynkin diagram, and a map
$$\Res^{\cla,\nilp}:\on{Op}_\fg^{\cla,\nilp}\to {\mathbf O}_\CJ/B.$$
Let us denote by $\Isom^{\cla,\nilp}$ the groupoid
$$\Isom^{\cla,\nilp}:=\on{Op}_\fg^{\cla,\nilp}\underset{{\mathbf
O}_\CJ/B}\times \on{Op}_\fg^{\cla,\nilp},$$ and let
$\isom^{\cla,\nilp}$ be the corresponding algebroid on
$\on{Op}_\fg^{\cla,\nilp}$.

As in the case of $\cla$ there exists a natural closed embedding
$$\Isom^{\cla,\nilp}\hookrightarrow \Isom|_{\on{Op}_\fg^{\cla,\nilp}}.$$
Repeating the proofs in the $\cla=0$ case we obtain:

\begin{prop}  \label{algebroid for lambda}
The subscheme $\on{Op}_\fg^{\cla,\nilp}\subset \on{Op}_\fg(\D^\times)$
is co-isotropic. The map \eqref{abstract coisotrop} induces an
isomorphism
$$N^*_{\on{Op}_\fg^{\cla,\nilp}/\on{Op}_\fg(\D^\times)}\simeq
\isom^{\cla,\nilp}.$$
\end{prop}

Let us consider a particular example of $\cla=-\crho$. In this case
$\CJ=\CI$, and $${\mathbf O}_\CJ/B\simeq \on{pt}.$$ 
Therefore, the map
$$\isom^{\crho,\nilp}\to T(\on{Op}_\fg^{\crho,\nilp})$$
is surjective. Therefore by \propref{algebroid for lambda}, the map
$$N^*_{\on{Op}_\fg^{\cla,\nilp}/\on{Op}_\fg(\D^\times)}\to
T(\on{Op}_\fg^{\crho,\nilp}),$$ given by the Poisson structure, is
surjective as well. By Serre's duality, the map
$$\Omega^1(\on{Op}_\fg^{\cla,\nilp})\to
N_{\on{Op}_\fg^{\cla,\nilp}/\on{Op}_\fg(\D^\times)}$$ is
injective. This means that the map
$\isom^{-\crho,\nilp}\hookrightarrow
\isom|_{\on{Op}_\fg^{-\crho,\nilp}}$ is an isomorphism.  In fact, it
is easy to see that the map $\Isom^{-\crho,\nilp}\hookrightarrow
\Isom|_{\on{Op}_\fg^{-\crho,\nilp}}$ is an isomorphism.

\medskip

Finally, let us consider the case of $\cla$ which is integral and
dominant.  We have the subscheme $\on{Op}_\fg^{\cla,\reg}\subset
\on{Op}_\fg^{\cla,\nilp}$, and we claim that the behavior of the
groupoid $\Isom$ and the algebroid $\isom$, restricted to it, are the
same as in the $\cla=0$ case. In particular, the analogs of
Corollaries \ref{nilp dir} and \ref{nilpotent directions} hold, when
we replace $\on{Op}_\fg^\nilp$ by $\on{Op}_\fg^{\cla,\nilp}$.

\ssec{Grading on the Lie algebroid} \label{identification of
algebroids, grading}

Recall the action of the group-scheme $\on{Aut}(\D)$ on the scheme
$\on{Op}_\fg(\D)$ and the ind-scheme $\on{Op}_\fg(\D^\times)$. It is
is easy to see that this action lifts to a map from $\on{Aut}(\D)$ to
the groupoids $\Isom$ and $\Isom^\reg$, respectively. In particular,
we obtain a map
$$\on{Der}(\hCO)\simeq\Lie(\on{Aut}(\D))\to \isom^\reg.$$

We choose a coordinate on $\D$ and consider two distinguished elements
$L_0=t\pa_t$ and $L_{-1}=\pa_t$ in $\on{Der}(\hCO)$. The action of
$L_0$ integrates to an action of $\BG_m$, thus defining a grading on
$\isom^\reg$.  Recall also that this choice of a coordinate
trivializes the $B$-bundle $\CP_{B,\on{Op}_\fg^\reg}$ on
$\on{Op}_\fg^\reg$.

\begin{prop} \label{grading on algebroid}  \hfill

\smallskip

\noindent{\em (1)} The image of $L_{-1}$ under
$$\isom^\reg\simeq \on{At}(\CP_{G,\on{Op}_\fg^\reg})\to
\on{At}(\CP_{G,\on{Op}_\fg^\reg})/\on{At}(\CP_{B,\on{Op}_\fg^\reg})\simeq
(\fg/\fb)_{\on{Op}_\fg^\reg}$$ identifies, under the trivialization of
$(\fg/\fb)_{\on{Op}_\fg^\reg}$ corresponding to the above choice of a
coordinate, with the element $p_{-1}\in \fg/\fb$.

\smallskip

\noindent{\em (2)} Under the above trivialization of $\CP_B$, the
subspace $\fg\subset \fg_{\on{Op}_\fg^\reg}$ is $L_0$-stable, and
grading arising on it equals the one induced by $\on{ad}_{\crho}$.

\end{prop}

\begin{proof}

The proof is essentially borrowed from \cite{BD}, Prop. 3.5.18.

By definition, the action of $L_{-1}$ on $\on{Op}_\fg^\reg\wh\times
\D$ lifts onto the triple
$(\CF_{G,\on{Op}_\fg^\reg},\nabla_{\on{Op}_\fg^\reg},
\CF_{B,\on{Op}_\fg^\reg})$.  The lift of $L_{-1}$ onto the $G$-bundle
$\CP_{G,\on{Op}_\fg^\reg}$ is obtained via the identification of the
latter with the space of horizontal (with respect to $\nabla$)
sections of $\CF_{G,\on{Op}_\fg^\reg}$. This lift does not preserve
the reduction of $\CP_{G,\on{Op}_\fg^\reg}$ to $B$; the resulting
element in
$\on{At}(\CP_{G,\on{Op}_\fg^\reg})/\on{At}(\CP_{B,\on{Op}_\fg^\reg})$,
which is the element appearing in point (1) of the proposition,
equals, by definition, to the value of
$$\langle \nabla_{\on{Op}_\fg^\reg}/\CF_{B,\on{Op}_\fg^\reg}, \pa_t\rangle\in
(\fg/\fb)_{\CF_{B,\on{Op}_\fg^\reg}}$$
at $\on{Op}_\fg^\reg\times x\subset \on{Op}_\fg^\reg\hat\times \D$.

When the triple
$(\CF_{G,\on{Op}_\fg^\reg},\nabla_{\on{Op}_\fg^\reg},\CF_{B,\on{Op}_\fg^\reg})$
is written as a connection on the trivial $B$-bundle in the form
$\nabla^0+p_{-1}dt+\bq(t)dt,\bq(t)\in \fb[[t]]$, the above value
equals $p_{-1}$, as required.

\medskip

The second point of the proposition follows immediately from
\lemref{grading}(1).

\end{proof}

\newpage

\vspace*{10mm}

{\Large \part{Categories of representations}}

\vspace*{10mm}


This Part of the paper is devoted to the discussion of various 
categories of representations of affine Kac-Moody algebras of critical level.

In \secref{cat rep} we recall the results of \cite{FF,F:wak} about the
structure of the center of the completed universal enveloping algebra
of $\hg$ at the critical level. According to \cite{FF,F:wak}, the
spectrum of the center is identified with the space $\Op(\D^\times)$
of $\cg$-opers over the formal punctured disc. This means that the
category $\hg_\crit\mod$ "fibers" over the affine ind-scheme
$\Op(\D^\times)$. Next, we introduce the categories of representations
that we study in this project, and in \secref{statement of conj} we
formulate our \mainconjref{main} and \mainthmref{equiv of quot}.

In \secref{generalities} we collect some results concerning the
structure of the category $\hg_\crit\mod$ over its center. In
particular, we discuss the various incarnations of the renormalized
universal enveloping algebra at the critical level. This
renormalization is a phenomenon that has to do with the fact that we
are dealing with a one-parameter family of associative algebras (the
universal enveloping of the Kac-Moody Lie algebra, depending on the
level), which at some special point (the critical level) acquires a
large center.

In \secref{regular} we discuss the subcategory $\hg_\crit\mod_\reg$ of
representations at the critical level, whose support over
$\Op(\D^\times)$ belongs to the subscheme of regular opers. We study
its relation with the category of D-modules on the affine Grassmannian
$\Gr_G\simeq G\ppart/G[[t]]$, and this leads us to \mainconjref{conj
on Grassmann}. We prove \thmref{Hecke fully faithful} which states
that a natural functor in one direction is fully faithful at the
level of derived categories. The formalism of convolution action,
developed in \cite{BD}, Sect. 7, and reviewed in Part V below,
allows us to reduce this assertion to a comparison of self-Exts of a
certain basic object in both cases. On one side the required
computation of Exts had been performed in \cite{ABG}, and on the other
side it follows from the recent paper \cite{FT}.

\secref{manipulation} plays an auxiliary role: we give a proof of one
of the steps in the proof of \thmref{Hecke fully faithful} mentioned
above, by analyzing how the algebra of $G[[t]]$-equivariant self-Exts
of the vacuum module $\BV_\crit$ interacts with $G$-equivariant
cohomology of the point.

\bigskip

\section{Definition of categories}      \label{cat rep}

\ssec{}

Let $\fg$ be a simple finite-dimensional Lie algebra. For an invariant
inner product $\kappa$ on $\fg$ (which is unique up to a scalar) define
the central extension $\hg_\kappa$ of the formal loop algebra $\fg
\otimes \BC\ppart$ which fits into the short exact sequence
$$
0 \to \BC {\mb 1} \to \hg_\kappa \to \fg \otimes \BC\ppart \to 0.
$$
This sequence is split as a vector space, and the commutation
relations read
\begin{equation}    \label{KM rel}
[x \otimes f(t),y \otimes g(t)] = [x,y] \otimes f(t) g(t) + \kappa(x,y)\cdot 
\on{Res}(g\, df)\cdot {\mb 1},
\end{equation}
and ${\mb 1}$ is a central element. The Lie algebra $\hg_\kappa$ is
the {\em affine Kac-Moody algebra} associated to $\kappa$. We will
denote by $\hg_\kappa\mod$ the category of {\em discrete}
representations of $\hg_{\kappa}$ (i.e., such that any vector is
annihilated by $\fg \otimes t^n\BC[[t]]$ for sufficiently large $n$),
on which ${\mb 1}$ acts as the identity.

Let $U_{\kappa}(\hg)$ be the quotient of the universal enveloping
algebra $U(\hg_{\kappa})$ of $\hg_{\kappa}$ by the ideal generated by
$({\mb 1}-1)$. Define its completion $\wt{U}_{\kappa}(\hg)$ as follows:
$$
\wt{U}_{\kappa}(\hg) = \underset{\longleftarrow}\lim \;
U_{\kappa}(\hg)/U_{\kappa}(\hg) \cdot (\fg \otimes t^n\BC[[t]]).
$$ It is clear that $\wt{U}_{\kappa}(\hg)$ is a topological algebra,
whose discrete continuous representations are the same as objects
of $\hg_\kappa\mod$.

The following theorem, due to \cite{FF,F:wak}, describes the center
$Z_\kappa(\hg)$ of $\wt{U}_{\kappa}(\hg)$.

Let $\kappa_{\crit}$ be the {\em critical} inner product on $\fg$
defined by the formula
$$
\kappa_{\crit}(x,y) = - \frac{1}{2} \on{Tr} (\on{ad} (x) \circ \on{ad} (y)).
$$
Denote by $\check G$ the group of adjoint type whose Lie algebra
$\check \fg$ is Langlands dual to $\fg$ (i.e., the Cartan matrix of
$\check\fg$ is the transpose of that of $\fg$).

\begin{thm}   \label{FF} \hfill

\noindent{\em (1)} $Z_\kappa(\hg) = \BC$ if $\kappa \neq \kappa_{\crit}$.

\noindent{\em (2)} $Z_{\crit}(\hg)$ is isomorphic to the algebra
$\on{Fun} (\on{Op}_{\check G}(\D^\times))$ of functions on the space of 
$\check G$-opers on the punctured disc $\D^\times$.
\end{thm}

{}From now on we will denote $Z_{\crit}(\hg)$ simply by $\fZ_\fg$.

\ssec{}  \label{I-categories}

Let $I$ be the Iwahori subgroup of the group $G[[t]]$, i.e., the
preimage of a fixed Borel subgroup $B \subset G$ under the evaluation
homomorphism $G[[t]] \to G$.  Let $I^0\subset I$ be the pro-unipotent
radical of $I$. Non-canonically we have a splitting
$\Lie(I)=\Lie(I^0)\oplus \fh$, where $\fh$ is a Cartan subalgebra of $\fg$.

Recall that an object $\CM\in \fg_\kappa\mod$ is called $I$-integrable
(resp., $I^0$-integrable)
if the action of $\Lie(I)\subset \hg_\kappa$ (resp., $\Lie(I^0)$)
on $\CM$ integrates to an action of the pro-algebraic group $I$
(resp., $I^0$).
In the case of $I^0$ this condition is equivalent to saying that 
$\Lie(I^0)$ acts locally nilpotently, and in the case of $I$ that,
in addition, $\fh$ acts semi-simply with eigenvalues corresponding 
to integral weights. (The latter condition is easily seen to be 
independent of the choice of the splitting $\fh\to \Lie(I)$).

Following the conventions of \secref{equiv cat}, we will denote the
corresponding subcategories of $\hg_\crit\mod$ by $\hg_\crit\mod^I$
and $\hg_\crit\mod^{I^0}$, respectively. We will denote by
$D(\hg_\kappa\mod)^I$ and $D(\hg_\kappa\mod)^{I^0}$ the corresponding
triangulated categories, see  \secref{equiv der cat}.  
Moreover, the functor $D^+(\hg_\kappa\mod)^{I^0}\to D^+(\hg_\kappa\mod)$ is
fully-faithful.

Recall also that an object 
$\CM\in \fg_\kappa\mod$ is called $I$-monodromic if it is
$\Lie(I^0)$-integrable and $\fh$ acts locally
finitely with generalized eigenvalues corresponding to integral weights.
It is evident that a module $\CM$ is $I$-monodromic if and
only if it has an increasing filtration  with successive quotients 
being $I$-integrable. We will denote the subcategory 
of monodromic modules by $\fg_\kappa\mod^{I,m}$.
We will denote by  $D(\hg_\kappa\mod)^{I,m}$ the full subcategory
of $D(\hg_\kappa\mod)$ consisting of complexes with $I$-monodromic
cohomology.

Let us note that the above notions make sense more generally
for an arbitrary category $\CC$ endowed with a Harish-Chandra action 
of $I$ (see \secref{HCh action of groups}). Namely, we have the full
subcategories $$\CC^I\subset \CC^{I^0}\subset \CC$$
along with the equivariant categories $D(\CC)^I$, $D(\CC)^{I,m}$,
$D(\CC)^{I_0}$. Since the group $I^0$ is pro-unipotent, 
the functor $$D^+(\CC)^{I_0}\to D^+(\CC)$$ is fully faithful
and its image consists of complexes, whose cohomologies
are $I^0$-equivariant. We also introduce 
the $I$-monodromic category $\CC^{I,m}$ as the full subcategory of 
$\CC$, consisting of objects that admit a filtration, whose subquotients 
belong to $\CC^I$; we let $D(\CC)^{I,m}$ to be the full subcategory of 
$D(\CC)$, which consists of complexes, whose cohomologies 
belong to $\CC^{I,m}$.

\medskip

{}From now on let us take $\kappa=\kappa_{\crit}$.
Recall the subscheme $\nOp\subset \Op(\D^\times)$, see
\secref{nilpotent opers}.
Let $\hg_\crit\mod_\nilp\subset \fg_\crit\mod$ be the subcategory
consisting of modules, on which the action of the center 
$\fZ_\fg\simeq \on{Fun} (\Op(\D^\times))$
factors through the quotient $\fZ^\nilp_\fg:=\Fun(\nOp)$.
This is a category endowed with an action of $G((t))$,
and in particular, of $I$. 

Our main object of study is the category $\Cat$, where we follow
the above conventions regarding the notion of the $I$-monodromic
subcategory. In other words,
$$\Cat=\hg_\crit\mod^{I,m}\cap \hg_\crit\mod_\nilp.$$
The following will be established in \secref{induction from O}:

\begin{lem}   \label{mon=N-eq}
The inclusion functor
$$\hg_\crit\mod^{I,m}\cap \hg_\crit\mod_\nilp \to
\hg_\crit\mod^{I^0}\cap \hg_\crit\mod_\nilp$$
is an equivalence.
\end{lem}
(In other words, any module in $\hg_\crit\mod$, which
is $I^0$-integrable, and on which the center acts via
$\fZ^\nilp_\fg$, is automatically $I$-monodromic.)

By the above lemma, the inclusion
$$D^+(\hg_\crit\mod_\nilp)^{I,m}\hookrightarrow
D^+(\hg_\crit\mod_\nilp)^{I^0}$$ is in fact an equivalence, and both
these categories identify with the full subcategory of
$D^+(\hg_\crit\mod_\nilp)$, consisting of complexes, whose
cohomologies belong to $\Cat$.

The following assertion seems quite plausible, but we are unable to
prove it at the moment:
\begin{conj}   \label{naive category OK}
The natural functor $D(\Cat)\to \DCat$ 
is an equivalence.
\end{conj}

We will not need it in what follows.

\section{The main conjecture}   \label{statement of conj}

\ssec{}
We shall now formulate our main conjecture. Recall the scheme $\nMOp$, see
\secref{Miura with nilp sing}. Let $D^b\left(\on{QCoh}(\nMOp)\right)$
be the bounded derived category of quasi-coherent sheaves on $\nMOp$.

Our main conjecture is as follows:

\begin{mainconj} \label{main}
We have an equivalence of  triangulated categories 
$$\DCat\simeq D^b\left(\on{QCoh}(\nMOp)\right).$$
\end{mainconj}

In what follows we will provide some motivation for this conjecture.
We will denote a functor establishing the conjectural equivalence
$\DCat\to D^b\left(\on{QCoh}(\nMOp)\right)$ by $\sF$.

\medskip

Note that both categories $D^b(\on{QCoh}(\nMOp))$ and $\DCat$ come
equipped with natural t-structures. The functor $\sF$ will not be
exact, but we expect it to be of bounded cohomological amplitude, and
hence to extend to an equivalence of the corresponding unbounded
derived categories.

\medskip

Recall the ind-scheme $\on{MOp}_{\cg,\gen}(D^\times)$ from \secref{tr
Miura}.  Following \cite{FF:si,F:wak}, to a quasi-coherent sheaf $\CR$
on $\on{MOp}_{\cg,\gen}(D^\times)\simeq
\on{Conn}_{\check{H}}(\omega_{\D^\times}^\rho)$ one can attach a
Wakimoto module $\BW^{w_0}_\crit(\CR)\in \hg_\crit\mod$ (see
\secref{wakimodules} for a review of this construction).

It turns out that if $\CR$ is supported on the closed subscheme
$\on{MOp}_{\cg,\gen}(D^\times)\underset{\Op(\D^\times)}\times
\Op^\nilp$, then $\BW^{w_0}_\crit(\CR)$ belongs to the subcategory
$\Cat$. The main compatibility property that we expect from the
functor $\sF$ is that $\sF(\BW^{w_0}_\crit(\CR))$ will be isomorphic
to the direct image of $\CR$ under the morphism
\begin{equation} \label{split into comp}
\on{MOp}_{\cg,\gen}(D^\times)\underset{\Op(\D^\times)}\times
\Op^\nilp\to \nMOp
\end{equation}
of \thmref{four versions of Miura}.

\medskip

In view of this requirement, the functor $\sF^{-1}$, inverse to $\sF$,
should be characterized by the property that it extends the Wakimoto
module construction from quasi-coherent sheaves on
$\on{MOp}_{\cg,\gen}(D^\times)\underset{\Op(\D^\times)}\times
\Op^\nilp$ to those on $\nMOp$. This was in fact the main motivation
for \mainconjref{main}.

\ssec{}   \label{localization motiv}

In this subsection we would like to explain a point of view on
\conjref{main} as a localization-type statement for affine algebras at
the critical level that connects D-modules on the affine flag variety
to $\hg_\crit\mod$.

This material will not be used in what follows, and for that reason we
shall allow ourselves to appeal to some results and constructions that
are not available in the published literature. One set of such results
is Bezrukavnikov's theory of perverse sheaves on the affine flag
scheme (see \cite{Bez}) and another the formalism of triangulated
categories over stacks (to be developed in \cite{Ga2}).

\medskip

Let $\Fl_G$ be the affine flag scheme corresponding to $G$, i.e.,
$\Fl_G\simeq G\ppart/I$. Let $\fD(\Fl_G)\mod$ denote the category of
right D-modules on $\Fl$.  Let $\fD(\Fl_G)\mod^I$,
$\fD(\Fl_G)\mod^{I^0}$ and $\fD(\Fl_G)\mod^{I,m}$ be the 
subcategories of $I$-equivariant, $I^0$-equivariant and $I$-monodromic
D-modules, respectively. One easily shows that the inclusion functor
$$\fD(\Fl_G)\mod^{I,m} \to \fD(\Fl_G)\mod^{I^0}$$
is in fact an equivalence of categories.

Let $D\left(\fD(\Fl_G)\mod\right)^I$ and
$D\left(\fD(\Fl_G)\mod\right)^{I^0}$ denote the corresponding
triangulated categories.

\medskip

Recall the Grothendieck alteration $\tg\to \cg$ from \secref{Miura
with nilp sing}. Let $\tN$ be the Springer resolution of the nilpotent
cone $\CN_{\cG}\subset \cg$. Let $\on{St}_{\cG}$ be the "thickened"
Steinberg variety $$\on{St}_{\cG}:= \tg\underset{\cg}\times \tN.$$
Note that the scheme $\tn:=\cn\underset{\cg}\times \tg$ introduced in
\secref{Miura with nilp sing} equals the preimage of $\cG/\cB^-\times
\{\check\fb\}$ under the natural map $\on{St}_{\cG}\to \cG/\cB^-\times
\cG/\cB^-$, and we have natural isomorphisms of stacks:
$$\tN/\cG\simeq \cn/\cB \text{ and } \on{St}_{\cG}/\cG\simeq \tn/\cB.$$

The next lemma insures that the definition of the scheme $\tn$ (and,
hence, of $\on{St}_{\cG}$) is not too naive, i.e., that we do not
neglect lower cohomology:

\begin{lem}
The derived tensor product 
$$\Fun(\tg)\overset{L}{\underset{\Fun(\cg)}\otimes} \Fun(\cn)\in
\on{QCoh}(\tg)$$ is concentrated in cohomological dimension $0$.
\end{lem}

\begin{proof}

Consider the vector space $\cg/\cn$. It is enough to show that the
composed map
$$\tg\to\cg\to \cg/\cn$$ is flat near $0\in \cg/\cn$. Since the
varieties we are dealing with are smooth, it is enough to check that
the dimension of the fibers is constant.  The latter is evident.

\end{proof}

\medskip

According to \cite{AB}, there exists a natural tensor functor
$$D^b\left(\on{Coh}(\tN/\cG)\right)\to D^b(\fD(\Fl_G)\mod)^I.$$
In particular, using the convolution action of the monoidal category
$D^b(\fD(\Fl_G)\mod)^I$ on the entire $D^b(\fD(\Fl_G)\mod)$, we obtain 
a monoidal action of $D^b\left(\on{Coh}(\tN/\cG)\right)$ 
on $D^b(\fD(\Fl_G)\mod)$. This construction can be upgraded to a
structure on $D^b(\fD(\Fl_G)\mod)$ of triangulated category over 
the stack $\tN/\cG$, see \cite{Ga2}. In particular, it makes sense
to consider the base-changed triangulated category
\begin{equation}  \label{cartesian product categories}
D^b\left( \fD(\Fl_G)\mod\right)
\underset{\tN/\cG}\times \nOp,
\end{equation}
where we are using the map $\Res^\nilp:\nOp\to \cn/\cB\simeq \tN/\cG$.

A far-reaching generalization of \conjref{main} is the following statement
in the spirit of the localization theorem of \cite{BB}:

\begin{conj}    \label{quasi}
There is an equivalence of triangulated categories
$$
D^b\left( \fD(\Fl_G)\mod\right)
\underset{\tN/\cG}\times \nOp \simeq D^b(\hg_\crit\mod_\nilp).
$$
\end{conj}

A version of this conjecture concerning $\hg_\crit\mod_\reg$, rather
than $\hg_\crit\mod_\nilp$, can be made precise without the machinery
of categories over stacks, and it will be discussed in
\secref{regular}.

Let us explain the convection between the above \conjref{quasi} and
\conjref{main}.  Namely, we claim that the latter is obtained from the
former by passing to the corresponding $I^0$-equivariant categories on
both sides. In order to explain this, we recall the main result of
Bezrukavnikov's theory:

\begin{thm}  \label{bezr}
There is a natural equivalence 
$$D^b\left(\fD(\Fl_G)\mod\right)^{I^0}\simeq 
D^b\Bigl(\on{Coh}\left(\on{St}_{\cG}/\cG\right)\Bigr).$$
\end{thm}

This theorem implies that the base-changed category 
$$D^b\left( \fD(\Fl_G)\mod\right)^{I_0}
\underset{\tN/\cG}\times \nOp$$ is equivalent to
$$D^b\left(\QCoh(\on{St}_{\cG}/\cG\underset{\tN/\cG}\times \nOp)\right),$$
which by \corref{nilp Miura opers as cart} is the same as
$D^b(\QCoh(\nMOp))$.

\ssec{}   \label{partially integrable}

We shall now formulate one of the main results of this paper, which
amounts to an equivalence as in \mainconjref{main}, but at the level
of certain quotient categories. This result provides us with the main
supporting evidence for the validity of \mainconjref{main}. Before
stating the theorem, let us give some motivation along the lines of
\thmref{bezr}.

Let $\CF$ be an $I^0$-integrable D-module on $\Fl_G$. We will say that
it is partially integrable if $\CF$ admits a filtration
$\CF=\underset{k\geq 0}\cup\, \CF_k$ such that each successive
quotient $\CF_k/\CF_{k-1}$ is equivariant with respect to a parahoric
subalgebra $\fp^\iota=\Lie(I)+{\mathfrak {sl}}_2^\iota$ for some
vertex of the Dynkin graph $\iota\in \CI$.

Similarly, we will call an object $\CM$ of $\hg_\crit\mod^{I^0}$
partially integrable if there exists a filtration $\CM=\underset{k\geq
0}\cup\, \CM_k$ such that for each successive quotient
$\CM_k/\CM_{k-1}$ there exists a parahoric subalgebra $\fp^\iota$ as
above such that its action integrates to an action of the
corresponding pro-algebraic group.  More generally, the notion of
partial integrability makes sense in any category equipped with a
Harish-Chandra action of $G\ppart$ (see 
\secref{conv section}, where the latter notion is introduced).

In both cases it is easy to see that partially integrable objects form
a Serre subcategory. Let ${}^f\fD(\Fl_G)\mod^{I^0}$ (resp., $\Catf$)
denote the quotient category of $\fD(\Fl_G)\mod^{I^0}$ (resp., $\Cat$)
by the subcategory of partially integrable objects.  We will denote by
${}^f D\left(\fD(\Fl_G)\mod\right){}^{I^0}$ (resp., $\DCatf$)
the triangulated quotient categories by the subcategories consisting
of objects whose cohomologies are partially integrable.

\medskip

Let us recall the statement from \cite{Bez} that describes the
category ${}^f D^b\left(\fD(\Fl_G)\mod\right){}^{I^0}$ in terms of
quasi-coherent sheaves.

Let $h_0$ denote the algebra of functions on the scheme
$\varpi^{-1}(0)$, where $\varpi$ is the natural
projection $\fh^*\to \fh^*\qu W$. This is a nilpotent algebra of length 
$|W|$.

Recall also that $\fh^*\simeq \check \fh$. We have a natural map
$$\on{St}_{\cG}\simeq \tg\underset{\cg}\times \tN\to
\check \fh \underset{\check \fh\qu W}\times \tN\simeq
\on{Spec}(h_0)\times \tN.$$

\begin{thm}  \label{bezr quot}
There is a canonical  equivalence
$${}^f D^b\left(\fD(\Fl_G)\mod\right){}^{I^0}\simeq
D^b\left(\on{QCoh}\left(\on{Spec}(h_0)\times \tN/\cG\right)\right),$$
so that under the equivalence of \thmref{bezr} the functor
$$D^b\left(\fD(\Fl_G)\mod\right)^{I^0}\to {}^f
D^b\left(\fD(\Fl_G)\mod\right)^{I^0}$$ corresponds to the direct image
under the projection $\on{St}_{\cG}/\cG\to \on{Spec}(h_0)\times
\tN/\cG$.
\end{thm}

Combining this with the \conjref{quasi}, we arrive at the following
statement, which is proved in Part IV of this paper and is
one of our main results.

\begin{mainthm}  \label{equiv of quot}
We have an equivalence:
$$^f\sF:\DCatf\to
D^b\left(\on{QCoh}\left(\on{Spec}(h_0)\times \nOp\right)\right).$$
Moreover, this functor is exact in the sense that it preserves the
natural t-structures on both sides.
\end{mainthm}

\section{Generalities on $\hg_\crit$-modules}   \label{generalities}

\ssec{}

Recall that the ind-scheme $\Op(\D^\times)$ contained the following
subschemes:
$$\Op^\reg\subset \Op^\nilp\subset \Op^\RS=\Op^{\ord_1}\subset
\Op^{\ord_k} \text{ for } k\geq 1.$$

Let us denote by 
$$\fZ_\fg^{\ord_k}\twoheadrightarrow \fZ^\RS_\fg \twoheadrightarrow
\fZ^\nilp_\fg \twoheadrightarrow \fZ^\reg_\fg,$$ respectively, the
corresponding quotients of $\fZ_\fg\simeq \Fun(\Op(\D^\times))$.

\medskip

Let us denote by $\fZ_\fg\mod$ the category of discrete
$\fZ_\fg$-modules.  By definition, any object of this category is a
union of subobjects, each of which is acted on by $\fZ_\fg$ via the
quotient $\fZ_\fg^{\ord_k}$ for some $k$.

\medskip

Let $\imath^\reg$ (resp., $\imath^\nilp$, $\imath^\RS$,
$\imath^{\ord_k}$) denote the closed embedding of
$\Spec(\fZ^\reg_\fg)$ (resp., $\Spec(\fZ_\fg^\nilp)$,
$\Spec(\fZ_\fg^\RS)$, $\Spec(\fZ_\fg^{\ord_k})$) into the ind-scheme
$\Spec(\fZ_\fg)$, and let $\imath^\reg_!$ (resp., $\imath^\nilp_!$,
$\imath^\RS_!$, $\imath^{\ord_k}_!$) denote the corresponding direct
image functor on the category of modules.

It is easy to see that at the level of derived categories we have
well-defined right adjoint functors from $D^+(\fZ_\fg\mod)$ to
$D^+(\fZ^\reg_\fg\mod)$, $D^+(\fZ^\nilp_\fg\mod)$,
$D^+(\fZ^\RS_\fg\mod)$ and $D^+(\fZ^{\ord_k}_\fg\mod)$, denoted
$\imath^\reg{}^!$, $\imath^\nilp{}^!$, $\imath^\nilp{}^!$ and
$\imath^{\ord_k}{}^!$, respectively.

\ssec{}

Let $\hg_\crit\mod_\reg$ (resp., $\hg_\crit\mod_\nilp$,
$\hg_\crit\mod_\RS$, $\hg_\crit\mod_{\ord_k}$) denote the subcategory
of $\hg_\crit\mod$ whose objects are modules on which $\fZ_\fg$ acts
through the corresponding quotient.

The following basic result was established in \cite{BD}, Theorem 3.7.9.

\begin{thm}   \label{induced modules}
The induced module $\on{Ind}^{\hg_\crit}_{t^k\fg[[t]]\oplus
\BC\one}(\BC)$ belongs to $\hg_\crit\mod_{\ord_k}$.
\end{thm}

Here and below, when considering the induced modules such as
$\on{Ind}^{\hg_\crit}_{t^k\fg[[t]]\oplus \BC\one}(\BC)$, we will
assume that $\one$ acts as the identity. We will also need the
following:

\begin{lem}   \label{induction flat}
The module $\on{Ind}^{\hg_\crit}_{t^k\fg[[t]]\oplus \BC\one}(\BC)$
is flat over $\fZ^{\ord_k}_\fg$.
\end{lem}

\begin{proof}

By construction, the PBW filtration on $\fZ_\fg$ induces 
a filtration on $\fZ^{\ord_k}_\fg$ such that
$$\on{gr}(\fZ^{\ord_k}_\fg)\simeq 
\left(\Sym\left(\fg\ppart/t^k\fg[[t]]\right)\right)^{G[[t]]}.$$
This filtration is compatible with the natural filtration on 
$\on{Ind}^{\hg_\crit}_{t^k\fg[[t]]\oplus \BC\one}(\BC)$, 
and it suffices to check
the flatness on the associated graded level.

This reduces the assertion to showing that the algebra
$\Sym\left(\fg\ppart/t^k\fg[[t]]\right)$ is flat over
$\left(\Sym\left(\fg\ppart/t^k\fg[[t]]\right)\right)^{G[[t]]}$.
However, the multiplication by $t^{-k}$ reduces us to the situation
when $k=0$, in which case the required assertion is proved in
\cite{EF}.

\end{proof}

Let us denote by $\imath^\reg_!$ (resp., $\imath^\nilp_!$,
$\imath^\RS_!$, $\imath^{\ord_k}_!$) the evident functor from
$\hg_\crit\mod_\reg$ (resp., $\hg_\crit\mod_\nilp$,
$\hg_\crit\mod_\RS$, $\hg_\crit\mod_{\ord_k}$) to $\hg_\crit\mod$. It
is easy to show that each of these functors admits adjoint, denoted
$\imath^\reg{}^!$ (resp., $\imath^\nilp{}^!$, $\imath^\RS{}^!$,
$\imath^{\ord_k}{}^!$), defined on $D^+(\hg_\crit\mod)$.

{}From \lemref{induction flat} and \lemref{! ind restrictions}
we obtain the following:

\begin{lem}  \label{! as needed}
The functor $\imath^\reg{}^!:D^+(\hg_\crit\mod)\to D^+(\hg_\crit\mod_\reg)$
commutes in the natural sense with the forgetful functors 
$D^+(\hg_\crit\mod)\to D^+(\fZ_\fg\mod)$ and $D^+(\hg_\crit\mod_\reg)\to 
D^+(\fZ_\fg^\reg\mod)$,
and similarly for the $\nilp,\RS$ and $\ord_k$ versions.
\end{lem}

\ssec{}   \label{equivariant restriction}

Let now $K$ be a group-subscheme of $G[[t]]$. Following our conventions,
we will denote by
$\hg_\crit\mod^K$ (resp., $\hg_\crit\mod^K_\reg$,
$\hg_\crit\mod^K_\nilp$, $\hg_\crit\mod^K_\RS$,
$\hg_\crit\mod^K_{\ord_k}$) the corresponding abelian categories of
$K$-equivariant objects, see \secref{HCh action of groups}. We will denote by
$D(\hg_\crit\mod)^K$ (resp., $D(\hg_\crit\mod_\reg)^K$,
$D(\hg_\crit\mod_\nilp)^K$, $D(\hg_\crit\mod_\RS)^K$,
$D(\hg_\crit\mod_{\ord_k})^K$) the corresponding triangulated
categories.

The functors $\imath^\reg_!$ (resp., $\imath^\nilp_!$, $\imath^\RS_!$,
$\imath^{\ord_k}_!$) extend to the $K$-equivariant setting in a
straightforward way. By \propref{! ind res equiv}, we have:

\begin{lem} \label{existence of !, equiv}
There exist functors $\imath^\reg{}^!:D(\hg_\crit\mod)^K\to
D(\hg_\crit\mod_\reg)^K$ (resp.,
$\imath^\nilp{}^!:D^+(\hg_\crit\mod)^K\to D^+(\hg_\crit\mod_\nilp)^K$,
$\imath^\RS{}^!:D^+(\hg_\crit\mod)^K\to D^+(\hg_\crit\mod_\RS)^K$,
$\imath^{\ord_k}{}^!:D^+(\hg_\crit\mod)^K\to
D^+(\hg_\crit\mod_{\ord_k})^K$), that are right adjoint to the
functors $\imath^\reg_!$ (resp., $\imath^\nilp_!$, $\imath^\RS_!$,
$\imath^{\ord_k}_!$), and which commute with the forgetful functors to
the corresponding 
derived categories
$D^+(\hg_\crit\mod)$, $D^+(\hg_\crit\mod_\reg)$,
$D^+(\hg_\crit\mod_\nilp)$, $D^+(\hg_\crit\mod_\RS)$ and
$D^+(\hg_\crit\mod_{\ord_k})$.
\end{lem}

This lemma implies that if $\CM_1,\CM_2$ are two objects of, say 
$\hg_\crit\mod_\reg^K$, then there exists a spectral sequence,
converging to $\on{Ext}^\bullet_{D(\hg_\crit\mod)^K}(\imath^\reg_!(\CM_1),
\imath^\reg_!(\CM_2))$, and whose second term $E_2^{p,q}$ is given by
\begin{equation} \label{spec seq term}
\on{Ext}^p_{D(\hg_\crit\mod_\reg)^K}(\CM_1,\CM_2)\underset{\fZ^\reg_\fg}
\otimes \Lambda^q(N_{\fZ^\reg_\fg/\fZ_\fg}),
\end{equation}
where $N_{\fZ^\reg_\fg/\fZ_\fg}$ denotes the normal bundle to
$\Spec(\fZ^\reg_\fg)$ inside $\Spec(\fZ_\fg)$. The same spectral
sequence exists when we replace the index $\reg$ by either of $\nilp$,
$\RS$ or $\ord_k$.

\ssec{} \label{action of renorm}

We shall now recall a construction related to that of the
{\it renormalized} universal enveloping algebra at the
critical level, following \cite{BD}, Sect. 5.6.

The main ingredient is the action of the algebra $\fZ_\fg$ on
$\wt{U}_\crit(\hg)$ by {\it outer} derivations. Let us recall
the construction:

\medskip

Let us pick a non-zero (symmetric, invariant) pairing
$\kappa:\fg\otimes \fg\to \BC$, and using it construct a $1$-parameter
deformation of the critical pairing:
$\kappa_\hslash=\kappa_\crit+\hslash\cdot \kappa$. We obtain a
$1$-parametric family of topological associative algebras
$\wt{U}_\hslash(\hg)$.  For an element $a\in \fZ_\fg$, and its lift
$a_\hslash\in \wt{U}_\hslash(\hg)$ and $b\in \wt{U}_\hslash(\hg)$, the
element $[a_\hslash,b]\in \wt{U}_\hslash(\hg)$ is $0$ modulo
$\hslash$.

Hence, the operation $b\mapsto \frac{[a_\hslash,b]}{\hslash}\,
\on{mod}\, \hslash$ is a derivation of $\wt{U}_\crit(\hg)$. It does
not depend on the choice of the lifting $a_\hslash$ up to inner
derivations. This construction has the following properties:

\begin{lem} \hfill

\smallskip

\noindent{\em (a)} The constructed map $\fZ_\fg\to
\on{Der}^{\on{out}}(\wt{U}_\crit(\hg))$ is a derivation, i.e., it
extends to a (continuous) map of (topological) $\fZ_\fg$-modules
$\Omega^1(\fZ_\fg)\to \on{Der}^{\on{out}}(\wt{U}_\crit(\hg))$.

\smallskip

\noindent{\em (b)} Each of the above derivations preserves the
subalgebra $\fZ_\fg\subset \wt{U}_\crit(\hg)$, i.e., $\fZ_\fg$ is a
topological Poisson algebra and $\Omega^1(\fZ_\fg)$ is an algebroid
over $\Spec(\fZ_\fg)$.

\end{lem}

\medskip

The following result, which relates the Poisson algebra structure on
$\fZ_\fg$ with Langlands duality, is crucial for this paper:

Recall from \secref{identification of algebroids} that $\dIsom$
denotes the groupoid over the ind-scheme $\Op(\D^\times)$, whose fiber
over $\chi,\chi'\in \Op(\D^\times)$ is the scheme of isomorphisms of
$\cG$-local systems on $\D^\times$, corresponding to $\chi$ and
$\chi'$, respectively, and $\disom$ denotes its algebroid. One of the
key properties of the the isomorphism \thmref{FF}, proved in
\cite{FF,F:wak}, is that it respects the Poisson structures. In other
words, in terms of the corresponding Lie algebroids (see
\secref{identification of algebroids}) we have the following

\begin{thm}  \label{FF for algebroids}
Under the isomorphism $\fZ_\fg\simeq \Fun(\Op(D^\times))$,
we have a canonical identification of the Lie algebroids
$\Omega^1(\fZ_\fg)\simeq \disom$.
\end{thm}

\medskip

Let us now derive some consequences from the construction
described above. By \lemref{regular coisotrop} and its variant
for the $\nilp$, $\RS$ and $\ord_k$ cases, we obtain the following

\begin{cor}
The ideal of each of the quotient algebras
$\fZ_\fg^\reg$ (resp., $\fZ_\fg^\nilp$, $\fZ_\fg^\RS$,
$\fZ_\fg^{\ord_k}$) is stable under the Poisson bracket, i.e.,
$N^*_{\fZ_\fg^\reg/\fZ_\fg}$ (resp., $N^*_{\fZ_\fg^\nilp/\fZ_\fg}$,
$N^*_{\fZ_\fg^\RS/\fZ_\fg}$, $N^*_{\fZ_\fg^{\ord_k}/\fZ_\fg}$) is an
algebroid over the corresponding algebra.
\end{cor}

Let observe that
for any $\hg_\crit$-module $\CM$ we obtain a map
\begin{equation} \label{action of algebroid on category}
\Omega^1(\fZ_\fg)\to \on{Ext}^1_{\hg_\crit\mod}(\CM,\CM).
\end{equation}
This map is functorial in the sense that for a morphism of
$\hg\crit$-modules $\CM\to \CM'$, the two compositions
$$\Omega^1(\fZ_\fg)\to \on{Ext}^1_{\hg_\crit\mod}(\CM,\CM)\to
 \on{Ext}^1_{\hg_\crit\mod}(\CM,\CM')$$
and
$$\Omega^1(\fZ_\fg)\to \on{Ext}^1_{\hg_\crit\mod}(\CM',\CM')\to
 \on{Ext}^1_{\hg_\crit\mod}(\CM,\CM')$$
coincide.

\medskip

The next series of remarks is stated for the subscheme
$\Spec(\fZ^\reg_\fg)\subset \Spec(\fZ_\fg)$; however, they equally
apply to the cases when $\reg$ is replaces by either of $\nilp$, $\RS$
or $\ord_k$.

Note that the Poisson structure, viewed as a map $\Omega^1(\fZ_\fg)\to
T(\fZ_\fg)$ gives rise to a commutative diagram:
$$
\CD 0 @>>> N^*_{\fZ_\fg^\reg/\fZ_\fg} @>>>
\Omega^1(\fZ_\fg)|_{\Spec(\fZ^\reg_\fg)} @>>> \Omega^1(\fZ^\reg_\fg)
@>>> 0 \\ & & @VVV @VVV @VVV & \\ 0 @>>> T(\fZ^\reg_\fg) @>>>
T(\fZ_\fg)|_{\Spec(\fZ^\reg_\fg)} @>>> N_{\fZ_\fg^\reg/\fZ_\fg} @>>>
0.  \endCD
$$

Let $\CM$ and $\CM'$ be two objects of $\hg_\crit\mod_\reg$. 
Note that we have an exact sequence
\begin{align*}
&0\to \on{Ext}^1_{\hg_\crit\mod_\reg}(\CM,\CM')\to
\on{Ext}^1_{\hg_\crit\mod}(\CM,\CM')\to
\Hom(\CM,\CM')\underset{\fZ^\reg_\fg}\otimes
N_{\fZ^\reg_\fg/\fZ_\fg}\to \\
&\to\on{Ext}^2_{\hg_\crit\mod_\reg}(\CM,\CM').
\end{align*}

It is easy to see that the composed map
$$\Hom(\CM,\CM')\underset{\fZ_\fg}\otimes \Omega^1(\fZ_\fg)\to
\on{Ext}^1_{\hg_\crit\mod}(\CM,\CM')\to
\Hom(\CM,\CM')\underset{\fZ^\reg_\fg}\otimes
N_{\fZ^\reg_\fg/\fZ_\fg}$$ comes from the map $\Omega^1(\fZ_\fg)\to
N_{\fZ^\reg_\fg/\fZ_\fg}$ from the above commutative diagram. Thus, we
obtain the following commutative diagram:
\begin{equation} \label{map to Ext 1}
\CD
\Hom(\CM,\CM')\underset{\fZ^\reg_\fg}\otimes N^*_{\fZ_\fg^\reg/\fZ_\fg} 
@>>> \on{Ext}^1_{\hg_\crit\mod_\reg}(\CM,\CM') \\ 
@VVV  @VVV \\
\Hom(\CM,\CM')\underset{\fZ^\reg_\fg}\otimes 
\Omega^1(\fZ_\fg)|_{\Spec(\fZ^\reg_\fg)}  @>>>
\on{Ext}^1_{\hg_\crit\mod}(\CM,\CM') \\
@VVV @VVV \\
\Hom(\CM,\CM')\underset{\fZ^\reg_\fg}\otimes \Omega^1(\fZ^\reg_\fg) @>>>
\Hom(\CM,\CM')\underset{\fZ^\reg_\fg}\otimes N_{\fZ^\reg_\fg/\fZ_\fg}
\endCD
\end{equation}
and a natural map
\begin{equation} \label{map to Ext 2}
\Hom(\CM,\CM')\underset{\fZ^\reg_\fg}\otimes
\left(N_{\fZ^\reg_\fg/\fZ_\fg}/\Omega^1(\fZ^\reg_\fg)\right)\to 
\on{Ext}^2_{\hg_\crit\mod_\reg}(\CM,\CM').
\end{equation}

\medskip

Let us now consider once again the family $\wt{U}_\hslash(\hg)$, and
inside $\wt{U}_\hslash(\hg)\underset{\BC[[\hslash]]}\otimes
\BC((\hslash))$ consider the $\BC[[\hslash]]$-subalgebra generated by
$\wt{U}_\hslash(\hg)$ and elements of the form
$$\frac{a_\hslash}{\hslash} \text{ for } a_\hslash \text{ mod }\hslash\in
\on{ker}(\fZ_\fg\to \fZ^\reg_\fg).$$  Taking this algebra modulo
$\hslash$, we obtain an algebra, denoted $U^{\ren,\reg}(\hg_\crit)$, and
called the renormalized enveloping algebra at the critical level. The
algebra $U^{\ren,\reg}(\hg_\crit)$ has a natural filtration, with the
$0$-th term isomorphic to $\wt{U}_\crit(\hg)\underset{\fZ_\fg}\otimes
\fZ^\reg_\fg$, and the first associated graded quotient isomorphic to
$$\left(\wt{U}_\crit(\hg)\underset{\fZ_\fg}\otimes
\fZ^\reg_\fg\right)\underset{\fZ^\reg_\fg} {\wh\otimes}
N^*_{\fZ_\fg^\reg/\fZ_\fg}.$$

Let us $U^{\ren,\reg}(\hg_\crit)\mod$ denote the category of (discrete)
$U^{\ren,\reg}(\hg_\crit)$-modules.  We have a tautological homomorphism
$\wt{U}_\crit(\hg)\to U^{\ren,\reg}(\hg_\crit)$, whose restriction to
$\fZ_\fg$ factors through $\fZ^\reg_\fg$; thus we have a restriction
functor $U^{\ren,\reg}(\hg_\crit)\mod\to \hg_\crit\mod_\reg$.  In addition,
the adjoint action of the algebra
$\wt{U}_\crit(\hg)\underset{\fZ_\fg}\otimes \fZ^\reg_\fg$ on itself
extends to an action of the first term of the above-mentioned
filtration on $U^{\ren,\reg}(\hg_\crit)$.

Let now $\CM_\hslash$ be an $\hslash$-family of modules over
$\hg_\hslash$ such that the action of $\fZ_\fg$ on
$\CM:=\CM_\hslash/\hslash\cdot \CM_\hslash$ factors through
$\fZ^\reg_\fg$. Then $\CM$ is naturally acted on by
$U^{\ren,\reg}(\hg_\crit)$.  This construction provides a supply of objects
of $U^{\ren,\reg}(\hg_\crit)\mod$.

\begin{lem}  \label{action of algebroid on exts}
Let $\CM,\CM'$ be $U^{\ren,\reg}(\hg_\crit)$-modules. Then

\smallskip

\noindent{\em (a)} The map
$\Hom_{\hg_\crit\mod}(\CM,\CM')\underset{\fZ^\reg_\fg}\otimes
N^*_{\fZ_\fg^\reg/\fZ_\fg} \to
\on{Ext}^1_{\hg_\crit\mod_\reg}(\CM,\CM')$ vanishes.

\smallskip

\noindent{\em (b)}
We have a natural action of the algebroid $N^*_{\fZ_\fg^\reg/\fZ_\fg}$ on
$\on{Ext}^\bullet_{\hg_\crit\mod_\reg}(\CM,\CM')$.
\end{lem}

\medskip

Finally, let us note that the category of
$U^{\ren,\reg}(\hg_\crit)$-modules carries a Harish-Chandra action of
$G\ppart$. In particular, if $K$ is a group-subscheme of $G[[t]]$, we
can introduce the categories $U^{\ren,\reg}(\hg_\crit)\mod^K$ and
$D(U^{\ren,\reg}(\hg_\crit))^K$. In addition, analogs of the diagrams
appearing above remain valid for
$\on{Ext}^\bullet_{\hg_\crit\mod_\reg}(\CM,\CM')$ and
$\on{Ext}^\bullet_{\hg_\crit\mod}(\CM,\CM')$ replaced by
$\on{Ext}^\bullet_{D(\hg_\crit\mod_\reg)^K}(\CM,\CM')$ and
$\on{Ext}^\bullet_{D(\hg_\crit\mod)^K}(\CM,\CM')$, respectively.

\ssec{}

For the rest of this section we will be concerned with the category
$\hg_\crit\mod_{\on{ord}_1}$, denoted also by $\hg_\crit\mod_\RS$.

Consider now the functor $\fg\mod\to\hg_\crit\mod$ given by
\begin{equation} \label{ind funct}
M\mapsto \on{Ind}^{\hg_\crit}_{\fg[[t]]\oplus \BC\one}(M),
\end{equation}
where $\fg[[t]]$ acts on $M$ via the evaluation map $\fg[[t]]\to \fg$
and $\one$ acts as identity.

By definition.
$$\on{Ind}^{\hg_\crit}_{\fg[[t]]\oplus \BC\one}(U(\fg))\simeq
\on{Ind}^{\hg_\crit}_{t\fg[[t]]\oplus \BC\one}(\BC),$$
and by \thmref{induced modules}, this module belongs to
$\hg_\crit\mod_\RS$. This implies that the module
$\on{Ind}^{\hg_\crit}_{\fg[[t]]\oplus \BC\one}(M)\in \hg_\crit\mod_\RS$
for any $M$.

\medskip

In what follows we will need the following technical assertions, in
which we use the notion of quasi-perfectness introduced in 
\secref{DG categories}.

\begin{prop}   \label{induction perfect} \hfill

\smallskip 

\noindent{\em (1)} Representations of the form
$\on{Ind}^{\hg_\crit}_{\fg[[t]]\oplus \BC\one}(M)$ for $M\in \fg\mod$
are quasi-perfect as objects of $D(\hg_\crit\mod)$.

\smallskip

\noindent{\em (2)} Any object $\CM\in \hg_\crit\mod_\RS$, which is
quasi-perfect in $\hg_\crit\mod$, is also quasi-perfect in
$\hg_\crit\mod_\RS$. The same is true when the RS condition
is replaced by any of $\on{ord}_k$, $\nilp$
or $\reg$.

\end{prop}

\begin{proof}

Since the induction functor is exact, by adjunction,
$$\Hom_{D(\fg_\crit\mod)}(\on{Ind}^{\hg_\crit}_{\fg[[t]]\oplus
\BC\one}(M),\CM_1^\bullet)\simeq
\Hom_{D(\fg[[t]]\mod)}(M,\CM_1^\bullet).$$ 
When $\CM_1^\bullet$ is bounded from below the latter
is computed by the standard cohomological complex of $\fg[[t]]$
(see \secref{standard complex}), which manifestly commutes with direct
sums.  This proves the first point of the proposition.

The second point follows from \propref{quasi-perf}.

\end{proof}

\ssec{}   \label{tau !}

Denote by $\BM_\lambda$ (resp., $\BM_\lambda^\vee$, $\BL_\lambda$) the
$\hg_\crit$-module induced from the Verma module $M_\lambda$ (resp.,
the contragredient Verma module $M_\lambda^\vee$, the irreducible
module $L_\lambda$) with highest weight $\lambda$ over $\fg$:
$$
\BM_\lambda = \on{Ind}^{\hg_\crit}_{\fg[[t]]\oplus
\BC\one}(M_\lambda), \quad \BM_\lambda^\vee =
\on{Ind}^{\hg_\crit}_{\fg[[t]]\oplus \BC\one}(M_\lambda^\vee), \quad
\BL_\lambda = \on{Ind}^{\hg_\crit}_{\fg[[t]]\oplus
\BC\one}(L_\lambda).
$$

Recall that we have the natural residue map $\Res^\RS:\Op^\RS\to 
\check \fh\qu W\simeq \fh^*\qu W$. At the level of algebras of functions
we have therefore a map
\begin{equation}  \label{residue}
\Res^\RS{}^*:\Sym(\fh)^W\to \fZ^\RS_\fg.
\end{equation}

Thus, for every $M\in \fg\mod$ we obtain two a priori different
actions of $\Sym(\fh)^W$ on
$\on{Ind}^{\hg_\crit}_{\fg[[t]]\oplus \BC\one}(M)$:

One action corresponds to the map $\Res^\RS{}^*$ and the action
of $\fZ^\RS_\fg$ on objects of $\hg_\crit\mod_\RS$. Another action comes
from the Harish-Chandra isomorphism 
\footnote{which we normalize so that the central character of $M_\lambda$ equals
$\varpi(\lambda+\rho)$}
\begin{equation} \label{HCh isomo}
\Sym(\fh)^W\simeq Z(U(\fg)),
\end{equation}
the action of $Z(U(\fg))$ by endomorphisms on $M$, and, hence, by functoriality,
its action on
$\on{Ind}^{\hg_\crit}_{\fg[[t]]\oplus \BC\one}(M)$.

Let $\tau$ be the involution of $Z(U(\fg))$, induced by the anti-involution 
$x\mapsto -x$ of $U(\fg)$. Alternatively, $\tau$ can be thought of as induced
by the outer involution of $\fg$ that acts on the weights as $\lambda\mapsto -w_0(\lambda)$.

\begin{prop}  \label{discrepancy}
The above two actions of $\Sym(\fh)^W$ on
$\on{Ind}^{\hg_\crit}_{\fg[[t]]\oplus \BC\one}(M)$ 
differ by $\tau$.
\end{prop}

\begin{proof}

It is enough to consider the universal example of
$M=U(\fg)$. In the course of the proof of the proposition
we will essentially reprove \thmref{induced modules}.

Consider the grading on $\hg_\crit$ induced by the $\BG_m$-action on
$\D$ by loop rotations. Then all our objects, such as $\fZ_\fg$,
$\fZ^\RS_\fg$ and $\on{Ind}^{\hg_\crit}_{t\fg[[t]]\oplus
\BC\one}(\BC)$ acquire a natural grading; the degree $i$ subspace will
be denoted by the subscript $i$, i.e., $(\cdot)_i$.

Consider the ideal $\fZ_\fg\cdot (\fZ_\fg)_{>0}$ in $\fZ_\fg$ generated 
by elements of positive degree. From \secref{grading}, we know that the 
quotient $\fZ_\fg/\fZ_\fg\cdot (\fZ_\fg)_{>0}$ is precisely $\fZ_\fg^\RS$.
Since the grading on $\on{Ind}^{\hg_\crit}_{t\fg[[t]]\oplus \BC\one}(\BC)$ is 
non-positive and the module is generated by the subspace of degree $0$, 
the above ideal annihilates this module.

\medskip

Consider now the subalgebra of degree $0$ elements
$(\fZ_\fg/\fZ_\fg\cdot (\fZ_\fg)_{>0})_0\subset \fZ_\fg^\RS$.
According to \secref{grading}, it is isomorphic to $\Sym(\fh)^W$ 
and the resulting embedding 
\begin{equation}  \label{inv deg zero}
\Sym(\fh)^W\to \fZ_\fg^\RS
\end{equation}
is the homomorphism $\Res^\RS{}^*$.

The action of $(\fZ_\fg/\fZ_\fg\cdot (\fZ_\fg)_{>0})_0$ on
$\on{Ind}^{\hg_\crit}_{t\fg[[t]]\oplus \BC\one}(\BC)$ 
preserves the subspace of degree $0$ elements.
However, the latter subspace is isomorphic to $U(\fg)$.
Therefore, $(\fZ_\fg/\fZ_\fg\cdot (\fZ_\fg)_{>0})_0$ acts on
$U(\fg)$ commuting with both left and right module 
structure; hence it comes from a homomorphism
$(\fZ_\fg/\fZ_\fg\cdot (\fZ_\fg)_{>0})_0\to Z(U(\fg))$.

\medskip

It remains to compare the resulting homomorphism
$$\Sym(\fh)^W\to (\fZ_\fg/\fZ_\fg\cdot (\fZ_\fg)_{>0})_0\to
Z(U(\fg))$$ with the Harish-Chandra isomorphism. This has been proved
in \cite{F:wak}, Sect. 12.6. Let us repeat the argument for
completeness:

It is enough to show that for any weight $\lambda\in \fh^*$, 
the two characters, corresponding to $\Sym(\fh)^W$ acting in the
two ways on the module $\BM^\vee_\lambda$, differ by $\tau$.

\medskip

Let $\BW^{w_0}_{\crit,\lambda}$ be the Wakimoto module corresponding
to the weight $\lambda$, as in \secref{wakimoto lambda}. By
\lemref{wakimoto with lambda are flat}, the character of
$\Sym(\fh)^W$, acting on $\BW^{w_0}_{\crit,\lambda}$ via 
\eqref{inv deg zero}, is given by $\varpi(-\lambda-\rho)$.

By \secref{wakimoto lambda}, we have have a non-trivial homomorphism
$\BM^\vee_\lambda\to \BW^{w_0}_{\crit,\lambda}$, and hence the
center $\fZ_\fg$ acts on both modules by the same character.

\end{proof}

Recall that for $\chi\in \fh^*//W\simeq \check\fh//W$ we have a
subscheme $\Op^{\RS,\chi}\subset \Op^\RS$; if $\mu\in \fh^*$ is
integral and anti-dominant, then $\Op^{\RS,\varpi(\mu)}\simeq
\Op^{-\mu-\rho,\nilp}$; if, moreover, $\mu+\rho$ is anti-dominant,
then the latter scheme contains the subscheme $\Op^{-\mu-\rho,\reg}$.

Let us denote by $\fZ_\fg^{\RS,\chi}$, $\fZ_\fg^{-\mu-\rho,\nilp}$ and
$\fZ_\fg^{-\mu-\rho,\reg}$, respectively, the corresponding quotients
of $\fZ_\fg$. Let $\hg_\crit\mod_{\RS,\chi}$,
$\hg_\crit\mod_{-\mu-\rho,\nilp}$, $\hg_\crit\mod_{-\mu-\rho,\reg}$ be
the corresponding subcategories of $\hg_\crit\mod$. The general
results stated in this section, concerning the behavior of
$\hg_\crit\mod_{\reg}$, $\hg_\crit\mod_{\nilp}$, $\hg_\crit\mod_{\RS}$
and $\hg_\crit\mod_{\ord_k}$, are equally applicable to
$\hg_\crit\mod_{\RS,\chi}$, $\hg_\crit\mod_{-\mu-\rho,\nilp}$ and
$\hg_\crit\mod_{-\mu-\rho,\reg}$.

\medskip

{}From \propref{discrepancy} we obtain the following:

\begin{cor}
The modules $\BM_\lambda$, $\BM^\vee_\lambda$ and $\BL_\lambda$
belong to $\hg_\crit\mod_{\RS,\varpi(-\lambda-\rho)}$.
\end{cor}

For a dominant integral weight $\lambda$, let $V^\lambda$ be the
corresponding irreducible finite-dimensional $\fg$-module. Let
$\BV^\lambda_\crit$ denote the corresponding induced module at the
critical level. In \secref{W & V} we will also establish the
following:

\begin{prop}  \label{support of V lambda}
The module $\BV^\lambda_\crit$ belongs to
$\hg_\crit\mod_{\lambda,\reg}$.
\end{prop}

\ssec{}    \label{induction from O}

Recall now that the subscheme $\nOp\subset \Op^\RS$ is the preimage of
$\varpi(-\rho)\in \fh^*\qu W$ under the map $\res:\Op^\RS\to \fh^*\qu W$.

In particular, if we denote by $\CO_0$ the subcategory of the usual
category $\CO$ corresponding to $\fg$-modules with central character
equal to $\varpi(-\rho)$, we obtain that the induction \eqref{ind
funct} defines a functor $\CO_0\to \Cat$. In particular, the modules
$\BM_{w(\rho)-\rho}$, $\BM^\vee_{w(\rho)-\rho}$ for $w\in W$ all
belong to $\Cat$.

\medskip

In what follows we will consider sections of right D-modules on the 
affine flag variety $\Fl_G$. Instead of ordinary right D-modules, we will 
consider the ones twisted by a line bundle, which is the tensor product 
of the critical line bundle on $\Gr_G$ and the $G\ppart$-equivariant 
line bundle, corresponding to the weight $2\rho$ (this choice
is such that the twisting induced on $G/B\subset \Fl_G$ corresponds
to {\it left} D-modules on $G/B$.)

We will denote the resulting category by $\fD(\Fl_G)_\crit\mod$,
and by a slight abuse of language we will continue to call its
objects D-modules. Of course, as an abstract category 
$\fD(\Fl_G)_\crit\mod$ is equivalent to $\fD(\Fl_G)\mod$, but
the functor of global sections is different. We have
$$\on{R}\Gamma:D^+\left(\fD(\Fl_G)_\crit\mod\right)\to
D^+(\hg_\crit\mod).$$

In particular, $\Gamma(\Fl_G,\delta_{1_{\Fl_G}})\simeq \BM_{-2\rho}$.

\medskip

As usual, if $K$ is a subgroup of $G[[t]]$, we
will denote by $\fD(\Fl_G)_\crit\mod^K$ the abelian category
of $K$-equivariant D-modules, and by 
$D\left(\fD(\Fl_G)_\crit\mod^K\right)^K$ the corresponding
triangulated category.

For $\CF^\bullet\in D^+\left(\fD(\Fl_G)_\crit\mod\right)$, we have
$$\Gamma(\Fl_G,\CF^\bullet)\simeq \CF^\bullet\star \BM_{-2\rho}.$$
Hence, we obtain the following

\begin{cor}
The functor of global sections gives rise to a functor
$$D^+\left(\fD(\Fl_G)_\crit\mod\right)\to D^+(\hg_\crit\mod_\nilp).$$
\end{cor}

\ssec{}

Let us now prove \lemref{mon=N-eq}.

\begin{proof}

Let $\CM$ be an $I^0$-integrable module. Then it admits a filtration
whose subquotients are quotients of modules of the form
$\on{Ind}^{\hg_\crit}_{\fg[[t]]\oplus \BC\one}(M)$, where
$M$ is an $N$-integrable $\fg$-module, 

If we impose the condition that $\CM\in \hg_\crit\mod_\nilp$, then by
\propref{discrepancy}, we can assume that the above $M$ 
has central character $\varpi(-\rho)$. But, as is well-known, this implies
that $M\in \CO_0$.

\end{proof}

In addition, we have the following result:

\begin{lem}  \label{hom from L}
Any object $\CM\in \Cat$ admits a non-zero map
$\BL_{w(\rho)-\rho}\to \CM$.
\end{lem}

\begin{proof}

By definition, any $\CM$ contains a vector, annihilated by $\Lie(I^0)$,
and which is an eigenvector of $\fh$.
Hence, we have a non-trivial map $\BM_\lambda\to \CM$. By
\propref{discrepancy}, $\lambda$ must be of the form
$w(\rho)-\rho$ for some $w\in W$.  

The Verma module $M_{w(\rho)-\rho}$ admits a filtration, whose
subquotients are the irreducibles $L_{w'(\rho)-\rho}$, $w'\geq w$.
Since the induction functor is exact, $\BM_{w(\rho)-\rho}$ admits
a filtration with subquotients isomorphic to $\BM_{w'(\rho)-\rho}$.

Let $w'$ be the maximal element such that corresponding term
of the filtration on $\BM_{w(\rho)-\rho}$ maps non-trivially to $\CM$.
This gives the desired map.

\end{proof}

\section{The case of regular opers}   \label{regular}

\ssec{}

Recall that the pre-image of $\on{pt}/\cB\hookrightarrow
\cn/\cB$ under $\nOp\to \cn/\cB$ is the scheme
$\Op^\reg$ of regular $\cG$-opers on the disc $\D$.  From the point of
view of representations, the algebra $\fZ^\reg_\fg\simeq\Fun(\Op^\reg)$ is
characterized as follows. Let
$$
\BV_{\crit} \simeq
\on{Ind}^{\hg_\crit}_{\fg[[t]]\oplus \BC\one}(\BC)
$$
be the vacuum Verma module of critical level. According to
\cite{FF,F:wak}, the action of the center $\fZ_\fg \simeq
\Fun(\Op(\D^\times))$ on $\BV_{\crit}$ factors
through its quotient $\Fun(\Op^\reg)$. Moreover, the latter algebra is
isomorphic to the algebra of endomorphisms of $\BV_{\crit}$.

In this section we will be concerned with the category
$\hg_\crit\mod_\reg$ and its derived category
$D\left(\hg_\crit\mod_\reg\right)$. We will see that there are many
parallels between the categories $D\left(\hg_\crit\mod_\reg\right)$
and $D\left(\hg_\crit\mod_\nilp\right)$, but the structure of the
former is considerably simpler.

Let $\Catr$ denote the full subcategory of $\hg_\crit\mod_\reg$ equal
to the intersection $\hg_\crit\mod_\reg\cap \hg_\crit\mod^{I,m}$; we
let $\DCatr:=D(\hg_\crit\mod_\reg)^{I_0}$ denote the corresponding
full triangulated category.

In this section we will formulate a conjecture that describes these
categories in terms of D-modules on the affine Grassmannian.

\ssec{}

Before stating the conjecture we would like to motivate it by
Bezrukavnikov's theory in the spirit of \secref{localization motiv}.
In this subsection the discussion will be informal.

Let $\Gr_G=G\ppart/G[[t]]$ be the affine Grassmannian of the group
$G$. We will consider right D-modules on $\Gr_G$ and denote this
category by $\fD(\Gr_G)\mod$. As before, we have the subcategories
$\fD(\Gr_G)\mod^I$, $\fD(\Gr_G)\mod^{I_0}\simeq \fD(\Gr_G)\mod^{I,m}$
and the corresponding triangulated categories
$D\left(\fD(\Gr_G)\mod\right)^I$,
$D\left(\fD(\Gr_G)\mod\right)^{I_0}\subset
D\left(\fD(\Gr_G)\mod\right)$.

\medskip

Consider the two categories appearing in \conjref{quasi}, and let us
apply a further base change with respect to the map $\Op^\reg\to
\nOp$. We obtain an equivalence:
\begin{equation} \label{quasi-conj for reg}
D^b(\fD(\Fl_G)\mod)\underset{\tN/\cG}\times (\nOp\underset{\nOp}\times
\Op^\reg)\simeq D^b(\hg_\crit\mod_\nilp)\underset{\nOp}\times
\Op^\reg.
\end{equation}

The right-hand side is by definition equivalent to 
$D^b(\hg_\crit\mod_\reg)$. The left-hand side can be rewritten as
$$\left(D^b(\fD(\Fl_G)\mod)\underset{\tN/\cG}\times \on{pt}/\cB\right)
\underset{\on{pt}/\cB}\times \Op^\reg.$$

\medskip

The theory of spherical sheaves on the affine Grassmannian implies
that $D^b\left(\fD(\Gr_G)\mod\right)$ is naturally {\it a category
over the stack} $\on{pt}/\cG$ in the sense explained in
\secref{Satake}. It follows from Bezrukavnikov's theory \cite{Bez}
that the categories $D^b\left(\fD(\Gr_G)\mod\right)$ and
$D^b(\fD(\Fl_G)\mod)$ are related as follows:
$$D^b(\fD(\Fl_G)\mod)\underset{\tN/\cG}\times \on{pt}/\cB \simeq
D^b\left(\fD(\Gr_G)\mod\right)\underset{\on{pt}/\cG}\times \on{pt}/\cB.$$

Hence, from \eqref{quasi-conj for reg} we obtain the following
conjecture:
$$D^b\left(\fD(\Gr_G)\mod\right)\underset{\on{pt}/\cG}\times
\Op^\reg \simeq D^b(\hg_\crit\mod_\reg).$$

Our \conjref{conj on Grassmann} below reformulates the last statement
in terms that do not require the formalism of categories over a stack.

\ssec{}    \label{discussion of groupoids}

Recall from \secref{identification of algebroids, reg} the groupoid
$\dIsom^\reg$ on $\Op^\reg$ and the corresponding Lie algebroid
$\disom^\reg$, which is the Atiyah algebroid of the principal
$\cG$-bundle $\CP_{\cG,\Op^\reg}$. For $V\in \Rep$ we will denote
by $V_{\Op^\reg}$ the corresponding (projective) module over
$\Fun(\Op^\reg)$.

Using \thmref{FF} we can transfer these objects to $\Spec(\fZ^\reg_\fg)$,
and we will denote them by $\Kdv$, $\kdv$, $\CP_{\cG,\fZ_\fg^\reg}$ and
$V_{\fZ^\reg_\fg}$, respectively. From \thmref{FF for algebroids} and 
\secref{identification of algebroids, reg}, we obtain the following

\begin{cor}  \label{FF for algebroids, reg} \hfill

\smallskip

\noindent{\em (a)}
Under the isomorphism $\fZ^\reg_\fg\simeq \Fun(\Op^\reg)$,
we have a canonical identification of Lie algebroids
$N^*_{\fZ^\reg_\fg/\fZ_\fg}\simeq \isom^\reg$.

\smallskip

\noindent{\em (b)}
We have a commutative diagram
$$
\CD
& & 0 & & 0 & & & & \\
&  & @VVV @VVV & &   \\
& & \cg_{\fZ^\reg_\fg} @>{\on{id}}>> \cg_{\fZ^\reg_\fg} @>>> 0 & & \\
& & @VVV @VVV @VVV & \\
0 @>>>  N^*_{\fZ^\reg_\fg/\fZ_\fg} @>>>
\Omega^1(\fZ_\fg)|_{\Spec(\fZ^\reg_\fg)} @>>>
\Omega^1(\fZ^\reg_\fg) @>>> 0 \\
& & @VVV @VVV @VVV & \\
 0 @>>> T(\fZ^\reg_\fg) @>>> T(\fZ_\fg)|_{\Spec(\fZ^\reg_\fg)} @>>>
 N_{\fZ^\reg_\fg/\fZ_\fg} @>>> 0 \\
 & & @VVV @VVV @VVV & \\
 & & 0 @>>> \cg_{\fZ^\reg_\fg} @>{\on{id}}>> \cg_{\fZ^\reg_\fg} & & \\
 &  & & & @VVV @VVV \\
 & & & & 0 & & 0 & &
 \endCD
 $$

\end{cor}

\ssec{}   \label{Dmod on Gr}

In what follows will work not with usual right D-modules on $\Gr_G$,
but rather with the D-modules twisted by the critical line bundle, as
in \secref{induction from O}.  We will denote the corresponding
category by $\fD_\crit(\Gr_G)\mod$.  We have the following result,
established in \cite{FG}:

\begin{thm}  \label{exactness on Grassmannian}
The functor of global sections
$\Gamma:\fD_\crit(\Gr_G)\mod\to\hg_\crit\mod$ is exact and
faithful. Moreover, it factors canonically through a functor
$\Gamma^{\ren}:\fD_\crit(\Gr_G)\mod\to U^{\ren,\reg}(\hg_\crit)\mod$,
and the latter functor is fully-faithful.
\end{thm}

Consider the category
$\on{Sph}_G:=\fD_\crit(\Gr_G)\mod^{G[[t]]}$. According to the results
of Lusztig, Drinfeld, Ginzburg and Mirkovi\'c-Vilonen (see \cite{MV}),
this is a tensor category under the convolution product, which is
equivalent to the category $\Rep$ of representations of the algebraic
group $\cG$. For $V\in \Rep$ we will denote by $\CF_V$ the
corresponding (critically twisted) $G[[t]]$-equivariant D-module on
$\Gr_G$.

\medskip

Let us recall the basic result of \cite{BD}, Sect. 5.5 and 5.6, that
describes global sections of the (critically twisted) D-modules
$\CF_V$.

\begin{thm}  \label{BD} \hfill
 
\smallskip
 
\noindent{\em (a)}
We have a canonical isomorphism of $\hg_\crit$-modules
$$\Gamma(\Gr_G,\CF_V)\simeq \BV_{\crit}\underset{\fZ^\reg_\fg}\otimes
V_{\fZ^\reg_\fg},$$ compatible with tensor product of representations.

\smallskip
 
\noindent{\em (b)} The isomorphism of (a) and that of \corref{FF for
algebroids, reg} are compatible in the sense that the
$N^*_{\fZ^\reg_\fg/\fZ_\fg}$-action on $\Hom(\BV_\crit,
\Gamma(\Gr_G,\CF_V))$, coming from \thmref{exactness on Grassmannian}
and \secref{action of renorm}, corresponds to the canonical
$\disom^\reg$-action on $V_{\Op^\reg}$.

\end{thm}

\medskip

We can take the convolution product of any D-module on $\Gr_G$ with a
spherical one. A priori, this will be a complex of D-modules on
$\Gr_G$, but as in \cite{Ga} one shows that this is a single
D-module. (Alternatively, this follows from the lemma below, using
\thmref{exactness on Grassmannian}).  Thus, we obtain an action of the
tensor category $\Rep$ on $\fD_\crit(\Gr_G)\mod$:
$$\CF,V\mapsto \CF\star \CF_V.$$ 

\begin{lem}   \label{sections of convolution}
For $\CF\in \fD_\crit(\Gr_G)\mod$ and $V\in \Rep$ we have a canonical 
isomorphism:
$$\Gamma(\Gr_G,\CF\star \CF_V)\simeq 
\Gamma(\Gr_G,\CF)\underset{\fZ_\fg^\reg}\otimes V_{\fZ^\reg_\fg}.$$
\end{lem}

\begin{proof}
Let us recall the formalism of the convolution action
(see \secref{convolution action}).  We have the functors
$$D^b(\fD_\crit(\Gr_G)\mod)\times D^b(\fD_\crit(\Gr_G)\mod)^{G[[t]]}\to
D^b(\fD_\crit(\Gr_G)\mod)$$ and 
$$D^b(\fD_\crit(\Gr_G)\mod)\times D^b(\hg_\crit\mod_\reg)^{G[[t]]}\to 
D^b(\hg_\crit\mod_\reg),$$
which are intertwined by the functor $\Gamma$.
Note that $\Gamma(\Gr_G,\delta_{1_{\Gr_G}})\simeq \BV_\crit$, and
$\Gamma(\Gr_G,\CF)\simeq \CF\star \BV_{\crit}$.
 
Hence, we have
\begin{align*}
&\Gamma(\Gr_G,\CF\star \CF_V)\simeq (\CF\star\CF_V)\star
\BV_\crit\simeq \CF\star (\CF_V\star \BV_\crit)\simeq \CF\star
\Gamma(\Gr_G,\CF_V)\simeq \\ &\simeq \CF\star
(\BV_{\crit}\underset{\fZ_\fg^\reg}\otimes V_{\fZ^\reg_\fg})\simeq
(\CF\star \BV_\crit)\underset{\fZ_\fg^\reg}\otimes
V_{\fZ^\reg_\fg}\simeq \Gamma(\Gr_G,\CF)\underset{\fZ_\fg^\reg}\otimes
V_{\fZ^\reg_\fg},
\end{align*}
where the second-to-last isomorphism is given by \thmref{BD}.

\end{proof}

\ssec{}

After these preparations we introduce the category
$\fD(\Gr_G)_\crit^{\on{Hecke}}\mod$ which is conjecturally equivalent to
$\hg_\crit\mod_\reg$.

Its objects are (critically twisted) D-modules $\CF$ on $\Gr_G$, endowed 
with an action of the algebra $\fZ_\fg^\reg$ by endomorphisms,
and a family of functorial isomorphisms
$$\alpha_V:\CF\star \CF_V\simeq V_{\fZ^\reg_\fg}
\underset{\fZ_\fg^\reg}\otimes \CF,\,\,\,\, V\in \Rep,$$
compatible with tensor products of
representations in the sense that for $U,V\in \Rep$ the diagram
$$
\CD
(\CF\star \CF_U)\star \CF_V @>>> \CF\star (\CF_U\star \CF_V) \\
@VVV   @VVV   \\
(U_{\fZ^\reg_\fg}\underset{\fZ_\fg^\reg}\otimes \CF)\star \CF_V & & 
(U_{\fZ^\reg_\fg}\underset{\fZ_\fg^\reg}\otimes V_{\fZ^\reg_\fg}) 
\underset{\Fun(\Op^\reg)}\otimes  \CF \\
@VVV @VVV \\
U_{\fZ^\reg_\fg}\underset{\fZ_\fg^\reg}\otimes (\CF\star \CF_V) @>>>
U_{\fZ^\reg_\fg}\underset{\fZ_\fg^\reg}\otimes (V_{\fZ^\reg_\fg} 
\underset{\fZ_\fg^\reg} \otimes \CF)
\endCD
$$
is commutative, and that $\alpha_V$, for $V$ being the trivial
representation, is the identity map.

In fact, one can show as in \cite{AG1} that it is sufficient to give a
family of {\it morphisms} $\{\alpha_V\}$ satisfying the above
conditions; the fact that they are isomorphisms is then
automatic. Morphisms in this category are maps of D-modules that
commute with the action of $\fZ_\fg^\reg$ and the data of $\alpha_V$.

\medskip

Note that the category $\fD(\Gr_G)_\crit^{\on{Hecke}}\mod$ is
precisely the category
$$
\fD(\Gr_G)\mod \underset{\on{pt}/\cG}\times \Op^\reg
$$
introduced above.

\medskip

Consider the groupoid $\Kdv$ and note that the algebra $\Fun(\Kdv)$
is isomorphic to
$$\underset{V\in \on{Irr}(\Rep)}\oplus \,
V_{\fZ^\reg_\fg}\underset{\BC}\otimes
V^*_{\fZ_\fg^\reg},$$
and the unit section corresponds to the map
$$V_{\fZ_\fg^\reg}\underset{\BC}\otimes V^*_{\fZ_\fg^\reg}\to
V_{\fZ_\fg^\reg}\underset{\Fun(\fZ_\fg^\reg)}\otimes
V^*_{\fZ_\fg^\reg}\to
\Fun(\fZ_\fg^\reg).$$

Let us consider the space of global sections of an object $\CF\in
\fD(\Gr_G)_\crit^{\on{Hecke}}\mod$. From \lemref{sections of convolution}
we obtain the following

\begin{lem}
For an object $\CF$ of $\fD(\Gr_G)_\crit^{\on{Hecke}}\mod$, the action of
$\fZ^\reg_\fg$ on $\Gamma(\Gr_G,\CF)$ by
$\hg_\crit$-endomorphisms, canonically extends to an action of 
$\Fun(\Kdv)$. 
\end{lem}

Consider the functor $\Gamma^{\on{Hecke}}:
\fD(\Gr_G)_\crit^{\on{Hecke}}\mod\to \hg_\crit\mod_\reg$ given by
$$\CF\mapsto \Gamma(\Gr,\CF) \underset{\Fun(\Kdv)}\otimes
\fZ_\fg^\reg,$$ where $\fZ_\fg^\reg$ is considered as a
$\Fun(\Kdv)$-algebra via the unit section.

We propose the following:

\begin{mainconj} \label{conj on Grassmann}
The above functor $\Gamma^{\on{Hecke}}$ is exact and defines an
equivalence of categories $\fD(\Gr_G)_\crit^{\on{Hecke}}\mod\to
\hg_\crit\mod_\reg$.
\end{mainconj}

Note that by definition, the category
$\fD(\Gr_G)_\crit^{\on{Hecke}}\mod$ carries a Harish-Chandra action of
$G\ppart$ at the critical level. By construction, the functor
$\Gamma^{\on{Hecke}}$ preserves this structure.  In particular, we can
consider the subcategories of $I^0$-equivariant objects on both
sides. As a consequence we obtain another conjecture:

\begin{mainconj} \label{I-monodromic on Grassmann}
The category $\Catr$ is equivalent to
$\fD(\Gr_G)_\crit^{\on{Hecke}}\mod^{I_0}$.
\end{mainconj}

\ssec{}

We now present another way of formulating \mainconjref{conj on
Grassmann}. Recall from \cite{Ga1}, that if $\CY$ is an affine
variety, $\CC$ is a $\Fun(\CY)$-linear abelian category and $\CG_\CY$
is an affine groupoid over $\CY$, it then makes sense to speak about a
lift of the $\CG_\CY$-action on $\CY$ to $\CC$.

We take $\CY=\Spec(\fZ^\reg_\fg)$, $\CG_\CY=\Kdv$ and
$\CC=\hg_\crit\mod_\reg$.  One can show that \mainconjref{conj on
Grassmann} is equivalent to the following one:

\begin{conj} \label{action of groupoid} 
The action of $\Kdv^\reg$ on $\Spec(\fZ^\reg_\fg)$ lifts to an action
on $\hg_\crit\mod_\reg$ in such a way that:

\smallskip

\noindent{\em(1)} This structure commutes in the natural
sense with the Harish-Chandra action of $G\ppart$ on 
$\hg_\crit\mod_\reg$.

\smallskip

\noindent{\em(2)}
The functor $\Gamma$ establishes an equivalence between the 
category $\fD(\Gr_G)_\crit\mod$ and the category
of $\Kdv$-equivariant objects in $\hg_\crit\mod_\reg$.
\end{conj}

\begin{remark} If we had an action of $\Kdv$ on
$\hg_\crit\mod_\reg$, as conjectured above, then at the infinitesimal
level we would have functorial maps
$$\kdv^\reg\to \on{Ext}^1_{\hg_\crit\mod_\reg}(\CM,\CM),$$ for any
$\CM\in \hg_\crit\mod_\reg$. However, the latter maps are known to
exist, as follows from \eqref{map to Ext 1} in 
\secref{action of renorm}.
\end{remark}

\ssec{}

Although we are unable to prove \mainconjref{conj on Grassmann} at the
moment, we will establish one result in its direction, which we will use
later on.

Let us denote by $\on{L}\Gamma^{\on{Hecke}}:
D^-\left(\fD(\Gr_G)_\crit^{\on{Hecke}}\mod\right) \to
D^-\left(\hg_\crit\mod_\reg\right)$ the functor given by
$$\CF\mapsto \Gamma(\Gr_G,\CF)\overset{L}{\underset{\Fun(\Kdv)}\otimes}
\fZ^\reg_\fg,$$ where $\overset{L}\otimes$ is defined using a left
resolution of $\fZ^\reg_\fg$ by projective $\Fun(\Kdv)$-modules.

One easily shows (and we will see this in the course of the proof of
the next theorem) that $\on{L}\Gamma^{\on{Hecke}}$ is in fact the left
derived functor of $\Gamma^{\on{Hecke}}$.

\begin{thm} \label{Hecke fully faithful}
The functor $\on{L}\Gamma^{\on{Hecke}}$, restricted to 
$D^b\left(\fD(\Gr_G)_\crit^{\on{Hecke}}\mod\right)$, is fully faithful.
\end{thm}

Before giving the proof, we need some preparations.

\ssec{}

Let us observe that the obvious forgetful functor
$\fD(\Gr_G)_\crit^{\on{Hecke}}\mod\to \fD(\Gr_G)_\crit\mod$ admits a left
adjoint, which we will denote by $\on{Ind}^{\on{Hecke}}$. Indeed, it is
given by $$\CF\mapsto \underset{V\in \on{Irr}(\Rep)}\oplus \, 
(\CF\star \CF_{V^*})\underset{\BC}\otimes V_{\fZ^\reg_\fg}.$$

Evidently, we have
\begin{lem}
$$\on{L}\Gamma^{\on{Hecke}}(\Gr_G,\on{Ind}^{\on{Hecke}}(\CF))\simeq
\Gamma^{\on{Hecke}}(\Gr_G,\on{Ind}^{\on{Hecke}}(\CF))\simeq
\Gamma(\Gr_G,\CF).$$
\end{lem}

Therefore, \thmref{Hecke fully faithful} implies the following:

\begin{thm} \label{Exts in reg}
For $\CF^\bullet_1,\CF^\bullet_2,\in D^b(\fD(\Gr_G)_\crit)$ and
$\CM^\bullet_i=\Gamma(\Gr_G,\CF^\bullet_i)\in D^b(\hg_\crit\mod_\reg)$
the map, given by the functor $\on{L}\Gamma^{\on{Hecke}}$,
\begin{align*}
&R\on{Hom}_{D\left(\fD(\Gr_G)_\crit\mod\right)}
\left(\CF^\bullet_1,\underset{V\in \on{Irr}(\Rep)}\oplus \,
(\CF^\bullet_2\star \CF_{V^*})\underset{\BC}\otimes
V_{\fZ^\reg_\fg}\right)\to \\
&\on{RHom}_{D(\hg_\crit\mod_\reg)}(\CM_1^\bullet,\CM_2^\bullet)
\end{align*}
is an isomorphism.
\end{thm}

{}From this theorem we obtain that all
$\on{R}^i{\Hom}_{D(\hg_\crit\mod_\reg)}
(\Gamma(\Gr_G,\CF_1^\bullet),\Gamma(\Gr_G,\CF_2^\bullet))$, viewed as
quasi-coherent sheaves on $\Spec(\fZ^\reg_\fg)$, are equivariant with
respect to the groupoid $\Kdv$.  We claim that we know a priori that
the above $\on{R}^i{\Hom}$ is acted on by the algebroid $\kdv\simeq
N^*_{\fZ^\reg_\fg/\fZ_\fg}$, and the map in \thmref{Exts in reg} is
compatible with the action of $N^*_{\fZ^\reg_\fg/\fZ_\fg}$. This
follows from \lemref{action of algebroid on exts} and \thmref{BD}.

\ssec{Proof of \thmref{Hecke fully faithful}}

It is clear that any object of $\fD(\Gr_G)_\crit^{\on{Hecke}}\mod$
admits a surjection from an object of the form
$\on{Ind}^{\on{Hecke}}(\CF)$ for some $\CF\in \fD(\Gr_G)_\crit\mod$.
Therefore, any bounded from above complex in
$\fD(\Gr_G)_\crit^{\on{Hecke}}\mod$ admits a left resolution by a
complex consisting of objects of this form.  Hence, it is sufficient
to show that for $\CF_1\in \fD(\Gr_G)_\crit\mod$ and $\CF^\bullet_2\in
D^+\left(\fD(\Gr_G)_\crit^{\on{Hecke}}\mod\right)$ the map
\begin{align} \label{map of Homs Hecke}
&\on{RHom}_{D\left(\fD(\Gr_G)^{\on{Hecke}}_\crit\mod\right)}
\left(\on{Ind}^{\on{Hecke}}(\CF_1),\CF^\bullet_2\right)\to \\
&\on{RHom}_{D(\hg_\crit\mod_\reg)}
\left(\Gamma(\Gr_G,\CF_1),\on{L}\Gamma^{\on{Hecke}}
(\Gr_G,\CF^\bullet_2)\right).
\end{align}
is an isomorphism. Note that by adjunction the LHS of the above
formula is isomorphic to
$\on{RHom}_{D\left(\fD(\Gr_G)_\crit\mod\right)}(\CF_1,\CF^\bullet_2)$,
where we regard $\CF_2^\bullet$ just as an object of
$D^+\left(\fD(\Gr_G)_\crit\mod\right)$.

\medskip

Without loss of generality we can assume that $\CF_1$ is finitely
generated, and is equivariant with respect to some congruence subgroup
$K\subset G[[t]]$.  By \secref{averaging}, we can replace
$\CF_2^\bullet$ by $\on{Av}_K(\CF_2^\bullet)$, i.e., without
restriction of generality, we can assume that $\CF_2^\bullet$ is also
$K$-equivariant.

We will use the Harish-Chandra action of $G\ppart$ on
$\fD(\Gr_G)_\crit^{\on{Hecke}}\mod$ and $\hg_\crit\mod_\reg$.  Namely,
we will interpret $\CF_1$ as $\CF_1\star \delta_{1_{\Gr_G}}\in
\fD(\Gr_G)_\crit\mod^K$, and hence
$$\on{Ind}^{\on{Hecke}}(\CF_1)\simeq \CF_1\star
(\on{Ind}^{\on{Hecke}}(\delta_{1_{\Gr_G}}))\in
\fD(\Gr_G)^{\on{Hecke}}_\crit\mod^K.$$
Similarly, $\Gamma(\Gr_G,\CF_1)\simeq \CF_1\star \BV_\crit$.

Let $\wt{\CF_1}$ be the dual D-module in
$\fD(G\ppart/K)_\crit\mod^{G[[t]]}$, see \secref{HCh adjunction}. Set
$$\CF^\bullet := \wt{\CF_1}\star \CF_2\in
D^+\left(\fD(\Gr_G)^{\on{Hecke}}_\crit\mod\right)^{G[[t]]}.$$

By \secref{HCh adjunction}, we have
$$\on{RHom}_{D\left(\fD(\Gr_G)_\crit\mod\right)}(\CF_1,\CF^\bullet_2)
\simeq \on{RHom}_{D\left(\fD(\Gr_G)_\crit\mod\right)^{G[[t]]}}
\left(\delta_{1_{\Gr_G}},\CF^\bullet\right)$$ and
$\on{RHom}_{D(\hg_\crit\mod_\reg)}
\left(\Gamma(\Gr_G,\CF_1),\on{L}\Gamma^{\on{Hecke}}
(\Gr_G,\CF^\bullet_2)\right)$ is isomorphic to
$$\on{RHom}_{D(\hg_\crit\mod_\reg)^{G[[t]]}}
\left(\BV_\crit,\on{L}\Gamma^{\on{Hecke}}(\Gr_G,\CF^\bullet)\right).$$

Evidently, we can assume that $\CF^\bullet$ is an object, denoted
$\CF$, of the abelian category $\fD(\Gr_G)^{\on{Hecke}}_\crit\mod^{G[[t]]}$.

\medskip

Since the category $\fD(\Gr_G)_\crit\mod^{G[[t]]}$ is equivalent to
$\Rep$, we obtain that the category
$\fD(\Gr_G)^{\on{Hecke}}_\crit\mod^{G[[t]]}$ is equivalent to the category
of $\fZ^\reg_\fg$-modules, with the functor being given by
$$\CL\mapsto \on{Ind}^{\on{Hecke}}(\delta_{1_{\Gr_G}})\otimes \CL.$$
Therefore, the D-module $\CF$ above has such a form for some
$\fZ^\reg_\fg$-module $\CL$.  

By a $\reg$- and $G[[t]]$-equivariant version of \propref{induction
perfect}, we can assume that $\CL$ is finitely presented. Since
$\fZ^\reg_\fg$ is a polynomial algebra, every finitely presented
module admits a finite resolution by projective ones. This reduces us
to the case when $\CL=\fZ^\reg_\fg$. Thus, we obtain that it is enough
to show the following:

\medskip

\noindent (*)\hskip1cm {\it The map
$$\on{Ext}^\bullet_{D\left(\fD(\Gr_G)_\crit\mod\right)^{G[[t]]}}
\left(\delta_{1_{\Gr_G}},\on{Ind}^{\on{Hecke}}(\delta_{1_{\Gr_G}})\right)\to 
\on{Ext}^\bullet_{D(\hg_\crit\mod_\reg)^{G[[t]]}}(\BV_\crit,\BV_\crit)$$

is an isomorphism.}

\medskip

To establish (*) we proceed as follows. It is known from \cite{ABG},
Theorem 7.6.1, that
$$\on{Ext}^\bullet_{D(\fD(\Gr)_\crit\mod)^{G[[t]]}}
\left(\delta_{1_{\Gr_G}}, \underset{V\in \on{Irr}(\Rep)}\oplus \, 
\CF_{V}\underset{\BC}\otimes V^*\right)\simeq 
\on{Sym}^{\bullet}(\cg),$$
viewed as a graded algebra with an action of $\cG$,
where the generators $\cg\subset \on{Sym}^\bullet(\cg)$
have degree $2$. 

Hence, the left-hand side in (*) is isomorphic 
to the graded algebra over $\fZ_\fg^\reg$ obtained from the $\cG$-torsor
$\CP_{\cG,\fZ^\reg_\fg}$ and the $\cG$-algebra $\on{Sym}^{\bullet}(\cg)$, i.e., 
\begin{equation}  \label{twisted symmetric algebra}
\CP_{\cG,\fZ^\reg_\fg}\overset{\cG}\times \on{Sym}^\bullet(\cg)\simeq 
\Sym^\bullet_{\fZ^\reg_\fg}(\cg_{\fZ^\reg_\fg}).
\end{equation}

Now we claim that the right-hand side in (*) is also isomorphic
to the algebra appearing in \eqref{twisted symmetric algebra}:

\begin{thm} \label{self ext of vac}
There exists a canonical isomorphism of algebras
$$\on{Ext}^\bullet_{D(\hg_\crit\mod_\reg)^{G[[t]]}}(\BV_\crit,\BV_\crit)
\simeq \Sym^\bullet_{\fZ^\reg_\fg}(\cg_{\fZ^\reg_\fg}),$$ compatible
with the action of $N^*_{\fZ^\reg_\fg/\fZ_\fg}$, where the generators
$\cg\subset \on{Sym}^\bullet(\cg)$ have degree $2$.
\end{thm}

\ssec{Proof of \thmref{self ext of vac}}

{}From \corref{FF for algebroids, reg} and \eqref{map to Ext 2} we
obtain a map
\begin{equation} \label{map on generators}
\cg_{\fZ^\reg_\fg}\to
\on{Ext}^2_{D(\hg_\crit\mod_\reg)^{G[[t]]}}(\BV_\crit,\BV_\crit),
\end{equation}
compatible with the action of $N^*_{\fZ^\reg_\fg/\fZ_\fg}$.  We are
going to show that this map induces the isomorphism stated in the
theorem.  We will do it by analyzing the spectral sequence of
\secref{equivariant restriction}.

\medskip

Since the $\hg_\crit$-action on $\BV_\crit$ can be canonically
extended to an action of $U^{\ren,\reg}(\hg_\crit)$, from
\lemref{action of algebroid on exts} and \eqref{map to Ext 1}, we
obtain a map
$$\Omega^1(\fZ^\reg_\fg)\to
\on{Ext}^1_{D(\hg_\crit\mod)^{G[[t]]}}(\BV_\crit,\BV_\crit).$$
We will use the following result of \cite{FT}:
\begin{thm} \label{vacuum FT}
The cup-product induces an isomorphism of algebras
$$\Omega^\bullet(\fZ^\reg_\fg)\to
\on{Ext}^\bullet_{D(\hg_\crit\mod)^{G[[t]]}}(\BV_\crit,\BV_\crit).$$
\end{thm}

\medskip

Recall that $\imath^\reg$ denotes the embedding
$\Spec(\fZ_\fg^\reg)\hookrightarrow \Spec(\fZ_\fg)$. Consider the
object $\imath^\reg{}^!(\BV_\crit)\in D(\hg_\crit\mod_\reg)^{G[[t]]}$,
see \secref{equivariant restriction}. By {\it loc.cit.}, the $j$-th
cohomology of this complex is isomorphic to
$\BV_\crit\underset{\fZ_\fg^\reg}\otimes
\Lambda^j(N_{\fZ^\reg_\fg/\fZ_\fg})$.

Consider the cohomological truncation of $\imath^!(\BV_\crit)$,
leaving the segment in the cohomological degrees $j$ and $j+1$. It
gives rise to a map in the derived category
\begin{equation} \label{2 step complex}
\phi_j:\BV_\crit\underset{\fZ^\reg_\fg}\otimes
\Lambda^{j+1}(N_{\fZ^\reg_\fg/\fZ_\fg})\to
\BV_\crit\underset{\fZ^\reg_\fg}\otimes
\Lambda^{j}(N_{\fZ^\reg_\fg/\fZ_\fg})[2].
\end{equation}

\begin{lem} \label{multiplicativity of diff}
The map $\phi_j$ equals the composition
$$\BV_\crit\underset{\fZ^\reg_\fg}\otimes
\Lambda^{j+1}(N_{\fZ^\reg_\fg/\fZ_\fg})\to
\BV_\crit\underset{\fZ^\reg_\fg}\otimes N_{\fZ^\reg_\fg/\fZ_\fg}
\underset{\fZ^\reg_\fg}\otimes
\Lambda^{j}(N_{\fZ^\reg_\fg/\fZ_\fg})\overset{\phi_1\otimes \on{id}}
\to\BV_\crit\underset{\fZ^\reg_\fg}\otimes
\Lambda^{j}(N_{\fZ^\reg_\fg/\fZ_\fg})[2].$$
\end{lem}

\medskip

By \secref{equivariant restriction}, we obtain a spectral sequence,
converging to
$\on{Ext}^\bullet_{D(\hg_\crit\mod)^{G[[t]]}}(\BV_\crit,\BV_\crit)$,
whose second term is given by
$$E_2^{i,j}=\on{Ext}^i_{D(\hg_\crit\mod_\reg)^{G[[t]]}}(\BV_\crit,\BV_\crit)
\underset{\fZ^\reg_\fg}\otimes \Lambda^j(N_{\fZ^\reg_\fg/\fZ_\fg}).$$

Note also that by \lemref{multiplicativity of diff}, the differential
in the above spectral sequence, which maps $E_2^{i-2,j+1}\to
E_2^{i,j}$ can be expressed through the case when $j=0$ as

\begin{align*}
&E_2^{i-2,j+1}\simeq E_2^{i-2,0}\underset{\fZ^\reg_\fg}\otimes
\Lambda^{j+1}(N_{\fZ^\reg_\fg/\fZ_\fg}) \to
E_2^{i-2,0}\underset{\fZ^\reg_\fg}\otimes N_{\fZ^\reg_\fg/\fZ_\fg}
\underset{\fZ^\reg_\fg}\otimes
\Lambda^j(N_{\fZ^\reg_\fg/\fZ_\fg})\simeq \\ &\simeq
E_2^{i-2,1}\underset{\fZ^\reg_\fg}\otimes
\Lambda^{j}(N_{\fZ^\reg_\fg/\fZ_\fg})\to
E_2^{i,0}\underset{\fZ^\reg_\fg}\otimes
\Lambda^{j}(N_{\fZ^\reg_\fg/\fZ_\fg})\simeq E_2^{i,j}.
\end{align*}

Let us observe that the canonical map 
\begin{equation} \label{inj in spec}
\on{Ext}^j_{D(\hg_\crit\mod)^{G[[t]]}}(\BV_\crit,\BV_\crit)\to E_2^{0,j}
\end{equation}
identifies by construction with the map $\Omega^j(\fZ^\reg_\fg)\to 
\Lambda^j(N_{\fZ^\reg_\fg/\fZ_\fg})$ coming from
of \eqref{map to Ext 1}; in particular, it is injective.

\medskip

We will prove by induction on $i=1,2,...$ the following statements:

\begin{itemize}

\item(i) $E_2^{2i-1,0}=0$,

\smallskip

\item(ii) $E_2^{2i,0}\simeq \on{Sym}^i(\cg_{\fZ^\reg_\fg})$ such that
the differential $E_2^{2i-2,1}\to E_2^{2i,0}$ identifies with the map
$\on{Sym}^{i-1}(\cg_{\fZ^\reg_\fg})\underset{\fZ^\reg_\fg}\otimes
N_{\fZ^\reg_\fg/\fZ_\fg}\to \on{Sym}^{i}(\cg_{\fZ^\reg_\fg})$.

\end{itemize}

Note that item (i) above implies that $E_2^{2i-1,j}=0$ for any $j$ and
that item (ii) implies that $E_2^{2i,j}\simeq
\on{Sym}^i(\cg_{\fZ^\reg_\fg}) \underset{\fZ^\reg_\fg}\otimes
\Lambda^j(N_{\fZ^\reg_\fg/\fZ_\fg})$ such that the differential
identifies with the Koszul differential
$\on{Sym}^{i-1}(\cg_{\fZ^\reg_\fg})\underset{\fZ^\reg_\fg}\otimes
\Lambda^j(N_{\fZ^\reg_\fg/\fZ_\fg})\to
\on{Sym}^{i}(\cg_{\fZ^\reg_\fg})\underset{\fZ^\reg_\fg}\otimes
\Lambda^{j-1}(N_{\fZ^\reg_\fg/\fZ_\fg})$.

\medskip

Consider first the base of the induction, i.e., the case $i=1$. In
this case we know a priori that $E_2^{1,0}=0$.  We obtain that
$\Omega^1(\fZ^\reg_\fg)$ maps isomorphically onto the kernel of the
map $N_{\fZ^\reg_\fg/\fZ_\fg}\simeq E_2^{0,1}\to E_2^{2,0}$. In
particular, the map of \eqref{map on generators}
$\cg_{\fZ^\reg_\fg}\hookrightarrow E_2^{2,0}$ is injective. We claim
that the latter map is surjective as well. Indeed, if it were not, the
map in \eqref{inj in spec} would not be injective for $j=2$.

Hence, the differential $E_2^{0,j}\to E_2^{2,j-1}$ does identify with
the corresponding term of the Koszul differential. In particular,
$\on{Ext}^i_{D(\hg_\crit\mod)^{G[[t]]}}(\BV_\crit,\BV_\crit)$ maps
isomorphically to $E_3^{0,j}\simeq \on{ker}(E_2^{0,i}\to
E_2^{2,i-1})$. This implies, in particular, that all the higher
differentials $E_k^{0,j}\to E_k^{k,j-k-1}$ for $k\geq 3$ vanish.

\medskip

Let us now perform the induction step. Observe that the induction
hypothesis, all the terms of the spectral sequence $E_k^{i',j}$ for
$0<i'\leq 2i-2$ with $k\geq 3$ vanish. Therefore, the term
$E_2^{2i+1,0}$ injects into
$$\on{ker}
\left(\on{Ext}^{2i+1}_{D(\hg_\crit\mod)^{G[[t]]}}(\BV_\crit,\BV_\crit)
\to E_2^{0,2i+1}\right),$$ and as the latter map is injective, we
obtain that $E_2^{2i+1,0}=0$.

By a similar argument we obtain that $E_3^{2i,1}=E_3^{2i+2,0}=0$.
Hence, $E_2^{2i+2,0}$ identifies with 
$$\on{coker}\left(\on{Sym}^{i-1}(\cg_{\fZ^\reg_\fg})
\underset{\fZ^\reg_\fg}\otimes \Lambda^2(N_{\fZ^\reg_\fg/\fZ_\fg})\to
\on{Sym}^i(\cg_{\fZ^\reg_\fg})
\underset{\fZ^\reg_\fg}\otimes N_{\fZ^\reg_\fg/\fZ_\fg}\right)\simeq \on{Sym}^{i+1}(\cg_{\fZ^\reg_\fg}).$$

\medskip

To finish the proof of theorem it remains to remark that, by
construction, the cup-product map
\begin{align*}
&\on{Ext}^{2}_{D(\hg_\crit\mod_\reg)^{G[[t]]}}(\BV_\crit,\BV_\crit)\otimes 
\on{Ext}^{2i}_{D(\hg_\crit\mod_\reg)^{G[[t]]}}(\BV_\crit,\BV_\crit)\to \\
&\on{Ext}^{2i+2}_{D(\hg_\crit\mod_\reg)^{G[[t]]}}(\BV_\crit,\BV_\crit)
\end{align*}
identifies with the multiplication map $\cg_{\fZ^\reg_\fg}\otimes
\on{Sym}^i(\cg_{\fZ^\reg_\fg})\to \on{Sym}^{i+1}(\cg_{\fZ^\reg_\fg})$.

\ssec{}

Thus, the two graded algebras appearing in (*) are abstractly
isomorphic to one-another.  It remains to see that the existing map
indeed induces an isomorphism. Since both algebras are freely
generated by their degree $2$ part, it is sufficient to show that the
map 
\begin{align} \label{degree 2}
&\cg_{\fZ^\reg_\fg}\simeq
\on{Ext}^2_{D\left(\fD(\Gr_G)_\crit\mod\right)^{G[[t]]}}
\left(\delta_{1_{\Gr_G}},\underset{V\in \on{Irr}(\Rep)}\oplus \,
\CF_{V}\underset{\BC}\otimes V^*_{\fZ^\reg_\fg}\right)\to \\
&\on{Ext}^2_{D(\hg_\crit\mod_\reg)^{G[[t]]}}(\BV_\crit,\BV_\crit)\simeq
\cg_{\fZ^\reg_\fg}
\end{align}
is an isomorphism. Since the map of \thmref{Exts in reg} is compatible
with the action of the algebroid $N^*_{\fZ^\reg_\fg/\fZ_\fg}$, and
since $\cg_{\fZ^\reg_\fg}$ is irreducible as a $N^*_{\fZ^\reg_\fg/\fZ_\fg}$-module,
if the map in \eqref{degree 2} were not an
isomorphism, it would be zero. We claim that this leads to a
contradiction:

\medskip

Consider the canonical maps of $H^\bullet(\on{pt}/G)\simeq
H^\bullet_{G[[t]]}(\on{pt})$ to both the LHS and RHS of (*).  Note
that we have a canonical identification
$$H^\bullet(\on{pt}/G)\simeq \on{Sym}^\bullet(\fh^*)^W\simeq  
\on{Sym}^\bullet(\check\fh)^W\simeq
\on{Sym}^\bullet(\cg)^{\cG}.$$ By the construction of the isomorphism 
in \cite{ABG}, Theorem 7.6.1, $$H^\bullet(\on{pt}/G)\to 
\on{Ext}^\bullet_{D(\fD(\Gr)_\crit\mod)^{G[[t]]}}\left(\delta_{1_{\Gr}},
\on{Ind}^{\on{Hecke}}(\delta_{1_{\Gr}})\right)$$
corresponds to the canonical embedding 
$\on{Sym}^\bullet(\cg)^{\cG}\to \on{Sym}^{\bullet}(\cg)$.
Therefore, if the map of (*) was $0$ on the generators, it would also
annihilate the augmentation ideal in $H^\bullet(\on{pt}/G)$. However, we
have the following assertion:

\begin{thm}   \label{equivariant cohomology}
The map $$H^\bullet(\on{pt}/G)\simeq H^\bullet_{G[[t]]}(\on{pt})\to
\on{Ext}^\bullet_{D(\hg_\crit\mod_\reg)^{G[[t]]}}(\BV_\crit,\BV_\crit)$$
corresponds under the isomorphism of \thmref{self ext of vac} to the
map
$$H^\bullet(\on{pt}/G)\simeq
\on{Sym}^\bullet(\cg)^{\cG}\overset{\tau}\simeq \on{Sym}^\bullet(\cg)^{\cG}
\to \CP_{\cG,\fZ^\reg_\fg}\overset{\cG}\times \on{Sym}^\bullet(\cg)
\simeq \Sym^\bullet(\cg_{\fZ^\reg_\fg}),$$
where $\tau$ is as in \secref{tau !}.
\end{thm}

The proof of this theorem will be given in the next section.

\section{A manipulation with equivariant cohomology: proof of \thmref{equivariant cohomology}}
\label{manipulation}

\ssec{}    

We will consider the algebra of self-Exts of $\BV_\crit$ in a category
bigger than $\hg_\crit\mod_\reg$, namely, in the category
$\hg_\crit\mod_\nilp$.

Let $\CP_{\cB,\Op^\nilp}$ be the canonical $\cB$-torsor on the scheme
$\Op^\nilp$, and let $\CP_{\cB,\Op^\reg}$ be its restriction to $\Op^\reg$.
We will denote by $\CP_{\cB,\fZ_\fg^\nilp}$ and $\CP_{\cB,\fZ_\fg^\reg}$
the corresponding $\cB$-torsors on $\Spec(\fZ^\nilp_\fg)$ and
$\Spec(\fZ^\reg_\fg)$, respectively. If $V$ is a representation of $\cB$
(in practice we will take $V=\cb,\cn,\cg/\cn$, etc.), we will denote by
$V_{\fZ^\nilp_\fg}$, $V_{\fZ^\reg_\fg}$ the corresponding modules
over $\fZ^\nilp_\fg$ and $\fZ^\reg_\fg$, respectively.

Recall that by \corref{nilpotent directions}, the image
of the normal $N_{\Op^\reg/\Op^\nilp}$ in the quotient
$$N_{\fZ^\reg_\fg/\fZ_\fg}/\Omega^1(\fZ^\reg_\fg)$$ identifies with
$\cn_{\fZ^\reg_\fg}\subset \cg_{\fZ^\reg_\fg}$. From the proof of
\thmref{self ext of vac} we obtain the following statement.

\begin{lem}
The natural map
$$\on{Ext}^\bullet_{D(\hg_\crit\mod_\reg)^{G[[t]]}}(\BV_\crit,\BV_\crit)\to
\on{Ext}^\bullet_{D(\hg_\crit\mod_\nilp)^{G[[t]]}}(\BV_\crit,\BV_\crit)$$
induces an isomorphism
$$\Sym^\bullet\left((\cg/\cn)_{\fZ^\reg_\fg}\right)\simeq
\on{Ext}^\bullet_{D(\hg_\crit\mod_\nilp)^{G[[t]]}}(\BV_\crit,\BV_\crit).$$
\end{lem}

By the equivariance of the map in (*) with respect to the algebroid
$N^*_{\fZ^\reg_\fg/\fZ_\fg}$, the image of $H^\bullet(\on{pt}/G)$ in
$\Sym^\bullet_{\fZ^\reg_\fg}(\cg_{\fZ^\reg_\fg})$ is {\it a priori}
contained in the subalgebra $\on{Sym}^\bullet(\cg)^{\cG}$. Hence, it
is sufficient to show that the composition
$$\Sym^\bullet(\check\fh)^W\simeq H^\bullet(\on{pt}/G)\to
\on{Ext}^\bullet_{D(\hg_\crit\mod_\nilp)^{G[[t]]}}(\BV_\crit,\BV_\crit)\simeq
\Sym^\bullet\left((\cg/\cn)_{\fZ^\reg_\fg}\right)$$
equals the natural map
$$\Sym^\bullet(\check\fh)^W\overset{\tau}\simeq
\Sym^\bullet(\check\fh)^W\to \Sym^\bullet(\check\fh)\to
\Sym^\bullet\left((\cb/\cn)_{\fZ^\reg_\fg}\right)\hookrightarrow
\Sym^\bullet\left((\cg/\cn)_{\fZ^\reg_\fg}\right).$$

\ssec{}

Consider now the module $\BM_0\in \hg_\crit\mod_\nilp^I$. Since
$\on{Av}_{G[[t]]/I}(\BM_0)\simeq \BV_\crit$, by \secref{proper conv},
we obtain an isomorphism
$$\on{RHom}_{D(\hg_\crit\mod_\nilp)^{G[[t]]}}(\BV_\crit,\BV_\crit)\simeq
\on{RHom}_{D(\hg_\crit\mod_\nilp)^I}(\BM_0,\BV_\crit).$$

It is easy to see that the composition
$$H^\bullet(\on{pt}/G[[t]])\to
\on{Ext}^\bullet_{D(\hg_\crit\mod_\nilp)^{G[[t]]}}(\BV_\crit,\BV_\crit)\to
\on{Ext}^\bullet_{D(\hg_\crit\mod_\nilp)^I}(\BM_0,\BV_\crit)$$ equals
the map
$$H^\bullet(\on{pt}/G[[t]])\to H^\bullet(\on{pt}/I)\to
\on{Ext}^\bullet_{D(\hg_\crit\mod_\nilp)^I}(\BM_0,\BM_0)\to
\on{Ext}^\bullet_{D(\hg_\crit\mod_\nilp)^I}(\BM_0,\BV_\crit).$$

By \corref{Verma is flat}, the module $\BM_0$ is flat over
$\fZ^\nilp_\fg$. Hence, by \lemref{* res}
$$\on{RHom}^\bullet_{D(\hg_\crit\mod_\nilp)^I}(\BM_0,\CM)\simeq
\on{RHom}^\bullet_{D(\hg_\crit\mod_\reg)^I}(\BM_{0,\reg},\CM)$$ for
any $\CM\in \hg_\crit\mod_\reg^I$, where 
$\BM_{0,\reg}:=\BM_0\underset{\fZ^\nilp_\fg} \otimes \fZ^\reg_\fg$.

Moreover, the map $$H^\bullet(\on{pt}/I)\to
\on{Ext}^\bullet_{D(\hg_\crit\mod_\nilp)^I}(\BM_0,\BM_0)\to
\on{Ext}^\bullet_{D(\hg_\crit\mod_\nilp)^I}(\BM_0,\BV_\crit)$$
that appears above, equals the map
\begin{align*}
&H^\bullet(\on{pt}/I) \to
\on{Ext}^\bullet_{D(\hg_\crit\mod_\reg)^I}(\BM_{0,\reg},\BM_{0,\reg})\to
\\
&\to\on{Ext}^\bullet_{D(\hg_\crit\mod_\reg)^I}(\BM_{0,\reg},\BV_\crit)
\simeq \on{Ext}^\bullet_{D(\hg_\crit\mod_\nilp)^I}(\BM_0,\BV_\crit).
\end{align*}

Thus, we obtain a commutative diagram
$$
\CD
H^\bullet(\on{pt}/G[[t]]) @>>> H^\bullet(\on{pt}/I) \\
@VVV   @VVV  \\
\Sym^\bullet\left((\cg/\cn)_{\fZ^\reg_\fg}\right) @>{\sim}>>
\on{Ext}^\bullet_{D(\hg_\crit\mod_\reg)^I}(\BM_{0,\reg},\BV_\crit),
\endCD
$$
and it is easy to see that the resulting map
$$\Sym^\bullet(\check\fh)\simeq H^\bullet(\on{pt}/I)\to 
\Sym^\bullet\left((\cg/\cn)_{\fZ^\reg_\fg}\right)$$ is
a homomorphism of algebras.

Therefore,  it suffices to show that the map
\begin{align*} \label{H 2 map}
&\check\fh\simeq H^2(\on{pt}/I)\to
\on{Ext}^2_{D(\hg_\crit\mod_\reg)^I}(\BM_{0,\reg},\BM_{0,\reg})\to \\
&\to \on{Ext}^2_{D(\hg_\crit\mod_\reg)^I}(\BM_{0,\reg},\BV_\crit)
\simeq (\cg/\cn)_{\fZ^\reg_\fg}
\end{align*}
equals the negative of the tautological map $\check\fh\to  (\cg/\cn)_{\fZ^\reg_\fg}$.

\medskip

Let $\CM$ be an arbitrary $I$-equivariant object of
$\hg_\crit\mod_\reg$. One easily establishes the following
compatibility of spectral sequences:

\begin{lem}
The composition
$$\on{Ext}^1_{\hg_\crit\mod}(\CM,\CM)\to \on{Hom}(\CM,\CM)
\underset{\fZ^\reg_\fg}\otimes N_{\fZ^\reg_\fg/\fZ_\fg}
\to \on{Ext}^2_{D(\hg_\crit\mod_\reg)^I}(\CM,\CM)$$
equals the composition
$$\on{Ext}^1_{\hg_\crit\mod}(\CM,\CM)\to \on{Hom}(\CM,\CM)\otimes
H^2(\on{pt}/I) \to \on{Ext}^2_{D(\hg_\crit\mod_\reg)^I}(\CM,\CM).$$
\end{lem}

Therefore, to complete the proof of \thmref{equivariant cohomology},
it is sufficient to construct a map $\fh^*\to
\on{Ext}^1_{\hg_\crit\mod}(\BM_{0,\reg},\BM_{0,\reg})$ such that the
composition
$$\fh^*\to \on{Ext}^1_{\hg_\crit\mod}(\BM_{0,\reg},\BM_{0,\reg})\to
\on{Hom}(\BM_{0,\reg},\BM_{0,\reg})\otimes H^2(\on{pt}/I)$$ comes from
the natural isomorphism $\fh^*\to H^2(\on{pt}/I)$, and the composition
\begin{align*}
&\check\fh\simeq \fh^*\to
\on{Ext}^1_{\hg_\crit\mod}(\BM_{0,\reg},\BM_{0,\reg})\to
\on{Hom}(\BM_{0,\reg},\BM_{0,\reg})\underset{\fZ^\reg_\fg}\otimes
N_{\fZ^\reg_\fg/\fZ_\fg}\to \\ &\to
\on{Hom}(\BM_{0,\reg},\BM_{0,\reg}) \underset{\fZ^\reg_\fg}\otimes
(\cg/\cn)_{\fZ^\reg_\fg}
\end{align*}
equals the negative of the embedding $\check\fh\to (\cb/\cn)_{\fZ^\reg_\fg}$.

\ssec{}

The required map $\fh^*\to
\on{Ext}^1_{\hg_\crit\mod}(\BM_{0,\reg},\BM_{0,\reg})$ is constructed
as follows. By deforming the highest weight, we obtain the "universal"
Verma module $U(\fg)\underset{U(\fn)}\otimes \BC=:M_{\on{univ}}$, and
the corresponding induced module $\BM_{\on{univ}}$ over
$\hg_\crit$. In particular, we have a map $\fh^*\to
\on{Ext}^1(\BM_0,\BM_0)$.

Evidently, the composition $\fh^*\to \on{Ext}^1(\BM_0,\BM_0)\to
\on{Ext}^1(\BM_0,\BM_{0,\reg})$ factors canonically through
$\on{Ext}^1_{\hg_\crit\mod}(\BM_{0,\reg},\BM_{0,\reg})$.

The fact that the composition $\fh^*\to
\on{Ext}^1_{\hg_\crit\mod}(\BM_{0,\reg},\BM_{0,\reg})\to
\on{Hom}(\BM_{0,\reg},\BM_{0,\reg})\otimes H^2(\on{pt}/I)$ comes from
$\fh^*\to H^2(\on{pt}/I)$ follows from the corresponding property of
the composition $\fh^*\to \on{Ext}^1_{\fg\mod}(M_0,M_0)\to
\on{Hom}(M_0,M_0)\otimes H^2(\on{pt}/I)$.

\medskip

Consider the composition $\check\fh\simeq \fh^*\to
\on{Hom}(\BM_{0,\reg},\BM_{0,\reg}) \underset{\fZ^\reg_\fg}\otimes
(\cg/\cn)_{\fZ^\reg_\fg}$. This map is equivariant with respect to the
group $\on{Aut}(\D)$. In particular, if we choose a coordinate on
$\D$, the above map has degree $0$ with respect to the action of
$\BG_m$ by loop rotations. Since $\check\fh$ equals the degree $0$
subspace of
$\on{Hom}(\BM_{0,\reg},\BM_{0,\reg})\underset{\fZ^\reg_\fg}\otimes
(\cg/\cn)_{\fZ^\reg_\fg}$ (see \secref{identification of algebroids,
grading}), we obtain that the map in question factors through 
{\it some} map $\check\fh\to \check\fh$.

To prove that the latter map is in fact the negative of the identity, we proceed as
follows.  By \secref{nilpotent opers},
we have an identification $\fh^*\otimes \fZ^\nilp_\fg\simeq
N_{\fZ^\nilp_\fg/\fZ^{\on{RS}}_\fg}$.  Moreover, by \lemref{h-directions} and \propref{discrepancy}
the composition
$$\fh^*\otimes \fZ^\reg_\fg\simeq
N_{\fZ^\nilp_\fg/\fZ^{\on{RS}}_\fg}|_{\Spec(\fZ^\reg_\fg)}\to
N_{\fZ^\nilp_\fg/\fZ_\fg}|_{\Spec(\fZ^\reg_\fg)}\to
(\cg/\cn)_{\fZ^\reg_\fg},$$ maps identically onto $\fh^*\otimes
\fZ^\reg_\fg\subset (\cg/\cn)_{\fZ^\reg_\fg}$.

Now, our assertion follows from the fact, that the map
$$\fh^*\otimes \fZ^\nilp_\fg\to \on{Ext}^1(\BM_0,\BM_0)\to
\on{Hom}(\BM_0,\BM_0) \underset{\fZ^\nilp_\fg}\otimes
N_{\fZ^\nilp_\fg/\fZ_\fg}$$ equals the negative of
$$\fh^*\otimes \fZ^\nilp_\fg\simeq
N_{\fZ^\nilp_\fg/\fZ^{\on{RS}}_\fg} \overset{1\otimes \on{id}}\longrightarrow
\on{Hom}(\BM_0,\BM_0) \underset{\fZ^\nilp_\fg}\otimes
N_{\fZ^\nilp_\fg/\fZ^\RS_\fg}\to 
\on{Hom}(\BM_0,\BM_0) \underset{\fZ^\nilp_\fg}\otimes
N_{\fZ^\nilp_\fg/\fZ_\fg}$$ by \propref{discrepancy}.

\newpage

\vspace*{10mm}

{\Large \part{Wakimoto modules}}

\vspace*{10mm}


In this Part we review the Wakimoto modules which were introduced for
an arbitrary affine Kac-Moody algebra $\ghat$ in
\cite{FF:usp,FF:si,F:wak} following the work of Wakimoto \cite{Wak} in
the case of $\su$. On the intuitive level, Wakimoto modules are
sections of certain D-modules on the Iwahori orbits on the
semi-infinite flag manifold $G\ppart/B\ppart$. The construction of
\cite{FF:usp,FF:si,F:wak} may be phrased in terms of a kind of
semi-infinite induction functor, as we explain below. This approach to
the Wakimoto modules is similar to the one discussed in
\cite{Ar,Vor,GMS}. It uses the formalism of chiral algebras, and in
particular, the chiral algebra of differential operators on the group
$G$. It also uses the language of semi-infinite cohomology, which was
introduced by Feigin \cite{Fei} and, in the setting of chiral
algebras, by Beilinson and Drinfeld \cite{CHA}.

Let $\oG$ be the big cell $B\cdot w_0\cdot B\subset G$, and for an
arbitrary level $\kappa$ we consider the chiral algebra
$\fD^{\on{ch}}(\oG)_\kappa$ of chiral differential operators on it.  In
\secref{bose} we define the chiral algebra $\fD^{\on{ch}}(\oGN)_\kappa$ as
a BRST reduction of $\fD^{\on{ch}}(\oG)_\kappa$ with respect to
$\fn\ppart$. This chiral algebra can be thought of as governing
D-modules on the big cell in $G\ppart/N\ppart$; we show that the
natural homomorphism to it from the chiral algebra, corresponding to
the Kac-Moody Lie algebra, $\hg_\kappa$ coincides with the free field
realization homomorphism of Feigin and Frenkel.

By construction, any chiral module over $\fD^{\on{ch}}(\oGN)_\kappa$
is a bi-module over $\hg_\kappa$ and the Heisenberg algebra
$\wh\fh_{-\kappa+\kappa_\crit}$. In \secref{sec Wak} for any such
module, we define the induction functor from the category
$\wh\fh_{\kappa-\kappa_\crit}\mod$ to $\hg_\kappa\mod$. The resulting
$\hg_\kappa$-modules are by definition the Wakimoto modules. Thus,
Wakimoto modules can be viewed as induced from
$\wh\fh_{\kappa-\kappa_\crit}$ to $\hg_\kappa$,  
using certain bi-modules.

In \secref{conv of wak} we study cohomological properties of Wakimoto
modules and, in particular, their behavior with respect to the
convolution functors. The crucial result that we need below is
\propref{lambda equivariance} that states that Wakimoto modules are
essentially invariant under convolution with "lattice" elements in the
Iwahori-Hecke algebra.
  
In \secref{Wak crit} we specialize to the case
$\kappa=\kappa_\crit$. The crucial result here, due to
\cite{F:wak}, is that certain Wakimoto modules are isomorphic to Verma
modules over $\hg_\crit$. This fact will allow us to obtain
information about the structure of Verma modules that will be used in
the subsequent sections.

\bigskip

\section{Free field realization}    \label{bose}

In what follows we will use the language of chiral algebras on a curve
$X$, developed in \cite{CHA}. We will fix a point $x\in X$ and
identify $\fD_X$-modules supported at this point with underlying
vector spaces. We will identify the formal disc $\D$ with the formal
neighborhood of $x$ in $X$.

\ssec{}   \label{Tate extension of b}

Let $L$ be a Lie-* algebra, which we assume to be projective and
finitely generated as a $\fD_X$-module. Recall that there exists a
canonical Tate central extension of $L$, which is a Lie-* algebra
$\wh{L}^{\on{Tate}}$
$$0\to \omega_X\to \wh{L}^{\on{Tate}}\to L\to 0,$$ see \cite{CHA},
Sect. 2.7. The key property of $\wh{L}^{\on{Tate}}$ is that if $\CM$ is a
chiral module over $\wh{L}^{-\on{Tate}}$ (here "$-$" signifies the Baer
negative central extension), then we have a well-defined complex of
$\fD_X$-modules, denoted $\fC^\semiinf(L,\CM)$, which we will refer to
as the semi-infinite complex of $\CM$ with respect to $L$. We will
denote by $H^\semiinf(L,\CM)$ (resp., $H^\semiinfi(L,\CM)$) the $0$-th
(resp., $i$-th) cohomology of this complex.

If $\CM$ is supported at the point $x\in X$, by definition,
$\fC^\semiinf(L,\CM)$ is given by the semi-infinite complex of the
Tate Lie algebra $H^0_{DR}(\D^\times,L)$ with respect to the lattice
$H^0_{DR}(\D,L)\subset H^0_{DR}(\D^\times,L)$.

If $\CA$ is a chiral algebra with a homomorphism $\wh{L}^{-\on{Tate}}\to \CA$, then 
$\fC^\semiinf(L,\CA)$ has a natural structure of a DG chiral algebra

\medskip

Let now $L'$ and $L''$ be two central extensions of $L$ by $\omega_X$,
whose Baer sum is identified with $\wh{L}^{-\on{Tate}}$, and let $\CM$
and $\CM'$ be $L'$- and $L''$-modules, respectively. Then $\CM\otimes
\CM'$ is a module over $\wh{L}^{-\on{Tate}}$, and in this case we will
use the notation
$$\CM\underset{L}\torsemiinf\CM' \text{ or }
\CM\underset{H^0_{DR}(\D^\times,L),H^0_{DR}(\D,L)}\torsemiinf\CM' $$
instead of $\fC^\semiinf(L,\CM\otimes \CM')$. If the latter is acyclic
away from cohomological degree $0$ we will denote by the same symbol
the corresponding $0$-th cohomology.

Finally, let $\fh$ be a finite-dimensional subspace in
$H^0_{DR}(X,L)$.  In this case, $\fC^\semiinf(L,\CM)$ admits a
subcomplex $\fC^\semiinf(L;\fh,\CM)$ of relative cochains. We will
sometimes also use the notation $\fC^\semiinf(L;\fh,\cdot)$ and
$\cdot\underset{L;\fh}\torsemiinf\cdot$.

\ssec{}

Let $L_\fg$, $L_\fb$ and $L_\fn$ be the Lie-* algebras, corresponding
to the Lie algebras $\fg$, $\fb$ and $\fn$, respectively. For a level
$\kappa$, we will denote by $L_{\fg,\kappa}$ the corresponding
Kac-Moody extension of $L_\fg$ by $\omega_X$, and by $L_{\fb,\kappa}$
the induced central extension of $L_\fb$.  Let $\wh{L}^{\on{Tate}}_\fb$ be
the Tate extension of $L_\fb$, and let $\wh{L}'_{\fb,\kappa}$ be the
Baer sum of $\wh{L}^{\on{Tate}}_\fb$ and $L_{\fb,\kappa'}$, where
$\kappa'=-\kappa-2\kappa_\crit$; let $\wh{L}_{\fb,\kappa}$ be the Baer
negative of $\wh{L}'_{\fb,\kappa}$. 

Since $\kappa'|_\fn=0$, the extension induced by $L_{\fb,\kappa'}$ on
$L_\fn$ is canonically trivialized. The extension induced by
$\wh{L}^{\on{Tate}}_\fb$ is also canonically trivialized, since $\fn$ is
nilpotent. Hence, $\wh{L}_{\fb,\kappa}$ comes from a well-defined
central extension $\wh{L}_{\fh,\kappa}$ of the commutative Lie-*
algebra $L_\fh$. We will denote by $\wh{L}'_{\fh,\kappa}$ the Baer
negative of $\wh{L}_{\fh,\kappa}$.

Note that when $\kappa$ is integral, the above central extensions of
Lie algebras $H^0_{DR}(\D^\times,?)$ all come from the corresponding
central extensions of loop groups.

We will denote by $\fH_\kappa$ (resp., $\fH'_\kappa$) the reduced 
universal enveloping chiral algebra of $\wh{L}_{\fh,\kappa}$ 
(resp., $\wh{L}'_{\fh,\kappa}$).
We will denote by $\CA_{\fg,\kappa}$ the reduced universal enveloping 
chiral algebra of $L_{\fg,\kappa}$.

\medskip

Let $\CM$ be a chiral $L_{\fb,\kappa'}$-module.
Since the Tate extension of $L_\fb$, induced by the adjoint action
equals the extension induced by the adjoint action on $L_\fn$, the
complex $\fC^\semiinf(L_\fn,\CM)$ carries a chiral action of
$\wh{L}'_{\fb,\kappa}$. The resulting action of $L_\fn\subset
\wh{L}'_{\fb,\kappa}$ on the individual semi-infinite cohomologies
$H^\semiinfi(L_\fn,\CM)$ is trivial. Hence, we obtain that each
$H^\semiinfi(L_\fn,\CM)$ is a chiral $\fH'_\kappa$-module. If $\CR$ is
an $\fH_\kappa$-module, regarded as a $\wh{L}_{\fb,\kappa}$-module,
$\fC^\semiinf(L_\fb,\CM\otimes \CR)$ makes sense. If we suppose,
moreover, that $\fC^\semiinf(L_\fn,\CM)$ is acyclic away from degree
$0$, then
$$\fC^\semiinf(L_\fb,\CM\otimes \CR)\simeq
H^\semiinf(L_\fn,\CM)\underset{L_\fh}\torsemiinf \CR.$$

\ssec{}

Recall now that for any level $\kappa$ we can introduce the
chiral algebra of differential operators (CADO) $\fD^{\on{ch}}(G)_\kappa$,
which admits two mutually commuting homomorphisms
$$\fl_\fg:\CA_{\fg,\kappa}\to \fD^{\on{ch}}(G)_\kappa \leftarrow
\CA_{\fg,\kappa'}:\fr_\fg.$$

Let $\oG$ denote the open Bruhat cell $B\cdot w_0\cdot B\subset G$,
where $w_0$ is the longest element of the Weyl group. 
We will denote by $\oGN$, $\oGB$ the corresponding open
subsets in $G/N$ and $G/B$, respectively.

Let $\fD^{\on{ch}}(\oG)_\kappa$ be the induced CADO on $\oG$.
Consider the chiral DG algebra 
$\fC^\semiinf(L_\fn,\fD^{\on{ch}}(\oG)_\kappa)$,
where we take $L_\fn$ mapping to $\fD^{\on{ch}}(\oG)_\kappa$ via
$$\fD^{\on{ch}}(\oG)_\kappa\leftarrow \fD^{\on{ch}}(G)_\kappa
\overset{\fr_\fg} \longleftarrow \CA_{\fg,\kappa'}\leftarrow 
L_{\fg,\kappa'} \leftarrow L_\fn.$$

Since $\oG\to \oGN$ is a principal $N$-bundle, from 
\cite{CHA}, Sect. 2.8.16 we obtain the following

\begin{lem}   \label{descent of CADO}
The complex $\fC^\semiinf(L_\fn,\fD^{\on{ch}}(\oG)_\kappa)$ is acyclic away
from degree zero, and the resulting chiral algebra is a CADO on
$\oGN$.
\end{lem}

Let us denote $H^\semiinf(L_\fn,\fD^{\on{ch}}(\oG)_\kappa)$ by
$\fD^{\on{ch}}(\oGN)_\kappa$. By construction, we have a homomorphism of
chiral algebras
$$\fD^{\on{ch}}(\oGN)_\kappa\leftarrow \fH'_\kappa,$$
which we will denote by $\fr_\fh$.
We define the chiral algebra $\fD^{\on{ch}}(\oGB)_\kappa$ as the
Lie-* centralizer of $\fH'_\kappa$ in $\fD^{\on{ch}}(\oGN)_\kappa$.
The map $\fl_\fg:\CA_{\fg,\kappa}\to
\fD^{\on{ch}}(G)_\kappa$ induces a homomorphism
\begin{equation} \label{bozonization}
\fl_\fg:\CA_{\fg,\kappa}\to \fD^{\on{ch}}(\oGB)_\kappa.
\end{equation}

Again, by construction, we have a canonical map
\begin{equation}  \label{H and H'}
\fD^{\on{ch}}(\oGB)_\kappa\to
\fD^{\on{ch}}(\oGN)_\kappa\underset{L_\fh;\fh}\torsemiinf \fH_\kappa.
\end{equation}

\begin{lem}   \label{G/B as semiinf}
The map in \eqref{H and H'} is an isomorphism.
\end{lem}

The proof will become clear from the discussion in the next section.

\ssec{}   \label{descr of CADO on G/B}

Note that $\fD^{\on{ch}}(\oGB)_\kappa$ is not a CADO on $\oGB$. We will now
give a more explicit, even if less canonical, description of the chiral
algebras $\fD^{\on{ch}}(\oGN)_\kappa$, $\fD^{\on{ch}}(\oGB)_\kappa$ and the
free field realization homomorphism.

Let us choose a representative of $w_0$ in $W$ and identify the
variety $\oG\simeq N\cdot w_0\cdot B$ with the product $N\times B$,
endowed with the action on $N$ on the left and of $B$ on the right.
Then $\fD^{\on{ch}}(\oG)_\kappa$ becomes a CADO on this group, isomorphic
to $\fD^{\on{ch}}(N)\otimes \fD^{\on{ch}}(B)_{\kappa'}$. We will denote the
existing maps
$$\CA_\fn\to \fD^{\on{ch}}(N)\subset \fD^{\on{ch}}(\oG)_\kappa \text{ and }
\CA_{\fb,\kappa'}\to \fD^{\on{ch}}(B)_{\kappa'}\subset
\fD^{\on{ch}}(\oG)_\kappa$$ by $\fl_\fn$ and $\fr_\fb$, respectively, and
the "new" maps, as in \cite{AG},
$$\CA_\fn\to \fD^{\on{ch}}(N)\subset \fD^{\on{ch}}(\oG)_\kappa \text{ and }
\wh{\CA}_{\fb,\kappa}\to \fD^{\on{ch}}(B)_{\kappa'}\subset
\fD^{\on{ch}}(\oG)_\kappa$$ by $\fr_\fn$, $\fl_\fb$, respectively, where
$\wh{\CA}_{\fb,\kappa}$ is the reduced chiral universal envelope of
$\wh{L}_{\fb,\kappa}$.

Then 
\begin{equation} \label{CADO on G/N}
\fD^{\on{ch}}(\oGN)_\kappa\simeq \fD^{\on{ch}}(N)\otimes
\wh{\fD}^{\on{ch}}(H)_\kappa,
\end{equation}
where $\wh{\fD}^{\on{ch}}(H)_\kappa$ is a CADO on $H$ with the maps
$$\fH_\kappa\overset{\fl_\fh}\longrightarrow
\wh{\fD}^{\on{ch}}(H)_\kappa\overset{\fr_\fh} \longleftarrow \fH'_\kappa.$$
As usual, the centralizer of $\fH'_\kappa$ in
$\wh{\fD}^{\on{ch}}(H)_\kappa$ is $\fH_\kappa$, and we obtain that
$$\fD^{\on{ch}}(\oGB)_\kappa\simeq \fD^{\on{ch}}(N)\otimes \fH_\kappa.$$
The above isomorphism makes the assertion of \lemref{G/B as semiinf}
manifest: indeed, it follows from the fact that
$\wh{\fD}^{\on{ch}}(H)_\kappa \underset{L_\fh;\fh}\torsemiinf
\fH_\kappa\simeq \fH_\kappa$, see \secref{semijective}.

Homomorphism \eqref{bozonization} therefore gives rise to a
homomorphism from the affine Kac-Moody algebra to the tensor product
of the chiral algebras $\fD^{\on{ch}}(N)$ and $\fH_\kappa$:
\begin {equation}    \label{free field}
\CA_{\fg,\kappa} \to \fD^{\on{ch}}(N) \otimes \fH_\kappa.
\end{equation}
This is the {\em free field realization} homomorphism of
\cite{FF:si,F:wak}.

The CADO $\fD^{\on{ch}}(N)$ may be identified with what physicists
call the free field $\beta\gamma$ system, and $\fH_\kappa$ is a
twisted form of a Heisenberg algebra, which is also related to a free
bosonic system. That is why the homomorphism \eqref{free field} is
referred to as free field realization.

\ssec{}   \label{chiralization}

Let us now explain in what sense the homomorphism
\begin{equation} \label{explicit bozonization}
\CA_{\fg,\kappa}\to \fD^{\on{ch}}(\oGB)_\kappa\simeq \fD^{\on{ch}}(N)\otimes
\fH_\kappa
\end{equation}
above is an affine analog (i.e., chiralization) of a well-known
phenomenon for finite-dimensional Lie algebras. We will appeal to
notations introduced in \cite{AG}.

Consider the variety $\oG\simeq N\times B$ with an action of the Lie
algebra $\fg$ on the left. This action defines a map $\fg\to T(\oGN)$,
whose image consists of vector fields, that are invariant with respect
to the action of $H\simeq B/N$ on the right. Since the Lie algebra of
such vector fields is isomorphic to $T(N)\oplus \left(\Fun(N)\otimes
\fh\right)$, we obtain a map
\begin{equation} \label{fd map}
\fg\to T(N)\oplus \left(\Fun(N)\otimes\fh\right).
\end{equation}
The restriction of this map to $\fn\subset \fg$ is the homomorphism
$\fl_n\to T(N)$. The restriction to $\fh\subset \fg$ is the sum of two
maps: one is $\fh\to T(N)$, corresponding to the natural adjoint of
$H$ on $N$, and the other is the identity map $\fh\to \fh\subset
\Fun(N)\otimes\fh$, twisted by $w_0$.

\medskip

The map of \eqref{fd map} can be chiralized in a straightforward way,
and we obtain a map of Lie-* algebras
\begin{equation}  \label{pre-bosonization}
L_\fg\to \Theta(N)\oplus \bigl(\Fun(\on{Jets}(N))\otimes L_\fh\bigr),
\end{equation}
where for an affine scheme $Y$, we denote by $\on{Jets}(Y)$ the
$\fD_X$-scheme of jets into $Y$, and $\Theta(Y)$ denotes the tangent
algebroid on this $\fD_X$-scheme. By construction, we have

\begin{lem}
The image of $L_{\fg,\kappa}\subset \CA_{\fg,\kappa}$ under
\eqref{explicit bozonization} belongs to $$\fD^{\on{ch}}(N)^{\leq 1}\oplus
\bigl(\Fun(\on{Jets}(N))\otimes (\fH_\kappa)^{\leq 1}\bigr),$$ where
$(\cdot)^{\leq i}$ denotes the PBW filtration.  The composition
\begin{align*}
&L_{\fg,\kappa}\to \fD^{\on{ch}}(N)^{\leq 1}\oplus
  \bigl(\Fun(\on{Jets}(N))\otimes (\fH_\kappa)^{\leq 1}\bigr)\to \\ &
  \bigl(\fD^{\on{ch}}(N)^{\leq 1}/\fD^{\on{ch}}(N)^{\leq 0}\bigr)\oplus
  \bigl(\Fun(\on{Jets}(N)\otimes \left((\fH_\kappa)^{\leq
  1}/(\fH_\kappa)^{\leq 0}\right)\bigr) \simeq \\ & \simeq
  \Theta(N)\oplus \bigl(\Fun(\on{Jets}(N))\otimes L_\fh\bigr)
\end{align*}
factors through $L_\fg$ and equals the map of \eqref{pre-bosonization}.
\end{lem}

\ssec{}

For the remainder of this section we will specialize to the case when
$\kappa=\kappa_\crit$. The following basic fact is established in
\cite{CHA}, Sect. 2.8.17.

\begin{prop}   \label{Miura=CADO}
The Lie-* algebra $\wh{L}'_{\fh,\kappa}$ is commutative if
and only if $\kappa=\kappa_\crit$. In this case there is
a canonical isomorphism 
$$\Spec(\fH_\crit)\simeq \ConHX^\fD,$$
respecting the torsor structure on both sides with
respect to the D-scheme of $\fh^*$-values $1$-forms 
on $X$.
\end{prop}

Since $\fH'_\crit$ is commutative, it is contained as a chiral
subalgebra in $\fD^{\on{ch}}(\oGB)_\crit$, moreover from \eqref{H and H'}
we infer:
\begin{equation} 
\fH'_\crit\simeq \fz(\fD^{\on{ch}}(\oGB)_\crit).
\end{equation}

Since $\fH_\crit$ is commutative as well, from the isomorphism
\eqref{H and H'} we obtain that there exists a homomorphism
(which is easily seen to be an isomorphism) from 
$\fH_\crit$ to $\fz(\fD^{\on{ch}}(\oGB)_\crit)$.

\begin{lem} \label{H symmetry}
The resulting homomorphism $\fH_\crit\to \fH'_\crit$ comes
from the sign-inversion isomorphism $L_{\fh,\crit}\to L'_{\fh,\crit}$
of commutative Lie-* algebras.
\end{lem}

\ssec{}

We will now study the homomorphism
\begin{equation} \label{critical bosonization}
\fl_\fg:\CA_{\fg,\crit}\to \fD^{\on{ch}}(\oGB)_\crit.
\end{equation}

\begin{prop}   \label{centralizer on G/B}
The centralizer of $\CA_{\fg,crit}$ in $\fD^{\on{ch}}(\oGB)_\crit$ equals
$\fH'_\crit$.
\end{prop}

\begin{proof}

Since $\fH'_\crit$ is the center of the chiral algebra
$\fD^{\on{ch}}(\oGB)_\crit$, the fact that it centralizes the image of
$\CA_{\fg,crit}$ is evident.

To prove the inclusion in the opposite direction, we will establish a
stronger fact. Namely, that the centralizer in $\fD^{\on{ch}}(\oGB)_\crit$
of the image of $L_{\fn}+\fh$ is already contained in $\fH'_\crit$.

Using the description of $\fD^{\on{ch}}(\oGB)_\crit$ given in \secref{descr
of CADO on G/B}, we obtain that the centralizer of $\fl_\fn(\CA_\fn)$
in it equals $\fr_\fn(\CA_\fn)\otimes \fH'_\crit$, in the notation of
{\it loc.cit.}

Consider now the action of $\fh\in \Gamma(X,\CA_{\fg,\crit})$ on
$\fD^{\on{ch}}(\oGN)_\crit\simeq \fD^{\on{ch}}(N)\otimes
\wh{\fD}^{\on{ch}}(H)_\crit$.  By \secref{chiralization}, this action
decomposes as a tensor product of the natural adjoint action on
$\fD^{\on{ch}}(N)$, and the action on $\wh{\fD}^{\on{ch}}(H)_\kappa$ given by
$\fl_\fh$, twisted by $w_0$.  Since $\wh{\fD}^{\on{ch}}(H)_\kappa$ is
commutative, the resulting action of $\fh$ on
$\fD^{\on{ch}}(\oGB)_\crit\simeq \fD^{\on{ch}}(N)\otimes \fH'_\crit$ is the
adjoint action along the first factor.

This implies our assertion since $(\CA_\fn)^\fh\simeq \BC$, as $\fh$
acts on $\fn$, and hence on $\CA_\fn$, by characters, which belong to
the positive span of $\Delta^+$.

\end{proof}

\ssec{}

Consider now the composition
$$\fz_\fg=\fz(\CA_{\fg,\crit})\to \fD^{\on{ch}}(\oGB)_\crit.$$

{}From \propref{centralizer on G/B} we immediately obtain the following
result.

\begin{cor} \label{center-to-center}
The image of $\fz_\fg$ is contained in
$\fz(\fD^{\on{ch}}(\oGB)_\crit)=\fH'_\crit$.
\end{cor}

Thus, we obtain a homomorphism of commutative chiral algebras 
\begin{equation} \label{RT Miura}
\fz_\fg\to \fH'_\crit\simeq \fH_\crit.
\end{equation}

Let us now recall that ultimate form of the isomorphism statement of
\cite{FF,F:wak} (see Theorem 11.3 of \cite{F:wak}):

\begin{thm} \label{FF with Miura}
There exists a canonical isomorphism of commutative chiral algebras
$\fz_\fg\simeq \Fun(\Op(X))$ such that the diagram
$$
\CD
\fz_\fg @>{\sim}>> \Fun(\Op(X)^\fD) \\
@VVV   @V{\on{MT}^*}VV  \\
\fH_\crit @>{\sim}>> \Fun(\ConHX^\fD) 
\endCD
$$
is commutative, where the left vertical arrow is the map of \eqref{RT
Miura}, the right vertical arrow is the Miura transformation of
\eqref{Miura transformation}, and the bottom horizontal arrow is the
isomorphism of \propref{Miura=CADO}, composed with the
automorphism, induced by the automorphism $\tau:=\cla\mapsto -w_0(\cla)$
of $\check{H}$.

\end{thm}

\ssec{}

To conclude this section let us return to the set-up of
\secref{chiralization}.  Consider the map $\fg\to T(N)$, obtained by
composing the map of \eqref{fd map} with the projection on the
$T(N)$-factor.

It is well-known that lifts of this map to a homomorphism of Lie
algebras $\fg\to \fD(N)^{\leq 1}$, which on $\fn\subset \fg$ induce
the map $\fl_\fn:\fn\to T(N)$, are classified by characters of $\fh$
(and correspond to $G$-equivariant twistings on $G/B$). We would like
to establish an affine analog of this statement.

The analog of characters of $\fh$ will be played by the set of chiral
algebra homomorphisms $\psi:\fH_\crit\to \CO_X$. For any such $\psi$,
the composition
$$\phi:L_{\fg,\crit}\to  \fD^{\on{ch}}(N)\otimes \fH_\crit\to  \fD^{\on{ch}}(N)$$
is a Lie-* algebra homomorphism, satisfying:

\begin{itemize}

\item
The image of $\phi$ belongs to $\fD^{\on{ch}}(N)^{\leq 1}$,

\item
The composition $L_{\fg,\crit}\to  \fD^{\on{ch}}(N)\to \Theta(N)$
equals the composition of the map \eqref{pre-bosonization},
followed by the projection on the $\Theta(N)$-factor,

\item
The restriction of $\phi$ to $L_\fn$ equals $\fl_\fn$.

\end{itemize}

\begin{prop}
Let $L'_\fg$ be a central extension of $L_\fg$ by means of $\omega_X$, 
split over $L_\fn$, and let $\phi:L'_\fg\to \fD^{\on{ch}}(N)$ be a 
homomorphism of Lie-* algebras, satisfying the three properties above.
Then $L'_\fg\simeq L_{\fg,\crit}$ and $\phi$ is obtained from some 
$\psi:\fH_\crit\to \CO_X$ in the manner described above.
\end{prop}

\begin{proof}

First, since $L'_\fg$ splits over $L_\fn$, we obtain that as a
$\fD_X$-module $L'_\fg\simeq L_\fg\oplus \omega_X$. Let us show that the
bracket on $L'_\fg$ corresponds to the critical pairing. For this, it
is sufficient to calculate the bracket on $L_\fh\subset
L_\fg$. However, since $L_\fh$ is commutative, the latter bracket is
independent of the choice of a pair $(L'_\fg, \phi)$. Hence, we may
choose the pair $L'_\fg=L_{\fg,\crit}$ and a homomorphism
corresponding to some homomorphism $\psi:\fH_\crit\to \CO_X$. In the
latter case, our assertion is clear.

Consider the set of all homomorphisms of chiral algebras
$\psi:\fH_\crit\to \CO_X$. By definition, this is a torsor over
$\Gamma(X,\omega_X\otimes \fh^*)$.  Consider now the space of
homomorphisms $\phi_\fb:L_{\fb,\crit}\to \fD^{\on{ch}}(N)^{\leq 1}$,
satisfying the same three conditions as $\phi$. This set is also a
torsor over $\Gamma(X,\omega_X\otimes \fh^*)$. Moreover, it is easy to
see that the map $\psi\mapsto \phi\mapsto \phi|_{L_{\fb}}=:\phi_\fb$
is a map of torsors.

Hence, for any $\phi$ as in the proposition, there exists a $\psi$,
such that the two homomorphism $L_{\fg,\crit}\to \fD^{\on{ch}}(N)^{\leq 1}$
coincide, when restricted to $L_{\fb,\crit}$.  We claim that in this
case the two homomorphisms in question coincide on the entire of
$L_{\fg,\crit}$.

Indeed, let $\phi_1$ and $\phi_2$ be two such homomorphisms. Then
$\phi_1-\phi_2$ is a map $L_\fg/L_\fb\to \Fun(\on{Jets}(N))\otimes
\omega_X$.  Let $\ssf$ be a section on $L_{\fn^-}$, and let $\sse$ be
a section of $L_\fn$ such that $[\sse,\ssf]\in L_\fh$. We obtain that
$[\phi_1(\ssf)-\phi_2(\ssf),\phi_\fb(\sse)]=0$. Hence, the image of
$\phi_1-\phi_2$ consists of $L_\fn$-invariant sections of
$\Fun(\on{Jets}(N))\otimes \omega_X$, and the latter subspace is
$\omega_X$.

Again, for $\ssf$ above, let $\ssh$ be a section of $L_\fh$ such that
$[\ssf,\ssh]=c\cdot \ssf$, where $c$ is a non-zero scalar. We obtain
$[\phi_1(\ssf)-\phi_2(\ssf),\phi_\fb(\ssh)]=c\cdot
(\phi_1(\ssf)-\phi_2(\ssf))$. However, by the above,
$\phi_1(\ssf)-\phi_2(\ssf)$ is central.  Hence, $c\cdot
(\phi_1(\ssf)-\phi_2(\ssf))=0$, implying our assertion.

\end{proof}

\section{Construction of Wakimoto modules}    \label{sec Wak}

\ssec{}   \label{constr of Wak}

Homomorphism \eqref{bozonization} allows us to produce representations
of $\CA_{\fg,\kappa}$, i.e., $\hg_\kappa$-modules, by restricting
modules of $\fD^{\on{ch}}(\oGB)_\kappa$. This should be regarded
as a chiral analog of the construction of $\fg$-modules by
taking sections of twisted D-modules on the big Schubert cell
$\oGB$.

\medskip

In the applications, modules over $\fD^{\on{ch}}(\oGB)_\kappa$ that we will
consider are obtained using \eqref{H and H'}, from pairs of modules:
$\CM\in \fD^{\on{ch}}(\oGN)_\kappa\mod$, and $\CR\in \fH_\kappa\mod$
by taking $\CM\underset{L_\fh}\torsemiinf \CR$. Let us describe 
the examples of $\fD^{\on{ch}}(\oGB)_\kappa$-modules that we will consider.

\ssec{}

First, note that if $\fD^{\on{ch}}(Y)$ is a CADO on (the scheme of jets
corresponding to) a smooth affine $X$-scheme $Y$, any left D-module on
the scheme $Y[[t]]$ gives rise to a chiral module over $\fD^{\on{ch}}(Y)$.

Indeed, if $\CF$ is such a D-module, it (or, rather, the space of its
global sections) is naturally a chiral module over
$\Fun(\on{Jets}(Y))$ and a Lie-* module over $\Theta_Y$. In this case
we can induce it and obtain a chiral module over $\fD^{\on{ch}}(Y)$.

If $Y'\subset Y$ is a smooth locally closed subvariety, let us denote
by $\Dist_Y(Y')$ the left D-module of distributions on $Y'$ (i.e., the
$*$-extension of the D-module $\Fun(Y')$), and let
$\Dist_{Y[[t]]}(\on{ev}^{-1}(Y'))$ denote the corresponding left
D-module on $Y[[t]]$, i.e.,
$$\Dist_{Y[[t]]}(\on{ev}^{-1}(Y'))\simeq \on{ev}^*(\Dist_Y(Y')).$$
Finally, let $\Dist_Y^{\on{ch}}(\on{ev}^{-1}(Y'))$ denote the resulting 
$\fD^{\on{ch}}(Y)$-module. 

\medskip

Let us take $Y=\oG$ and for each element $w\in W$ consider
$Y'=\on{Ad}_{w_0w^{-1}}(N)\cdot w_0\cdot N \subset \oG$.  For example,
if $w=w_0$ we get the D-module of functions on $N\cdot w_0\cdot N$,
and if $w=1$ we get the $\delta$-function at $w_0\cdot N$.

\medskip

Note that $\on{ev}^{-1}(N)=I^0$, so we obtain a left D-module
$\Dist_{G[[t]]}(\on{Ad}_{w_0w^{-1}}(I^0)\cdot w_0\cdot I^0)$ on
$G[[t]]$ and $\Dist_G^{\on{ch}}(\on{Ad}_{w_0w^{-1}}(I^0)\cdot w_0\cdot
I^0)_\kappa\in \fD^{\on{ch}}(\oG)_\kappa\mod$. Consider the chiral
$\fD^{\on{ch}}(\oGN)_\kappa$-module
\begin{equation} \label{descended D-module}
\Dist_{\oGN}^{\on{ch}}(\on{ev}^{-1}(\on{Ad}_{w_0w^{-1}}(N)\cdot
w_0))_\kappa:=
H^\semiinf(L_\fn,\on{Dist}_G^{\on{ch}}(\on{Ad}_{w_0w^{-1}}(I^0)\cdot
w_0\cdot I^0)_\kappa).
\end{equation}
In other words,
$\Dist_{\oGN}^{\on{ch}}(\on{ev}^{-1}(\on{Ad}_{w_0w^{-1}}(N)\cdot
w_0))_\kappa$ is obtained by the above construction for $Y=\oGN$ and
$Y'=\on{Ad}_{w_0w^{-1}}(N)\cdot w_0\subset \oGN$.  From \secref{descr
of CADO on G/B} we obtain that
$\Dist_{\oGN}^{\on{ch}}(\on{ev}^{-1}(\on{Ad}_{w_0w^{-1}}(N)\cdot
w_0))_\kappa$ is indeed acyclic away from degree $0$.

Moreover, as a module over $H^0_{DR}(\D^\times,L_\fn\oplus
\wh{L}'_{\fh,\kappa})$, it is isomorphic to
$$\Dist^{\on{ch}}_N(\on{Ad}_{w_0w^{-1}}(N)\cap N)\otimes 
\on{Ind}^{H^0_{DR}(\D^\times,\wh{L}'_{\fh,\kappa})}_{t\fh[[t]]\oplus \BC}
\Bigl(\Fun\left(H(t\BC[[t]]])\right)\Bigr).$$

In particular, as a $H^0_{DR}(\D^\times,\wh{L}'_{\fh,\kappa})$-module,
it is $H(t\BC[[t]]])$-integrable, and injective as a
$H(t\BC[[t]]])$-representation.  Furthermore, it is free over over
$\fh[t^{-1}]$ for any choice of a splitting $\fh[t^{-1}]\to
H^0_{DR}(\D^\times,\wh{L}'_{\fh,\kappa})$.

\ssec{}

Now, for $w\in W$ and an $\fH_\kappa$-module $\CR$ we define the
(complex of) $\fD^{\on{ch}}(\oGB)_\kappa$-modules 
\begin{equation} \label{definition of twisted wakimoto}
'\BW^w_\kappa(\CR) :=
\Dist_{\oGN}^{\on{ch}}(\on{ev}^{-1}(\on{Ad}_{w_0w^{-1}}(N)\cdot w_0))_\kappa 
\underset{\fh\ppart,t\fh[[t]]}\torsemiinf \CR.
\end{equation}

Note that by the $H(t\BC[[t]]])$-integrability of
$\Dist_{\oGN}^{\on{ch}}(\on{ev}^{-1}(\on{Ad}_{w_0w^{-1}}(N)\cdot
w_0))_\kappa$, we have
\begin{equation} \label{R integrable}
'\BW^w_\kappa(\CR) \simeq {}'\BW^w_\kappa(\on{Av}_{H(t\BC[[t]]])}(\CR)),
\end{equation}
where $\on{Av}_{H(t\BC[[t]]])}$ denotes the averaging functor with respect
to $H(t\BC[[t]]])$, see \secref{averaging}.  Therefore, with no
restriction of generality, we can (and will) assume that $\CR$ is
$H(t\BC[[t]]])$-integrable. Under this assumption, as a $\fn\ppart)$-module
\begin{equation}  \label{descr of 'W(R)}
'\BW^w_\kappa(\CR)\simeq 
\Dist^{\on{ch}}_N(\on{ev}^{-1}(\on{Ad}_{w_0w^{-1}}(N)\cap N))\otimes \CR.
\end{equation}
In particular, it is acyclic away from degree $0$.

\medskip

We restrict $'\BW^w_\kappa(\CR)$ to $\CA_{\fg,\kappa}$ and obtain an
object of $\hg_\kappa\mod$, denoted by the same symbol.  Note that
when defining $'\BW^w_\kappa(\CR)$, we can avoid mentioning the chiral
algebra $\fD^{\on{ch}}(\oGB)_\kappa$. Namely,
$$'\BW^w_\kappa(\CR)\simeq \left(
\Dist_G^{\on{ch}}(\on{Ad}_{w_0w^{-1}}(I^0)\cdot w_0\cdot I^0)_\kappa\right)
\underset{\fb\ppart,\fn[[t]]+t\fh[[t]]}\torsemiinf \CR.$$

\medskip

The $\fD^{\on{ch}}(G)_\kappa$-module
$\Dist_G^{\on{ch}}(\on{Ad}_{w_0w^{-1}}(I^0)\cdot w_0\cdot I^0)_\kappa$ is
by construction equivariant with respect to the group
$\on{Ad}_{w_0w^{-1}}(I^0)$, when we think of the action on $G\ppart$
on itself by left multiplication. Let $\Dist_G^{\on{ch}}(I^0\cdot w\cdot
I^0)_\kappa$ be the chiral $\CA_{\fg,\kappa}$-module, obtained from
the module $\Dist_G^{\on{ch}}(\on{Ad}_{w_0w^{-1}}(I^0)\cdot w_0\cdot
I^0)_\kappa$ by applying the left shift by $w\cdot w_0$.

Set
\begin{equation} \label{definition of general wakimoto}
\BW^w_\kappa(\CR)
:=\left(\Dist_G^{\on{ch}}(I^0\cdot w\cdot
I^0)_\kappa\right)\underset{\fb\ppart,\fn[[t]]+t\fh[[t]]}\torsemiinf
\CR.
\end{equation}
This is what we will call the Wakimoto module of type $w$ corresponding
to the $\fH_\kappa$-module $\CR$.

\medskip

Tautologically, as a $\hg_\kappa$-module, $\BW^w_\kappa(\CR)$ is
obtained from $'\BW^w_\kappa(\CR)$ by the automorphism
$\on{Ad}_{ww_0}$ of $\fg$, and it is $I^0$-equivariant.  Note,
however, that $\BW^w_\kappa(\CR)$ does not come by restriction from a
$\fD^{\on{ch}}(\oGB)_\kappa$-module, unless $w=w_0$.

We have a description of $\BW^w_\kappa(\CR)$ similar to \eqref{descr
of 'W(R)}, but with respect to the subalgebra $\fn^{ww_0}\ppart$, where
we set $\fn^w:=\on{Ad}_{w}(\fn)$, $N^w=\on{Ad}_{w}(N)$. Namely,
\begin{equation}  \label{descr of W(R)}
\BW^w_\kappa(\CR)\simeq
\Dist^{\on{ch}}_{N^{ww_0}}\left(\on{ev}^{-1}(N^{ww_0}\cap N)\right)\otimes
\CR
\end{equation}

\ssec{}

Assume now that the $\fH_\kappa$-module $\CR$ is $H[[t]]$-integrable.
(Having already the assumption that it is $H(t\BC[[t]]])$-integrable, this
amounts to requiring that $\fh$ acts semi-simply with eigenvalues
corresponding to integral weights.) We claim that in this case the
module $\BW^w_\kappa(\CR)$ is $I$-integrable.

\medskip

Indeed, let us instead of $\Dist_G^{\on{ch}}(I^0\cdot w\cdot I^0)_\kappa$
take the chiral $\fD^{\on{ch}}(G)_\kappa$-module
$$\Dist_G^{\on{ch}}(I^0\cdot w\cdot I)_\kappa.$$ 
As an object of $\fD^{\on{ch}}(G)_\kappa\mod$, it is clearly $I$-integrable with
respect to both left and right action of $G\ppart$ on itself.

Consider $H^\semiinf(L_\fn,\on{Dist}^{\on{ch}}(I^0\cdot w\cdot I)_\kappa)$. 
This is a $\fH'_\kappa$-module, which is
$H[[t]]$-integrable and injective as an $H[[t]]$-module.

One easily checks that $\BW^w_\kappa(\CR)$ is isomorphic to 
$$H^\semiinf\left(L_\fn,\on{Dist}^{\on{ch}}(I^0\cdot w\cdot I)_\kappa\right)
\underset{L_\fh;\fh}\torsemiinf \CR,$$
which is manifestly $I$-integrable.

\ssec{}    \label{wakimoto lambda}

For a weight $\lambda\in \fh^*$ consider the $1$-dimensional Lie-*
module over $\wh{L}_{\fh,\kappa}$ corresponding to the character
$\lambda$. Let us denote by $\pi_\lambda$ the induced chiral module
over $\fH_\kappa$.

For future use we introduce the notation
\begin{equation} \label{def wak lambda}
\BW^w_{\kappa,\lambda}:=\BW^w_\kappa(\pi_{w^{-1}(\lambda+\rho)+\rho}).
\end{equation}

Observe that the definition of $\BW^w_{\kappa,\lambda}$ can be
rewritten as
$$\left( H^\semiinf\left(\fn\ppart,\fn[[t]],\Dist_G^{\on{ch}}(I^0\cdot
w\cdot I^0)_\kappa \right)^{t\fh[[t]]}\right)
\underset{\fh}\otimes \BC^{w^{-1}(\lambda+\rho)+\rho}.$$

\medskip

Let $M^w_\lambda$ be the $\fg$-module equal to $\Dist_G(N\cdot w\cdot N)
\underset{\fb,\fn}\otimes \BC^{w^{-1}(\lambda+\rho)+\rho}$.  Note
that when $w=1$, $M^w_\lambda$ is the Verma module $M_\lambda$, and
when $w=w_0$, $M^w_\lambda$ is the dual Verma $M^\vee_\lambda$. In
general, $M^w_\lambda$ has always highest weight $\lambda$, and it is
characterized by the property that it is free with respect to the Lie
subalgebra $\fn^{ww_0}\cap \fn^-$ and co-free with respect to
$\fn^{ww_0}\cap \fn$.

Set $\BM^w_{\kappa,\lambda}:=\on{Ind}^{\hg_\kappa}_{\fg[[t]]\oplus
\BC\one}(M^w_\lambda)$ be the induced $\hg_\kappa$-module. We claim
that we always have a map
\begin{equation} \label{Verma to wakimoto}
\BM^w_{\kappa,\lambda}\to \BW^w_{\kappa,\lambda}.
\end{equation}

This amounts to constructing a map of $\fg[[t]]$-modules
$M^w_\lambda\to \BW^w_{\kappa,\lambda}$. We have
\begin{align*}
&M^w_\lambda\hookrightarrow \Dist_{G[[t]]}(I^0\cdot w\cdot I^0)
\underset{\fb[[t]],\fn[[t]]+t\fh[[t]]}\otimes
\BC^{w^{-1}(\lambda+\rho)+\rho}\hookrightarrow \\ &
\left(\Dist_G^{\on{ch}}(I^0\cdot w\cdot
I^0)_\kappa\right)^{\fn[[t]]+t\fh[[t]]} \underset{\fh}\otimes
\BC^{w^{-1}(\lambda+\rho)+\rho},
\end{align*}
which maps to the required semi-infinite cohomology. 

\section{Convolution action on Wakimoto modules}     \label{conv of wak}

\ssec{}

In this section we will apply the formalism of convolution functors
$$\star:\fD(G/K)_\kappa\mod^{K'}\times D(\hg_\kappa\mod)^K\to
D(\hg_\kappa\mod)^{K'},$$ where $K,K'$ are subgroups of $G[[t]]$ to
derive some additional properties of Wakimoto modules.

The subgroups that we will use will be either $I^0$ or $G^{(1)}$, the
first congruence sub group in $G[[t]]$, and if $\kappa$ is integral,
also $I$. When a confusion is likely to occur, we will use the
notation $\cdot \underset{K}\star\cdot$ to emphasize which equivariant
derived category we are working in, see \secref{convolution
action}. We will identify D-modules on $G$ (resp., $G/N$, $G/B$) with
the corresponding $\kappa$-twisted D-modules on $G\ppart/G^{(1)}$
(resp., $\wt{\Fl}_G=G\ppart/I^0$, $\Fl_G=G\ppart/I$).

\medskip

Another two pieces of notation that we will need are as follows.  If
${\mathbf g}$ is a point of $G\ppart$, and $\CF$ an object of an
arbitrary category with a Harish-Chandra action of $G\ppart$ at level
$\kappa$, we will denote by $\delta_{{\mathbf g}}\star \CF$ the twist
of $\CF$ by ${\mathbf g}$.  If $\CF$ is equivariant with respect to a
congruence subgroup $K\subset G[[t]]$, then
$$\delta_{{\mathbf g}}\star \CF\simeq \delta_{{\mathbf
g}_{G\ppart/K}}\underset{K}\star \CF,$$ where $\delta_{{\mathbf
g}_{G\ppart/K}}$ is the unique $\kappa$-twisted D-module on $G\ppart/K$,
whose $!$-fiber at the point ${\mathbf g}_{G\ppart/K}\in G\ppart/K$ is
$\BC$.

Let $\bU$ be a pro-unipotent subgroup such that
$\kappa|_{\on{Lie}(\bU)}$ is trivial. Then for $\CF$ as above,
$\uBC_{\bU}\star \CF$ will denote the same thing as
$\on{Av}_{\bU}(\CF)$. In other words, if $\CF$ is equivariant with
respect to some unipotent $K\subset G[[t]]$ containing a congruence
subgroup, and $\bU'=\bU\cap K$, then
$$\uBC_{\bU}\star \CF\simeq
\Dist_{G\ppart/K}(\bU/\bU')_\kappa\underset{K}\star \CF\otimes 
\det\left(\on{Lie}(\bU)/\on{Lie}(\bU')[1]\right)^{\otimes -1},$$
where $\uBC_{\bU/\bU'}$ denotes the cohomologically shifted
D-module on $\bU/\bU'$,
corresponding via Riemann-Hilbert to the {\it constant sheaf}
on $\bU/\bU'$, and $\Dist_{G\ppart/K}(\bU/\bU')_\kappa$ is the unique
$\kappa$-twisted D-module on $G\ppart/K$, supported on 
$\bU/\bU'\subset G\ppart/K$, and whose $!$-restriction to this subscheme
is $\Fun(\bU/\bU')$, see \secref{D-mod on groups and quotients}.

We will use the following observation:
\begin{lem} \label{product of subgroups}
Suppose that $\bU$ contains two subgroups 
$\bU_1$ and $\bU_2$ such that the multiplication map
defines an isomorphism $\bU_1\times \bU_2\to \bU$.
Then
$$\uBC_{\bU}\star \CF\simeq \uBC_{\bU_1}\star (\uBC_{\bU_2}\star \CF).$$
\end{lem}

\medskip

For $\wt{w}\in W_{\on{aff}}$ we will denote by $\wt{j}_{\kappa,\wt{w}}$ the
unique $\kappa$-twisted $I^0$-equivariant D-module on $\wt{\Fl}_G$,
supported on $I^0\cdot \wt{w}\subset \wt{\Fl}_G$, whose
$!$-restriction to this subscheme is isomorphic to $\Fun(I^0\cdot
\wt{w})$, as an $I^0$-equivariant quasi-coherent sheaf. Of course, the
isomorphism class of this D-module depends on the choice of a
representative of $\wt{w}$ in $G\ppart$.

Since $\wt{j}_{\kappa,\wt{w}}\simeq \uBC_{I^0}\star
\delta_{\wt{w}_{G\ppart/I^0}} \otimes
\det\left(\on{Lie}(I^0)/\on{Lie}(I^0)\cap
\on{Ad}_{\wt{w}}(\on{Lie}(I^0))[1]\right)$, from \lemref{product of
subgroups} we obtain the following:

\begin{lem}  \label{step-by-step convolution}
For $\wt{w}\in W_{\on{aff}}$ assume that $I^0$ can be written as a product
of subgroups $\bU_1\cdot \bU_2$ such that
$\on{Ad}_{\wt{w}{}^{-1}}(\bU_2)\subset I^0$.  Then for an
$I^0$-equivariant object $\CF$ of a category with a Harish-Chandra
action of $G\ppart$, we have a canonical isomorphism
$$\wt{j}_{\kappa,\wt{w}}\underset{I^0}\star \CF\simeq
\on{Av}_{\bU_1}\left(\delta_{\wt{w}}\star \CF\right)\otimes
\det\left(\on{Lie}(\bU_1)/\on{Lie}(\bU_1)\cap
\on{Ad}_{\wt{w}}(\on{Lie}(\bU_1))[1]\right).$$
\end{lem}

Suppose that $\kappa$ is integral, i.e., comes from a group ind-scheme
extension $\wh{G\ppart}$ of $G\ppart$ split over $G[[t]]$; let us denote
by $\CP^\kappa$ the resulting line bundle on $\Gr_G=G\ppart/G[[t]]$. In
this case we will denote by $\wcosta$ (resp., $\wsta$) the
$I$-equivariant $\kappa$-twisted D-modules on $\Fl_G$ given by the
$*$-extension (resp., $!$-extension) of the twisted right D-module on
$I\cdot \wt{w}\subset \Fl_G$, corresponding to the restriction of the
line bundle $\wh{G\ppart}/I\to G\ppart/I$ to this subscheme. If $\kappa$
is not integral the above $I$-equivariant D-modules still make sense
for $w\in W$.

If $l(\wt{w}_1)+l(\wt{w}_2)=l(\wt{w}_1\cdot \wt{w}_2)$, then
$$j_{\wt{w}_1,*}\underset{I}\star j_{\wt{w}_2,*}\simeq
j_{\wt{w}_1\cdot \wt{w_2},*} \text{ and }
j_{\wt{w}_1,!}\underset{I}\star j_{\wt{w}_2,!}\simeq j_{\wt{w}_1\cdot
\wt{w_2},!}.$$ Since the functor $\wcosta\underset{I}\star \cdot$ is
right exact, the above isomorphism implies that the functor
$\wsta\underset{I}\star \cdot$, being its quasi-inverse, is left
exact.

\medskip

Let us observe that the definition of $\wcosta$ (resp., $\wsta$) is
evidently independent of the choice of representatives $\wt{w}$ in
$G\ppart$. The direct image of $\wt{j}_{\kappa,\wt{w}}$ under
$\wt{\Fl}_G\to \Fl_G$ is isomorphic to $\wcosta\otimes
(\fl_\kappa^{\wt{w}})^{\otimes -1}$, where $\fl_\kappa^{\wt{w}}$ is
the line defined as
\begin{equation} \label{determinant line}
\fl_\kappa^{\wt{w}}:= \Gamma\left(I\cdot
\wt{w},\Omega^{\on{top}}(I^0\cdot \wt{w})\otimes \CP^\kappa|_{I^0\cdot
\wt{w}}\right)^{I_0}.
\end{equation}

\ssec{}

Let us first observe that for $w\in W$
$$\Dist_G^{\on{ch}}(I^0\cdot w\cdot I^0)_\kappa\simeq 
\wt{j}_{\kappa,w}  \underset{I^0}\star \Dist_G^{\on{ch}}(I^0)_\kappa.$$
Hence, we obtain the following

\begin{lem}  \label{wakimoto from one-another}
$\BW^w_\kappa(\CR)\simeq \wt{j}_{\kappa,w} \underset{I^0}\star
  \BW_\kappa^1(\CR)$.
\end{lem}

If $\CR$ is integrable with respect to $H[[t]]$, the above lemma
implies that
\begin{equation}  \label{reform of wakimoto from one-another integral}
\BW^w_\kappa(\CR)\simeq j_{w,*}\star \BW^1(\CR),
\end{equation}
which, in turn, implies that 
\begin{equation} \label{more waimoto from one-another}
j_{w,!}\star \BW^{w^{-1}}(\CR)\simeq \BW^1(\CR),
\end{equation}
and if $l(w_1\cdot w_2)=l(w_1)+l(w_2)$, then
\begin{equation} \label{even more wakimoto from one-another}
j_{w_1,*}\star \BW^{w_2}(\CR)\simeq \BW^{w_1\cdot w_2}(\CR).
\end{equation}

\ssec{}

For $w\in W$ recall that $\fn^w$ (resp., $\fb^w$) denotes the subalgebra 
$\on{Ad}_w(\fn)\subset \fg$ (resp., $\on{Ad}_w(\fb)\subset \fg$). Note that
the Cartan quotient of $\fb^w$ is still canonically identified with $\fh$. 
For $w=w_0$ we will sometimes also write $\fn^-$, $\fb^-$.

\begin{prop}   \label{wakimoto as reduction}
For any chiral $\fH_\kappa$-module $\CR$,
$$\BW^w_{\kappa}(\CR)\simeq \left(\Dist_G^{\on{ch}}(I^0)_\kappa\right)
\underset{\fb^w\ppart,t\fb^w[[t]]+\fn\cap\fn^w}\torsemiinf \CR.$$
\end{prop}

\begin{proof} 

It is enough to show that 
$$H^\semiinf(\fn\ppart,\fn[[t]],\Dist_G^{\on{ch}}(I^0\cdot w\cdot
I^0)_\kappa)\simeq H^\semiinf(\fn^w\ppart,t\fn^w[[t]] +
\fn^w\cap\fn,\Dist_G^{\on{ch}}(I^0)_\kappa),$$ in a way compatible with the
$\fH_\kappa$-actions.

Again, we have $\Dist_G^{\on{ch}}(I^0\cdot w\cdot I^0)_\kappa\simeq
\Dist_G^{\on{ch}}(I^0)_\kappa\underset{I^0}\star \wt{j}_{\kappa,w}$, where
we are using the action of $G\ppart$ on itself by right translations.
We have $I^0=(I^0\cap B^-[[t]])\cdot (I^0\cap N[[t]])$. By
\lemref{step-by-step convolution}, we obtain that
$$\Dist_G^{\on{ch}}(I^0\cdot w\cdot I^0)_\kappa\simeq 
\Dist_G^{\on{ch}}(I^0)_\kappa\star \delta_{w}\star \uBC_{N[[t]]}
\otimes \det\left(\fn/\fn^{w^{-1}}\cap \fn[1]\right).$$
Hence, by \secref{semiinf and de rham}
\begin{align*}
&H^\semiinf(\fn\ppart,\fn[[t]],\Dist_G^{\on{ch}}(I^0\cdot w\cdot
  I^0)_\kappa)\simeq \\ &\simeq
H^\semiinf(\fn\ppart,\fn[[t]],\Dist_G^{\on{ch}}(I^0)_\kappa\star \delta_{w})
\otimes \det\left(\fn/\fn^{w^{-1}}\cap \fn[1]\right),
\end{align*}
which, in turn, is isomorphic to
$H^\semiinf(\fn^w\ppart,\fn^w[[t]],\Dist_G^{\on{ch}}(I^0)_\kappa)\otimes
\det\left(\fn^w/\fn^w\cap \fn[1]\right)$. The determinant line exactly
accounts for the change of the lattice $\fn^w[[t]]\mapsto
t\fn^w[[t]]+\fn\cap\fn^w$.

\end{proof}

As a corollary, we obtain the following characterization of the
Wakimoto modules $\BW^w_{\kappa}(\CR)$:

\begin{cor}  \label{cohomological characterization}
For an $I^0$-integrable $\hg_{\kappa'}$-module $\CM$ and a chiral
$\fH_\kappa$-module $\CR$, we have a quasi-isomorphism
$$\CM\underset{\fg\ppart,\Lie(I^0)}\torsemiinf \BW^w_\kappa(\CR)\simeq
\fC^\semiinf(\fb^w\ppart,t\fb^w[[t]]+\fn\cap\fn^w,\CM\otimes \CR).$$
\end{cor}

\begin{proof}

In view of \propref{wakimoto as reduction}, it suffices to show that
for any $\CM$ as in the proposition,
$$\CM\underset{\fg\ppart,\Lie(I^0)}\torsemiinf
\left(\Dist_G^{\on{ch}}(I^0)_\kappa\right)\simeq \CM$$ as
$\hg_\kappa$-modules. But this follows from \secref{semijective}.

\end{proof}

\ssec{}

We will now show that Wakimoto modules of type $w\cdot w_0$ are
well-behaved with respect to the functor of semi-infinite cohomology
of the algebra $\fn^w\ppart$. This is in fact a fundamental property
of Wakimoto modules which was found in \cite{FF:si}.

Namely, let $\CL$ be a module over $\fn^w\ppart$, on which the subalgebra
$t\fn^w[[t]]+\fn\cap\fn^w$ acts locally nilpotently. Let $\CR$ be
an $\fH_\kappa$-module, on which $t\fh[[t]]$ acts locally nilpotently.
(By \eqref{R integrable} the latter is not really restrictive.)

\begin{prop}   \label{acyclicity of wakimoto}
Under the above circumstances,
$$\CL\underset{\fn^w\ppart,t\fn^w[[t]]+\fn\cap\fn^w}\torsemiinf
\BW_\kappa^{ww_0}(\CR)$$ is canonically isomorphic to $\CL\otimes
\CR$.
\end{prop}

\begin{proof}

By \eqref{descr of W(R)}, it suffices to show that
$$\CL\underset{\fn^w\ppart,t\fn^w[[t]]+\fn\cap\fn^w}\torsemiinf
 \Dist_{N^w}^{\on{ch}} \left(\on{ev}^{-1}(N^w\cap N)\right) \simeq
 \CL.$$

However, this readily follows from \corref{endomorphisms of
semijective}(2).

\end{proof}

\ssec{}

For an integral coweight $\check\lambda$ let us consider the
corresponding point $t^{\check\lambda}\in G\ppart$. We will also think
of $\check\lambda$ as an element of $W_{\on{aff}}$ corresponding to this
orbit. Note that if $\lambda$ is dominant, the orbit of $I\cdot
t^{\check\lambda}\subset \Fl_G$ has the property that under the
projection $\Fl_G\to \Gr_G$ it maps one-to-one.

We have already established the transformation property of Wakimoto modules
with respect to convolution with $\wt{j}_{\kappa,w}$ for $w\in W$, see 
\lemref{wakimoto from one-another}. Now we would like to 
study their behavior with respect to convolution with 
$\wt{j}_{\kappa,\check\lambda}$. 

\medskip

Note that we have a natural adjoint action of $H\ppart$ on
$H^0_{DR}(\D^\times,\wh{L}_{\fh,\kappa})$, and 
similarly for the Baer negative extension. Thus,
we obtain that $H\ppart$ acts on the categories
$\fH_\kappa\mod$ and $\fH'_\kappa\mod$. For
$t^{\check\lambda}\in H\ppart$  we will denote
the corresponding functor by $\CR\mapsto t^{\check\lambda}\star \CR$.

The following property of Wakimoto modules will play a crucial role:

\begin{prop}   \label{lambda equivariance}
For a dominant $\lambda$ we have
$$\wt{j}_{\kappa,\check\lambda}\underset{I^0}\star
\BW^{w_0}_\kappa(\CR)\simeq \BW^{w_0}_\kappa(t^{w_0(\check\lambda)}\star
\CR).$$
\end{prop}

\ssec{Proof of \propref{lambda equivariance}}

Consider the subscheme $I^0\cdot t^{\check\lambda}\cdot I^0\subset
G\ppart$. Clearly, there exists a unique irreducible object of
$\fD^{\on{ch}}(G)_\kappa\mod^{I^0,I^0}$, supported on this subset. Let us
denote it by $\Dist_G^{\on{ch}}(I^0 t^{\check\lambda} I^0)_\kappa$. In
particular, for $\check\lambda=0$ we recover
$\Dist^{\on{ch}}_G(I^0)_\kappa$.

We have
$$\wt{j}_{\kappa,\check\lambda}\underset{I^0}\star
\Dist^{\on{ch}}_G(I^0)_\kappa\simeq \Dist_G^{\on{ch}}(I^0t^{\check\lambda}
I^0)_\kappa\simeq \Dist^{\on{ch}}_G(I^0)_\kappa \underset{I^0}\star
\wt{j}_{\kappa',\check\lambda}.$$

Therefore, by \propref{wakimoto as reduction}, we have to show that
$$\left(\Dist^{\on{ch}}_G(I^0)_\kappa\underset{I^0}\star
\wt{j}_{\kappa',\check\lambda}\right)
\underset{\fb^-\ppart,t\fb^-[[t]]}\torsemiinf \CR \simeq
\left(\Dist^{\on{ch}}_G(I^0)_\kappa\right)
\underset{\fb^-\ppart,t\fb^-[[t]]}\torsemiinf (t^{w_0(\check\lambda)}
\star \CR).$$

\medskip

Let us write $I^0_+=I^0\cap B[[t]]$ and $I^0_-=I^0\cap N^-[[t]]$, and
recall that the product map defines an isomorphism $I^0=I^0_+\cdot
I^0_-$. Note also that $\on{Ad}_{t^{\check\lambda}}(I^0_+)\subset
I^0$.

Hence,
$$\Dist^{\on{ch}}_G(I^0)_\kappa\underset{I^0}\star
\wt{j}_{\kappa',\check\lambda}\simeq \Dist^{\on{ch}}_G(I^0)_\kappa\star
\delta_{t^{\check\lambda}}\star \uBC_{I^0_-}\otimes
\det\left(t\fn^-[[t]]/t\fn^-[[t]]\cap
\on{Ad}_{t^{-\check\lambda}}(t\fn^-[[t]])[1]\right).$$

\medskip

Therefore, by \lemref{step-by-step convolution} we obtain that
\begin{align*}
&\left(\Dist^{\on{ch}}_G(I^0)_\kappa\star
\wt{j}_{\kappa',\lambda}\right)\underset{\fb^-\ppart,
t\fb^-[[t]]}\torsemiinf \CR \simeq \\
&H^\semiinf\left(\fn^-\ppart,t\fn^-[[t]],
\left(\Dist^{\on{ch}}_G(I^0)_\kappa\star
\delta_{t^{\check\lambda}}\right)\right)
\underset{\fh\ppart,t\fh[[t]]}\torsemiinf \CR \otimes
\det\left(t\fn^-[[t]]/\on{Ad}_{t^{-\check\lambda}}(t\fn^-[[t]])[1]\right).
\end{align*}

For a $L'_{\fb^-,\kappa}$-module $\CM$ we have
$$H^\semiinf(\fn^-\ppart,t\fn^-[[t]],t^{\check\lambda}\star \CM)\simeq
t^{\check\lambda}\star H^\semiinf
(\fn^-\ppart,\on{Ad}_{t^{\check\lambda}}(t\fn^-[[t]]),\CM),$$ as
$\fH'_\kappa$-modules.  

Hence, the expression above can be rewritten as
$$\Biggl(H^\semiinf\left(\fn^-\ppart,t\fn^-[[t]],
\Dist^{\on{ch}}_G(I^0)_\kappa \right)\star t^{\check\lambda}\Biggr)
\underset{\fh\ppart,t\fh[[t]]}\torsemiinf \CR,$$ where we have
absorbed the determinant line into changing the lattice
$\on{Ad}_{t^{-\check\lambda}}(t\fn^-[[t]])\mapsto t\fn^-[[t]]$.  The
latter can, in turn, be rewritten as
$$H^\semiinf\left(\fn^-\ppart,t\fn^-[[t]],
\Dist^{\on{ch}}_G(I^0)_\kappa\right)
\underset{\fh\ppart,t\fh[[t]]}\torsemiinf (t^{w_0(\check\lambda)}\star
\CR),$$ which is what we had to show.
 \footnote{The replacement of $\cla$ by
$w_0(\cla)$ comes from the fact that the identifications
$B/N\simeq H\simeq B^-/N^-$ differ by the automorphism $w_0$
of $H$.}

\section{Wakimoto modules at the critical level}     \label{Wak crit}

\ssec{}     \label{wakimodules}

In this section we will consider in more detail Wakimoto modules at
the critical level. By \propref{Miura=CADO} and using the isomorphism
$\fH'_\crit\simeq \fH_\crit$, if $\CR$ is a quasi-coherent sheaf on
$\ConHDt$, we can define Wakimoto modules $\BW^w_\crit(\CR)$
for $w\in W$. 

\medskip

Note that for any $\CR$, the Wakimoto module $\BW^w_\crit(\CR)$
carries an action of $\fH_\crit$ by transport of structure, and the
isomorphism
\begin{equation}  \label{descr of critical wakimoto}
\BW^w_\crit(\CR)\simeq
\Dist^{\on{ch}}_{N^{ww_0}}\left(\on{ev}^{-1}(N^{ww_0}\cap N)\right)\otimes
\CR
\end{equation}
of \eqref{descr of W(R)} is compatible with the $\fH_\crit$-actions. 

\medskip

In particular, from \thmref{FF with Miura} combined with 
\propref{MT and residue} and \eqref{even more wakimoto from one-another},
we obtain the following result:

Recall that a weight
$\lambda$ is called anti-dominant, if $\langle \lambda,\check
\alpha\rangle\notin \BZ^{>0}$ for any $\alpha\in \Delta^+$, or
equivalently, if the intersection of the two sets
$\{\lambda-\on{Span}_{\BZ^+}(\Delta^+)\}$ and $\{w(\lambda),\, w\in
W\}$ consists only of the element $\alpha$.

\begin{cor} \label{wakimoto with lambda are flat} 
The action of the center $\fZ_\fg\simeq \Fun(\Op(\D^\times))$ on
$\BW^w_{\crit,\lambda}$ factors through
$\fZ^{\on{RS},\varpi(-\lambda-\rho)}_\fg\simeq
\Fun(\Op^{\RS,\varpi(-\lambda-\rho)})$. Moreover, if $w^{-1}(\lambda+\rho)$ is
dominant, then $\BW^w_{\crit,\lambda}$ is flat over
$\fZ^{\on{RS},\varpi(-\lambda-\rho)}_\fg$.
\end{cor}

Another useful observation is the following:

\begin{prop} \label{endomorphisms of wakimoto}
Let $\CR_1,\CR_2$ be two $\fH_\crit$-modules, on which $\fh
\subset\Gamma(X,\fH_\crit)$ acts by the same scalar. Then for any
$w\in W$ the map
$$\Hom_{\fH_\crit}(\CR_1,\CR_2)\to
\Hom_{\hg_\crit}(\BW^w_\crit(\CR_1),\BW^w_\crit(\CR_2))$$
is an isomorphism.
\end{prop}

\begin{proof}

Since $\CR$ is a subspace of $\BW^w_\crit(\CR)$, the fact that the map
in question is injective is evident. Let us prove the surjectivity.

It will be more convenient to work with $'\BW^w_\crit(\CR_i)$ instead
of $\BW^w_\crit(\CR_i)$, $i=1,2$. As in \eqref{descr of critical
wakimoto}, we have an identification
$$'\BW^w_\crit(\CR_i)\simeq \Dist^{\on{ch}}_{N}\left(\on{ev}^{-1}(N\cap
N^{w_0w^{-1}})\right)\otimes \CR_i,$$ respecting the actions of
$\fn\ppart$ and $\fH_\crit$.

Let us first analyze the space of endomorphisms
$\Dist^{\on{ch}}_{N}\left(\on{ev}^{-1}(N\cap N^{w_0w^{-1}})\right)$ as a
$\fn\ppart$-module. We obtain, as in \secref{endomorphisms of
semijective}, that the map $\fr_\fn:L_\fn\to \fD^{\on{ch}}(N)$ has the
property that the image of $U\left(\fn\ppart\right)$ is dense in
$\End_{\fn\ppart}\Bigl(\Dist^{\on{ch}}_{N}\left(\on{ev}^{-1}(N\cap
N^{w_0w^{-1}})\right)\Bigr)$.

By the assumption on the $\fh$-action, and arguing as in the proof of
\propref{centralizer on G/B}, we obtain that any map of vector spaces
$'\BW^w_\crit(\CR_1)\to {}'\BW^w_\crit(\CR_2)$, compatible with the
action of $\fn\ppart$ and $\fh\subset \Gamma(X,\CA_{\fg,crit})$ has the
form
$$\on{Id}_{\Dist^{\on{ch}}_{N}\left(\on{ev}^{-1}(N\cap
N^{w_0w^{-1}})\right)}\otimes \varphi,$$ where $\varphi$ is some map
$\CR_1\to \CR_2$ as vector spaces. To prove that $\varphi$ is a map of
$\fH_\crit$-modules, we argue as follows:

\medskip

Recall that for a $\hg_\crit$-module $\CM$, the semi-infinite cohomology
$$H^\semiinf(\fn\ppart,t\fn[[t]]+\fn\cap \fn^{w_0w^{-1}},\CM)$$ is
naturally a $\fH'_\crit$-module. We will regard it as a
$\fH_\crit$-module via the isomorphism $\fH_\crit\simeq
\fH'_\crit$. Recall the isomorphism
\begin{equation} \label{cohomology of wakimoto}
\CR\simeq H^\semiinf(\fn\ppart,t\fn[[t]]+\fn\cap
\fn^{w_0w^{-1}},{}'\BW^w_\crit(\CR))
\end{equation}
given by \propref{acyclicity of wakimoto}. From \lemref{H symmetry} we
obtain the following

\begin{lem}
The isomorphism \eqref{cohomology of wakimoto}
respects the $\fH_\crit$-module structures.
\end{lem}

{}From the construction of the isomorphism of \propref{acyclicity of
wakimoto} it is easy to see that any map $'\BW^w_\crit(\CR_1)\to
{}'\BW^w_\crit(\CR_2)$ of the form $\on{Id}\otimes \phi$ induces on
the left-hand side of \eqref{cohomology of wakimoto} the endomorphism
equal to $\varphi$. Hence, the above lemma implies that $\phi$
respects the $\fH_\crit$-actions.

\end{proof}

\ssec{}    

We will now recall a crucial result of \cite{F:wak} (see Proposition
6.3 and Remark 6.4) that establishes isomorphisms between Wakimoto
modules and Verma modules. Note that in \cite{F:wak} the module
$\BW^{w_0}_{\crit,\lambda}$ is denoted by $W_{\la,\kappa_c}$.

\begin{prop}   \label{Wakimoto=Verma}
Let $\lambda$ be such that $\lambda+\rho$ is anti-dominant.
Then $\BW^{w_0}_{\crit,\lambda}\simeq \BM_\lambda$.
\end{prop}

\begin{proof}

First, we claim that when $\lambda$ is anti-dominant,
the $\fg$-module $M^{w_0}_\lambda$ (see \secref{wakimoto lambda})
is in fact isomorphic to the Verma module $M_\lambda$.

Indeed, let us note first that $M^{w_0}_\lambda$ has a vector of
highest weight $\lambda$, i.e., there is a morphism $M_\lambda\to
M^{w_0}_\lambda$. Now, it is well-known that the anti-dominance
condition on $\lambda+\rho$ implies that $M_\lambda$ is irreducible,
hence the above map is injective. The assertion follows now from the
fact that the two modules have the same formal character.  This is a
prototype of the argument proving the proposition.

\medskip

By \secref{wakimoto lambda} we have a map in one direction 
\begin{equation} \label{map from Verma to Wakimoto}
\BM_\lambda\to \BW^{w_0}_{\crit,\lambda}
\end{equation}
and we claim that it is an isomorphism. We will regard both sides
of \eqref{map from Verma to Wakimoto} as modules over the
{\it Kac-Moody} algebra $\BC\cdot t\partial_t\ltimes \hg_\crit$,
where we normalize the action of $t\partial_t$ so that it
annihilates the generating vector in $\BM_\lambda$, and
it acts on $\BW^{w_0}_{\crit,\lambda}$ by loop rotation.

The map \eqref{map from Verma to Wakimoto} clearly
respects this action. Moreover, both sides have 
well-defined formal characters with respect to the extended
Cartan subalgebra $\BC\cdot t\partial_t\oplus \fh\oplus \BC\one$,
and a computation shows that these characters are equal.
Therefore, the map \eqref{map from Verma to Wakimoto} 
is surjective if and only if it is injective.

\medskip

Suppose that the kernel of the map in question is non-zero.
Let $v\in \BM_\lambda$ be a vector of a {\it highest weight}
with respect to $\BC\cdot t\partial_t\oplus \fh\oplus \BC\one$;
let us denote this weight by $\wh{\mu}$. 

Then the quotient $\BW^{w_0}_{\crit,\lambda}/\on{Im} \BM_\lambda$ also
contains a vector, call it $v'$, of weight $\wh{\mu}$.  Moreover, by
assumption, $v'$ projects non-trivially to the space of coinvariants
$\left(\BW^{w_0}_{\crit,\lambda}\right)_{\fn^-[t^{-1}]\oplus
t^{-1}\fb[t^{-1}]}$

However, from \eqref{descr of critical wakimoto}, it follows that the
projection
$$\Fun(N[[t]]) \to
\left(\BW^{w_0}_{\crit,\lambda}\right)_{t^{-1}\fn[t^{-1}]\oplus
t^{-1}\fh[t^{-1}]}$$ is an
isomorphism. Therefore, $\wh{\mu}$ must be of the form

\begin{equation} \label{form of mu}
\wh{\mu}:=(-n,\lambda-\beta,-\check h), \,\,n\in \BZ^{\geq 0},\,\,
\beta\in \on{Span}^+(\Delta^+).
\end{equation}

\medskip

We will now use the Kac-Kazhdan theorem \cite{KK} that describes the
possible highest weights of submodules of a Verma module. This theorem
says that there must exist a sequence of weights
$$(0,\lambda,-\check{h})=\wh{\mu_1},\, \wh{\mu_2},...,\wh{\mu}_{n-1},
\,\wh{\mu}_n=\wh{\mu}$$
and a sequence of positive affine roots
$\alpha_{\on{aff},k}$ such that
\begin{equation} \label{tras cond}
\wh{\mu}_{k+1}=\wh{\mu}_k-b_k\cdot \alpha_{\on{aff},k}
\end{equation}
with $b_k\in \BZ^{>0}$ and such that 
$$b_k\cdot (\alpha_{\on{aff},k},\alpha_{\on{aff},k})=2\cdot
(\alpha_{\on{aff},k},\wh{\mu}_k+\rho_{\on{\on{aff}}}),$$ where $(\cdot,\cdot)$
is the invariant inner product on the Kac-Moody algebra.

Let us write $\wh{\mu}_k=(n_k,\mu_k,-\check h)$ and
$$
\begin{cases}
& \alpha_{\on{aff},k}=(m_k,\epsilon_k\cdot \alpha_k,0),\,\, m_k\geq 0,\,\,
\alpha_k\in \Delta^+,\,\,\epsilon_k=\pm 1 \text{ if } \alpha_{\on{aff},k}
\text{ is real,} \\ & \alpha_{\on{aff},k}=(m_k,0,0) \text{ if }
\alpha_{\on{aff},k} \text{ is imaginary} .
\end{cases}
$$

In the latter case we obtain $\mu_{k+1}=\mu_k$. In the former case we have
$$b_k=\langle\wh{\mu}_k+\rho_{\on{aff}},\check
\alpha_{\on{aff},k}\rangle,$$ and since $\rho_{\on{aff}}=(0,\rho,\check h)$
we obtain that $b_k=\epsilon_k\cdot \langle
\mu_k+\rho,\alpha_k\rangle$, implying that
$$(\mu_{k+1}+\rho)=s_{\alpha_k}(\mu_k+\rho),$$
regardless of the sign of $\epsilon_k$.

In particular, we obtain that $(\lambda+\rho)-\beta$ belongs to the
$W$-orbit of $\lambda+\rho$, but this contradicts the anti-dominance
of $\lambda+\rho$.

\end{proof}

\ssec{}  \label{W & V}

We will use the above proposition to derive information about
the structure of other Wakimoto and Verma modules.

\begin{cor}  \label{other Wakimoto=Verma}
For $\lambda$ such that $\lambda+\rho$ is anti-dominant and $w\in W$
we have an isomorphism
$$\BW^{ww_0}_{\crit,w(\lambda+\rho)-\rho}\simeq
\BM_{w(\lambda+\rho)-\rho}.$$
\end{cor}

\begin{proof}

Let us assume that $\lambda$ is integral. In this case all
$\BW^w_{\crit,\lambda}$ and $\BM_{w(\lambda+\rho)-\rho}$ are
$I$-integrable, and we can use the convolution action of D-modules on
$G/B\subset \Fl_G$ to pass from one-another.

(If $\lambda$ is not integral, the proof is essentially the same, when
instead of $B$-equivariant D-modules on $G/B$, we will use
$\lambda$-twisted D-modules and replace the $B$-equivariant category
by a $\lambda$-twisted version.)

It is known that for $\lambda$ anti-dominant,
$j_{w,!}\underset{B}\star M_\lambda=M_{w(\lambda+\rho)-\rho}$.
Hence,
$j_{w,!}\underset{I}\star \BM_\lambda=\BM_{w(\lambda+\rho)-\rho}$.
This implies the corollary in view of \eqref{more waimoto from one-another}.

\end{proof}

Since for every weight $\lambda'$ there exists an element of the Weyl
group such that $\lambda'=w(\lambda+\rho)-\rho$ with $\lambda$
anti-dominant, every Verma module $\BM_{\lambda'}$ is isomorphic to an
appropriate Wakimoto module. By combining this with \propref{MT and
residue} and \propref{endomorphisms of wakimoto} we obtain the
following statement.

\begin{cor} \label{Verma is flat}
The module $\BM_{\lambda}$ is flat over
$\fZ^{\RS,\varpi(-\lambda-\rho)}_\fg$.  The map
$$\fZ^{\RS,\varpi(-\lambda-\rho)}_\fg\to
\on{End}_{\hg_\crit}(\BM_{\lambda})$$ is an isomorphism.
\end{cor}

Let us give an additional proof of the second assertion:

\begin{proof}

As above, we can reduce the statement to the case when
$\lambda+\rho$ is itself anti-dominant, and $\BM_\lambda\simeq
\BW^{w_0}_{\crit,\lambda}$.

In the latter case, we have to show that the embedding of
$\pi_{w_0(\lambda)}$ into the subspace of vectors of weight $\lambda$
in $\left(\BW^{w_0}_{\crit,\lambda}\right){}^{\Lie(I^0)}$ is an
isomorphism.

Consider the bigger subspace
$\left(\BW^{w_0}_{\crit,\lambda}\right){}^{\fn[[t]]}$.  As in the
proof of \propref{centralizer on G/B}, this subspace is isomorphic to
$\on{Ind}^{\fn\ppart}_{\fn[[t]]}(\BC)\otimes \pi_{w_0(\lambda)}$,
which implies that the vectors of weight $\lambda$ belong to $1\otimes
\pi_{w_0(\lambda)}$.

\end{proof}

Next, we shall prove \propref{support of V lambda}:

\begin{proof}

Consider the Wakimoto module $\BW^{w_0}_{\crit,\lambda}$ when
$\lambda$ is dominant. We claim that $\fZ_\fg$ acts on it via
$\fZ_\fg^{\lambda,\reg}$. This follows by combining
\propref{regular Miura} with \thmref{FF with Miura} and the fact that the
isomorphism of \propref{Miura=CADO} sends the chiral
$\fH_\crit$-module $\pi_\mu$ to
$\Fun\bigl(\on{Conn}_{\check{H}}(\omega_X^\rho)^{\RS,\mu}\bigr)$.

\medskip

Composing the map \eqref{Verma to wakimoto} with the natural embedding
$\BV^\lambda_\crit\to \BM^\vee_{\crit,\lambda}$ we obtain a map
$\BV^\lambda_\crit \to \BW^{w_0}_{\crit,\lambda}$,
which can be shown to be injective. \footnote{We will supply a proof
in the next paper in the series.} We obtain that the ideal in the center
that annihilates $\BW^{w_0}_{\crit,\lambda}$, annihilates
$\BV_\crit^\lambda$ as well, which is what we had to show.

\end{proof}

\ssec{}

Finally, let us derive a corollary of \propref{lambda equivariance} at
the critical level.  In this case the adjoint action of $H\ppart$ on
$H^0_{DR}(\D^\times,\wh{L}_{\fh,\crit})$ is trivial, and hence we
obtain the following

\begin{cor}  \label{lambda invariance}
For a dominant coweight $\check\lambda$ and a $\fH_\crit$-module $\CR$
we have an isomorphism
$$\wt{j}_{\check\lambda}\underset{I^0}\star \BW^{w_0}_\crit(\CR)\simeq
\BW^{w_0}_\crit(\CR).$$
\end{cor}

Suppose now that $\CR$ in $H[[t]]$-integrable. In this case we obtain
that $$\wt{j}_{\check\lambda}\underset{I^0}\star
\BW^{w_0}_\crit(\CR)\simeq j_{\check\lambda,*}\underset{I}\star
\BW^{w_0}_\crit(\CR)\simeq \BW^{w_0}_\crit(\CR),$$ where both the LHS
and RHS are canonically defined, i.e., are independent of the choice
of a representative $t^{\check\lambda}\in G\ppart$. However, the
isomorphism between them, given by \corref{lambda invariance}, does
depend on this choice. In what follows we will need a more precise
version of the above result:

\begin{cor}  \label{lambda equivariance of Verma}
We have an isomorphism
$$j_{\check\lambda,*}\underset{I}\star \BM_{-2\rho}\simeq
\BM_{-2\rho}\otimes \omega_x^{\langle -\rho,\check\lambda\rangle},$$
where $\omega_x$ is the fiber of $\omega_X$ at $x\in X$, compatible
with the natural actions of $\on{Lie}(\on{Aut}(\D))$ on both sides.
\end{cor}

Note that $\langle -\rho,\check\lambda\rangle$, appearing in the
corollary, may be a half-integer. In the above formula the expression
$\omega_x^{\langle -\rho,\check\lambda\rangle}$ involves a choice of a
square root of $\omega_X$, as does the construction of the critical
line bundle on $\Gr_{G_{ad}}$.  However, the character of
$\on{Lie}(\on{Aut}(\D))$ on $\omega_x^{\langle
-\rho,\check\lambda\rangle}$ is, of course, independent of this
choice.

\medskip

\begin{remark} By considering the action of the renormalized
universal enveloping algebra as in \secref{action of renorm}, one
shows that, more generally, there is an isomorphism
$$j_{\check\lambda,*}\underset{I}\star \BW^{w_0}_{\crit,\mu}\simeq
\BW^{w_0}_{\crit,\mu}\otimes \omega_x^{\langle
\mu+\rho,\check\lambda\rangle},$$ compatible with the
$\on{Lie}(\on{Aut}(\D))$-actions.
\end{remark}

\begin{proof}

The existence of an isomorphism stated in the corollary follows by
combining \propref{Wakimoto=Verma} and \corref{lambda invariance}. By
\propref{endomorphisms of wakimoto}, we obtain that there exists a
line, acted on by (a double cover of) $\on{Aut}(\D)$, and a canonical
isomorphism
$$j_{\check\lambda,*}\underset{I}\star \BM_{-2\rho}\simeq
\BM_{-2\rho}\otimes \fl,$$ compatible with the
$\on{Lie}(\on{Aut}(\D))$-actions.

We have to show that the character of $\on{Lie}(\on{Aut}(\D))$,
corresponding to $\fl$, equals that of $\omega_x^{\langle
-\rho,\check\lambda\rangle}$.  Let $t\partial_t\in
\on{Lie}(\on{Aut}(\D))$ be the Euler vector field, corresponding to
the coordinate $t$ on $\D$. It suffices to show that the highest
weight of $j_{\check\lambda,*}\underset{I}\star \BM_{-2\rho}$ with
respect to $\BC\cdot t\partial_t\oplus \fh$ equals $(-\langle
\rho,\check\lambda\rangle,-2\rho)$.

The module in question identifies with
$\Gamma(\Fl_G,j_{\check\lambda,*})$.  The highest weight line in
$\Gamma(\Fl_G,j_{\check\lambda,*})$ consists of $I^0$-invariant
sections of this D-module, that are supported on the $I$-orbit of
$t^{\check\lambda}$. Now the fact that $t\partial_t$ acts on this line
by the character equal to $-\langle \rho,\check\lambda\rangle$ is a
straightforward calculation, as in \cite{BD}, Sect. 9.1.

\end{proof}

\bigskip

\newpage

\vspace*{10mm}

{\Large \part{Proof of \mainthmref{equiv of quot}}}

\vspace*{10mm}


The goal of this Part is to prove \mainthmref{equiv of quot}:

\medskip
\noindent {\em There is an equivalence of categories}
$$^f\sF:\Catf \simeq \on{QCoh}\left(\on{Spec}(h_0)\times
\nOp\right).$$

\medskip

In \secref{module Ppi} we introduce the module $\Ppi$ as induced from
the big projective module $\Pi$ over the finite-dimensional Lie
algebra $\fg$. We first review the properties of $\Pi$, and the
corresponding properties of $\Ppi$, related to the notion of partial
integrability. The functor $^f\sF$ is then defined as
$$
^f\sF(M) = \Hom(\Ppi,M),
$$
and we state \mainthmref{RHom from Pi} which asserts that this functor
is exact.

As we shall see later (see \secref{proof of main thms}), both
\mainthmref{RHom from Pi} and \mainthmref{equiv of quot} follow once
we can compute $\on{RHom}_{\DCat}(\Ppi,\BM_{w_0})$, where $\BM_{w_0}$
is the corresponding Verma module over $\hg_\crit$. Also in
\secref{proof of main thms} (see \propref{Ext from Pi to Verma}) we
show that it is sufficient to compute
$\on{RHom}_{\DCatr}(\Ppi_\reg,\BM_{w_0,\reg})$, where $\Ppi_\reg$ and
$\BM_{w_0,\reg}$ are the restrictions of the corresponding modules to
the subscheme $\Op^\reg\subset \Op^\nilp$.

The computation of $\on{RHom}_{\DCatr}(\Ppi_\reg,\BM_{w_0,\reg})$ is
carried out in \secref{Pi via Gr}. We reduce it to a calculation
involving D-modules on the affine Grassmannian once we can identify
$\Ppi_\reg$ as sections of some specific critically twisted D-module
on $\Gr$. The latter identification is given by \thmref{Wakimoto via
Grassmannian}.  This theorem is proved in \secref{proof of W v G} by a
rather explicit argument.

Having proved \mainthmref{equiv of quot}, we compare in \secref{comp
with semiinf} the functor $\Hom_{\Cat}(\Ppi,\cdot)$ with the one given
by semi-infinite cohomology with respect to the Lie algebra
$\fn^-\ppart$ against a non-degenerate character. We show that the two
functors are isomorphic. We also express the semi-infinite cohomology
of $\fn\ppart$ with coefficients in a $\hg_\crit$-module of the form
$\Gamma(\Gr_G,\CF)$, where $\CF$ is a critically twisted D-module on
$\Gr_G$, in terms of the de Rham cohomologies of the restrictions of
$\CF$ to $N\ppart$-orbits in $\Gr_G$.

\bigskip

\section{The module $\Ppi$}     \label{module Ppi}

\ssec{}

Recall from \secref{induction from O} that $\CO_0$ denotes the
subcategory of the category $\CO$ of $\fg$-modules, whose objects are
modules with central character $\varpi(\rho)$.  According to
\cite{BB}, the functor of global sections induces an equivalence
between the category of $N$-equivariant (or, equivalently,
$B$-monodromic) left D-modules on $G/B$ and $\CO_0$.

To simplify our notation slightly, we will use the notation $M_w$
instead of $M_{w(\rho)-\rho}$ and $M^\vee_w$ instead of
$M_{w(\rho)-\rho}^\vee$. We will denote by $L_w$ the irreducible
quotient of $M_w$. By $\BM_w$, $\BM^\vee_w$ and $\BL_w$ we will denote
the corresponding induced representations of $\ghat$ at the critical
level.

By definition, objects of $\CO_0$ are $N$-integrable, and the
condition on the central character implies that they are in fact
$B$-monodromic. Hence, every object $M\in \CO_0$ carries an action of
the commutative algebra $\fh$. This is the obstruction to being
$B$-equivariant.  (The notions of $B$-integrability (equivalently,
$B$-equivariance), $N$-integrability and $B$-monodromicity are
defined as their $I$- and $I^0$-counterparts and make sense in any
category $\CC$ with a Harish-Chandra action of $G$.)

\begin{lem} \label{h_0}
For every $M\in \CO_0$, the action of $\Sym(\fh)$ on $M$ factors through
$\Sym(\fh)\to h_0$.
\end{lem}

The lemma follows, e.g., from the localization theorem of \cite{BB}.
Thus, we obtain that the algebra $h_0$ maps to the center of $\CO_0$.
In fact, it follows from \cite{Be} that $h_0$ is isomorphic to the
center of $\CO_0$.

\medskip

As in \secref{partially integrable}, we will call an object $M\in \CO_0$ 
partially integrable if it admits a filtration such that each successive
quotient is integrable with respect to a parabolic subalgebra
$\fb+{\mathfrak {sl}}_2^\imath$ for some $\iota\in \CI$. 
This notion makes sense in an arbitrary category with a Harish-Chandra
action of $G$.

We will denote by $^f\CO_0$ the quotient abelian category of $\CO_0$
by the subcategory of partially integrable objects. We will denote by
$M\mapsto {}^fM$ the projection functor $\CO_0\to {}^f\CO_0$.

\medskip

Let $\Pi$ be a "longest" indecomposable projective in $\CO_0$.
By definition:
$$\Hom(\Pi,L_{w_0})=\BC,\,\, \Hom(\Pi,L_w)=0 \text{ if } w\neq w_0.$$
Moreover, $\Pi$ is known to be isomorphic (non-canonically) to
its contragredient dual. We have the following result (see \cite{BG}):

\begin{lem}  \label{Pi in quotient}
\hfil

\smallskip

\noindent{\em (1)}
The map $h_0\to \End(\Pi)$ is an isomorphism.

\noindent{\em (2)}
The functor $M\mapsto \Hom(\Pi,M)$ induces an
equivalence $^f\CO_0\to h_0\mod$.
\end{lem}

By construction, the image $^f\Pi$ of $\Pi$ in $^f\CO_0$ identifies
with the free $h_0$-module with one generator. The maps
$M_{w_0}\to L_{w_0}\to M^\vee_{w_0}$ induce isomorphisms
$^f M_{w_0}\to {}^f L_{w_0}\to {}^f M{}^\vee_{w_0}$, and all
identify with the trivial $h_0$-module $\BC$.
 
 \medskip

We will now recall the construction of $\Pi$ as 
\begin{equation} \label{Pi=Xi}
\Pi=\Gamma(G/B,\Xi),
\end{equation}
where $\Xi$ is a certain left D-module on $G/B$.

\ssec{}   \label{character psi}

To describe $\Xi$ we need to introduce some notation, which will
also be used in the sequel. Let $\psi:N^-\to \BG_a$ be a 
non-degenerate character. By a slight abuse of notation, we will
denote also by $\psi$ its differential: $\fn^-\to \BC$.

Let $\epsi$ denote the pull-back of the "$e^x$" D-module from $\BG_a$
to $N^-$. This is a "character sheaf" in the sense of \secref{N, psi
equiv}.

If $N^-$ acts (in the Harish-Chandra sense) 
on a category $\CC$, we will denote by $\CC^{N^-,\psi}$
the corresponding $(N^-,\psi)$-equivariant category 
(see \secref{N, psi equiv}), and by $D(\CC)^{N^-,\psi}$ the corresponding
triangulated category. Since $N^-$ is unipotent, the natural forgetful
functor $D(\CC)^{N^-,\psi}\to D(\CC)$ is fully faithful,
see \secref{N, psi equiv}.

Following \secref{N, psi equiv}, we will denote denote by
$\on{Av}_{N^-,\psi}$ the functor
$$\CM\mapsto \epsi\star \CM \otimes\det(\fn^-[1])^{-1}:D(\CC)\to
D(\CC)^{N^-,\psi}.$$ This functor is the right adjoint and a left
quasi-inverse to $D(\CC)^{N^-,\psi}\to D(\CC)$.

\begin{lem}  \label{psi kills int}
Suppose that $\CC$ is endowed with a Harish-Chandra action of $G$, and
let $\CM\in \CC^{B,m}$ be partially integrable, then
$\on{Av}_{N^-,\psi}(\CM)=0$.
\end{lem}

\begin{proof}

We can assume that $\CM$ is an object of $\CC$ integrable with respect
to a parabolic subgroup $P^\iota$ for some $\iota\in\CI$.  Then the
convolution $\epsi\star \CF$ factors through the direct image of
$\epsi$ under $N^-\hookrightarrow G\twoheadrightarrow G/P^\iota$, and
the latter is clearly $0$.

\end{proof}

\medskip

For example, if $N^-$ acts on a scheme $Y$, in this way we obtain the
category of $(N^-,\psi)$-equivariant D-modules on $Y$. In other words,
its objects are D-modules $\CF$ on $Y$, together with an isomorphism
$$\on{act}^*(\CF)\simeq \epsi\boxtimes \CF\in \fD(N^-\times Y)\mod,$$
compatible with the restriction to the unit section and associative
in the natural sense. One can show that in this case the functor
$$D(\fD(Y)\mod^{N^-,\psi})\to D(\fD(Y)\mod)^{N^-,\psi}$$ is
an equivalence.

If we restrict ourselves to holonomic D-modules, or, rather, if we take
the corresponding triangulated category (which, by definition, is the
full subcategory of $D(\fD(Y)\mod)$, consisting of complexes with
holonomic cohomologies), then in addition to the functor $\CF\mapsto
\epsi\star \CF$ we also have a functor
$$\CF\mapsto \epsi\overset{!}\star\CM: D(\fD(Y)_{\on{hol}}\mod)\to
D(\fD(Y)_{\on{hol}}\mod)^{N^-,\psi},$$ corresponding to taking direct
image with compact supports. This functor, tensored with
$\det(\fn^-[1])$, is the left adjoint and a left quasi-inverse to the
tautological functor $D(\fD(Y)_{\on{hol}}\mod)^{N^-,-\psi}\to
D(\fD(Y)_{\on{hol}}\mod)$.

\begin{prop}  \label{psi averaging} 
Suppose that $Y$ is acted on by $G$. Then for $\CF\in
D(\fD(Y)_{\on{hol}}\mod)^{B,m}$ the canonical arrow:
$\epsi\overset{!}\star \CF\to \epsi\star \CF$ is an isomorphism. In
particular, the functor
$$D(\fD(Y)_{\on{hol}}\mod)^{B,m}\to D(\fD(Y)_{\on{hol}}\mod)
\overset{\epsi\star \cdot}\longrightarrow
D(\fD(Y)_{\on{hol}}\mod)^{N^-,-\psi}$$ is exact.
\end{prop}

\begin{proof}

It is enough to analyze the functor $\CF\mapsto \epsi\star \CF$ on the
subcategory $\fD(Y)_{\on{hol}}\mod^B$.

The basic observation is that the D-module
$$\Dist_{G/B}(N^-,\psi):=\epsi\star \delta_{1_{G/B}}\in
\fD(G/B)\mod^{B^-,\psi},$$ which is by definition the $*$-extension of
$\epsi$ under $N^-\cdot 1_{G/B}\hookrightarrow G/B$, is clean. This
means that the $*$-extension coincides with the $!$-extension, or,
which is the same, that the arrow $\epsi\overset{!}\star
\delta_{1_{G/B}}\to \epsi\star \delta_{1_{G/B}}$ is an
isomorphism. (One easily shows that by observing that for any $\bg\in
G/B\setminus N^-$, the restriction of $\psi$ to its stabilizer in $N^-$ is
non-trivial.)  In particular, $\Dist_{G/B}(N^-,\psi)$ is the Verdier
dual of $\Dist_{G/B}(N^-,-\psi)$.

Note that for $\CF\in \fD(Y)_{\on{hol}}\mod^B$,
\begin{equation} \label{conv with psi as clean}
\epsi\star \CF\simeq \Dist_{G/B}(N^-,\psi)\underset{B} \star \CF,
\end{equation} 
and similarly for $\epsi\overset{!}\star \CF$.  This establishes the
assertion of the proposition.

\end{proof}

\ssec{}

After these preliminaries, we are ready to introduce $\Xi$:
$$\Xi:=\uBC_N\overset{!}\star \Dist_{G/B}(N^-,\psi)\otimes \det(\fn[1])
\in D(\fD(G/B)\mod)^N,$$
where 
$$\CF\mapsto \uBC_N\overset{!}\star\CF\otimes \det(\fn[1])^{\otimes
2}$$ is the functor $D(\fD(G/B)_{\on{hol}}\mod)\to
D(\fD(G/B)_{\on{hol}}\mod)^N$, left adjoint to the tautological functor
$D(\fD(G/B)_{\on{hol}}\mod)^N\to D(\fD(G/B)_{\on{hol}}\mod)$.  Explicitly,
$\uBC_N\overset{!}\star\cdot\otimes \det(\fn[1])$ is given by
convolution with compact supports with the constant D-module on $N$.

\begin{prop} \label{properties of Xi} \hfill

\smallskip

\noindent{\em (1)}   \label{ppts of Xi}
The complex $\Xi$ is concentrated in cohomological degree $0$.

\smallskip

\noindent{\em (2)}
$\Xi$ is projective as an object of $\fD(G/B)\mod^N$

\smallskip

\noindent{\em (3)}
$\Xi$ is non-canonically Verdier self-dual, i.e.
$\Xi\simeq \uBC_N\star \Dist_{G/B}(N^-,-\psi)\otimes \det(\fn[1])^{-1}$.

\smallskip

\noindent{\em (4)}
$\Xi$ is canonically independent of the choice of $\psi$.

\end{prop}

\begin{proof}

Consider the functor $\CF\mapsto \on{RHom}(\Xi,\CF)$ on the category 
$D(\fD(G/B)\mod)^N$. We have
$$\on{RHom}_{D(\fD(G/B)\mod)^N}(\Xi,\CF)\simeq
\on{RHom}_{D(\fD(G/B)\mod)}(\epsi\overset{!}\star
\delta_{1_{G/B}},\CF)\otimes \det(\fn[1]),$$ which, in turn, is
isomorphic to
$$\on{RHom}_{D(\fD(G/B)\mod)}(\delta_{1_{G/B}},\epsip\star \CF)\otimes
\det(\fn[1]),$$ where $\psi'=-\psi$.

By \lemref{psi averaging}, $\epsip\star \CF$ is concentrated in
cohomological degree $0$. Moreover, it is lisse near $1_{G/B}$. Hence,
the above $\on{RHom}$ is concentrated in cohomological degree $0$.

Now, we will use the fact that $D(\fD(G/B)\mod)^N$ is equivalent to the
derived category of the abelian category $\fD(G/B)\mod$. Then the above
property of $\on{RHom}$ implies simultaneously assertions (1) and (2) of 
the proposition.

\medskip

The above expression for $\on{RHom}(\Xi,\CF)$ also implies that it is
$0$ if $\CF$ is partially integrable, and
$\on{RHom}(\Xi,\delta_{1_{G/B}})$ is one-dimensional. This implies
that $\Xi$ corresponds to a projective cover of $\delta_{1_{G/B}}\in
\fD(G/B)\mod^N$, i.e., $\Gamma(G/B,\Xi)\simeq \Pi$.

Since it is known that contravariant duality on $\CO_0$ goes over to
Verdier duality on $\fD(G/B)\mod^N$, assertion (3) of the proposition
holds.

The fact that $\Xi$ is {\it non-canonically} independent of the choice
of $\psi$ also follows, since we have shown that \eqref{Pi=Xi} is
valid for any choice of $\psi$.  To establish that it is canonically
independent, we argue as follows:

\medskip

Let $\psi'$ be another non-degenerate character of $N^-$. Then there
exists an element $\bh\in H$, which, under the adjoint action of $H$
on $N^-$, transforms $\psi$ to $\psi'$.

Since $\Xi$ is $B$-monodromic, we have a canonical isomorphism of
D-modules $\bh^*(\Xi)\simeq \Xi$. However, from the construction of
$\Xi$, we have $\bh^*(\Xi)\simeq \Xi'$, where the latter is the D-module
constructed starting from $\psi'$.

\end{proof}

For any category $\CC$ with a Harish-Chandra action of $G$ we 
can consider the functor 
$$\CF\mapsto \Xi\underset{B}\star \CF:D(\CC)^B\to D(\CC)^N.$$

\begin{prop} \label{convolution with Xi} 
\hfill

\smallskip

\noindent{\em (1)} The above functor is exact, and it annihilates an
object $\CF\in\CC$ if and only if $\CF$ is partially integrable.

\smallskip

\noindent{\em (2)} For $\CM_1,\CM_2\in D(\CC)^B$ we have a
non-canonical but functorial isomorphism $$\on{RHom}_{D(\CC)}(\Xi\star
\CM_1,\CM_2)\simeq \on{RHom}_{D(\CC)}(\CM_1,\Xi\star\CM_2).$$

\end{prop}

\begin{proof}

Using \propref{ppts of Xi}(3), we can rewrite the functor in question as
$$\CF\mapsto \uBC_N\star \epsi\star \CF\otimes \det(\fn[1])^{-1}.$$
Hence, the fact that it annihilates partially integrable objects
follows from \lemref{psi kills int}.

Recall that the object $\Pi\in \CO_0$ is tilting, i.e., it admits two
filtrations: one, whose successive quotients are isomorphic to Verma
modules, and another, whose successive quotients are dual Vermas.
Hence, $\Xi$ also admits such filtrations, with subquotients being
$j_{w,!}$ and $j_{w,*}$, respectively.  It is clear that convolution
with the latter is right exact. The convolution with $j_{w,!}$, being
a quasi-inverse of the convolution with $j_{w^{-1},*}$, is therefore
left exact.This proves the exactness assertion of the proposition.

Finally, let us show that if $\CF$ is not partially integrable, then
$\Xi\star \CF\neq 0$. Let ${}^f\CC$ be the quotient category
of $\CC$ by the Serre subcategory of partially integrable objects.
Let $^f\CF$ be the image of $\CF$ in ${}^f\CC$.

We claim that the image of $\Xi\star \CF$ in ${}^f\CC$ is endowed
with an increasing filtration of length $|W|$, whose subquotients are
all isomorphic to $^f\CF$. This follows from the existence of the
filtration on $\Xi$ by $j_{w,!}$: Indeed, the cokernel of the map
$\delta_{1_{G/B}}\to j_{w,!}$ is partially integrable, hence
$^f\CF\to {}^f (j_{w,!}\underset{B}\star \CF)$ is an isomorphism.

\medskip

Now let us prove assertion (2) of the proposition. By \secref{HCh
adjunction},
$$\on{RHom}_{D(\CC)}(\Xi\underset{B}\star \CM_1,\CM_2)\simeq
\on{RHom}_{D(\CC)^B}(\CM_1,\wt{\Xi}\star \CM_2),$$
where $\wt{\Xi}$ is the corresponding dual D-module on $B\backslash G$.
Since $\CM_2$ was assumed $B$-equivariant,
$$\wt{\Xi}\star \CM_2\simeq (\wt{\Xi}\star \uBC_B)\underset{B}\star
\CM_2\otimes \det(\fb[1]).$$ Similarly, by the $B$-equivariance of
$\CM_1$,
$$\on{RHom}_{D(\CC)}(\CM_1,\Xi\star\CM_2)\simeq
\on{RHom}_{D(\CC)^B}(\CM_1,(\uBC_B\star \Xi)\star \CM_2).$$

Hence, it remains to see that 
$$\uBC_B\star \Xi
\simeq \wt{\Xi}\star \uBC_B\otimes\det(\fb[1])
\in D(\fD(G/B)\mod)^{B}\simeq D(\fD(G)\mod)^{B\times B}.$$

Using \propref{ppts of Xi}(3), the left-hand side is isomorphic to
$\on{Av}_{B\times B}(\Dist_G(N^-,\psi))$, and using \propref{ppts of
Xi}(4), the right-hand side is isomorphic to the same thing.

\end{proof}

\ssec{}

Let us return to representations of affine algebras at the critical
level. We define the module $\Ppi\in \hg_\crit\mod$ as
\begin{equation}  \label{defn of Ppi}
\Ppi=\on{Ind}^{\hg_\crit}_{\fg[[t]]\oplus \BC\one}(\Pi).
\end{equation}
By \secref{induction from O}, $\Ppi$ belongs to $\Cat$.

{}From the tilting property of $\Pi$, we obtain that $\Ppi$ admits two
filtrations: one whose subquotients are modules of the form $\BM_w$,
and another, whose subquotients are of the form $\BM_w^\vee$. Together
with \corref{Verma is flat} this implies:

\begin{cor} \label{Pi is flat}
The module $\Ppi$ is flat over $\fZ^\nilp_\fg$.
\end{cor}

Using our conventions concerning twisted D-modules on $\Fl_G$, we can
rewrite the definition of $\Ppi$ as
$$\Gamma(\Fl_G,\Xi),$$ where we think of $\Xi$ as living on $\Fl_G$
via $G/B\hookrightarrow \Fl_G$.

Note that for $\CM^\bullet\in D(\hg_\crit\mod)^I$, the convolution
$\Xi\underset{I}\star \CM^\bullet$ is tautologically the same as
$\Xi\underset{B}\star \CM^\bullet$, when we think of $\CM^\bullet$ is a $\fg$-module
via $\fg\hookrightarrow \hg_\crit$.

\begin{prop}  \label{RHom from Pi to part int}
If an object $\CM$ of $\Cat$ is partially integrable, then
$$\on{RHom}_{\DCat}(\Ppi,\CM)=0.$$
\end{prop}

\begin{proof}

By \propref{induction perfect}, we can assume that $\CM$ is
$I$-integrable.  In this case the assertion follows readily from
\propref{convolution with Xi}.

\end{proof}

Obviously, the induction functor $\CO_0\to \Cat$ descends to a
well-defined functor $^f\CO_0\to \Catf$. Let $$\CM\mapsto {}^f\CM$$
denote the projection functor $\Cat\to \Catf$. In
particular, we obtain the modules ${}^f\BM_w$ and ${}^f\Ppi$ in
${}^f\Cat$.

{}From \propref{RHom from Pi to part int} we obtain the following

\begin{cor} \label{map from Pi in quot}
The map
$$\on{RHom}_{\DCat}(\Ppi,\CM)\to
\on{RHom}_{\DCatf}({}^f\Ppi,{}^f\CM)$$ is an isomorphism.
\end{cor}

Since we have a surjection $\Ppi\to \BM_{w_0}$, we also obtain the
following

\begin{cor}   \label{hom from Verma to part int}
If an object $\CM$ of $\Cat$ is partially integrable, then
$\Hom(\BM_{w_0},\CM)=0$. 
\end{cor}

The main theorem in Part IV, from which we will derive
\mainthmref{equiv of quot} is the following:

\begin{mainthm} \label{RHom from Pi}
For any object $\CM$ of $\Cat$ we have
$$\on{R}^i\Hom_{\DCat}(\Ppi,\CM)=0 \text{ for } i>0.$$ 
\end{mainthm}

\section{The module $\Ppi_\reg$ via the affine Grassmannian}
\label{Pi via Gr}

We proceed with the proof of \mainthmref{RHom from Pi}.

\ssec{}

Consider the quotient $\fZ_\fg^\reg$ of $\fZ^\nilp_\fg$. We will
denote by $\Ppi_\reg$ and $\BM_{w_0,\reg}$ the modules
$\Ppi\underset{\fZ^\nilp_\fg}\otimes \fZ_\fg^\reg$ and
$\BM\underset{\fZ^\nilp_\fg}\otimes \fZ_\fg^\reg$, respectively.  The
goal of this section is to express these $\hg_\crit$-modules as
sections of critically twisted D-modules on the affine Grassmannian.

\medskip

Consider the element of the extended affine Weyl group equal to
$w_0\cdot \check\rho= -\check\rho \cdot w_0$. Let $j_{w_0\cdot
  \check\rho,*}$ and $j_{w_0\cdot \check\rho,!}$ denote the
corresponding critically twisted D-modules on $\Fl_G$.

Note that $w_0\cdot \check\rho$ is minimal in its coset in
$W\backslash W_{\on{aff}}/W$, in particular, the orbit $I\cdot
(w_0\cdot\check\rho)\subset \Fl_G$ projects one-to-one under $\Fl_G\to
\Gr_G$. Hence, $j_{w_0\cdot \check\rho,!}\underset{I}\star
\delta_{1_{Gr_G}}$ it the D-module on $\Gr_G$, obtained as the
extension by $0$ from the Iwahori orbit of the element
$t^{-\check\rho}\in \Gr_G$. Let us denote by $\on{IC}_{w_0\cdot
\check\rho,\Gr_G}$ the intersection cohomology D-module corresponding
to the above $I$-orbit.

We have the maps
\begin{equation} \label{IC of w_0 rho}
j_{w_0\cdot \check\rho,!}\underset{I}\star
\delta_{1_{Gr_G}}\twoheadrightarrow \on{IC}_{w_0\cdot
\check\rho,\Gr_G}\hookrightarrow j_{w_0\cdot
\check\rho,*}\underset{I}\star \delta_{1_{Gr_G}},
\end{equation}
such that the kernel of the first map and cokernel of the second map
are supported on the closed subset 
$\ol\Gr_G{}^{\check\rho}-\Gr_G^{\check\rho}\subset \Gr_G$.

\begin{prop}  \label{Xi clean}
The maps
$$\Xi\underset{I}\star j_{w_0\cdot \check\rho,!}\underset{I}\star
\delta_{1_{Gr_G}} \to \Xi\underset{I}\star \on{IC}_{w_0\cdot
\check\rho,\Gr_G}\to \Xi\underset{I}\star j_{w_0\cdot
\check\rho,*}\underset{I}\star \delta_{1_{Gr_G}}$$ are isomorphisms.
\end{prop}

The proposition follows from \propref{convolution with Xi},
using the following lemma:

\begin{lem}  \label{non-regular weight part int}
Any $I$-monodromic D-module on $\Gr_G$, supported on 
$\ol\Gr_G{}^{\check\rho}-\Gr_G^{\check\rho}$ is partially
integrable.
\end{lem}

\begin{proof}

Recall that the $G[[t]]$--orbits on $\Gr_G$ are labeled by the set of
dominant coweights of $G$; for a coweight $\check\lambda$ we will
denote by
$$\Gr_G^{\check\lambda}\overset{\on{emb}_{\check\lambda}}\hookrightarrow
\Gr_G$$ the embedding of the corresponding orbit.  The quotient
$G^{(1)}\backslash\Gr_G^{\check\lambda}$ is a $G$-homogeneous space,
isomorphic to a partial flag variety. We identify it with with $G/P$
by requiring that the point $w_0\cdot t^{\check\lambda}\in G\ppart$
project to $1_{G/P}\subset G/P$; we have $P=B$ if and only if
$\check\lambda$ is regular.

Note that the $G[[t]]$-orbits appearing in
$\ol\Gr_G{}^{\check\rho}-\Gr_G^{\check\rho}$ all correspond to
irregular $\check\lambda$. Therefore, it is enough to show that an
irreducible $I$-equivariant D-module on $\Gr_G^{\check\lambda}$ with
irregular $\check\lambda$ is partially integrable.

Any such D-module comes as a pull-back from a an irreducible
$B$-equivariant D-module on $G/P$ for some parabolic $P$, strictly
larger than $B$.  By definition, irreducible $B$-equivariant D-modules
on $G/P$ are IC-sheaves of closures of $B$-orbits on $G/P$. So, it is
enough to show that any such closure is stable under $SL_2^\iota$ for
some $\iota\in \CI$. But this is nearly evident:

The orbit of $1_{G/P}$ is clearly $P$-stable. Any other orbit
corresponds to some element $w\in W$ of length more than $1$.  Hence
there exists a simple reflection $s_\iota$ such that $s_\iota\cdot
w<w$. Then the orbit, corresponding to $s_\iota\cdot w$ is contained
in the closure of the one corresponding to $w$, and their union is
$SL_2^\iota$ stable.

\end{proof}

Let us denote by $\pi_{\check\lambda}$ the map from $\Gr_G^{\check\lambda}$
to the corresponding partial flag variety $G/P$. The following lemma
follows directly from definitions.

\begin{lem}    \label{swap}
Let $\wt{w}$ be an element of $W_{\on{aff}}$ which is minimal in its double
coset $W\backslash W_{\on{aff}}/W$, and $\check\lambda$ the corresponding
dominant coweight. Assume that $\check\lambda$ is regular, and let
$\CF$ be a D-module on $G/B$. We have
$$\CF \star j_{\wt{w}, !} \star \delta_{1_{\Gr_G}} \simeq
(\on{emb}_{\check\lambda})_!\circ \pi_{\check\lambda}^*(\CF) \text{
and } \CF \star j_{\wt{w}, *} \star \delta_{1_{\Gr_G}} \simeq
(\on{emb}_{\check\lambda})_*\circ \pi_{\check\lambda}^*(\CF).$$
\end{lem}

Therefore, we can rewrite
$$\Xi\underset{I}\star j_{w_0\cdot \check\rho,!}\underset{I}\star
\delta_{1_{Gr_G}} \simeq (\on{emb}_{\check\rho})_!\circ
\pi_{\check\rho}^*(\Xi), \text{ and } \Xi\underset{I}\star j_{w_0\cdot
\check\rho,*}\underset{I}\star \delta_{1_{Gr_G}} \simeq
(\on{emb}_{\check\rho})_*\circ \pi_{\check\rho}^*(\Xi)$$ Hence, the
assertion of \propref{Xi clean} can be reformulated as cleanness of
the perverse sheaf $\pi_{\check\rho}^*(\Xi)$ on $\Gr_G^{\check\rho}$,
i.e., that the map
$$(\on{emb}_{\check\rho})_!\circ 
\pi_{\check\rho}^*(\Xi) \to (\on{emb}_{\check\rho})_*\circ 
\pi_{\check\rho}^*(\Xi)$$
is an isomorphism.

\ssec{}

Set $$\bL_{w_0}:=\IC_{w_0\cdot\check\rho,\Gr_G}\underset{G}\star
\BV_\crit= \Gamma(\Gr_G, \IC_{w_0\cdot\check\rho,\Gr_G}).$$

A key result, from which we will derive the main theorem is the
following:

\begin{thm} \label{Wakimoto via Grassmannian}
There exists a canonically defined map $\BM_{w_0,\reg} \otimes
\omega_x^{\langle \rho,\check\rho\rangle} \to \bL_{w_0}$ such that:

\smallskip

\noindent{\em (a)}
The above map is surjective and its kernel is partially integrable.

\noindent{\em (b)}
The induced map
$$\Ppi_\reg\otimes \omega_x^{\langle \rho, \check\rho\rangle} \simeq
\Xi\underset{I}\star \BM_{w_0,\reg} \otimes \omega_x^{\langle \rho,
\check\rho\rangle} \to \Xi\underset{I}\star
\IC_{w_0\cdot\check\rho,\Gr_G}\underset{G}\star \BV_\crit\simeq
\Gamma(\Gr_G,\Xi\underset{I}\star \IC_{w_0\cdot\check\rho,\Gr_G})$$ is
an isomorphism.
\end{thm}

This theorem will be proved in \secref{proof of W v G}. Let us now
state a corollary of \thmref{Wakimoto via Grassmannian} that will be
used in the proof of \mainthmref{RHom from Pi}.

\begin{cor}   \label{RHom from Pi reg}
For any $i>0$,
$\on{R}^i\Hom_{\DCatr}(\Ppi_\reg,\BM_{w_0,\reg})=0$, and the
natural map $\fZ^\reg_\fg\to \Hom(\Ppi_\reg,\BM_{w_0,\reg})$ is an
isomorphism.
\end{cor}

Let us prove this corollary. By \propref{RHom from Pi to part int}, it
is sufficient to compute
$$\on{RHom}_{\DCat} \left(\Gamma(\Gr_G,\Xi\underset{I}\star
\IC_{w_0\cdot\check\rho,\Gr_G}), \Gamma(\Gr_G,
\IC_{w_0\cdot\check\rho,\Gr_G})\right).$$ By \thmref{Exts in reg}, the
latter $\on{RHom}$ is isomorphic to
$$\on{RHom}_{D\left(\fD(\Gr_G)_\crit\mod\right)}
\left(\Xi\underset{I}\star \IC_{w_0\cdot\check\rho,\Gr_G},
\underset{V\in \on{Irr}(\Rep)}\oplus \,
\IC_{w_0\cdot\check\rho,\Gr_G}\star \CF_{V^*}\underset{\BC}\otimes
V_{\fZ^\reg_\fg}\right).$$

Let $I^{-,0}$ be the subgroup of $G[[t]]$ equal to the preimage of
$N^-\subset G$ under the evaluation map. By composing with
$\psi:N^-\to \BG_a$, we obtain a character on $I^{-,0}$, denoted
in the same way, and we can consider the category
$\fD(\Gr_G)_\crit\mod^{I^{-,0},\psi}$ of $(I^{-,0},\psi)$-equivariant
D-modules, and the corresponding triangulated category.

As in \secref{character psi}, the forgetful functor
$D(\fD(\Gr_G)_\crit\mod^{I^{-,0},\psi})\hookrightarrow
D(\fD(\Gr_G)_\crit\mod)$ admits a right adjoint, which we will denote
by $\on{Av}_{I^{-,0},\psi}$, given by convolution with the corresponding
D-module on $I^{-,0}$.

{}From \propref{psi averaging} we obtain that the composition
$$D(\fD(\Gr_G)_\crit\mod)^{I,m}\to D(\fD(\Gr_G)_\crit\mod)\to
D(\fD(\Gr_G)_\crit\mod)^{I^{-,0},\psi},$$ where the last arrow is the
functor $\CF\mapsto \on{Av}_{I^{-,0},\psi}(\CF)\otimes
\det(\fn^-[1])^{-1}$ is exact, and essentially commutes with the
Verdier duality on the holonomic subcategory.

\medskip

By the construction of $\Xi$, for $\CF_1\in \fD(\Gr_G)_\crit\mod^I$ and
$\CF_2\in \fD(\Gr_G)_\crit\mod^{I,m}$
\begin{align*}
&\on{RHom}_{D(\fD(\Gr_G)_\crit\mod)}(\Xi\underset{I}\star
\CF_1,\CF_2)\simeq \\
&\on{RHom}_{D(\fD(\Gr_G)_\crit\mod)^{I^{-,0},\psi}}
\left(\on{Av}_{I^{-,0},\psi}(\CF_1),
\on{Av}_{I^{-,0},\psi}(\CF_2)\right).
\end{align*}

Using the exactness property of $\on{Av}_{I^{-,0},\psi}$ mentioned
above, \corref{RHom from Pi reg} follows from the next general
result:

\begin{thm} \label{Cass-Shal}
For any two $\CF'_1,\CF'_2\in \fD(\Gr_G)_\crit\mod^{I^{-,0},\psi}$ and
$i>0$
$$\on{R}^i\Hom_{D(\fD(\Gr_G)_\crit\mod)}(\CF'_1,\CF'_2)=0.$$ The
functor $\fD(\Gr_G)_\crit\mod^{G[[t]]}\to
\fD(\Gr_G)_\crit\mod^{I^{-,0},\psi}$, given by
$$\CF\mapsto \on{Av}_{I^{-,0},\psi} (\IC_{w_0\cdot \check\rho}\star
\CF),$$ is an equivalence of abelian categories.
\end{thm}

The proof of \thmref{Cass-Shal} is a word-for-word repetition of the
proof of the main theorem of \cite{FGV}, using the fact that the
combinatorics of $I^0$ (resp., $I^{-,0}$) orbits on $\Gr_G$ is the
same as that of $N\ppart$ (resp., $N^-\ppart$) orbits. The main point is
that any irreducible object of $\fD(\Gr_G)_\crit\mod^{I^{-,0},\psi}$
is a clean extension from a character sheaf on an orbit.

\section{Proofs of the main theorems}      \label{proof of main thms}

In this section we will prove \mainthmref{RHom from Pi} and 
derive from it \mainthmref{equiv of quot}, assuming 
\thmref{Wakimoto via Grassmannian}  
(which is proved in the next section).

\ssec{}

In \corref{RHom from Pi reg} we computed the extensions between
$\Ppi_\reg$ and $\BM_{w_0,\reg}$ in the category $\DCatr$. Now we
use this result to compute the extensions between $\Ppi$ and
$\BM_{w_0}$ in the category $\DCat$.

\begin{prop}  \label{Ext from Pi to Verma}
The morphism $\fZ^\nilp_\fg\to \on{Hom}_{\DCat}(\Ppi,\BM_{w_0})$
is an isomorphism and $R^i\on{Hom}_{\DCat}(\Ppi,\BM_{w_0})=0$
for $i>0$.
\end{prop}

\begin{proof}

Let us note that for any two objects $\CM^\bullet_1,\CM^\bullet_2\in
D^b(\hg_\crit\mod_\nilp)$ the complex
$\on{RHom}_{D(\hg_\crit\mod_\nilp)}(\CM_1^\bullet,\CM_2^\bullet)$ is
naturally an object of $D^+(\fZ^\nilp_\fg\mod)$. Recall also that
$\fZ^\reg_\fg$, as a module over $\fZ^\nilp_\fg$, admits a finite
resolution by finitely generated projective modules.  Therefore, the
functor
$$\CM^\bullet\mapsto
\CM^\bullet\overset{L}{\underset{\fZ^\nilp_\fg}\otimes} \fZ^\reg_\fg$$
is well-defined as a functor $D^b(\hg_\crit\mod_\nilp)\to
D^b(\hg_\crit\mod_\nilp)$.

Almost by definition we obtain the following result.

\begin{lem}
$$\left(\on{RHom}_{D(\hg_\crit\mod_\nilp)}(\CM_1^\bullet,
\CM_2^\bullet)\right) \overset{L}{\underset{\fZ^\nilp_\fg}\otimes}
\fZ^\reg_\fg$$ is isomorphic to
$$\on{RHom}_{D(\hg_\crit\mod_\nilp)}\left(\CM_1^\bullet,
(\CM_2^\bullet \overset{L}{\underset{\fZ^\nilp_\fg}\otimes}
\fZ^\reg_\fg)\right).$$
\end{lem}

\medskip

Since $\Ppi$ is flat over $\fZ^\nilp_\fg$, by \lemref{* res},
we obtain that for any $\CM\in \hg_\crit\mod_\reg$,
\begin{lem}  \label{Ext reg ou pas}
$$\on{RHom}_{D(\hg_\crit\mod_\nilp)}(\Ppi,\CM)\simeq
\on{RHom}_{D(\hg_\crit\mod_\reg)}(\Ppi_\reg,\CM).$$
\end{lem}
By combining this with \corref{RHom from Pi reg}, we obtain that the 
natural map
$$\fZ_\fg^\reg\to
\left(\on{RHom}_{D(\hg_\crit\mod_\nilp)}(\Ppi,\BM_{w_0})\right)
\overset{L}{\underset{\fZ^\nilp_\fg}\otimes} \fZ_\fg^\reg$$ is a
quasi-isomorphism. We will now derive the assertion of \propref{Ext from Pi to Verma}
by a Nakayama lemma type argument.

\medskip

Consider the $\BG_m$-action on $\hg_\crit$, coming from
$\BG_m\hookrightarrow \on{Aut}(\D)$. We obtain that $\BG_m$ acts
weakly on the categories $\hg_\crit\mod_\nilp$ and
$\hg_\crit\mod$. Since the objects $\Ppi$ and $\BM_{w_0}$ are
$\BG_m$-equivariant, the Ext groups
$$\on{Ext}^i_{\hg_\crit\mod}(\Ppi,\BM_{w_0}) \text{ and }
\on{Ext}^i_{\hg_\crit\mod_\nilp}(\Ppi,\BM_{w_0})$$
acquire an action of $\BG_m$ by \lemref {inner hom weak}.

We claim that the grading arising on
$\on{Ext}^i_{\hg_\crit\mod_\nilp}(\Ppi,\BM_{w_0})$ is bounded from
above. First, let us note that the grading on
$\on{Ext}^i_{\hg_\crit\mod}(\Ppi,\BM_{w_0})$ is non-positive.  This is
evident since $\on{Ext}^i_{\hg_\crit\mod}(\Ppi,\BM_{w_0})$ are
computed by the standard complex
$\fC^\bullet\left(\fg[[t]],\Hom_{\BC}(\Pi,\BM_{w_0})\right)$, whose
terms are non-positively graded. Note also that the algebra
$\fZ^\nilp_\fg$ is non-positively graded, and the grading on
$N_{\fZ^\nilp_\fg/\fZ_\fg}$ is such that only finitely many free
generators have positive degrees. Now, the spectral sequence of
\secref{equivariant restriction} implies by induction on $i$ that the
grading on $\on{Ext}^i_{\hg_\crit\mod_\nilp}(\Ppi,\BM_{w_0})$ is
bounded from above.

\medskip

Since the algebra $\Fun(\nOp)$ is itself non-positively graded,
$\on{RHom}_{D(\hg_\crit\mod_\nilp)}(\Ppi,\BM_{w_0})$ can be
represented by a complex of graded modules such that the grading on
each term is bounded from above, and which lives in non-negative
cohomological degrees.  Recall again that the ideal of
$\Fun(\Op^\reg)$ in $\Fun(\nOp)$ is generated by a regular sequence of
homogeneous negatively graded elements.  The proof is concluded by the
following observation:

\begin{lem}
Let 
$$
Q^\bullet:=Q^0\to Q^1\to...\to Q^n\to...
$$
be a complex of graded modules over a graded algebra
$A=\BC[x_1,...,x_n]$, where $\deg(x_i)<0$ such that the grading on
each $Q^i$ is bounded from above.  Assume that
$Q^\bullet\overset{L}{\underset{A}\otimes}\BC$ is acyclic away from
cohomological degree $0$. Then $Q^\bullet$ is itself acyclic away from
cohomological degree $0$.
\end{lem}

\end{proof}

\begin{cor}  \label{RHom to any L}
For any $\fZ^\nilp_\fg$-module $\CL$ 
$$\on{R}^i\Hom_{\DCat}(\Ppi,\BM_{w_0}
\underset{\fZ^\nilp_\fg}\otimes \CL)=0$$
for $i>0$ and is isomorphic to $\CL$ for $i=0$.
\end{cor}

\begin{proof}

Since any module $\CL$ is a direct limit of finitely presented ones,
by \propref{induction perfect}, we may assume that $\CL$ is finitely
presented. Since $\fZ^\nilp_\fg$ is isomorphic to a polynomial
algebra, any finitely presented module admits a finite resolution by
projective ones:
$$\CL_n\to...\to \CL_1\to \CL_0\to \CL.$$ Since $\BM_{w_0}$ is flat
over $\fZ^\nilp_\fg$ (cf, \corref{Verma is flat}), we obtain a
resolution
$$\BM_{w_0}\underset{\fZ^\nilp_\fg}\otimes \CL_n\to...\to
\BM_{w_0}\underset{\fZ^\nilp_\fg}\otimes \CL_1\to 
\BM_{w_0}\underset{\fZ^\nilp_\fg}\otimes \CL_0\to
\BM_{w_0}\underset{\fZ^\nilp_\fg}\otimes \CL.$$

Hence, we obtain a spectral sequence, converging to
$$\on{R}^i\Hom_{D(\hg_\crit\mod)^I}(\Ppi,\BM_{w_0}
\underset{\fZ^\nilp_\fg}\otimes \CL),$$ whose first term $E^{i,j}_1$
is given by
$$\on{R}^i\Hom_{D(\hg_\crit\mod)^I}(\Ppi,\BM_{w_0}
\underset{\fZ^\nilp_\fg}\otimes \CL_{-j}).$$

Since $\CL_\bullet$ are projective, by \propref{Ext from Pi to Verma},
we obtain that $E^{i,j}_1=0$ unless $i=0$, and in the latter case, it
is isomorphic to $\CL_{-j}$, implying the assertion of the corollary.

\end{proof}

\begin{cor}   \label{kernel part int}
For any object $L$ of $\Cat$ and the $\fZ^\nilp_\fg$-module
$\CL:=\Hom(\BM_{w_0},L)$ the kernel of the natural map
$$\BM_{w_0}\underset{\fZ^\nilp_\fg}\otimes \CL\to L$$ is
partially integrable.
\end{cor}

\begin{proof}

Let $\CM$ be the kernel of $\BM_{w_0}\underset{\fZ^\nilp_\fg}\otimes
\CL\to L$, and suppose that it is not partially integrable. Let
$\CM'\subset \CM$ be the maximal partially integrable
submodule. Consider the short exact sequence
$$0\to \CM'\to \CM\to \CM''\to 0.$$
By \lemref{hom from L}, we have a non-zero map $\BL_w\to \CM''$ for some
$w\in W$. We claim that $w$ necessarily equals $w_0$.

Indeed, all modules $\BL_{w}$ with $w\neq w_0$ are partially
integrable, and we would obtain that the preimage in $\CM$ of
$\on{Im}(\BL_w)$ is again integrable, and is strictly bigger than
$\CM$.

\medskip

Hence, we have a map $\BM_{w_0}\to \CM''$, and by composing, we obtain a map
$\Ppi\to \CM''$. Now, by \propref{RHom from Pi to part int}, this maps lifts
to a map $\Ppi\to \CM$, i.e., $\Hom(\Ppi,\CM)\neq 0$.

Consider now the exact sequence
$$0\to \Hom(\Ppi,\CM)\to
\Hom(\Ppi,\BM_{w_0}\underset{\fZ^\nilp_\fg}\otimes \CL)\to
\Hom(\Ppi,L).$$

By \propref{RHom to any L}, the middle term is isomorphic to $\CL$,
and it maps injectively to $\Hom(\Ppi,L)$, since
$$\CL\simeq \Hom(\BM_{w_0},L)\hookrightarrow \Hom(\Ppi,L),$$
which is a contradiction.

\end{proof}

Now we are able to prove \mainthmref{RHom from Pi}:

\begin{proof}

Let $\CM$ be an object of $\Cat$. It admits a filtration
$0=\CM_0\subset \CM_1\subset \CM_2....$, whose subquotients
$\CM_j/\CM_{j-1}$ have the property that each is a quotient of the
module $\BL_w$ for some $w\in W$. By \propref{induction perfect}, to
prove that $\on{Ext}^i_{\hg_\crit\mod_\nilp}(\Ppi,\CM)=0$ for $i=0$,
by devissage, we can assume that $\CM$ itself is a quotient of some
$\BL_w$.

If $w\neq w_0$, then $\CM$ is partially integrable and the vanishing
of Exts follows from \propref{RHom from Pi to part int}. Hence, we can
assume that $\CM$ is a quotient of $\BL_{w_0}=\BM_{w_0}$. In this
case, the assertion of the theorem follows from \corref{RHom to any L}
combined with \corref{kernel part int}.

\end{proof}

\ssec{Proof of \mainthmref{equiv of quot}}

Now we derive \mainthmref{equiv of quot} from \mainthmref{RHom from
Pi}.  We define the functor $\hg_\crit\mod\to \fZ^\nilp_\fg\otimes
h_0\mod$ by
$$\CM\mapsto \Hom(\Ppi,\CM).$$

Composing with the forgetful functor $\Cat\to \hg_\crit\mod$ we obtain a
functor
$$\Cat\to \fZ^\nilp_\fg\otimes h_0\mod.$$
By \mainthmref{RHom from Pi}, the latter functor is exact, and by 
\propref{RHom from Pi to part int} it factors through
$\Catf$. This defines the desired functor
$$^f\sF:\DCatf\to D^b(\fZ^\nilp_\fg\otimes h_0\mod).$$

\medskip

We define a functor $^f\sG:\fZ^\nilp_\fg\otimes h_0\mod\to
\Catf$ by
$$\CL\mapsto {}^f\Ppi\underset{\fZ^\nilp_\fg\otimes h_0}\otimes \CL.$$
{}From \lemref{Pi in quotient} and \corref{Verma is flat} it follows that
this functor is exact. We will denote by the same character the
resulting functor
$$D^b(\fZ^\nilp_\fg\otimes h_0)\to \DCatf.$$

By \corref{map from Pi in quot}, for $\CL^\bullet\in
D^-(\fZ^\nilp_\fg\otimes h_0\mod)$ and $\CM^\bullet\in
\DCatf$ we have a natural isomorphism
$$\Hom_{D^b(\fZ^\nilp_\fg\otimes h_0\mod)}(\CL,{}^f\sF(\CM^\bullet))\simeq
\Hom_{\DCatf}({}^f\sG(\CL),\CM^\bullet).$$

Hence, $\sG$ and $\sF$ are mutually adjoint. Let us show that they are in
fact mutually quasi-inverse.

\medskip

Let us first show that the adjunction morphism $\on{Id}\to
{}^f\sF\circ {}^f\sG$ is an isomorphism. By exactness, it suffices
to show that for a $\fZ^\nilp_\fg\otimes h_0$-module $\CL$, on which
the action of $h_0$ is trivial, the map
\begin{equation} \label{F circ G}
\CL\mapsto \Hom\left({}^f\Ppi,{}^f\Ppi\underset{\fZ^\nilp_\fg\otimes
h_0}\otimes \CL\right)
\end{equation}
is an isomorphism.

We have $\Pi\underset{h_0}\otimes \BC\simeq M_1^\vee$, and hence
$\Ppi\underset{h_0}\otimes \BC\simeq \BM_{1}^\vee$. Since the kernel
of $\BM_1^\vee\to \BM_{w_0}$ is partially integrable, we obtain that
$${}^f\Ppi\underset{\fZ^\nilp_\fg\otimes h_0}\otimes \CL\simeq
{}^f\BM_{w_0} \underset{\fZ^\nilp_\fg}\otimes \CL,$$ and the
assertion follows from \corref{RHom to any L}.

\medskip

To show that the adjunction $^f\sG\circ {}^f\sF\to \on{Id}$ is an
isomorphism, by exactness, it is again sufficient to evaluate it on a
single module $\CM$.  Since the functor $^f\sF$ is faithful, it is
enough to show that
$${}^f\sF\circ {}^f\sG\circ {}^f\sF(\CM)\to {}^f\sF(\CM)$$ is an
isomorphism. But we already know that $^f\sF(\CM)\to {}^f\sF\circ
{}^f\sG\circ  {}^f\sF(\CM)$ is an isomorphism, and our assertion follows.

This completes the proof of \mainthmref{equiv of quot} modulo
\thmref{Wakimoto via Grassmannian}. \qed

\section{Proof of \thmref{Wakimoto via Grassmannian}}   \label{proof
  of W v G}

\ssec{}   \label{one map}

Let us first  construct the map
\begin{equation} \label{map from Wak to irr}
\BM_{w_0,\reg} \otimes \omega_x^{\langle \rho,\check\rho\rangle}\to
\Gamma(\Gr_G, \IC_{w_0\cdot\check\rho,\Gr_G}),
\end{equation}
whose existence is stated in \thmref{Wakimoto via Grassmannian}.

\medskip

Consider the $\hg_\crit$-module $\Gamma(\Gr_G, j_{w_0\cdot \check\rho,*}
\underset{I}\star\delta_{1,\Gr_G})$; it is equivariant with respect 
to the action of $\BG_m$ acting by loop rotations. This module contains
a unique line, corresponding to those sections of the twisted D-module 
$j_{w_0\cdot \check\rho,*}\underset{I}\star\delta_{1,\Gr_G}$, which are 
scheme-theoretically supported on the closure of the $I$-orbit of the 
element $t^{-\check\rho}\in \Gr_G$, and which are $I^0$-invariant. 

This line has weight $-2\check\rho$ with respect to $\fh$, and has the
highest degree with respect to the $\BG_m$-action. Moreover, a
straightforward calculation (see \cite{BD}, Sect. 9.1.13) shows 
that this line can be canonically identified with
$\omega_x^{\langle \rho,\check\rho\rangle}$.  This defines a map
$$\BM_{w_0,\reg} \otimes \omega_x^{\langle \rho,\check\rho\rangle}
\to \Gamma(\Gr_G, j_{w_0\cdot \check\rho,*}
\underset{I}\star\delta_{1,\Gr_G}).$$

We claim that the above map factors through $\Gamma(\Gr_G,
\IC_{w_0\cdot\check\rho,\Gr_G}) \subset \Gamma(\Gr_G, j_{w_0\cdot
\check\rho,*} \underset{I}\star\delta_{1,\Gr_G})$.  Indeed, by
\lemref{non-regular weight part int}, the quotient module is partially
integrable, and from \corref{hom from Verma to part int} we obtain
that it cannot receive a non-zero map from $\BM_{w_0,\reg}$.

\begin{prop} \label{surj of Wak on irr}
The map $\BM_{w_0,\reg} \otimes \omega_x^{\langle
\rho,\check\rho\rangle}\to \Gamma(\Gr_G,
\IC_{w_0\cdot\check\rho,\Gr_G})$ constructed above is surjective.
\end{prop}

The proof will be given at the end of this section. We will now
proceed with the proof of \thmref{Wakimoto via Grassmannian}.

\ssec{}   \label{another map}

We shall now construct a map
\begin{equation} \label{map from irr to Wak}
\Gamma(\Gr_G, \IC_{w_0\cdot\check\rho,\Gr_G})\to
\BM_{1,\reg} \otimes \omega_x^{\langle \rho,\check\rho\rangle}.
\end{equation}

First, by \secref{averaging}, for any $I$-equivariant $\hg_\crit$-module $\CM$,
$$\Hom(\BV_\crit,\CM)\simeq 
R^0\Hom_{\hg\mod^{G[[t]]}}(\BV_\crit,\on{Av}_{G[[t]]/I}(\CM)).$$

Applying this to $\CM=\BM_{w_0,\reg}$, we calculate:
\begin{equation} \label{averaging of Wakimoto}
\on{Av}_{G[[t]]/I}(\BM_{w_0,\reg})\simeq
\on{Av}_{G[[t]]/I}(\BM_{w_0})\underset{\fZ^\nilp_\fg}{\overset{L}\otimes} 
\fZ^\reg_\fg\simeq 
\BV_\crit[-\dim(G/B)]\underset{\fZ^\nilp_\fg}{\overset{L}\otimes} 
\fZ^\reg_\fg.
\end{equation}
Hence, the $0$-th cohomology of $\on{Av}_{G[[t]]/I}(\BM_{w_0,\reg})$
is isomorphic to
$$\on{Tor}_{\dim(G/B)}^{\fZ^\nilp_\fg}(\BV_\crit,\fZ^\reg_\fg)\simeq
\BV_\crit\otimes
\on{Tor}_{\dim(G/B)}^{\fZ^\nilp_\fg}(\fZ^\reg_\fg,\fZ^\reg_\fg).$$
However, from \corref{nilpotent directions} and \propref{grading on algebroid}(2), 
it follows that
$$\on{Tor}_{\dim(G/B)}^{\fZ^\nilp_\fg}(\fZ^\reg_\fg,\fZ^\reg_\fg)\simeq
\Lambda^{\dim(G/B)}(N^*_{\fZ^\reg_\fg/\fZ^\nilp_\fg})\simeq
\fZ^\reg_\fg\otimes \omega_x^{-\langle 2\rho,\check\rho\rangle}.$$
Hence, the above $0$-th cohomology is isomorphic to $\BV_\crit\otimes
\omega_x^{-\langle 2\rho,\check\rho\rangle}$, and we obtain a map
$$\BV_\crit\to  \BM_{w_0,\reg} \otimes 
\omega_x^{\langle 2\rho,\check\rho\rangle}.$$

By applying the convolution $j_{w_0\cdot
\check\rho,*}\underset{I}\star\cdot$ to both sides we obtain a map
\begin{equation} \label{crucial map}
\Gamma(\Gr_G, j_{w_0\cdot\check\rho,*}\underset{I}\star \delta_{1,\Gr_G})\to
j_{w_0\cdot\check\rho,*}\underset{I}\star \BM_{w_0,\reg} 
\otimes \omega_x^{\langle 2\rho,\check\rho\rangle}.
\end{equation}

However, by \eqref{lambda equivariance of Verma},
\begin{align*}
&j_{w_0\cdot\check\rho,*}\underset{I}\star \BM_{w_0,\reg} \simeq
j_{w_0,!}\underset{I}\star j_{w_0,*}\underset{I}\star
j_{w_0\cdot\check\rho,*}\underset{I}\star \BM_{w_0,\reg} \simeq \\
&j_{w_0,!}\underset{I}\star j_{\check\rho,*}\underset{I}\star
\BM_{w_0,\reg} \simeq j_{w_0,!}\underset{I}\star \BM_{w_0,\reg}\otimes
\omega_x^{-\langle \rho,\check\rho\rangle} \simeq \BM_{1,\reg}\otimes
\omega_x^{-\langle \rho,\check\rho\rangle},
\end{align*}
and by composing with the embedding $\Gamma(\Gr_G,
\IC_{w_0\cdot\check\rho,\Gr_G})\hookrightarrow \Gamma(\Gr_G,
j_{w_0\cdot\check\rho,*}\underset{I}\star \delta_{1,\Gr_G})$ we obtain
the map of \eqref{map from irr to Wak}. By constriction, this map
respects the $\BG_m$-action.

\ssec{}

Consider now the composition
\begin{equation} \label{composed map}
\BM_{w_0,\reg} \otimes \omega_x^{\langle \rho,\check\rho\rangle}\to
\Gamma(\Gr_G, \IC_{w_0\cdot\check\rho,\Gr_G}) \to
\BM_{1,\reg} \otimes \omega_x^{\langle \rho,\check\rho\rangle}.
\end{equation}

\begin{lem}  \label{map between Vermas}
The resulting map $\BM_{w_0,\reg}\to \BM_{1,\reg}$ is a non-zero
multiple of the canonical map, coming from the 
embedding $M_{w_0}\to M_1$.
\end{lem}

\begin{proof}

First, the map in question is non-zero by \propref{surj of Wak on
irr}. Secondly, our map $\BM_{w_0,\reg}\to \BM_{1,\reg}$ respects the
$\BG_m$-action by loop rotations. Since $\BM_{w_0,\reg}$ is generated
by a vector of degree $0$, and the subspace in $\BM_{1,\reg}$,
consisting of elements of degree $0$ is isomorphic to the Verma module
$M_{0}$, any map $\BM_{w_0,\reg}\to \BM_{1,\reg}$, compatible with the
grading, is a scalar multiple of the canonical map.

\end{proof}

\ssec{}

Let us now derive \thmref{Wakimoto via Grassmannian} from
\lemref{map between Vermas}.

Let us apply the convolution $\Xi\underset{I}\star\cdot$ to the three
terms appearing in \eqref{composed map}. We obtain the maps
\begin{equation}  \label{convolved composed map}
\Ppi \otimes \omega_x^{\langle \rho,\check\rho\rangle}\simeq
\Xi\underset{I}\star \BM_{w_0,\reg} \otimes \omega_x^{\langle
\rho,\check\rho\rangle}\to \Gamma(\Gr_G,\Xi\underset{I}\star
\IC_{w_0\cdot\check\rho,\Gr_G})\to \Xi\underset{I}\star \BM_{1,\reg}
\otimes \omega_x^{\langle \rho,\check\rho\rangle}.
\end{equation}

However, the canonical map $M_{w_0}\to M_0$ has the property that
its cokernel is partially integrable. Hence, the cone of the resulting map
$\BM_{w_0,\reg}\to \BM_{1,\reg}$ is also partially integrable. 

Hence, by \thmref{map between Vermas} and \propref{convolution with Xi},
the composed map in \eqref{convolved composed map} is an isomorphism.
In particular, we obtain that $\Ppi$ is a direct summand of 
$\Gamma(\Gr_G,\Xi\underset{I}\star \IC_{w_0\cdot\check\rho,\Gr_G})$.

\begin{lem}
The map
$$\fZ^\reg_\fg\otimes h_0\to 
\on{End}(\Gamma(\Gr_G,\Xi\underset{I}\star
\IC_{w_0\cdot\check\rho,\Gr_G}))$$
is an isomorphism. 
\end{lem}

\begin{proof}

By \thmref{Exts in reg} the assertion of the lemma is equivalent to
the fact that $h_0\simeq \on{End}(\Xi\underset{I}\star
\IC_{w_0\cdot\check\rho,\Gr_G})$, and $\on{Hom}(\Xi\underset{I}\star
\IC_{w_0\cdot\check\rho,\Gr_G}, \Xi\underset{I}\star
\IC_{w_0\cdot\check\rho,\Gr_G}\underset{G[[t]]}\star
\CF_{V^{\check\lambda}})=0$ for $\check\lambda\neq 0$.

The former isomorphism follows from the fact that $h_0\simeq
\on{End}(\Xi)$, combined with \propref{Xi clean} and the fact that the
projection $\Gr_G^{\check\rho}\to G/B$ is smooth with connected
fibers.

\medskip

To prove the vanishing for $\check\lambda\neq 0$, it is enough to show that
$$\on{Hom}(\Xi\underset{I}\star \IC_{w_0\cdot\check\rho,\Gr_G},
\IC_{w_0\cdot\check\rho,\Gr_G}\underset{G[[t]]}\star
\CF_{V^{\check\lambda}})=0,$$ because modulo partially integrable
objects, $\Xi$ appearing in the right-hand side is an extension of
several copies of $\delta_{1_{G/B}}$.

As in the proof of \corref{RHom from Pi reg}, the latter Hom is
isomorphic to
$$R^0\on{Hom}_{D(\fD(\Gr_G)_\crit\mod)^{I^{-,0},\psi}}
\left(\on{Av}_{I^{-,0},\psi}(\delta_{1_{\Gr_G}}),
\on{Av}_{I^{-,0},\psi}(\delta_{1_{\Gr_G}})\underset{G[[t]]}\star
\CF_{V^{\check\lambda}}\right),$$ and the latter vanishes, according
to \thmref{Cass-Shal}.

\end{proof}

Thus, we obtain that the ring
$\on{End}(\Gamma(\Gr_G,\Xi\underset{I}\star
\IC_{w_0\cdot\check\rho,\Gr_G}))$ has no idempotents. In particular,
the map $\Ppi \otimes \omega_x^{\langle \rho,\check\rho\rangle}\to
\Gamma(\Gr_G,\Xi\underset{I}\star \IC_{w_0\cdot\check\rho,\Gr_G})$ is
an isomorphism, establishing point (b) of \thmref{Wakimoto via
Grassmannian}.

\medskip

\propref{surj of Wak on irr} states that the map \eqref{map from Wak
to irr} is surjective. Thus, it remains to show that the kernel of the
map \eqref{map from Wak to irr} is partially integrable. But this
follows from point (b) and \propref{convolution with Xi}. Therefore we
obtain point (a) of \thmref{Wakimoto via Grassmannian}. This completes
the proof of \thmref{Wakimoto via Grassmannian} modulo \propref{surj
  of Wak on irr}, which is proved in the next section. \qed

\ssec{Proof of \propref{surj of Wak on irr}}

The crucial fact used in the proof of this proposition is that the
module $\bL_{w_0}$ carries an action of the renormalized algebra
$U^{\ren,\reg}(\hg_\crit)$, see \secref{action of renorm}.  Moreover,
as an object of the category $U^{\ren,\reg}(\hg_\crit)\mod$, the
module $\bL_{w_0}$ is irreducible, because the D-module
$\on{IC}_{w_0\cdot \check\rho,\Gr_G}$ is irreducible, and the global
sections functor $\fD(\Gr_G)\mod \to U^{\ren,\reg}(\hg_\crit)\mod$ is
fully faithful, according to \cite{FG}.

Recall that the algebra $U^{\ren,\reg}(\hg_\crit)\mod$ is naturally
filtered,
$$U^{\ren,\reg}(\hg_\crit)\mod=\underset{i}\cup\,
(U^{\ren,\reg}(\hg_\crit)\mod)^i, \text{ so that}$$
$(U^{\ren,\reg}(\hg_\crit)\mod)^0\simeq
\wt{U}_\crit(\hg)\underset{\fZ_\fg}\otimes \fZ^\reg_\fg$, and
$(U^{\ren,\reg}(\hg_\crit)\mod)^1/(U^{\ren,\reg}(\hg_\crit)\mod)^0$
being a free $(U^{\ren,\reg}(\hg_\crit)\mod)_0$-module, generated by
the algebroid $N^*_{\fZ^\reg_\fg/\fZ_\fg}$.

Let us denote by $(\bL_{w_0})^0\subset \bL_{w_0}$ the image of the map
$\BM_{w_0,\reg}\otimes \omega_x^{\langle \rho,\check\rho\rangle}\to
\bL_{w_0}$, and we define the submodule $(\bL_{w_0})^i$ inductively as
the image of $(\bL_{w_0})^{i-1}$ under the action of
$(U^{\ren,\reg}(\hg_\crit)\mod)^1$.  In particular, we have surjective
maps
$$N^*_{\fZ^\reg_\fg/\fZ_\fg}\otimes (\bL_{w_0})^i/(\bL_{w_0})^{i-1}\to
(\bL_{w_0})^{i+1}/(\bL_{w_0})^i,$$
and, hence, also surjective maps
$$\left(N^*_{\fZ^\reg_\fg/\fZ_\fg}\right)^{\otimes i}\otimes
\Bigl(\BM_{w_0,\reg}\otimes \omega_x^{\langle
\rho,\check\rho\rangle}\Bigr) \twoheadrightarrow
(\bL_{w_0})^i/(\bL_{w_0})^{i-1},$$ and $\underset{i}\cup\,
(\bL_{w_0})^i=\bL_{w_0}$. Our task is to show that
$(\bL_{w_0})^0=(\bL_{w_0})^1$, i.e., that $(\bL_{w_0})^0$ is stable
under the action of $(U^{\ren,\reg}(\hg_\crit)\mod)^1$.

\medskip

\begin{lem}   \label{no part int}
The module $\bL_{w_0}$ has no partially integrable subquotients.
\end{lem}

\begin{proof}

First, let us show first that $\bL_{w_0}$ has no partially integrable
quotient modules.  Suppose that $\CM$ is such a quotient module. Let
$i$ be the minimal integer such that the projection $(\bL_{w_0})^i\to
\CM$ is non-zero; by definition this projection factors through
$(\bL_{w_0})^i/(\bL_{w_0})^{i-1}$. Hence, some element of
$\left(N^*_{\fZ^\reg_\fg/\fZ_\fg}\right)^{\otimes i}$ gives rise to a
non-trivial map $\BM_{w_0,\reg}\to \CM$. But this is a contradiction,
since $\BM_{w_0}$, and hence $\BM_{w_0,\reg}$, cannot map to any
partially integrable module.

\medskip

Consider now $\bL_{w_0}$ as a graded module, i.e., as a module over
$\BC\cdot t\partial_t\ltimes \hg_\crit$. It is easy to see that a
graded module admits no partially integrable subquotients as a
$\hg_\crit$-module if and only if it has the same property with
respect to $\BC\cdot t\partial_t\ltimes \hg_\crit$.

As was remarked earlier, the commutative Lie algebra $\BC\cdot
t\partial_t\oplus \fh$ has finite-dimensional eigenspaces on the
module $\Gamma(\Gr_G,j_{w_0\cdot \check\rho,*}\underset{I}\star
\delta_{1_{Gr_G}})$; hence the same will be true for $\bL_{w_0}$.

\medskip

Consider the maximal $\BC\cdot t\partial_t\ltimes \hg_\crit$-stable
submodule of $\bL_{w_0}$, that does not contain the highest weight
line, and take the quotient. Since this quotient is generated by a
vector of weight $-2\rho$ with respect to $\fh$, it is not partially
integrable.

Let $\CM''$ be maximal $\BC\cdot t\partial_t\ltimes \hg_\crit$-stable
quotient of $\bL_{w_0}$, which admits no partially integrable
subquotients. It is well-defined due to the above
finite-dimensionality property.  It is non-zero, because we have just
exhibited one such quotient.

Let $\CM':=\on{ker}(\bL_{w_0}\to \CM'')$, and assume that $\CM'\neq
0$.  As above, some section of $N^*_{\fZ^\reg_\fg/\fZ_\fg}$ induces
a non-zero map of $\hg_\crit$-modules $\CM'\to \CM''$. Therefore,
$\CM'$ also admits a quotient, which has no partially integrable
subquotients.  This contradicts the definition of $\CM''$. Hence,
$\bL_{w_0}$ has no partially integrable subquotients.

\end{proof} 
 
Let us continue viewing $\bL_{w_0}$ as a graded module. For an integer
$n$ we will denote by $(\cdot)_n$ the subspace of elements of degree
$n$. By \secref{one map}, $(\bL_{w_0})_{n+\langle
\rho,\check\rho\rangle}=0$ if $n>0$, and $(\bL_{w_0})_{\langle
\rho,\check\rho\rangle}$ identifies with the Verma module $M_{w_0}$;
in particular, it is contained in $(\bL_{w_0})^0$.

\begin{lem}  \label{deg -1}
The subspace $(\bL_{w_0})_{\langle \rho,\check\rho\rangle-1}$ is also
contained in $(\bL_{w_0})^0$.
\end{lem}

\begin{proof}

Let $\sV\subset \Gamma(\Gr_G,j_{w_0\cdot
\check\rho,*}\underset{I}\star \delta_{1_{Gr_G}})$ be the subspace of
sections, scheme-theoretically supported on the $I$-orbit $I\cdot
t^{-\check\rho}\subset \Gr_G$.

Let $\on{Lie}(I^0)^-\subset \fg\ppart$ be the subalgebra, opposite to
$\on{Lie}(I^0)$, i.e., the one spanned by $t^{-1}\fg([t^{-1}])$ and
$\fn^{-1}\subset \fg$. The module $\Gamma(\Gr_G,j_{w_0\cdot
\check\rho,*}\underset{I}\star \delta_{1_{Gr_G}})$ is generated from
$\sV$ by means of $\on{Lie}(I^0)^-\subset \fg\ppart$.

Hence, the subspace $\Gamma(\Gr_G,j_{w_0\cdot \check\rho,*}
\underset{I}\star \delta_{1_{Gr_G}})_{\langle
\rho,\check\rho\rangle-1}$ is the direct sum of $(\fg\otimes
t^{-1})\cdot U(\fn^-)\cdot \sV_{\langle \rho,\check\rho\rangle}$ and
$U(\fn^-)\cdot \sV_{\langle \rho,\check\rho\rangle-1}$.  Note that
$\sV_{\langle \rho,\check\rho\rangle}$ is the highest weight line in
$\bL_{w_0}$.  Hence, $$(\fg\otimes t^{-1})\cdot U(\fn^-)\cdot
\sV_{\langle \rho,\check\rho\rangle} \subset (\bL_{w_0})^0.$$

Therefore,, it remains to show that $\sV_{\langle
\rho,\check\rho\rangle-1}\cap \bL_{w_0}$ is contained in
$(\bL_{w_0})^0$. Suppose not, and consider the image of $\sV_{\langle
\rho,\check\rho\rangle-1}\cap \bL_{w_0}$ in
$(\bL_{w_0})^1/(\bL_{w_0})^0$. This is a subspace annihilated by
$\fg(t\BC[[t]]])$, and stable under the $\fb$-action. Take some highest
weight vector. It gives rise to a map $\BM_w\to
(\bL_{w_0})^1/(\bL_{w_0})^0$ for some element $w\in W$; moreover
$w=w_0$ if and only if the above highest weight is $-2\rho$.

However, the algebra of functions on $I\cdot t^{-\check\rho}$ is
generated by elements, whose weights with respect to $\fh$ are in
$\on{Span}^+(\alpha_i)-0$. Therefore, the above highest weight is
different from $-2\rho$. Thus, we obtain a non-zero map $\BM_w\to
(\bL_{w_0})^1/(\bL_{w_0})^0$ for $w\neq w_0$, where $\BM_w$ is endowed
with a $\BG_m$-action such that its generating vector has degree
$\langle \rho,\check\rho\rangle-1$. But this leads to a contradiction:

By \lemref{no part int}, the image of $\BM_w$ in
$(\bL_{w_0})^1/(\bL_{w_0})^0$ equals the image of the submodule
$\BM_{w_0}\subset \BM_w$, as the quotient is partially
integrable. Hence, $\BM_w$ admits a quotient, which is simultaneously
a quotient module of $\BM_{w_0}$. However, this is impossible, since
we are working with the {\it Kac-Moody} algebra $\BC\cdot
t\partial_t\ltimes \hg_\crit$, and it is known that for Kac-Moody
algebras, Verma modules have simple and mutually non-isomorphic
co-socles.

\end{proof}

Now we are ready to finish the proof of \propref{surj of Wak on irr}. 
Consider the $\nilp$-version of the renormalized universal
enveloping algebra at the critical level, $U^{\ren,\nilp}(\hg_\crit)$,
see \secref{action of renorm}.  We have a natural homomorphism
$U^{\ren,\nilp}(\hg_\crit)\to U^{\ren,\reg}(\hg_\crit)$.

Consider the $\hslash$-family of $\hg_\hslash$-modules equal to
$\BM_{-2\rho+\kappa_\hslash(\check\rho,\cdot)}$. Its specialization at
$\hslash=0$ is the module $\BM_{w_0}$; 
and hence it acquires a $U^{\ren,\nilp}(\hg_\crit)$-action.

\begin{lem}
The map $\BM_{w_0}\otimes \omega_x^{\langle \rho,\check\rho\rangle}\to \bL_{w_0}$ 
is compatible with the $U^{\ren,\nilp}(\hg_\crit)$-actions.
\end{lem}

\begin{proof}
This follows from the fact that the map $\BM_{w_0}\otimes 
\omega_x^{\langle \rho,\check\rho\rangle}\to
\Gamma(\Gr_G,j_{w_0\cdot \check\rho,*}\underset{I}\star
\delta_{1_{Gr_G}})$, constructed in \secref{one map}, deforms away
from the critical level.
\end{proof}

\medskip

By \thmref{FF for algebroids} and 
\corref{nilpotent directions}, we have a short exact sequence
$$0\to N^*_{\fZ^\nilp_\fg/\fZ_\fg}|_{\Spec(\fZ^\reg_\fg)}\to
N^*_{\fZ^\reg_\fg/\fZ_\fg}\to (\cg/\cb)_{\fZ^\reg_\fg}\to 0.$$

Let $L_{-1}=\partial_t$ be the renormalized Sugawara operator, which
we view as an element of $(U^{\ren,\reg}(\hg_\crit))^1$.  By
\propref{grading on algebroid}, the image of $L_{-1}$ in
$N^*_{\fZ^\reg_\fg/\fZ_\fg}\twoheadrightarrow
(\cg/\cb)_{\fZ^\reg_\fg}$ is a principal nilpotent element. Hence,
$N^*_{\fZ^\nilp_\fg/\fZ_\fg}|_{\Spec(\fZ^\reg_\fg)}$ and $L_{-1}$
generate $N^*_{\fZ^\reg_\fg/\fZ_\fg}$ as an algebroid. This, in turn,
implies that $L_{-1}$ and
$(U^{\ren,\nilp}(\hg_\crit))^1\underset{\fZ_\fg^\nilp}\otimes
\fZ_\fg^\reg$ generate $(U^{\ren,\reg}(\hg_\crit))^1$ as an algebroid
over $(U^{\ren,\reg}(\hg_\crit))^0$.

\medskip

Thus, to prove \propref{surj of Wak on irr}, it remains to check that
$L_{-1}$ preserves $(\bL_{w_0})^0$. Since $L_{-1}$ normalizes
$(U^{\ren,\reg}(\hg_\crit))^0$, and since $(\bL_{w_0})^0$ is generated
over $(U^{\ren,\reg}(\hg_\crit))^0$ by its highest weight line, it
suffices to show that $L_{-1}$ maps this highest weight line to
$(\bL_{w_0})^0$.

However, the image of the highest weight line under $L_{-1}$ has
degree $\langle \rho,\check\rho\rangle-1$, and our assertion follows
from \lemref{deg -1}. This completes the proof of \propref{surj of Wak
  on irr}. \qed

\ssec{}

We conclude this section by the following observation:

\begin{prop}
For every $\chi\in \Spec(\fZ^\reg_\fg)$, the module 
$\bL_{w_0}\underset{\fZ^\reg_\fg}\otimes \BC_{\chi}$
is irreducible.
\end{prop}

\begin{proof}

Let us observe that, on the one hand, \corref{kernel part int} 
implies that the module $\BM_{w_0}\underset{\fZ^\reg_\fg}\otimes
\BC_{\chi}$ has a unique irreducible quotient, denoted 
$\bL_{w_0,\chi}$, such that the kernel of the projection
$$\BM_{w_0}\underset{\fZ^\reg_\fg}\otimes\BC_{\chi}\to \bL_{w_0,\chi}$$
is partially integrable. 

On the other hand, by \thmref{Wakimoto via Grassmannian}, the
above projection factors through
$$\BM_{w_0}\underset{\fZ_\fg^\reg}\otimes \BC_{\chi}\twoheadrightarrow 
\bL_{w_0}\underset{\fZ_\fg^\reg}\otimes \BC_{\chi}\to \bL_{w_0,\chi}.$$

Thus, we obtain a surjective map 
$\bL_{w_0}\underset{\fZ_\fg^\reg}\otimes \BC_{\chi}\to \bL_{w_0,\chi}$,
whose kernel is partially integrable. However, by \lemref{no part int},
we conclude that this map must be an isomorphism.

\end{proof}

\section{Comparison with semi-infinite cohomology}     \label{comp
  with semiinf}

\ssec{}

Consider the group ind-scheme $N^-\ppart$, and let $\Psi_0$ denote a
non-degenerate character $N^-\ppart\to \BG_a$ of conductor $0$. This
means that the restriction of $\Psi_0$ to $N^-[[t]]$ is trivial, and
its restriction to $\on{Ad}_{t^{\check\alpha_\iota}}(N^-[[t]])\subset
N^-\ppart$ is non-trivial for each $\iota\in \CI$. Note that to specify
$\Psi_0$ one needs to make a choice: e.g., of a non-vanishing $1$-form
on $\D$, in addition to a choice of $\psi:N^-\to \BG_a$.

For a coweight $\check\lambda$, let $\Psi_{\check\lambda}$ denote the
character obtained as a composition
$$N^-\ppart\overset{\on{Ad}(t^{\check\lambda})} \longrightarrow
N^-\ppart\overset{\Psi_0} \longrightarrow \BG_a.$$ Note that for
$\check\lambda=-\check\rho$, this character is canonical, modulo a
choice of $\psi$ (the latter we will consider fixed).

We will identify $N^-$ and $N$ by means of conjugation by a chosen
left of the element $w_0\in W$; and denote by the same symbol
$\Psi_{\check\lambda}$ the corresponding character on $N\ppart$. We
will also use the same notation for the corresponding characters on
the Lie algebras.

In this section we will study semi-infinite cohomology of
$\fn^-\ppart$ twisted by the characters $\Psi_{\check\lambda}$ with
coefficients in $\hg_\crit$-modules. The complex computing
semi-infinite cohomology was introduced by Feigin \cite{Fei}; the
construction is recalled in \secref{standard complex}. We denote it by
$$\CM\mapsto \fC^\semiinf(\fn^-\ppart,?,\CM\otimes \Psi_{\check\lambda}),$$
where $?$ stands for a choice of a lattice in $\fn^-\ppart$. Its
cohomology will be denoted by
$$
H^{\semiinf+\bullet}(\fn^-\ppart,\fn^-[[t]],\CM\otimes
\Psi_{\check\lambda}).
$$

\begin{prop}  \label{passage between characters}
For $\CM\in D(\hg_\crit\mod)^{I^0}$,
$$\fC^\semiinf(\fn\ppart,\fn[[t]],\CM\otimes \Psi_0)\simeq
\fC^\semiinf\left(\fn^-\ppart,t\fn^-[[t]],(\wt{j}_{w_0\cdot
\check\rho,*}\underset{I^0}\star\CM) \otimes
\Psi_{-\check\rho}\right).$$
\end{prop}

We do not give the proof, since it essentially repeats the proof of
\propref{lambda equivariance}. (In particular, the assertion is valid
at any level $\kappa$.)

Another important observation (also valid at any level) is the following:

\begin{lem} \label{semiinf kills part int}
If $\CM\in \hg_\crit\mod^{I^0}$ is partially integrable, then
$$H^{\semiinf+\bullet}(\fn^-\ppart,\fn^-[[t]],\CM\otimes
\Psi_{-\check\rho})=0.$$
\end{lem}

\begin{proof}

We can assume that $\CM$ is integrable with respect to ${\mathfrak
{sl}}_2^\iota$ for some $\iota\in \CI$. Let $f_\iota\in \fn^-\subset
\fn^-\ppart$ be the corresponding Chevalley generator. With no
restriction of generality, we can assume that
$\Psi_{-\check\rho}(f_\iota)=1$.

Consider the complex $\fC^\semiinf(\fn^-\ppart,\fn^-[[t]],\CM\otimes
\Psi_{-\check\rho})$, and recall (see \secref{standard complex}) that
we have an action of $\fn^-\ppart[1]$ on it by "annihilation
operators", $x\mapsto i(x)$, and the action of $\fn^-\ppart$ by Lie
derivatives $x\mapsto \Lie_x$ such that
$$[d,i(x)]=\Lie_x+\on{Id}\cdot \Psi_{-\check\rho}(x).$$

Hence, $i(f_\iota)$ defines a homotopy between the identity map on
$\fC^\semiinf(\fn^-\ppart,\fn^-[[t]],\CM\otimes \Psi_{-\check\rho})$ and
the map given by $\Lie_{f_\iota}$. However, by assumption, the latter
acts locally nilpotently, implying the assertion of the lemma.

\end{proof}

\ssec{}

In the rest of this section we will collect several additional facts
concerning the semi-infinite cohomology functor
$H^\semiinfi(\fn^-\ppart,t\fn^-[[t]],?\otimes \Psi_{-\check\rho})$.  By
\lemref{semiinf kills part int}, this functor, when restricted to
$\DCat$, factors through $\DCatf$.

\begin{thm} \label{semiinf is exact}
The two functors $\DCatf\to D(\Vect)$
$$\CM^\bullet\mapsto H^{\semiinf+\bullet}(\fn^-\ppart,t\fn^-[[t]],\
\CM^\bullet\otimes \Psi_{-\check\rho}) \text{ and } \CM^\bullet\mapsto
\on{Hom}(\Ppi,\CM^\bullet)$$ are isomorphic. In particular, for $0\neq
\CM\in \Catf$, we have
$$H^\semiinfi(\fn^-\ppart,t\fn^-[[t]],\CM\otimes \Psi_{-\check\rho})=0
\text{ for } i\neq 0 \text{ and }
H^\semiinf(\fn^-\ppart,t\fn^-[[t]],\CM\otimes \Psi_{-\check\rho})\neq 0.$$
\end{thm}

The proof of the theorem is based on the following observation:

\begin{lem}  \label{semiinf of Verma}
For any $\fZ^\nilp_\fg$-module $\CL$ and $w\in W$, we have:
$$
H^\semiinfi\bigl(\fn^-\ppart,t\fn^-[[t]],(\BM_{w}
  \underset{\fZ^\nilp_\fg}\otimes \CL) \otimes
  \Psi_{-\check\rho}\bigr)\simeq \begin{cases}
&  \CL, \,\, i=0, \\ & 0,\,\, i\neq 0.
\end{cases}
$$
\end{lem}

\begin{proof}
Since the quotients $\BM_1/\BM_w$ are all partially integrable, we can
assume that the element $w\in W$, appearing in the lemma, equals
$1$. In the latter case, the assertion follows from
\propref{acyclicity of wakimoto} and \corref{other Wakimoto=Verma}.
\end{proof}

Let us now prove \thmref{semiinf is exact}.

\begin{proof}

In view of \mainthmref{equiv of quot}, to prove the theorem we have to
establish an isomorphism
$$H^\semiinf(\fn^-\ppart,t\fn^-[[t]],\Ppi\otimes \Psi_{-\check\rho})\simeq
\fZ^\nilp_\fg\otimes h_0.$$

Consider the filtration on $\Ppi$, induced by the tilting filtration
on $\Pi$ with quotients $M_w$.  By \lemref{semiinf of Verma},
we obtain that
$H^\semiinfi(\fn^-\ppart,t\fn^-[[t]],\Ppi\otimes \Psi_{-\check\rho})=0$
for $i\neq 0$, and that $H^\semiinf(\fn^-\ppart,t\fn^-[[t]],\Ppi\otimes
\Psi_{-\check\rho})$ has a filtration, with subquotients isomorphic to
$\fZ^\nilp_\fg$.

Hence it remains to show that
$H^\semiinf(\fn^-\ppart,t\fn^-[[t]],\Ppi\otimes \Psi_{-\check\rho})$ is
flat as an $h_0$-module, where the action of $h_0$ is induced from the
identification $h_0\simeq \on{End}(\Pi)$.

It suffices to check that
$$\left(H^\semiinf(\fn^-\ppart,t\fn^-[[t]],\Ppi\otimes
\Psi_{-\check\rho})\right)\overset{L} {\underset{h_0}\otimes}
\BC\simeq \fZ^\nilp_\fg.$$

By \lemref{Pi in quotient},
${}^f\Pi\overset{L}{\underset{h_0}\otimes} \BC \simeq {}^f M_{w_0}$.
Hence,
$$\left(H^\semiinf(\fn^-\ppart,t\fn^-[[t]],\Ppi\otimes
\Psi_{-\check\rho})\right)\overset{L} {\underset{h_0}\otimes}
\BC\simeq H^\semiinf(\fn^-\ppart,t\fn^-[[t]],\BM_{w_0}
\otimes\Psi_{-\check\rho}),$$ and the assertion follows from
\lemref{semiinf of Verma}.

\end{proof}

\begin{cor}  \label{semiinf with n plus is right exact}
For any object $\CM$ of $\Cat$ and a dominant coweight $\check\lambda$
$$H^\semiinfi(\fn\ppart,\fn[[t]],\CM\otimes \Psi_{\check\lambda})=0
\text{ for } i>0.$$
\end{cor}

\begin{proof}
We have
$$H^\semiinfi(\fn\ppart,\fn[[t]],\CM\otimes \Psi_{\check\lambda})\simeq
H^\semiinf(\fn^-\ppart,t\fn^-[[t]],
\wt{j}_{\check\lambda,*}\underset{I^0}\star \wt{j}_{w_0\cdot
\check\rho,*}\underset{I^0}\star\CM \otimes \Psi_{-\check\rho}),$$ as
in \propref{passage between characters}.

Now the assertion of the corollary follows from the fact that the
functor $\CM\mapsto \wt{j}_{\check\lambda,*}\underset{I^0}\star
\wt{j}_{w_0\cdot \check\rho,*}\underset{I^0}\star\CM$ is right exact.

\end{proof}

As another application, we give an alternative proof of the following
result of \cite{FB} (see Theorem 15.1.9).

\begin{thm}   \label{semiinf of vacuum}
The natural map $\fZ^\reg_\fg\to H^\semiinf\left(\fn\ppart,\fn[[t]],
\BV_\crit \otimes \Psi_0\right)$ is an isomorphism, and all other
cohomologies $H^{\semiinf+i}\left(\fn\ppart,\fn[[t]], \BV_\crit
\otimes \Psi_0\right), i \neq 0$, vanish.
\end{thm}

\begin{proof}

Consider the map $$\BM_{w_0,\reg}\otimes
\omega_x^{\langle\rho,\crho\rangle}\to j_{w_0\cdot
\crho,*}\underset{I}\star \BV_{\crit}$$ of \secref{one map}. Its
kernel and cokernel are partially integrable; hence it induces
isomorphisms
$$H^\semiinfi(\fn^-\ppart,t\fn^-[[t]],\BM_{w_0,\reg}
\otimes\Psi_{-\check\rho})\to
H^\semiinfi(\fn^-\ppart,t\fn^-[[t]],j_{w_0\cdot
\crho,*}\underset{I}\star \BV_{\crit}\otimes\Psi_{-\check\rho}).$$

By \lemref{semiinf of Verma},
$$
H^\semiinfi(\fn^-\ppart,t\fn^-[[t]],\BM_{w_0,\reg}
\otimes\Psi_{-\check\rho})\simeq \begin{cases}
& \fZ^\reg_\fg, \,\, i=0, \\
& 0,\,\, i\neq 0,
\end{cases}
$$
we obtain that
$$
H^\semiinfi(\fn^-\ppart,t\fn^-[[t]],j_{w_0\cdot \crho,*}\underset{I}
\star \BV_{\crit}\otimes\Psi_{-\check\rho})\simeq \begin{cases}
& \fZ^\reg_\fg, \,\,
i=0, \\ & 0,\,\,
i\neq 0.
\end{cases}
$$

Applying \propref{passage between characters} for $\CM=\BV_\crit$, we
obtain that $H^{\semiinf+i}\left(\fn\ppart,\fn[[t]], \BV_\crit\otimes
\Psi_0\right)=0$ for $i\neq 0$, and
$$H^\semiinf\left(\fn\ppart,\fn[[t]],\BV_\crit \otimes
\Psi_0\right)\simeq \fZ^\reg_\fg.$$

Moreover, by unraveling the isomorphism of \propref{passage between
characters}, we obtain that the above isomorphism coincides with the
one appearing in the statement of the theorem.

\end{proof}

\ssec{}

Let $\CF$ be a critically twisted D-module on $\Gr_G$. In this
subsection we will express the semi-infinite cohomology
\begin{equation}  \label{semiinf of D-module}
H^\semiinfb\left(\fn\ppart,\fn[[t]], \Gamma(\Gr_G,\CF)\otimes \Psi_0\right)
\end{equation}
in terms of the de Rham cohomologies of $\CF$ along the
$N\ppart$-orbits in $\Gr_G$.

For a coweight $\check\lambda$, consider the $N\ppart$-orbit of the
point $t^{\check\lambda}$ on $\Gr_G$; by pulling back $\CF$,
by \secref{D-mod on groups and quotients}, we obtain a D-module on
$N\ppart$. We will denote it by $\CF|_{N\ppart\cdot t^{\check\lambda}}$.
If $\Psi_0$ is a non-degenerate character of conductor $0$, we will
denote by $H^\bullet(N\ppart,\CF|_{N\ppart\cdot
t^{\check\lambda}}\otimes \Psi_0)$ the resulting de Rham
cohomology. Note that this cohomology vanishes automatically unless
$\check\lambda$ is dominant, since otherwise $\Psi_0$ would be
non-trivial on the stabilizer of $t^{\check\lambda}\in \Gr_G$.

By decomposing $\CF$ in the derived category with respect to the
stratification of $\Gr_G$ by $N\ppart\cdot t^{\check\lambda}$, using
\secref{semiinf and de rham}, we obtain that, as an object of the
derived category of $\fZ_\fg^\reg$-modules,
$\fC^\semiinfb\left(\fn\ppart,\fn[[t]], \Gamma(\Gr_G,\CF)\otimes
\Psi_0\right)$ is a successive extension of complexes
\begin{equation}  \label{lambda contribution}
H^\semiinfb\left(\fn\ppart,\fn[[t]],\Gamma(\Gr_G,
\delta_{t^{\check\lambda}}) \otimes \Psi_0\right) \otimes
H^\bullet(N\ppart,\CF|_{N\ppart\cdot t^{\check\lambda}}\otimes
\Psi_0).
\end{equation}

Note also that $\Gamma(\Gr_G,\delta_{t^{\check\lambda}})$ is
isomorphic to the vacuum module, twisted by $t^{\check\lambda}\in
T\ppart$. Hence,
\begin{equation}  \label{twist module}
H^\semiinfb\left(\fn\ppart,\fn[[t]],
\Gamma(\Gr_G,\delta_{t^{\check\lambda}}) \otimes \Psi_0\right) \simeq
H^\semiinfb\left(\fn\ppart,\on{Ad}_{t^{\check\lambda}}(\fn[[t]]),
\BV_\crit\otimes \Psi_{\check\lambda}\right).
\end{equation}

We will prove the following:

\begin{thm} \label{spec seq splits}
\hfill

\smallskip

\noindent{\em (1)} For $\CF\in \fD(\Gr_G)_\crit\mod$ there is a
canonical direct sum decomposition
\begin{align*}
&H^\semiinfb\left(\fn\ppart,\fn[[t]], \Gamma(\Gr_G,\CF)\otimes
\Psi_0\right)\simeq \\ &\underset{\check\lambda}\oplus
H^\semiinfb\left(
\fn\ppart,\fn[[t]],\Gamma(\Gr_G,\delta_{t^{\check\lambda}})
\otimes \Psi_0\right) \otimes H^\bullet(N\ppart,\CF|_{N\ppart\cdot
t^{\check\lambda}}\otimes \Psi_0).
\end{align*}

\smallskip

\noindent{\em (2)}
The cohomology $H^\semiinfi\left(\fn\ppart,\fn[[t]], \BV_\crit
\otimes \Psi_{\check\lambda}\right)$ vanishes unless $\check\lambda$ is 
dominant and $i=0$, and in the latter case, it is canonically isomorphic to
$V^{\check\lambda}_{\fZ^\reg_\fg}$.

\end{thm}

The rest of this section is devoted to the proof of this theorem. Let
us first prove point (2). The fact that the semi-infinite cohomology
in question vanishes unless $\check\lambda$ is dominant, follows by
the same argument as in \lemref{semiinf kills part int}. Therefore,
let us assume that $\check\lambda$ is dominant and consider the
D-module $\CF_{V^{\check\lambda}}$, see \secref{Dmod on Gr}.

By the geometric Casselman-Shalika formula, see \cite{FGV},
$$H^\bullet(N\ppart,\CF_{V^{\check\lambda}}|_{N\ppart\cdot
t^{\check\mu}}\otimes \Psi_0)=0$$ unless $\mu=\lambda$. Therefore, all
terms with $\mu\neq \lambda$ in the spectral sequence \eqref{lambda
contribution} vanish. We obtain, therefore,
$$H^\semiinfb\left(\fn\ppart,\fn[[t]],
\Gamma(\Gr_G,\CF_{V^{\check\lambda}})\otimes \Psi_0\right) \simeq
H^\semiinfb\left(\fn\ppart,\fn[[t]], \BV_\crit \otimes
\Psi_{\check\lambda}\right).$$ But by \thmref{BD},
$\Gamma(\Gr_G,\CF_{V^{\check\lambda}})\simeq
\BV_\crit\underset{\fZ^\reg_\fg}\otimes
V^{\check\lambda}_{\fZ^\reg_\fg}$. By combining this with
\thmref{semiinf of vacuum}, we obtain
\begin{equation}  \label{character lambda}
H^\semiinf\left(\fn\ppart,\fn[[t]], \BV_\crit \otimes
\Psi_{\check\lambda}\right) \simeq V^{\check\lambda}_{\fZ^\reg_\fg},
\quad H^\semiinfi\left(\fn\ppart,\fn[[t]], \BV_\crit \otimes
\Psi_{\check\lambda}\right)=0, i\neq 0,
\end{equation}
as required.

\medskip

To prove point (1) we need some preparations.

\begin{prop} \label{action of algebroid on semiinf cohomology}
Suppose that $\CM$ is an object of $\hg_\crit\mod_\reg$ that
comes by restriction from a $U^{\ren,\reg}(\hg_\crit)$-module.
Then all $H^\semiinfi(\fn\ppart,\fn[[t]],\CM\otimes \Psi_0)$ are 
naturally modules over the algebroid $N^*_{\fZ^\reg_\fg/\fZ_\fg}$.
\end{prop}

\begin{proof}

We will assume that $\Psi_0$ comes from a character of the Lie-*
algebra $L_\fn$. In this case, the BRST complex
$\fC^\semiinf(L_\fn,\CA_{\fg,\crit}\otimes \Psi_0)$ is itself a
DG-chiral algebra.

Let $\CA_{\fg,\hslash}$ be a 1-st order deformation of
$\CA_{\fg,\crit}$ away from the critical level; i.e.,
$\CA_{\fg,\hslash}$ is flat over $\BC[\hslash]/\hslash^2$, and
$\CA_{\fg,\hslash}/\hslash\simeq \CA_{\fg,\crit}$.

\medskip

Let us consider the DG-chiral algebra
$\fC^\semiinf(L_\fn,\CA_{\fg,\hslash}\otimes \Psi_0)$.  From
\thmref{semiinf of vacuum}, it follows that it is acyclic off
cohomological degree $0$; in particular, its $0$-th cohomology is
$\BC[\hslash]/\hslash^2$-flat.

This implies that any section $a\in \fz_\fg$, which we think of as a
$0$-cocycle in $\fC^\semiinf(L_\fn,\CA_{\fg,\crit}\otimes \Psi_0)$, can
be lifted to a $0$-cocycle $a_\hslash\in
\fC^\semiinf(L_\fn,\CA_{\fg,\hslash}\otimes \Psi_0)$.

We will think of $\frac{a_\hslash}{\hslash}$ an element of the Lie-*
algebra $\CA^\sharp_\fg\otimes \on{Cliff}(L_\fn)$, where
$\CA^\sharp_\fg$ is as in \cite{FG}, and $\on{Cliff}(L_\fn)$ is the
Clifford chiral algebra, used in the definition of the BRST complex.

\medskip

By the construction of $\CA^\sharp_\fg$, for $\CM$ satisfying the
properties of the proposition, we have an action of $\CA^\sharp_\fg$
on $\CM$, and hence, an action of the Lie-* algebra
$\CA^\sharp_\fg\otimes \on{Cliff}(L_\fn)$ on the complex
$\fC^\semiinf(L_\fn,\CM\otimes \Psi_0)$. By taking the Lie-* bracket
with the above element $\frac{a_\hslash}{\hslash}$ we obtain an
endomorphism of $\fC^\semiinf(L_\fn,\CM\otimes \Psi_0)$, which
commutes with the differential.

It is easy to see that for a different choice of $a_\hslash$ the
corresponding endomorphisms of $\fC^\semiinf(L_\fn,\CM\otimes \Psi_0)$
will differ by a coboundary. Thus, we obtain a Lie-* action of
$\fz_\fg$ on each $H^\semiinfi(L_\fn,\CM\otimes \Psi_0)$. One easily
checks that this action satisfies the Leibniz rule with respect to the
$\fz_\fg$-module structure on $H^\semiinfi(L_\fn,\CM\otimes \Psi_0)$,
and hence extends to an action of the Lie-* algebroid
$\Omega^1(\fz_\fg)$. The latter is the same as an action of the
$\fZ^\reg_\fg$-algebroid $N^*_{\fZ^\reg_\fg/\fZ_\fg}$.

\end{proof}

We are now ready to finish the proof of \thmref{spec seq splits}.  By
\secref{action of renorm} and \propref{action of algebroid on semiinf
cohomology}, the terms of the spectral sequence \eqref{lambda
contribution} are acted on by the algebroid
$N^*_{\fZ^\reg_\fg/\fZ_\fg}$.

It is easy to see that the $N^*_{\fZ^\reg_\fg/\fZ_\fg}$-action on
$H^\semiinf\left(\fn\ppart,\fn[[t]], \BV_\crit \otimes \Psi_0\right)$
identifies via \thmref{semiinf of vacuum} with the canonical
$N^*_{\fZ^\reg_\fg/\fZ_\fg}$-action on $\fZ^\reg_\fg$. Moreover, from
\thmref{BD}(a) we obtain that the isomorphisms of \eqref{character
lambda} are compatible with the $N^*_{\fZ^\reg_\fg/\fZ_\fg}$-action.

This implies the canonical splitting of the spectral sequence. Indeed,
from \thmref{BD}(b) it is easy to derive that there are no non-trivial
Hom's and $\on{Ext}^1$'s between different
$V^{\check\lambda}_{\fZ^\reg_\fg}$, regarded as
$N^*_{\fZ^\reg_\fg/\fZ_\fg}$-modules.

\newpage

\vspace*{10mm}

{\Large \part{Appendix}}

\vspace*{10mm}


This Part, which may be viewed as a user's guide to Sect. 7 of
\cite{BD}, reviews some technical material that we need in the main
body of this paper.

In \secref{misc} we review some background material:
the three monoidal structures on the category of topological
vector spaces, the notion of a family of objects of an abelian 
category over a scheme or an ind-scheme, and the formalism
of DG-categories.

In \secref{groups act on cat} we introduce the notion of action of
a group-scheme on an abelian category. In fact, there are two such
notions that correspond to weak and strong actions, respectively. A
typical example of a weak action is when a group $H$ acts on a scheme
$S$, and we obtain an action of $H$ on the category $\QCoh_H$ of
quasi-coherent sheaves.  A typical example of a strong (equivalently,
infinitesimally trivial or Harish-Chandra type) action is when in the
above situation we consider the action of $H$ on the category
$\fD(S)\mod$ of D-modules on $S$. We also discuss various notions
related to equivariant objects and the corresponding derived
categories.

In \secref{D-mod on group ind-schemes} we make a digression and
discuss the notion of D-module over a group ind-scheme. The
approach taken here is different, but equivalent, to the one developed
in \cite{AG} via chiral algebras.

In \secref{conv section} we generalize the discussion of
\secref{groups act on cat} to the case of group ind-schemes. The goal
of this section is to show that if $\CC$ is a category that carries a
Harish-Chandra action of some group ind-scheme $G$, then at the level
of derived categories we have an action of the monoidal category of
D-modules over $G$ on $\CC$.  This formalism was developed in Sect. 7
of \cite{BD}, and in this section we essentially repeat it.

Finally, \secref{cat over top alg} serves a purely auxiliary role: we
prove some technical assertions concerning the behavior of an abelian
category over its center provided that a certain flatness assumption
is satisfied.

\bigskip

\section{Miscellanea}   \label{misc}

Unless specified otherwise, the notation in this part will be
independent of that of Parts I--IV. We will work over the ground field
$\BC$, and all additive categories will be assumed
$\BC$-linear. Unless specified otherwise, by tensor product, we will
mean tensor product over $\BC$.

If $\CC$ is a category, and $X_i$ is a directed system of objects in it,
then following the notation of SGA 4(I) notation, we write $"\underset{\longrightarrow}{\lim}"\, X_i$ 
for the resulting object in $\on{Ind}(\CC)$, thought of as a
contravariant functor on $\CC$. In contrast, $\underset{\longrightarrow}{\lim}\, X_i$ will denote the
object of $\CC$ representing the functor $\underset{\longleftarrow}{\lim}\,
\on{Hom}(\lim X_i,?)$ on
$\CC$, provided that it exists.

\ssec{Topological vector spaces and algebras}  \label{topological
  vector spaces}

In this subsection we will briefly review the material of
\cite{CHA:top}. By a topological vector space we will mean a vector
space over $\BC$ equipped with a linear topology, assumed complete and
separated. We will denote this category by $\CTop$; it is closed under
projective and inductive limits (note that the projective limits
commute with the forgetful functor to vector spaces, and inductive
limits do not). Every such topological vector space $\bV$ can be
represented as $\underset{\longleftarrow}{\lim}\, V^i$, where $V^i$
are usual (i.e., discrete) vector spaces and the transition maps
$V^j\to V^i$ are surjective.

For a topological vector space $\bV$ represented as projective limit as
above, its dual $\bV^*$ is by definition the object of $\CTop$ equal
to
$$\underset{\longrightarrow}{\lim}\, (V_i)^*,$$ where each $(V_i)^*$
is dual of the corresponding discrete vector space $V_i$, endowed with
the natural (pro-finite dimensional) topology. It is easy to see that
$\bV^*$ is well-defined, i.e., independent of the presentation of
$\bV$ as a projective limit.

A topological vector space $\bV$ is said to be of Tate type if it can
be written in the form $\bV_1\oplus \bV_2$, where $\bV_1$ is discrete
and $\bV_2$ is pro-finite dimensional.  In this case $\bV^*$ is also
of Tate type, and the natural map $(\bV^*)^*\to \bV$ is an
isomorphism.

\medskip

Following \cite{CHA:top}, we endow the category $\CTop$ with three
different monoidal structures:
$$\bV_1, \bV_2 \mapsto \bV_1\startimes \bV_2, \bV_1\arrowtimes\bV_2
\text{ and } \bV_1\shriektimes \bV_2.$$ They are constructed as
follows. Let us write $\bV_1=\underset{\longleftarrow}{\lim}\, V_1^i$,
$\bV_2=\underset{\longleftarrow}{\lim}\, V_2^j$. Then
$$\bV_1\shriektimes
\bV_2=\underset{i,j}{\underset{\longleftarrow}{\lim}}\, V_1^i\otimes
V_2^j.$$ It is easy to see that this monoidal structure is in fact a
tensor one.

To define $\bV_1\arrowtimes\bV_2$ we proceed in 2 steps. If $\bV_2=V$
is discrete and equal to $\underset{k}\cup\, V_k$, where $V_k$ are
finite-dimensional, we set
$$\bV_1\arrowtimes V= \underset{k}{\underset{\longrightarrow}{\lim}}\,
\bV_1\otimes V_k,$$ where the inductive limit is taken in $\CTop$.
For an arbitrary $\bV_2$ written as
$\bV_2=\underset{\longleftarrow}{\lim}\, V_2^j$, we set
$$\bV_1\arrowtimes\bV_2=\underset{j}{\underset{\longleftarrow}{\lim}}\,
(\bV_1\arrowtimes V^j_2).$$

Finally, $\bV_1\startimes \bV_2$ is characterized by the property that
$\Hom(\bV_1\startimes \bV_2,V)$, where $V$ is discrete, is the set of
bi-linear continuous maps $\bV_1\times \bV_2\to W$. This monoidal
structure is also tensor in a natural way.

We have natural maps
$$\bV_1\startimes \bV_2\to \bV_1\arrowtimes \bV_2\to \bV_1\shriektimes
\bV_2,$$ where the first arrow is an isomorphism if $\bV_2$ is
discrete and the second one is an isomorphism if $\bV_1$ is discrete.

Note also that for three objects $\bV_1,\bV_2,\bV_3\in \CTop$ we have
natural maps
$$(\bV_1\shriektimes \bV_2)\arrowtimes \bV_3\to \bV_1\shriektimes
(\bV_2\arrowtimes \bV_3) \text{ and } \bV_1\arrowtimes
(\bV_2\shriektimes \bV_3) \to (\bV_1\arrowtimes \bV_2)\shriektimes
\bV_3$$ and hence the map
\begin{equation} \label{exchange map}
(\bV_1\shriektimes \bV_2)\startimes \bV_3\to \bV_1\shriektimes
(\bV_2\startimes \bV_3).
\end{equation}

By an action of a topological vector space $\bV$ from a discrete
vector space $W_1$ to a discrete vector space $W_2$ we will mean a map
$$\bV\arrowtimes W_1\simeq \bV\startimes W_1\to W_2.$$ The latter
amounts to a compatible system of maps $V'\otimes W'_1\to W_2$,
defined for every finite-dimensional subspace $W'_1\subset W_1$ for
some sufficiently large discrete quotient $V'$ of
$\bV$.

\medskip

By definition, a topological associative algebra is an object $\bA\in
\CTop$ endowed with an associative algebra structure with respect to
the $\arrowtimes$ product.  By construction, any such $\bA$ can be
represented as $\underset{\bI}{\underset{\longleftarrow}{\lim}}\,
\bA/\bI$, where $\bI\subset \bA$ are open left ideals. A discrete
module over a topological associative algebra $\bA$ is a vector space
$V$ endowed with an associative action map $\bA\arrowtimes V\to V$;
we shall denote the category of discrete $\bA$-modules by $\bA\mod$.

\medskip

A topological associative algebra is called commutative if the
operation $\bA\arrowtimes \bA\to \bA$ factors through $\bA\shriektimes
\bA\to \bA$ and the latter map is commutative (in the sense of the
commutativity constraint for the $\shriektimes$ product). In this case
$\bA$ can be represented as
$\underset{i}{\underset{\longleftarrow}{\lim}}\, A_i$, where $A_i$ are
discrete commutative quotients of $\bA$.

For a commutative associative topological algebra, by a topological
$\bA$-module we shall mean a topological vector space $\bV$, endowed
with an associative map
$\bA\shriektimes\bV\to \bV$ such that $\bV$ is separated and complete
in the topology defined by open $\bA$-submodules.  Any such $\bV$ can
be represented as $\underset{\longleftarrow}{\lim}\, V^i$, with $V_i$
being discrete $\bA$-modules, on each of which $\bA$ acts through
a discrete quotient. If $f:\bA\to \bB$ is a homomorphism, we
define $f^*(\bV)$ as $\underset{\longleftarrow}{\lim}\,
\bB\underset{\bA}\otimes V^i$.

Note that if we regard $\bA$ as an associative topological algebra,
a discrete $\bA$-module is a topological $\bA$-module in the above
sense if and only if $\bA$ acts on it though some discrete quotient.

\medskip

A topological Lie algebra $\fg$ is a Lie algebra in the sense of the
$\startimes$ structure.  A discrete module over such $\fg$ is a vector
space $V$ endowed with a map $\fg\startimes V\to V$, which is
compatible with the bracket on $\fg$ in a natural way.

\medskip

Let $\bA$ be a commutative associative topological algebra. A Lie
algebroid over $\bA$ is a topological Lie algebra $\fg$ endowed with a
topological $\bA$-module structure $\bA\shriektimes \fg\to \fg$ and a
Lie algebra action map $\fg\startimes \bA\to \bA$, which satisfy the
usual compatibility conditions via \eqref{exchange map}.

\ssec{}   \label{topological algebroids}

Here we shall recall some notions related to infinite-dimensional
vector bundles and ind-schemes, borrowed from \cite{BD} and
\cite{Dr1}.

By an ind-scheme we will understand an ind-object in the category of
schemes, which can be represented as $\CY:=\underset{i\in
I}{\underset{\longrightarrow}{"\lim"}}\, \CY_i$, where the
transition maps $f_{i,j}:\CY_i\to \CY_j$ are closed embeddings. We
will always assume that the indexing set $I$ is countable.

A closed subscheme $Z$ of $\CY$ is called reasonable if for every $i$,
the ideal of the subscheme $Z\cap \CY_i$ of $\CY_i$ is locally
finitely generated. The ind-scheme $\CY$ is called reasonable if it
can be represented as an inductive limit of its reasonable subschemes
(or, in other words, one can choose a presentation such that the
ideal of $\CY_i$ in $\CY_j$ is locally finitely generated).

We shall say that $\CY$ is ind-affine if all the schemes $\CY_i$ are
affine. In this case, if we denote by $A_i$ the algebra of functions
of $\CY_i$, we will write $\CY=\Spec(\bA)$, where
$\bA=\underset{\longleftarrow}{\lim}\, A_i$ and $\bA=\CO_\CY$.

\medskip

Assume that $\CY=G$ is ind-affine and is endowed with a structure of
group ind-scheme.  This amounts to a co-associative co-unital map
$\CO_G\to \CO_G\shriektimes \CO_G$.  By definition, an action of $G$
on a topological vector space $\bV$ is a map
$$\bV\to \CO_G\shriektimes \bV,$$ such that the two morphisms
$$\bV\rightrightarrows \CO_G\shriektimes \CO_G\shriektimes \bV$$
coincide.

If $\bV$ is an associative or Lie topological algebra, we define in an
evident way what it means for an action to be compatible with the
operation of product on $\bV$.

\medskip

Assume now that $G$ is a group-scheme $H=\Spec(\CO_H)$.

\begin{lem}
Every $\bV$, acted on by $H$, can be written as 
$\underset{\longleftarrow}{\lim}\, V_i$, where 
$V_i\in \RepH$ are quotients of $\bV$.
\end{lem}

\begin{proof}

Let $V$ be some discrete quotient of $\bV$. We must show that we can
find an $H$-stable quotient $V'$ such that $\bV\twoheadrightarrow
V'\twoheadrightarrow V$.  Consider the map
$$\bV\to \CO_H\shriektimes \bV\to \CO_H\otimes V,$$
Let $\bV'$ be the kernel of this map; this is an open subspace in $\bV$.
The associativity of the action implies that $\bV'$ is $H$-stable. Hence,
$\bV/\bV'$ satisfies our requirements.

\end{proof}

Let $\CY$ be an ind-scheme. A topological *-sheaf on $\CY$ is a rule
that assigns to a commutative algebra $R$ and an $R$-point $y$ of
$\CY$ a topological $R$-module $\CF_{y}$, and for a morphism of
algebras $f:R\to R'$ an isomorphism $\CF_{y}\simeq f^*(\CF_{y'})$,
where $y'$ is the induced $R'$-point of $\CY$, compatible with
two-fold compositions.  Morphisms between topological *-sheaves are
defined in an evident manner and we will denote the resulting category
by $\QCoh_\CY^{\on{top},*}$. The cotangent sheaf $\Omega^1(\CY)$ is an
example of an object of $\QCoh_\CY^{\on{top},*}$.

We let $\on{\on{Tate}}_\CY$ denote the full subcategory of
$\QCoh_\CY^{\on{top},*}$ formed by Tate vector bundles (i.e., those,
for which each $\CF_{y}$ is an $R$-module of Tate type), see
\cite{Dr1}, Sect. 6.3.2. The following basic result was established in
\cite{Dr1}, Theorem 6.2:

\begin{thm}  \label{Tate families}
Let $\CY_1\to \CY_2$ be a formally smooth morphism between ind-schemes
with $\CY_1$ being reasonable. Then the topological *-sheaf of
relative differentials $\Omega^1(\CY_1/\CY_2)$ is a Tate vector bundle
on $\CY_1$.
\end{thm}

\medskip

Assume now that $\CY$ is affine and isomorphic to $\Spec(\bA)$ for a
commutative associative topological algebra $\bA$. In this case, the
category $\QCoh_\CY^{\on{top},*}$ is tautologically equivalent to that
of topological $\bA$-modules. We have the notion of Lie
algebroid over $\CY$ (which is the same as a topological Lie algebroid
over $\bA$).

Let now $\CG$ be an ind-groupoid over an ind-affine ind-scheme $\CY$,
such that both (r, equivalently, one of the) projections
$l,r:\CG\rightrightarrows \CY$ is formally smooth. Then, by the above
theorem, the normal to $\CY$ in $\CG$, denoted $N_{\CY/\CG}$, which is
by definition the dual of the restriction to $\CY$ of
$\Omega^1(\CG/\CY)$ with respect to either of the projections, is a
Tate vector bundle. The standard construction endows it with a
structure of Lie algebroid.

\ssec{A class of categories}   \label{class of cat}

Let $\CC$ be an abelian category, and let $\Ind(\CC)$ denote its
ind-completion.  We will assume that $\CC$ is closed under inductive
limits, i.e., that the tautological embedding $\CC\to \Ind(\CC)$
admits a right adjoint $\on{limInd}:\Ind(\CC)\to \CC$, and that the
latter functor is exact. In particular, it makes sense to tensor
objects of $\CC$ by vector spaces.

We shall say that an object $X\in \CC$ is finitely generated (or
compact) if the functor $\Hom(X,\cdot):\CC\to \Vect$ commutes with
direct sums. Let us denote by $\CC^c$ the full subcategory of $\CC$
formed by compact objects. We will assume that $\CC^c$ is equivalent
to a small category (i.e., that isomorphism classes of compact objects
in $\CC$ form a set).

\medskip

We shall say that $\CC$ satisfies (*) if every object of $\CC$ is
isomorphic to the inductive limit of its compact subobjects.

\begin{lem}
Assume that $\CC$ satisfies (*), and let $\sG$ be a left exact
contravariant functor $\CC^c\to \Vect$. The following conditions are
equivalent:

\smallskip

\noindent{\em(1)}
$\sG$ is representable by $X\in \CC$.

\smallskip

\noindent{\em(2)} For an inductive system $\{X_i\}\in \CC^c$, whenever
$X:=\underset{\longrightarrow}{\lim}\, X_i$ belongs to $\CC^c$, the
natural map
$$\sG(X)\to \underset{\longleftarrow}{\lim}\, \sG(X_i)$$
is an isomorphism.

\smallskip

\noindent{\em(3)} $\sG$ extends to a functor $\CC\to \Vect$ such that
for any inductive system $\{X_i\}\in \CC$, the map
$$\sG\left(\underset{\longrightarrow}{\lim}\, X_i\right)\to
 \underset{\longleftarrow}{\lim}\, \sG(X_i)$$
is an isomorphism.

\end{lem}

\medskip

We shall say that $\CC$ satisfies (**) if there exists an exact and
faithful covariant functor $\sF:\CC\to \Vect$, which commutes with
inductive limits.

The following is standard:
\begin{lem}  Assume that $\CC$ satisfies (*) and (**). Then:

\smallskip

\noindent{\em (1)} $\sF$ is representable by some
$\underset{i}{\underset{\longleftarrow}{"\lim"}}\, X_i\in
\Pro(\CC^c)$.

\smallskip

\noindent{\em (2)} Assume that $\sF$ has the following additional
property. Whenever a system of maps $\alpha_k:X\to Y_k$ is such that
for any non-zero subobject $X'\subset X$ not all maps $\alpha_k|_{X'}$
are zero, then the map
$$\sF(X)\to \underset{k}\Pi\, \sF(Y_k)$$ is injective. 

Then the projective system $\{X_i\}$ as above can be chosen so that
all the transition maps $X_{i'}\to X_i$ are surjective.

\smallskip

\noindent{\em (3)} Under the assumption of (2), the functor $\sF$
gives rise an to an equivalence $\CC\to \bA\mod$, where $\bA$ is the
topological associative algebra
$\underset{i}{\underset{\longleftarrow}{"\lim"}}\, \sF(X_i)\simeq
\End(\sF)$.

\end{lem}

\ssec{}    \label{objects over algebras}

If $A$ is an associative algebra, we will denote by $A\mod\otimes\CC$
the category, whose objects are objects of $\CC$, endowed with an
action of $A$ by endomorphisms, and morphisms being $\CC$-morphisms,
compatible with $A$-actions. This is evidently an abelian category.

If $M$ is a left $A$-module and $X\in \CC$, we produce an example of
an object of $A\mod\otimes\CC$ by taking $M\otimes X$.

\medskip

Let $M$ be a right $A$-module. We have a naturally defined right exact
functor
$$A\mod\otimes \CC\to \CC:X\mapsto M\underset{A}\otimes X.$$

\begin{lem}  \label{ff}
If $M$ is flat (resp., faithfully-flat $A$-algebra), then the above
functor is exact (resp., exact and faithful).
\end{lem}

For the proof see \cite{Ga1}, Lemma 4 and Proposition
5. \footnote{Whereas the first of the assertion of the lemma is
obvious from Lazard's lemma, the second is less so, and it was pointed
out to us by Drinfeld.}

\medskip

We will say that $X\in A\mod\otimes\CC$ is $A$-flat if the functor
$M\mapsto M\underset{A}\otimes X:A^{\on{op}}\mod\to \CC$ is exact. The
functor of tensor product can be derived in either (or both) arguments
and we obtain a functor
$$D^-(A^{\on{op}}\mod)\times D^-(A\mod\otimes \CC)\to D^-(\CC).$$

\medskip

If $M$ is a left $A$-module and $X\in A\mod\otimes\CC$, we define a
contravariant functor on $\CC$ by
$$Y\mapsto \on{Hom}_{\CC\otimes A\mod}(Y\otimes M,X).$$

This functor is representable by an object that we will denote by
$\Hom_A(M,X)$. If $M$ is finitely presented as an $A$-module, the
functor $X\mapsto \Hom_A(M,X)$ commutes with inductive limits.

\medskip

Let $\phi:A\to B$ be a homomorphism of algebras. We have a natural
forgetful functor $\phi_*:B\mod\otimes\CC\to A\mod\otimes\CC$, and its
left adjoint $\phi^*$, given by tensor product with $B$, viewed as a
right $A$-module. The right adjoint to $\phi_*$, denoted $\phi^!$, is
given by $X\mapsto \on{Hom}_A(B,X)$.

\medskip

\ssec{Objects parameterized by a scheme}

Assume now that $A$ is commutative and set $S=\Spec(A)$. In this case
we will use the notation $\QCoh_S\otimes \CC$ for
$A\mod\otimes\CC$. We will think of objects of $\QCoh_S\otimes\CC$ as
families of objects of $\CC$ over $S$.

For a morphism of affine scheme $f:S_1\to S_2$ we have the direct and
inverse image functors $f_*,f^*:\QCoh_{S_1}\otimes \CC\rightleftarrows
\QCoh_{S_2}\otimes \CC$, with $f^*$ being exact (resp., exact and
faithful) if $f$ is, by \lemref{ff}.

The usual descent argument shows:

\begin{lem}  \label{descent}
Let $S'\to S$ be a faithfully flat map. Then the category
$\QCoh_S\otimes \CC$ is equivalent to the category of descent data on
$\QCoh_{S'}\otimes \CC$ with respect to $S'\underset{S}\times
S'\rightrightarrows S'$.
\end{lem}

\medskip

This allows us to define the category $\QCoh_S\otimes \CC$ for any
separated scheme $S$. Namely, let $S'$ be an affine scheme covering
$S$. We introduce $\QCoh_S\otimes \CC$ as the category of descent data
on $\QCoh_{S'}\otimes \CC$ with respect to $S'\underset{S}\times
S'\rightrightarrows S'$.  \lemref{descent} above ensures that
$\QCoh_S\otimes \CC$ is well-defined, i.e., is independent of the
choice of $S'$ up to a unique equivalence.  (In fact, the same
definition extends more generally to stacks algebraic in the
faithfully-flat topology, for which the diagonal map is affine.)  For
a morphism of schemes $f:S_1\to S_2$ we have the evidently defined
direct and inverse image functors. If $f$ is a closed embedding and
the ideal of $S_1$ in $S_2$ is locally finitely generated, then we
also have the functor $f^!:\QCoh_{S_2}\otimes \CC\to
\QCoh_{S_1}\otimes \CC$, right adjoint to $f_*$.

\medskip

If $S_1$ is a closed subscheme of $S_2$ we say that an object $X\in
\QCoh_{S_2}\otimes \CC$ is set-theoretically supported on $S_1$, if
$X$ can be represented as an inductive limit of its subobjects, each
of which is the direct image of an object in some $\QCoh_{S'_1}\otimes
\CC$, where $S'_1$ is a nilpotent thickening of $S_1$ inside $S_2$.

\medskip

Suppose now that $S$ is of finite type over $\BC$. We will denote by
$\fD(S)\mod$ the category of right D-modules on $S$. We define the
category $\fD(S)\mod\otimes \CC$ as follows:

First, we assume that $S$ is affine and smooth. Then
$\fD(S)\mod\otimes \CC$ is by definition the category
$\Gamma(S,\fD(S))^{\on{op}}\mod\otimes \CC$.

If $S_1\to S_2$ is a closed embedding of affine smooth schemes, we
have an analog of Kashiwara's theorem, saying that
$\fD(S_1)\mod\otimes \CC$ is equivalent to the subcategory of
$\fD(S_2)\mod\otimes \CC$, consisting of objects set-theoretically
supported on $S_1$, when considered as objects of $\QCoh_{S_2}\otimes
\CC$.

This allows to define $\fD(S)\mod\otimes \CC$ for any affine scheme of
finite type, by embedding it into a smooth scheme. Finally, for an
arbitrary $S$, we define $\fD(S)\mod\otimes \CC$ using a cover by
affine schemes, as above.

\ssec{}

In this subsection we will assume that $\CC$ satisfies (*). Let $\bV$
be a topological vector space, and $X,Y\in \CC$. An action $\bV\times
X\to Y$ is a map
$$\bV\otimes X\to Y,$$ satisfying the following continuity condition:
For every compact subobject $X'\subset X$, the induced map $\bV\otimes
X'\to Y$ factors through $V\otimes X'\to Y$, where $V$ is a discrete
quotient of $\bV$.

\medskip

If $X'\to X$ (resp., $Y\to Y'$, $\bV'\to \bV$) is a map, and we have
an action $\bV\times X\to Y$, we produce an action $\bV\times X'\to Y$
(resp., $\bV\times X\to Y'$, $\bV'\times X\to Y$).

Note that if $\bV$ is pro-finite-dimensional, with the dual $\bV^*\in
\Vect$, an action $\bV\times X\to Y$ is the same as a map $X\to
\bV^*\otimes Y$.

\begin{lemconstr}
Let $\bV_2\times X\to Y$ and $\bV_1\times Y\to Z$ be actions. Then we
have an action
$$(\bV_1\arrowtimes \bV_2)\times X\to Z.$$
\end{lemconstr}

\begin{proof}

The construction immediately reduces to the case when $\bV_2=V_2$ is discrete,
$X$ is compact, and we have an action map $V_2\otimes X\to Y$. 

Then for every finite-dimensional subspace $V_2^k\subset V_2$ we can
find a compact subobject $Y^k\subset Y$, such that $V^k_2\otimes X\to
Y$ maps to $Y^k$ and the action $\bV_1\otimes Y^k\to Z$ factors
through a discrete quotient $V_1^k$ of $\bV_1$. Then
$$\underset{k}\cup\, \on{ker}(\bV_1\to \bV_1^k)\otimes V_2^k\subset 
\bV_1\arrowtimes V_2$$
is an open neighborhood of $0$,
and we have an action map
$$\left(\bV_1\arrowtimes V_2/\underset{k}\cup\, \on{ker}(\bV_1\to
\bV_1^k)\otimes V_2^k\right) \otimes X\simeq
\underset{\longrightarrow}{\lim} (V^k_1\otimes V^k_2)\otimes X\to Z.$$

\end{proof}

We shall say that $\bV$ acts on $X$, if we are given a map $\bV\times
X\to X$.  Objects of $\CC$, acted on by $\bV$ naturally form category,
which is abelian.

\medskip

Let $\bA$ be an associative topological algebra. We shall say that an
object $X\in \CC$ is acted on by $\bA$ if we are given an action map
$\bA\times X\to X$ such that the two resulting action maps
$(\bA\arrowtimes \bA)\times X\rightrightarrows X$ coincide. Objects of
$\CC$ acted on by $\bA$ form a category, denoted $\bA\mod\otimes \CC$.

\ssec{Objects of a category parameterized by an ind-scheme}
\label{objects over indscheme}

In this subsection we retain the assumption that $\CC$ satisfies
(*). Let $\CY$ be an ind-scheme, $\CY=\underset{i}\cup\, \CY_i$. We
introduce the category $\QCoh^*_\CY\otimes \CC$ to have as objects
collections $\{X_i\in \QCoh^*_{\CY_i}\otimes \CC\}$ together with a
compatible system if isomorphisms $f_{i,j}^*(X_j)\simeq X_i$, where
$f_{i,j}$ is the map $\CY_i\to \CY_j$.  Morphisms in the category are
evident.

It is easy to see that this category is independent of the
presentation of $\CY$ as an inductive limit. However,
$\QCoh^*_\CY\otimes \CC$ is, in general, not abelian.

Given an object of $\bX\in \QCoh^*_\CY\otimes \CC$ and a
scheme $S$ mapping to $\CY$, we have a well-defined object
$\bX|_{S}\in \QCoh_S\otimes \CC$.

\medskip

Assume now that $\CY$ is strict and reasonable. I.e., the system
$\CY_i$ can be chosen so that the maps $f_{i,j}$ are closed
embeddings, and the ideal of $\CY_i$ inside $\CY_j$ is locally
finitely generated.

We introduce the category $\QCoh^!_\CY\otimes \CC$ as follows. Its
objects are collections $\bX:=\{X_i\in \QCoh^*_{\CY_i}\otimes \CC\}$
together with a compatible system if isomorphisms $X_i\simeq
f_{i,j}^!(X_j)$. The morphisms in this category are evident.

\begin{lem}
$\QCoh^!_\CY\otimes \CC$  is an abelian category.
\end{lem}

\begin{proof}

If $\alpha=\{\alpha_i:X_i\to X'_i\}$ is a morphism in
$\QCoh^!_\CY\otimes \CC$, its kernel is given by the system
$\{\on{ker}(\alpha_i)\}$. It is easy to see that 
the cokernel and image of this morphism are given by the systems that
assign to each $i$
$$\underset{\underset{j\geq i}{\longrightarrow}}{\lim}\,
f_{i,j}^!\left(\on{coker}(\alpha_j)\right),\,\, 
\underset{\underset{j\geq i}{\longrightarrow}}{\lim}\,
f_{i,j}^!\left(\on{Im}(\alpha_j)\right),$$
respectively. The fact that the axioms of an
abelian category are satisfied is shown in the
same way as in the case of $\CC=\on{Vect}$.

\end{proof}

\medskip

Let now $\bA$ be a commutative topological algebra. Then $\bA$ can be
represented as $\underset{\longleftarrow}{\lim}\, A_i$, where $A_i$
are discrete commutative algebras.  Assume, moreover, that we can find
such a presentation that the ideal of $A_i$ in each $A_j$, $j\geq i$
is finitely generated. Then $\CY:=\underset{\longrightarrow}{\lim}\,
\Spec(A_i)$ is reasonable.

\begin{lem}
Under the above circumstances, the categories $\QCoh^!_\CY\otimes \CC$
and $\bA\mod\otimes \CC$ are equivalent.
\end{lem}

\begin{proof}
The functor $\QCoh^!_\CY\otimes \CC\to \bA\mod\otimes \CC$ is
evident. Its right adjoint is defined as follows: given an object
$X\in \bA\mod\otimes \CC$, represented as $\underset{i}\cup\, X_i$
with $X_i\in A_i\mod\otimes \CC$, we define an object $\{X'_i\}$ in
$\QCoh^!_\CY\otimes \CC$ by setting
$$X'_i=\underset{\underset{j\geq i}{\longrightarrow}}{\lim}\,
f_{i,j}^!(X_j).$$

The fact that the adjunction morphisms are isomorphisms is
shown as in the case $\CC=\on{Vect}$.

\end{proof}

\medskip

Let now $\bX=\{X_i\}$ be an object of $\QCoh^*_\CY\otimes \CC$ such
that each $X_i$ is $\CY_i$-flat. Let $M$ be an object in
$\QCoh_\CY^!\otimes \bA\mod$, where $\bA$ is some topological
algebra. We then have a well-defined tensor product
$$\bX\underset{\CO_\CY}\otimes M\in \QCoh_\CY^!\otimes (\bA\mod\otimes
\CC).$$ The corresponding system assigns to every $\CY_i$ the object
$$\bX|_{\CY_i}\underset{\CO_{\CY_i}}\otimes M_i\in \QCoh_{\CY_i}\otimes 
(\bA\mod\otimes \CC).$$

\medskip

Finally, let $\CY$ be a strict ind-scheme of ind-finite
type. Proceeding as above, one defines the category
$\fD(\CY)^!\mod\otimes \CC$ as the category of systems $\{X_i\}\in
\fD(\CY_i)^!\mod\otimes \CC$ with isomorphisms $X_i\simeq
f_{i,j}^!(X_j)$.

If $\CY$ is formally smooth, we can also introduce the DG-category of
$\Omega_\CY^\bullet$-modules with coefficients in $\CC$, and we will
have an equivalence between the corresponding derived category of
$\Omega_\CY^\bullet\mod\otimes \CC$ and the derived category of
$\fD(\CY)^!\mod\otimes \CC$.

\ssec{BRST complex}    \label{standard complex}

If $\fg$ is a topological Lie algebra, an action of $\fg$ on $X\in
\CC$ is an action map $\fg\times X\to X$ such that the difference of
the two iterations
$$(\fg\startimes \fg)\times X\to (\fg\arrowtimes \fg)\times X\to X$$
equals the action induced by the Lie bracket $\fg\startimes \fg\to \fg$.

Assume now that $\fg\simeq \sk$ is pro-finite dimensional. Then its
action on $X$ is the same as a co-action of the Lie co-algebra
$\sk^*\in \Vect$ on $X$, i.e., a map $a:X\to \sk^*\otimes X$,
satisfying the suitable axioms. In this case we can form a complex of
objects of $\CC$, called the standard complex, $\fC(\sk,X)$:

As a graded object of $\CC$, it is isomorphic to
$\fC(\sk,\CM):=X\otimes \Lambda^\bullet(\sk^*)$.  Let us denote by $i$
(resp., $i^*$) the action of $\sk[1]$ (resp., $\sk^*[-1]$) on
$\fC(\sk,X)$ by the "annihilation" (resp., "creation" operators), and
by $\Lie$ the diagonal action of $\sk$.  Then the differential $d$ on
$\fC(\sk,X)$ is uniquely characterized by the property that
$[d,i]=\Lie$. We automatically obtain that

\begin{itemize}

\item
$d^2=0$, 

\item 
The map $i^*:\Lambda^\bullet(\sk^*)\otimes \fC(\sk,X)\to \fC(\sk,X)$
is a map of complexes, where $\Lambda^\bullet(\sk^*)$ is endowed with
a differential coming from the Lie co-bracket.
\end{itemize}

If $X^\bullet$ is a complex of objects of $\CC$, acted on by $\sk$, we
will denote by $\fC(\sk,X^\bullet)$ the complex associated to the
corresponding bi-complex. It is clear that if $X^\bullet$ is bounded
from below and acyclic, then $\fC(\sk,X^\bullet)$ is acyclic as
well. However, this would not be true if we dropped the boundedness
from below assumption.

\medskip

The above set-up can be generalized as follows. Let now $\fg$ be a
topological Lie algebra, which is of Tate type as a topological vector
space. Let $\sCliff(\fg,\fg^*)$ be the (topological) Clifford algebra,
constructed on $\fg$ and $\fg^*$; it is naturally graded, where the
"creation" operators (i.e., elements of $\fg^*$) have degree $1$, and
the annihilation operators (i.e., elements of $\fg$) have degree
$-1$. Let $\sSpin(\fg)$ be some fixed irreducible representation of
$\sCliff(\fg)$, equipped with a grading.  (Of course, up to a grading
shift and a non-canonical isomorphism, $\sSpin(\fg)$ is unique.)

Recall that the canonical (i.e. Tate's) central extension $\fg_{\can}$
of $\fg$ is characterized by the property that the adjoint action of
$\fg$ on $\sCliff(\fg,\fg^*)$ is inner via a homomorphism
$\fg_{\can}\to \sCliff(\fg,\fg^*)$.  We will denote by $\fg_{-\can}$
the Baer negative central extension.

\medskip

Let $\CC$ be as above, and let $X\in \CC$ be acted on by
$\fg_{-\can}$.  Consider the graded object of $\CC$ given by
$$\fC^{\frac{\infty}{2}}(\fg,X):=X\otimes \sSpin(\fg).$$

As in the case when $\CC=\Vect$, one shows that
$\fC^{\frac{\infty}{2}}(\fg,X)$ acquires a canonical differential $d$,
characterized by the property that $[d,i]=\on{Lie}$, where $i$ denotes
the action of $\fg$ on $\fC^{\frac{\infty}{2}}(\fg,X)$ via
$\sSpin(\fg)$ by creation operators, and $\Lie$ is the diagonal action
of $\fg$ on $X\otimes \sSpin(\fg)$. We have

\begin{itemize}

\item $d^2=0$

\item The action $i^*$ of $\fg^*$ is compatible with the differential
$\fg^*\to \fg^*\shriektimes\fg^*$ given by the bracket.

\end{itemize}

If $X^\bullet$ is a complex of objects of $\CC$, acted on by
$\fg_{-\can}$, we will denote by
$\fC^{\frac{\infty}{2}}(\fg,X^\bullet)$ the complex, associated to the
corresponding bi-complex.

\begin{lem}
Assume that $X^\bullet$ is bounded from below and is acyclic. Then
$\fC^{\frac{\infty}{2}}(\fg,X^\bullet)$ is also acyclic.
\end{lem}

\begin{proof}

Let us choose a lattice $\sk\subset \fg$; we can then realize
$\sSpin(\fg)$ as $\sSpin(\fg,\sk)$--the module generated by an
element, annihilated by both $\sk\subset \fg\subset
\Lambda^\bullet(\fg)$ and $(\fg/\sk)^*\subset \fg^* \subset
\Lambda^\bullet(\fg^*)$.

In this case the complex $\fC^{\frac{\infty}{2}}(\fg,X^\bullet)$
acquires a canonical increasing filtration, numbered by natural
numbers, so that
$$\on{gr}^i\left(\fC^{\frac{\infty}{2}}(\fg,X^\bullet)\right)\simeq
\fC(\sk,X^\bullet\otimes \Lambda^i(\fg/\sk))[i].$$

This readily implies the assertion of the lemma.

\end{proof}

In what follows we will need to consider the following situation. Let
$X^\bullet$ be a complex of objects of $\CC$, endowed with two actions
of $\fg_{-\can}$, denoted $a$ and $a'$, respectively. Then
$X^\bullet\otimes \sSpin(\fg)$ acquires two differentials, $d$ and
$d'$.

Assume that there exists a self-anticommuting action
$$i_h:\fg[1]\times X^\bullet\to X^\bullet,$$ 
such that $a'(x)-a(x)=[d_X,i_h(x)]$, $[a'(x),i_h(y)]=i_h([x,y])$,
$[a(x),i_h(y)]=0$, where $d_X$ is the differential on $X^\bullet$.

\begin{lem}  \label{homotopy lemma}
Under the above circumstances, there exists a graded automorphism 
of the complex $X^\bullet\otimes \sSpin(\fg)$ that intertwines $d$ and $d'$.
\end{lem}

\begin{proof} 

Let $\Lambda^i(\fg)$ and $\Lambda^i(\fg^*)$ denote the !-completed
exterior powers of $\fg$ and $\fg^*$, respectively.

For a natural number $i$ consider the canonical element $\on{id}_i\in
\Lambda^i(\fg)\shriektimes \Lambda^i(\fg^*)$. We define the operator
$$T:X^\bullet\otimes \sSpin(\fg)\to X^\bullet\otimes \sSpin(\fg)$$ by
$\underset{i\in \BN}\Sigma\, (i_h\otimes i^*)(\on{id}_i)$, where $i^*$ and
$i_h$ denotes the extension of the actions of $\fg^*[-1]$ and
$\fg[1]$, respectively, to the exterior powers.

Clearly, $T$ is a grading-preserving isomorphism, and
$$T\circ i\circ T^{-1}=i+i_h.$$
One easily shows that $d'=T^{-1}\circ d\circ T$.
 
\end{proof}

\ssec{DG categories}    \label{DG categories}

We will adopt the conventions regarding DG categories from \cite{Dr}. 
Let $\bC$ be a $\BC$-linear DG category,  which admits arbitrary direct 
sums.

For $X^\bullet,Y^\bullet\in \bC$ we will denote by
$\CHom_{\bC}(X^\bullet,Y^\bullet)$ the corresponding complex,
and by $\Hom_{\bC}(X^\bullet,Y^\bullet)$ its $0$-th cohomology. By
definition, the homotopy category $Ho(\bC)$ has the same objects as
$\bC$, with the Hom space being $\Hom_{\bC}(X^\bullet,Y^\bullet)$

We will assume that $\bC$ is strongly pre-triangulated, i.e., that 
it admits cones. In this case $Ho(\bC)$ is triangulated.

\medskip

We will assume that $\bC$ is equipped with a cohomological functor
$\sH$ to an abelian category $\CC'$. We will denote by $D(\bC)$ the
corresponding localized triangulated category, and we will assume that
$\sH$ defines a t-structure on $D(\bC)$. We will denote by
$\RHom_{D(\bC)}(\cdot,\cdot)$ the resulting functor $D(\bC)^{\on{op}}\times
D(\bC)\to D(\Vect)$.

We will denote by $D^b(\bC)$ (resp., $D^+(\bC)$, $D^-(\bC)$,) the
subcategory consisting of objects $X^\bullet$, such that
$\sH(X^\bullet[i])=0$ for $i$ away from a bounded interval (resp.,
$i<<0$, $i>>0$.)

In what follows we will also use the following notion: we shall say
that an object $X^\bullet\in D(\bC)$ is quasi-perfect if it belongs to
$D^-(\bC)$, and the functor $Y \mapsto \Hom_{D(\bC)}(X^\bullet, Y[i])$
commutes with direct sums in the core of $\bC$ (i.e., those objects
$Y\in \bC$ for which $\sH(Y[j])=0$ for $j\neq 0$).

\begin{lem}
Let $X^\bullet \in D(\bC)$ be quasi-perfect and $Y^\bullet\in
D^+(\bC)$. Let $\CK^\bullet$ be a bounded from below complex of vector
spaces. Then $\RHom_{D(\bC)}(X^\bullet,Y^\bullet\otimes \CK^\bullet)$
is quasi-isomorphic to $\RHom_{D(\bC)}(X^\bullet,Y^\bullet)\otimes
\CK^\bullet$ in $D(\Vect)$.
\end{lem}

\medskip

The most typical example of this situation is, of course, when
$\bC=\bC(\CC)$ is the category of complexes of objects of an abelian
category $\CC$, and $\sH$ comes from an exact functor $\CC\to
\CC'$. If $\CC'\simeq \CC$, then $D(\bC)$ will be denoted $D(\CC)$;
this is the usual derived category of $\CC$.

An example of a quasi-perfect object of $D(\CC)$ is provided by a
bounded from above complex consisting of projective finitely generated
objects of $\CC$.

\medskip

Let $\bC_1,\bC_2$ be two DG categories as above, and let $\sG:\bC_1\to
\bC_2$ be a DG functor. We shall say that $\sG$ is exact if it sends
acyclic objects (in the sense of $\sH_1$) to acyclic ones (in the
sense of $\sH_2$).

The following (evident) assertion will be used repeatedly:

\begin{lem} \label{der adj}
Let $\sG:\bC_1\to \bC_2$ and $\sG':\bC_2\to \bC_1$
be mutually adjoint exact functors. Then
the induced functors $\sG,\sG':D(\bC_1)\rightleftarrows D(\bC_2)$ are
also mutually adjoint. 
\end{lem}

\begin{proof}

Let $\sG$ be the left adjoint of $\sG'$. Note first of all that the
functors induced by $\sG$ and $\sG'$ between the homotopy categories
$Ho(\bC_1)$ and $Ho(\bC_2)$ are evidently mutually adjoint.

Then for $X^\bullet\in D(\bC_1),
Y^\bullet\in D(\bC_2)$
\begin{equation} \label{RHom 1}
\Hom_{D(\bC_1)}(X^\bullet,\sG'(Y^\bullet))=
\underset{\underset{X'{}^\bullet\to
X^\bullet}{\longrightarrow}}{\lim}\,
\Hom_{Ho(\bC_1)}(X'{}^\bullet,\sG'(Y^\bullet))
\end{equation}
and
\begin{equation} \label{RHom 2}
\Hom_{D(\bC_2)}(\sG(X^\bullet),Y^\bullet)=
\underset{\underset{Y^\bullet\to
Y'{}^\bullet}{\longrightarrow}}{\lim}\,
\Hom_{Ho(\bC_2)}(\sG(X^\bullet),Y'{}^\bullet),
\end{equation}
where in both cases the inductive limits are taken over
quasi-isomorphisms, i.e., morphisms in the homotopy category that
become isomorphisms in the quotient triangulated category.

By adjunction, we rewrite the expression in \eqref{RHom 2} as
$$\underset{\underset{Y^\bullet\to
Y'{}^\bullet}{\longrightarrow}}{\lim}\,
\Hom_{Ho(\bC_2)}(X^\bullet,\sG'(Y'{}^\bullet)),$$ and we map it to
\eqref{RHom 1} as follows. For a quasi-isomorphism $Y^\bullet\to
Y'{}^\bullet$ the map $\sG(Y^\bullet)\to \sG(Y'{}^\bullet)$ is a
quasi-isomorphism as well, and given a map $X^\bullet\to
\sG'(Y'{}^\bullet)$, we can find a quasi-isomorphism $X'{}^\bullet\to
X^\bullet$, so that the diagram
$$
\CD
\sG'(Y^\bullet) @>>> \sG'(Y'{}^\bullet)  \\
@AAA  @AAA  \\
X'{}^\bullet @>>> X^\bullet
\endCD
$$ commutes in $Ho(\bC_1)$. The above map $X'\to \sG'(Y^\bullet)$
defines an element in \eqref{RHom 1}.

One constructs the map from \eqref{RHom 1} to \eqref{RHom 2} in a
similar way, and it is straightforward to check that the two are
mutually inverse.

\end{proof}

\section{Action of a group on a category}    \label{groups act on cat}

\ssec{Weak action} \label{weak action of groups}

Let $\CC$ be an abelian category as in \secref{class of cat}, and let
$H$ be an affine group-scheme. We will say that $H$ acts weakly on
$\CC$ if we are given a functor
$$\act^*:\CC\to \QCoh_H\otimes \CC,$$ and two functorial isomorphisms
related to it. The first isomorphism is between the identity functor
on $\CC$ and the composition $\CC\overset{\on{act}^*}\to
\on{QCoh}_H\otimes \CC\to \CC$, where the second arrow corresponds to
the restriction to $1\in H$.

\medskip

To formulate the second isomorphism, note that from the existing data
we obtain a natural functor $\on{act}^*_S: \on{QCoh}_S\otimes \CC\to
\on{QCoh}_{S\times H}\otimes \CC$ for any affine scheme $S$.

The second isomorphism is between the two functors $\CC\to
\on{QCoh}_{H\times H}\otimes \CC$ that correspond to the two circuits
of the diagram
\begin{equation}  \label{associativity}
\CD \CC @>{\on{act}^*}>> \on{QCoh}_H\otimes \CC \\ @V{\on{act}^*}VV
@V{\on{act}^*_H}VV \\ \on{QCoh}_H\otimes \CC @>{\on{mult}^*}>>
\on{QCoh}_{H\times H}^*\otimes \CC, \endCD
\end{equation}
where $\on{mult}$ denotes the multiplication map $H\times H\to H$.

We assume that the above two isomorphisms of functors satisfy the
usual compatibility conditions. We will refer to these isomorphisms as
the unit and associativity constraint of the action, respectively.

\begin{lem}   \label{act is flat}
The functor $\act^*$ is exact and faithful. For $X\in \CC$, the
$\CO_H$-family $\act^*(X)$ is flat.
\end{lem}

\begin{proof}
First, the faithfulness of $\act^*$ is clear, since the fiber at $1\in
H$ provides a left quasi-inverse $\on{QCoh}_H\otimes \CC\to \CC$.

Let $S$ be a scheme equipped with a map $\phi:S\to H$. Note that we
have a self-functor $\act_\phi^*:\QCoh_S\otimes \CC\to \QCoh_S\otimes
\CC$ given by
$$\QCoh_S\otimes \CC\overset{\act^*_S}\to \QCoh_S\otimes
\QCoh_H\otimes \CC \overset{(\on{id}_S\times \phi)^*}\to
\QCoh_S\otimes \CC.$$ This is an equivalence of categories and its
quasi-inverse is given by $\act^*_{\phi^{-1}}$, where $\phi^{-1}:S\to
G$ is obtained from $\phi$ by applying the inversion on $H$.  Note
that $\act_\phi^*$ is $\CO_S$-linear.

Let us take $S=H$ and $\phi$ to be inversion map. Then the composition
$\act^*_\phi\circ \act^*:\CC\to \QCoh_H\times \CC$ is isomorphic to
the functor $X\mapsto \CO_H\otimes X$, which is evidently
exact. Hence, $\act^*$ is exact as well.

Similarly, to show that $\act^*(X)$ is $\CO_H$-flat, it suffices to
establish the corresponding fact for $\act^*_\phi\circ \act^*(X)$,
which is again evident.

\end{proof}

Here are some typical examples of weak actions:

\noindent 1) Let $H$ act on an ind-scheme $\CY$. Then the category
$\on{QCoh}^!_\CY$ carries a weak $H$-action.

\noindent 2) Let $H$ act on a topological associative algebra $\bA$
(see \secref{topological algebroids}).  Then the category $\bA\mod$ of
discrete $\bA$-modules carries a weak $H$-action.

\ssec{Weakly equivariant objects}    \label{equiv cat}

Let us denote by $p^*$ the tautological functor 
$$\CC\to \QCoh_H\otimes \CC: X\mapsto \CO_H\otimes X,$$
where $\CO_H$ is the algebra of functions on $H$.

We will say that an object $X\in \CC$ is weakly $H$-equivariant
if we are given an isomorphism 
\begin{equation} \label{weak equivariance}
\on{act}^*(X)\simeq p^*(X), 
\end{equation}
which is compatible with the associativity constraint of the
$H$-action on $\CC$.

Evidently, weakly $H$-equivariant objects of $\CC$ form an abelian
category, which we will denote by $\CC^{w,H}$.  For example, let us
take $\CC$ to be $\Vect$--the category of vector spaces with the
obvious, i.e., trivial, $H$-action. Then $\CC^{w,H}$ is the category
of $H$-modules, denoted $H\mod$, or ${\mathcal Rep}(H)$.

\medskip

Let $X$ be an object of $\CC^{w,H}$, and $V\in \RepH$. We define a new
object $V\ast X\in \CC^{w,H}$ to be $V\otimes X$ as an object of
$\CC$, but the isomorphism $\on{act}^*(V\otimes X)\to p^*(V\otimes X)$
is multiplied by the co-action map $V\to \CO_H\otimes V$.

\medskip

We have an obvious forgetful functor $\CC^{w,H}\to \CC$, and it admits
a right adjoint, denoted $\Av_H^w$, given by $X\mapsto
p_*(\act^*(X))$.  For $X\in \CC^{w,H}$,
$$\Av_H^w(X)\simeq \CO_H\ast X.$$

\medskip

For two objects $X_1,X_2$ of $\CC^{w,H}$ we define
a contravariant functor $\uHom_{\CC}(X_1,X_2)$ of $\RepH$ by
$$\Hom_{\RepH}(V,\uHom_{\CC}(X_1,X_2))=
\Hom_{\CC^{w,H}}(V\ast X_1,X_2).$$
It is easy to see that this functor is representable.

\begin{lem}  \label{inner hom}
Let $X_1$ is finitely generated as an object of $\CC$.  Then the
forgetful functor $\RepH\to \Vect$ maps $\uHom_{\CC}(X_1,X_2)$ to
$\Hom_{\CC}(X_1,X_2)$.
\end{lem}

\begin{proof}
We have the map
$$\Hom_{\CC}(X_1,X_2)\simeq \Hom_{\CC^{w,H}}(X_1,\on{Av}_H^w(X_2))\to
\Hom_{\CC}(X_1,\on{Av}_H^w(X_2))\simeq \Hom_{\CC}(X_1,\CO_H\otimes
X_2),$$ and the latter identifies with $\CO_H\otimes
\Hom_{\CC}(X_1,X_2)$, by the assumption on $X_1$.

This endows $\Hom_{\CC}(X_1,X_2)$ with a structure of $H$-module.
It is easy to see that it satisfies the required adjunction property. 

\end{proof}

Note that since
$$\Hom_{\CC}(X_1,X_2)\simeq \Hom_{\CC^{w,H}}(X_1,\CO_H\ast X_2) \text{
and } \Hom_{\CC^{w,H}}(X_1,X_2)\simeq
\Hom_{H}(\BC,\uHom_{\CC}(X_1,X_2)),$$ we obtain that $X_1$ is finitely
generated as an object of $\CC$ if and only if it is one in
$\CC^{w,H}$.

\ssec{}

Let $\bC(\CC^{w,H})$ denote the DG category of complexes of objects of
$\CC^{w,H}$, and let $D(\CC^{w,H})$ be the corresponding derived
category. By \lemref{der adj}, the forgetful functor $D(\CC^{w,H})\to
D(\CC)$ admits a right adjoint given by $X^\bullet\mapsto
\act^*(X^\bullet)$.

\medskip

For $X_1^\bullet,X_2^\bullet\in D(\CC^{w,H})$ we define a
contravariant cohomological functor
$\uRHom_{D(\CC)}(X_1^\bullet,X_2^\bullet)$ on $D(\RepH)$ by
$$V^\bullet\mapsto \Hom_{D(\CC^{w,H})}(V^\bullet\ast
X_1^\bullet,X_2^\bullet).$$ It is easy to see that this functor is
representable.

\begin{lem} \label{inner hom weak}
Assume that $X_1^\bullet$ is quasi-perfect and $X_2^\bullet\in
D^+(\CC^{w,H})$.  Then the image of
$\uRHom_{D(\CC)}(X_1^\bullet,X_2^\bullet)$ under the forgetful functor
$D(\RepH)\to D(\Vect)$ is quasi-isomorphic to
$\RHom_{D(\CC)}(X_1^\bullet,X_2^\bullet)$.
\end{lem}

The proof repeats that of \lemref{inner hom}.

\ssec{Infinitesimally trivial actions}   \label{HCh action of groups}

Let $H^{(i)}$ be the $i$-th infinitesimal neighborhood of $1$ in $H$,
so that $H^{(0)}=1$ and $H^{(1)}=\on{Spec}(\BC\oplus
\epsilon\cdot\fh^*)$, where $\epsilon^2=0$. Note that if $H$ weakly
acts on $\CC$, the restriction to $H^{(1)}$ yields for every object
$X\in \CC$ a short exact sequence in $\CC$.
$$0\to \fh^*\otimes X\to X^{(1)}\to X\to 0,$$
where $X^{(1)}:=\on{act}^*(X)|_{H^{(1)}}$.

We will say that the action of $H$ on $\CC$ is infinitesimally
trivial, or of Harish-Chandra type, if we are given a functorial
isomorphism
\begin{equation} \label{Lie}
\on{act}^*(X)|_{H^{(1)}}\simeq p^*(X)|_{H^{(1)}},
\end{equation}
such that two compatibility condition (see below) are satisfied. 

\medskip

The first condition is that the isomorphism \eqref{Lie} respects the
identification of the restrictions of both sides to $1\in H$ with
$X$. (In view of this condition, the data of \eqref{Lie} amounts to a
functorial splitting $X\to X^{(1)}$.)

To formulate the second condition, consider the map of schemes
$$(h,h_1)\overset{\on{Ad}}\mapsto \on{Ad}_h(h_1): H\times H^{(1)}\to
H^{(1)}.$$ From \eqref{associativity} and \eqref{Lie} we obtain two
{\it a priori} different identifications
$$\on{Ad}^*(X^{(1)})\rightrightarrows
\CO_{H\times H^{(1)}}\otimes X\in \QCoh_{H\times H^{(1)}}\otimes \CC.$$
Our condition is that these two identifications coincide.

\medskip

Here are some typical examples of Harish-Chandra actions:

\noindent 1) Let $\CY$ be an ind-scheme of ind-finite type acted on by
$H$. Then the category $\fD(\CY)\mod$ carries a $H$-action of
Harish-Chandra type.

\noindent 2) Let $\bA$ be a topological associative algebra, acted on
by $H$, and assume that the derived action of $\fh$ on $\bA$ is inner,
i.e., comes from a continuous map $\fh\to \bA$. Then the action of $H$
on $\bA\mod$ is of Harish-Chandra type.

\medskip

Let now $X$ be an object of $\CC^{w,H}$. Note that in this case we
have two identifications between $\on{act}^*(X)|_{H^{(1)}}$ and
$p^*(X)|_{H^{(1)}}$.  Their difference is a map
$$a^\sharp:X\mapsto \fh^*\otimes X,$$ compatible with the co-bracket
on $\fh^*$, i.e., an action of $\fh$ on $X$, see \secref{topological
vector spaces}. We will call this map "the obstruction to strong
equivariance".

We will say that an object $X\in \CC^{w,H}$ is strongly
$H$-equivariant (or simply $H$-equivariant) if the map $a^\sharp$ is
zero. Strongly equivariant objects form a full subcategory in
$\CC^{w,H}$, which we will denote by $\CC^H$.

\medskip

Let us consider the example, where $\CC=\fD(\CY)\mod$, where $X$ is an
ind-scheme of ind-finite type acted on by $H$. Then
$\fD(\CY)\mod^{w,H}$ is the usual category of weakly $H$-equivariant
D-modules, and $\fD(\CY)\mod^H$ is the category strongly
$H$-equivariant $\fD$-modules.

More generally, if $\CC=\bA\mod$, where $\bA$ is a topological
associative algebra, acted on by $H$, then $\bA\mod^{w,H}$ consists of
$\bA$-modules, endowed with an algebraic action of $H$, compatible
with the action of $H$ on $\bA$.  If the action of $H$ on $\bA$ is of
Harish-Chandra type, and $X\in A\mod^{w,H}$, the map $a^\sharp:X\to
\fh^*\otimes X$ corresponds to the difference of the two actions of
$\fh$ on $X$.

\ssec{}  \label{equiv der cat}

Let $\bC(\CC^H)$ denote the DG category of complexes of objects of
$\CC^H$, and $D(\CC^H)$ the corresponding derived category.  We have a
natural functor $D(\CC^H)\to D(\CC)$, but in general it does not
behave well. Following Beilinson, we will now introduce the "correct"
triangulated category, along with its DG model, that corresponds to
strongly $H$-equivariant objects of $\CC$.

\medskip

Let $\bC(\CC)^H$ be the category whose objects are complexes
$X^\bullet$ of objects of $\CC^{w,H}$, endowed with a map of complexes
$$i^\sharp:X^\bullet\to \fh^*[-1]\ast X^\bullet,$$
such that the following conditions are satisfied:

\begin{itemize}

\item
$i^\sharp$ is a map in $\CC^{w,H}$.

\item
The iteration of $i^\sharp$, viewed as a map 
$X^\bullet\to \Lambda^2(\fh^*)[-2]\ast X^\bullet$, vanishes.

\item
The map $[d,i^\sharp]:X^\bullet\to \fh^*\ast X^\bullet$ equals the map
$a^\sharp$.

\end{itemize}

For two objects $X_1^\bullet$ and $X_2^\bullet$ of $\bC(\CC)^H$ we
define ${\mathcal Hom}^k_{\bC(\CC)^H}(X_1^\bullet,X_2^\bullet)$ to be
the subcomplex of ${\mathcal
Hom}^k_{\bC(\CC^{w,H})}(X_1^\bullet,X_2^\bullet)$ consisting of graded
maps $X_1^\bullet\to X_2^\bullet[k]$ that preserve the data of
$i^\sharp$. This defines on $\bC(\CC)^H$ a structure of DG-category.

Note that the usual cohomology functor defines a cohomological functor
$\bC(\CC)^H\to \CC^H$. We will denote by $D(\CC)^H$ the resulting
localized triangulated category, which henceforth we will refer to as
the "$H$-equivariant derived category of $\CC$".

By construction, the truncation functors $\tau^{<0}$, $\tau^{>0}$ are
well-defined at the level of $\bC(\CC)^H$. Therefore, objects of the
subcategory $D^b(\CC)^H$ (resp., $D^+(\CC)^H$, $D^-(\CC)^H$) can be
realized by complexes in $\bC(\CC)^H$ that are concentrated in
finitely many cohomological degrees (resp., cohomological degrees
$>>-\infty$, $<<\infty$).

\ssec{Examples}

Take first $\CC$ to be $\Vect$, in which case $\CC^{w,H}$ identifies
with the category $\RepH$, and $\CC^H$ is the same as ${\mc R}ep 
(H/H^0)$, where $H^0\subset H$ is the neutral connected component of
$H$.

We will denote the resulting DG category by $\bC(\on{pt}/H)$ and the
triangulated category by $D(\on{pt}/H)$. Note that $\bC(\on{pt}/H)$ is the
standard, i.e., Cartan, DG-model for the $H$-equivariant derived
category of the point-scheme.

\medskip

Consider the de Rham complex on $H$, denoted $\bDR_H$. The
multiplication on $H$ endows $\bDR_H$ with a structure of a DG
co-algebra. The category $\bC(\on{pt}/H)$ is tautologically the same as
the category of DG co-modules over $\bDR_H$.  In particular, $\bDR_H$
itself is naturally an object of $\bC\bigl(\on{pt}/(H\times H)\bigr)$.

\medskip

More generally, let $\CC$ be $\fD(\CY)\mod$ for $\CY$ as above.  In
this case $\bC(\fD(\CY)\mod)^H$ is the DG-model for the
$H$-equivariant derived category on $\CY$ studied in \cite{BD}; in
{\it loc.cit.} it is shown that the corresponding derived category is
equivalent to the category of \cite{BL}.

\ssec{Averaging}    \label{averaging}

Note that for any $\CC$ with an infinitesimally trivial action of $H$
we have a natural tensor product functor
$$V^\bullet,X^\bullet\mapsto V^\bullet\ast X^\bullet:
\bC(\on{pt}/H)\times \bC(\CC)^H\to \bC(\CC)^H,$$ which extends to a
functor $D(\on{pt} /H)\times D(\CC)^H\to D(\CC)^H$.

\medskip

We have a tautological forgetful functor $\bC(\CC)^H\to
\bC(\CC^{w,H})$, and we claim that it admits a natural right adjoint,
described as follows.

We will regard $X^\bullet\in \bC(\CC^{w,H})$ as a complex of objects
of $\CC$, acted on by $\fh$ via $a^\sharp$, and we can form the
standard complex
$$\fC(\fh,X^\bullet):=\Lambda^\bullet(\fh^*)\ast X^\bullet,$$ see
\secref{standard complex}. It is naturally an object of
$\bC(\CC^{w,H})$.  The action of the annihilation operators defines on
$\fC(\fh,X^\bullet)$ a structure of an object in $\bC(\CC)^H$.

The resulting functor $\bC(\CC^{w,H})\to \bC(\CC)^H$ is exact when
restricted to $\bC^+(\CC^{w,H})$, and the corresponding functor
$D^+(\CC)^H\to D^+(\CC^{w,H})$ is the right adjoint to the
tautological forgetful functor, by \lemref{der adj}.

\medskip

We will denote the composed functor
$$\bC(\CC)\overset{\Av_H^w}\longrightarrow \bC(\CC^{w,H})\to
\bC(\CC)^H$$ (and the corresponding functor $D^+(\CC)\to D^+(\CC)^H$)
by $\Av_H$.  This functor is the right adjoint to the forgetful
functor $\bC(\CC)^H\to \bC(\CC)$ (resp., $D(\CC)^H\to D(\CC)$).

Let us consider two examples:

\medskip

\noindent 1) For $\CC=\Vect$, we have $\on{Av}_H(\BC)\simeq \bDR_H\in
\bC(\on{pt}/H)$.

\medskip

\noindent 2)
Let $\CC=\fD(\CY)\mod$, where $\CY$ is an ind-scheme of ind-finite
type, acted on by $H$. The resulting functor at the level of derived
categories $D^+(\fD(\CY)\mod)\to D^+(\fD(\CY)\mod)^H$ is the
corresponding *-averaging functor:
$$\CF\mapsto p_*\circ \on{act}^!(\CF),$$
where $p$ and $\on{act}$ are the two maps $H\times \CY\to \CY$.

\ssec{The unipotent case}   

Assume now that $H$ is connected. We claim that in this case $\CC^H$
is a full subcategory of $\CC$. Indeed, for a $\CC$-morphism
$\phi:X_1\to X_2$ between objects of $\CC^H$, in the diagram
$$
\CD
\act^*(X_1) @>{\act^*(\phi)}>>  \act^*(X_2)  \\
@A{\sim}AA   @A{\sim}AA  \\
\CO_H\otimes X_1  @>>> \CO_H\otimes X_2
\endCD
$$
the bottom arrow is necessarily of the form $\on{id}\otimes \phi'$,
since its derivative along $H$ is $0$, as follows from the condition
that $a^\sharp|_{X_1}=a^\sharp|_{X_2}=0$. Then the unit constraint
forces $\phi'=\phi$.

\medskip

For $H$ connected let us denote by $D(\CC)_{\CC^H}$ the full
subcategory of $D(\CC)$ consisting of objects, whose cohomologies
belong to $\CC^H$.

\begin{prop}    \label{unipotent averaging}
Suppose that the group-scheme $H$ is pro-unipotent. Then the functor
$D^+(\CC)^H\to D^+(\CC)$ is fully-faithful, and it induces an
equivalence $D^+(\CC)^H\simeq D^+(\CC)_{\CC^H}$.
\end{prop}

\begin{proof}

Since $\Av_H:D^+(\CC)\to D^+(\CC)^H$ is the right adjoint to the
functor in question, to prove fully-faithfulness it suffices to show
that the adjunction map gives rise to an isomorphism between the
composition
$$D^+(\CC)^H\to D^+(\CC)\overset{\on{Av}_H}\to D^+(\CC)^H$$
and the identity functor. 

For $X^\bullet\in \bC^+(\CC)^H$, the object $\Av_H(X^\bullet)$ is
isomorphic to the tensor product of complexes
$$\bDR_H \ast X^\bullet,$$ and the adjunction map in question
corresponds to the natural map $\BC\to \bDR_H$.  The latter is a
quasi-isomorphism since $H$ was assumed pro-unipotent.

\medskip

It remains to show that $D^+(\CC)^H$ maps essentially surjectively
onto $D^+(\CC)_{\CC^H}$. For that it is sufficient to show that for
$X^\bullet\in D^+(\CC)_{\CC^H}$, the second adjunction map
$\Av_H(X^\bullet)\to X^\bullet$ is a quasi-isomorphism.

By devissage, we can assume that $X^\bullet$ is concentrated in one
cohomological dimension. However such an object is quasi-isomorphic
(up to a shift) to an object from $\CC^H$, which makes the assertion
manifest.

\end{proof} 

\ssec{Equivariant cohomology}   \label{eq vs non-eq Exts}

For $X^\bullet_1,X^\bullet_2\in \bC(\CC)^H$ we define a contravariant
functor
$$\uCHom_{\CC}(X^\bullet_1,X_2^\bullet):\bC(\on{pt}/H)\to \bC(\Vect)$$ by
$$V^\bullet\mapsto \CHom_{\bC(\CC)^H}(V^\bullet\ast
X^\bullet_1,X^\bullet_2).$$

This functor is easily seen to be representable. When $X_1^\bullet$ is
bounded from above and consists of objects that are finitely generated, the
forgetful functor $\bC(\CC)^H\to \bC(\Vect)$ maps
$\uCHom_{\CC}(X^\bullet_1,X_2^\bullet)$ to
$\CHom_{\bC(\CC)}(X^\bullet_1,X_2^\bullet)$.

Similarly, for $X^\bullet_1\in,X^\bullet_2 \in D(\CC)^H$ the
cohomological functor
$$V^\bullet\mapsto \RHom_{D(\CC)^H}(V^\bullet\ast
X^\bullet_1,X^\bullet_2)$$ is representable by some
$\uRHom_{D(\CC)}(X^\bullet_1,X_2^\bullet)\in D^+(\CC)^H$.  We have the
following assertion, whose proof repeats that of \lemref{inner hom}:

\begin{lem}   \label{Leray}
If $X^\bullet_1$ quasi-perfect as an object of $D(\CC)$ and
$X^\bullet_2$ is bounded from below. Then the forgetful functor
$D^+(\on{pt}/H)\to D^+(\Vect)$ maps
$\uRHom_{D(\CC)}(X^\bullet_1,X_2^\bullet)$ to
$\RHom_{D(\CC)}(X^\bullet_1,X_2^\bullet)$.
\end{lem}

The last lemma gives rise to the Leray spectral sequence that
expresses Exts in the $H$-equivariant derived category as equivariant
cohomology with coefficients in usual Exts.

\medskip

We will now recall an explicit way of computing Exts in the category
$D(\on{pt}/H)$, in a slightly more general framework. For what follows
we will make the following additional assumption on $H$ (satisfied in
the examples of interest):

\medskip

\noindent {\it We will assume that the group-scheme $H$ is such that 
its unipotent radical $H_u$ is of finite codimension in $H$.
We will fix a splitting $H/H_u=:H_{red}\to H$.}

\medskip

Let $\CC$ be an abelian category with the trivial action of
$H$. We will denote the resulting equivariant DG category by
$\bC(\on{pt}/H\otimes \CC)$. It consists of complexes of objects of
$\CC$, endowed with an algebraic $\CO_H$-action, and an action of
$\fh[1]$, satisfying the usual axioms.

Consider the functor $X^\bullet\mapsto
\uCHom_{\bC(\on{pt}/H)}(\BC,X^\bullet):\bC(\on{pt}/H\otimes \CC)\to
\bC(\CC)$, given by $$X^\bullet\mapsto (X^\bullet)^{H,\fh[1]}.$$
Consider the corresponding derived functor
$$\uCRHom_{D(\on{pt}/H)}(\BC,?):D(\on{pt}/H\otimes \CC)\to D(\CC).$$ Let
us show how to compute it explicitly.

\medskip

Let $BH^\bullet$ (resp., $EH^\bullet$) be the standard simplicial
model for the classifying space of $H$ (resp., the principal
$H$-bundle over it).  Let us denote by $\bDR_{EH^\bullet}$ be the de
Rham complex of $EH^\bullet$.  The action of $H$ on $EH^\bullet$ makes
$\bDR_{EH^\bullet}$ a co-module over $\bDR_H$, i.e., an object of
$\bC(\on{pt}/H)$. Since $EH^\bullet$ is contractible,
$\bDR_{EH^\bullet}$ is quasi-isomorphic to $\BC$.

\begin{lem}
For $X^\bullet\in \bC(\on{pt}/H\otimes \CC)$, there is a natural
quasi-isomorphism
$$\uCHom_{\bC(\on{pt}/H)}(\BC,\bDR_{EH^\bullet}\ast X^\bullet)\simeq
\uCRHom_{D(\on{pt}/H)}(\BC,X^\bullet).$$
\end{lem}

\begin{proof}

We only have to check that whenever $X^\bullet\in \bC(\on{pt}/H\otimes
\CC)$ is acyclic, then $$(\bDR_{EH^\bullet}\ast
X^\bullet)^{H,\fh[1]}$$ is acyclic as well.

Note that the rows of the corresponding bi-complex are isomorphic to
$$(\bDR_{H^n}\ast X^\bullet)^{H,\fh[1]}\simeq \bDR_{H^{n-1}}\otimes
X^\bullet.$$ In particular, they are acyclic if $X^\bullet$ is. In
other words, we have to show that the corresponding spectral sequence
is convergent.

Consider the maps $\bDR_{H^n}\to \bDR_{H_{red}^n}$, corresponding
to the splitting $H_{red}\to H$. They induce a quasi-isomorphism
 $$(\bDR_{EH^\bullet}\ast X^\bullet)^{H,\fh[1]}\to 
(\bDR_{EH_{red}^\bullet}\ast X^\bullet)^{H_{red},\fh_{red}[1]}.$$
This reduces us to the case when $H$ is finite-dimensional, for which
the convergence of the spectral sequence is evident.

\end{proof}

As a corollary, we obtain that the functor
$\uCRHom_{D(\on{pt}/H)}(\BC,?)$ commutes with direct sums. We will
sometimes denote the functor $\uCRHom_{D(\on{pt}/H)}(\BC,?)$ by
$H^\bullet_{\bDR}(\on{pt}/H,?)$.

\ssec{Harish-Chandra modules}   \label{der cat HCh}

Let $\fg$ be a Tate Lie algebra, acted on by $H$ by endomorphisms,
and equipped with a homomorphism $\fh\to \fg$, so that $(\fg,H)$ is
a Harish-Chandra pair. Then the category $\fg\mod$ is a category with
an infinitesimally trivial action of $H$.

The abelian category $\fg\mod^H$ is the same as $(\fg,H)\mod$, i.e., the
category of Harish-Chandra modules. For
$M^\bullet\in \bC(\fg\mod)^H$ we will denote by $x\mapsto a(x)$
the action of $\fg$ on $M^\bullet$ and for $x\in \fh$, by $a^\flat(x)$
the action obtained by deriving the algebraic $H$-action on $M^\bullet$.
(Then, of course, $a^\sharp(x)=a^\flat(x)-a(x)=[d,i^\sharp(x)]$).

Let $D(\fg\mod)^H$ 
be the corresponding derived category, and $D((\fg,H)\mod)$ be the naive
derived category of the abelian category $(\fg,H)\mod$.

\begin{prop}
Assume that $\fg$ is finite-dimensional. Then the evident functor
$D(\fg,H)\mod\to D(\fg\mod)^H$ is an equivalence.
\end{prop}

\begin{proof}

We will construct a functor $\Phi:\bC(\fg\mod)^H\to \bC((\fg,H)\mod)$
that would be the quasi-inverse of the tautological embedding at the
level of derived categories.

For $M^\bullet\in \bC(\fg\mod)$ consider the tensor product
\begin{equation}  \label{naive st}
U(\fg)\otimes \Lambda^\bullet(\fg)\otimes M^\bullet
\end{equation}
with the standard differential, where $U(\fg)$ is the
universal enveloping algebra.

Assume now that $M^\bullet$ is in fact an object of $\bC(\fg\mod)^H$.
Consider an action $i^o$ of $\fh[1]$ on \eqref{naive st}, given by 
$i^o(x)\cdot(u\otimes \omega\otimes m)=u\otimes \omega\wedge x\otimes
m+ u\otimes \omega\otimes i^\sharp(x)\cdot m$. Consider also a
$\fh$-action $\Lie^o$, given by
$$\on{Lie}^o_x\cdot (u\otimes \omega\otimes m)=
-u\cdot x\otimes \omega\otimes m+u\otimes \on{ad}_x(\omega)\otimes m+
u\otimes \omega\otimes a^\flat(x)(m).$$

We have the usual relation $[d,i^o(x)]=\on{Lie}^o_x$, and set
$$\Phi(M^\bullet):=(U(\fg)\otimes \Lambda^\bullet(\fg)\otimes
M^\bullet)_{\fh,\fh[1]}\simeq U(\fg)\underset{U(\fh)}\otimes
(\Lambda^\bullet(\fg) \underset{\Lambda^\bullet(\fh)}\otimes
M^\bullet).$$

\medskip

This is a complex of $\fg$-modules via the $\fg$-action by the
left multiplication on $U(\fg)$, Moreover, we claim that the action of
$\fh\subset \fg$ on $\Phi(M^\bullet)$ integrates to a $H$-action.
This follows from the fact that the $a^\flat$-action of $\fh$ on
$M^\bullet$ is integrable, and that the adjoint of $\fh$ on $\fg$ is
integrable.  Therefore, $\Phi(M^\bullet)$ is an object of
$\bC((\fg,H)\mod)$.

It is easy to see that $\Phi:\bC(\fg\mod)^H\to \bC((\fg,H)\mod)$ is
exact, and hence, it gives rise to a functor at the level of derived
categories.

\medskip

Note that for any $M^\bullet\in \bC(\fg\mod)^H$ we have the natural
maps
$$M^\bullet\leftarrow U(\fg)\otimes \Lambda^\bullet(\fg)\otimes
M^\bullet \to \Phi(M^\bullet),$$
both being quasi-isomorphisms. This implies the statement of the
proposition.

\end{proof}

\ssec{Relative BRST complex}

Assume now that $H$ is such that the adjoint action of $\fh$ on
$\sCliff(\fg,\fg^*)$ can be lifted to an algebraic action of $H$ on
$\sSpin(\fg)$. In particular, the canonical extension $\fg_{-\can}$
splits over $\fh$, and the category $\fg_{-\can}\mod$ also acquires an
infinitesimally trivial $H$-action.

\medskip

For an object $M^\bullet\in \bC(\fg_{-\can}\mod)^H$, consider the
complex $\fC^{\frac{\infty}{2}}(\fg,M^\bullet)$, associated with the
corresponding bi-complex. We claim that it is naturally an object of
$\bC(\on{pt}/H)$:

As a complex of vector spaces, it carries the diagonal action of the
group-scheme $H$ (we will denote the action of its Lie algebra by
$\on{Lie}^\flat$) and an action, denoted $i^\flat$, of
$\Lambda^\bullet(\fh)$ defined as $i|_{\fh}+i^\sharp$. Let us show how
to compute $H^\bullet_{DR}\left(\on{pt}/H,
\fC^{\frac{\infty}{2}}(\fg,M^\bullet)\right)\in D(\Vect)$ (see
\secref{eq vs non-eq Exts}).

For $M^\bullet$ as above, let us denote by
$\fC^{\frac{\infty}{2}}(\fg;H_{red},M^\bullet)$ (resp.,
$\fC^{\frac{\infty}{2}}(\fg;H,M^\bullet)$ the subcomplex of
$\fC^{\frac{\infty}{2}}(\fg,M^\bullet)$, equal to
$\left(\fC^{\frac{\infty}{2}}(\fg,M^\bullet)\right)^{H_{red},\fh_{red}[1]}$
(resp.,
$\left(\fC^{\frac{\infty}{2}}(\fg,M^\bullet)\right)^{H,\fh[1]}$).

\begin{lem}  \label{equiv cohomology of st complex}   \hfill

\smallskip

\noindent{\em (1)} The complex $H^\bullet_{DR}\left(
\on{pt}/H,\fC^{\frac{\infty}{2}}(\fg,M^\bullet)\right)$ is
quasi-isomorphic to $\fC^{\frac{\infty}{2}}(\fg;H_{red},M^\bullet)$.

\smallskip

\noindent{\em (2)} If each $M^i$ as above is injective as a
$H_u$-module, then the embedding
$$\fC^{\frac{\infty}{2}}(\fg;H,M^\bullet)\hookrightarrow
\fC^{\frac{\infty}{2}}(\fg;H_{red},M^\bullet)$$
is a quasi-isomorphism.

\end{lem}

\begin{proof}

First, by \secref{eq vs non-eq Exts}, we can assume that $M^\bullet$
is bounded from below. Secondly, arguing as in \propref{unipotent
averaging}, we can replace the original complex $M^\bullet$ by one,
which consists of modules that are injective over $H_u$ (and hence
over $H$).

Hence, it is sufficient to check that in this case
$$H^\bullet_{DR}\left(\on{pt}/H,
\fC^{\frac{\infty}{2}}(\fg,M^\bullet)\right)\leftarrow
\fC^{\frac{\infty}{2}}(\fg;H,M^\bullet)\to
\fC^{\frac{\infty}{2}}(\fg;H_{red},M^\bullet)$$ are
quasi-isomorphisms.

\medskip

Consider $\fC^{\frac{\infty}{2}}(\fg,M^\bullet)$ as a module over the
Clifford algebra $\sCliff(\fh)$, where the annihilation operators act
by $i^*$, and the creation operators act by means of $i^\flat$. We
obtain that
$$\fC^{\frac{\infty}{2}}(\fg,M^\bullet)\simeq \fC(\fh,M_1^\bullet),$$
for some complex $M_1^\bullet$ of $H$-modules, which are moreover,
consists of injective objects.  

\medskip

Thus, we have reduced the original problem to the case when
$\fg=\fh$. In this case, by \lemref{homotopy lemma} and \lemref{Leray},
$$H^\bullet_{DR}\left(\on{pt}/H,\fC(\fh,M^\bullet)\right)\simeq 
\RHom_{D(H-mod)}(\BC,M^\bullet),$$
which is quasi-isomorphic to 
$$\fC(\fh;H,M^\bullet)\simeq M^\bullet,$$
if $M^\bullet$ consists of injective $H$-modules.

Moreover, by the Hochshild-Serre spectral sequence,
$$\RHom_{D(H-mod)}(\BC,M^\bullet)\simeq 
\left(\RHom_{D(H_u-mod)}(\BC,M^\bullet)\right)^{H_{red}}\simeq
\fC(\fh;H_{red},M^\bullet).$$

\end{proof}

\ssec{Variant: equivariance against a character}   \label{N, psi equiv}

Let now $\psi$ be a homomorphism $H\to \BG_a$; we will denote by the
same character the resulting character on $\fh$. For a category $\CC$
as above, we introduce the category $\CC^{H,\psi}$ to be the full
subcategory of $\CC^{w,H}$, consisting of objects, for which the map
$a^\flat$ is given by the character $\psi$.

\medskip

Let us consider the example when $\CC=\fD(\CY)\mod$. Let $\be^\psi$ be
the pull-back of the Artin-Schreier D-module on $\BG_a$ under
$\psi$. Its fiber at $1\in H$ is trivialized and it is a character
sheaf in the sense that we have a canonical isomorphism
$\on{mult}^*(\be^\psi)\simeq \be^\psi\boxtimes \be^\psi$, which is
associative in the natural sense.

The category $\fD(\CY)\mod^{H,\psi}$ consists of D-modules $\CF$ on
$\CY$, endowed with an isomorphism $\act^*(\CF)\simeq
\be^\psi\boxtimes \CF\in \fD(H\times \CY)\mod$, satisfying the
associativity and unit conditions.

\medskip

Returning to the general situation, we introduce the category
$\bC(\CC)^{H,\psi}$ in the same way as $\bC(\CC)^H$ with the only
difference that we require that $[d,i^\sharp]=a^\sharp+\psi$. This is
DG-category with a cohomological functor to $\CC^{H,\psi}$. We will
denote by $D(\CC)^{H,\psi}$ the resulting triangulated category.

Much of the discussion about $\bC(\CC)^H$ carries over to this
situation.  For example, we have the averaging functor
$\Av_{H,\psi}:\bC(\CC)\to \bC(\CC)^{H,\psi}$, right adjoint to the
forgetful functor.  It is constructed as the composition of $\Av_H^w$
and the functor
$$X^\bullet\mapsto \fC(\fh,X^\bullet\otimes \BC^\psi):
\bC(\CC^{w,H})\to \bC(\CC)^{H,\psi},$$ where $\BC^\psi$ is the
$1$-dimensional representation of $\fh$ corresponding to the character
$\psi$.

When $H$ is pro-unipotent, one shows in the same way as above that the
functor $D^+(\CC)^{H,\psi}\to D^+(\CC)$ is an equivalence onto the
full subcategory, consisting of objects, whose cohomologies belong to
$\CC^{H,\psi}$.

\section{D-modules on group ind-schemes}   \label{D-mod on group
  ind-schemes}

\ssec{}

Let $G$ be an affine reasonable group ind-scheme, as in \cite{BD}. In
particular, its Lie algebra $\fg$ is a Tate vector space. We will
denote by $\CO_G$ the topological commutative algebra of functions on
$G$.

The multiplication on $G$ defines a map $\Delta_G:\CO_G\to
\CO_G\shriektimes \CO_G$. We will denote by $\Lie_l$ and $\Lie_r$ the
two maps
$$\fg\startimes \CO_G\to \CO_G,$$ corresponding to the action of $G$
on itself by left (resp., right) translations.

In addition, we have the maps
$$\Delta_\fg:\fg\to \CO_G\shriektimes\fg,\, \Delta_{\fg^*}:\fg^*\to
\CO_G\shriektimes\fg^*$$ that correspond to the adjoint and co-adjoint
actions of $G$ on $\fg$ and $\fg^*$, respectively.

\medskip

Let us denote by $T(G)$ (resp., $T^*(G)$) the topological $\CO_G$-module
of vector fields (resp., $1$-forms) on $G$. It identifies in two ways
with $\CO_G\shriektimes\fg$ (resp., $\CO_G\shriektimes\fg^*$),
corresponding to the realization of $\fg$ (resp., $\fg^*$) as right or
left invariant vector fields (resp., $1$-forms).  Note that $T(G)$ is
a topological Lie algebra and $T^*(G)$ is a module over it.

\ssec{}

Following \cite{AG}, we introduce the category of D-modules on
$G$, denoted $\fD(G)\mod$, as follows: 

Its objects are (discrete) vector spaces $\CM$, endowed with an action
$$\CO\arrowtimes \CM\simeq \CO\startimes \CM\overset{m}\to \CM$$ and a
Lie algebra action
$$a_l:\fg\arrowtimes \CM \simeq \fg\startimes \CM\to \CM,$$
such that the two pieces of data are compatible in the sense of the action
of $\fg$ on $\CO_G$ by {\it left} translations in the following sense:

We need that the difference of the two arrows:
$$\fg\startimes \CO_G\startimes \CM\overset{\id_G\otimes
m}\longrightarrow \fg\startimes \CM \overset{a_l}\to \CM$$ and
$$\fg\startimes \CO_G\startimes \CM\simeq \CO_G\startimes
\fg\startimes \CM \overset{\on{id}_{\CO_G}\otimes a_l}\longrightarrow
\CO_G\startimes \CM\overset{m}\to \CM$$ equals
$$\fg\startimes \CO_G\startimes \CM\overset{\Lie_l}\longrightarrow
\CO_G\startimes \CM\overset{m}\to \CM.$$

Morphisms in $\fD(G)\mod$ are maps of vector spaces $\CM_1\to \CM_2$ that 
commute with the actions of $\fg$ and $\CO_G$.

\ssec{Action of the Tate canonical extension}

Following Beilinson, we will show now that if $\CM$ is an object of
$\fD(G)\mod$, then the underlying vector space carries a canonical
action of $\fg_{-\can}$, denoted $a_r$, which commutes with the
original action of $\fg$, and which is compatible with the action of
$\CO_G$ via the action of $\fg$ on $\CO_G$ by {\it right}
translations.

\medskip

Set $\CM^{DR}=\CM\otimes \sSpin(\fg)$. Let us denote by $i_r$ and $i^*_r$
the actions on it of $\Lambda^\bullet(\fg)$ and $\Lambda^\bullet(\fg^*)$,
both of which are subalgebras in $\sCliff(\fg,\fg^*)$.

{}From the definition of $\fD(G)\mod$ it follows that $i_r$ and $i^*_r$
on $\CM$ extend to actions of the odd topological vector spaces $T(G)$
and $T^*(G)$, identified with $\CO_G\shriektimes \fg$ and
$\CO_G\shriektimes \fg^*$ using left-invariant vector fields and
forms, respectively. We will denote the resulting actions simply by
$i$ and $i^*$.

Using the map 
\begin{equation}   \label{left through right}
\fg\overset{-\Delta_{\fg}}\longrightarrow \CO_G\shriektimes\fg
\overset{\gamma\otimes \on{id}_\fg}\longrightarrow
\CO_G\shriektimes\fg,
\end{equation}
(here $\gamma$ is the inversion on $G$), we obtain a new action $i_l$
of $\Lambda^\bullet(\fg)$ on $\CM^{DR}$. Similarly, we have a new
action $i^*_l$ of $\Lambda^\bullet(\fg^*)$ on $\CM^{DR}$. Altogether,
we obtain a new action of the Clifford algebra $\sCliff(\fg,\fg^*)$ on
$\CM^{DR}$.

\medskip

We will denote by the symbol $\on{Lie}_l$ the action of $\fg$ on
$\CM^{DR}$ coming from the action $a_l$ of $\fg$ on $\CM$. We claim
that this action extends to an action of the Lie algebra $T(G)$,
denoted simply by $\Lie$.

First we define an action of the non-completed tensor product
$\CO_G\otimes \fg$ on $\CM^{DR}$. Namely, for $x\in \fg$, $f\in \CO_G$
and $v\in \CM$ we set we set
\begin{equation}  \label{anomaly cancellation}
(f\otimes x)\cdot v=f\cdot \on{Lie}_{l}(x)\cdot v+i_l^*(df)\cdot
  i_l(x)\cdot v.
\end{equation}
Note that
$$f\cdot \on{Lie}_{l}(x)\cdot v+i_l^*(df)\cdot i_l(x)\cdot v=
\on{Lie}_{l}(x)\cdot f\cdot v+i_l(x)\cdot i_l^*(df)\cdot v.$$ This
property implies that the above action of $\CO_G\otimes \fg$ extends
to the action of $\CO_G\shriektimes\fg\simeq T(G)$.  Indeed, when $x$
is contained in a deep enough neighborhood of zero, then both
$\on{Lie}_{l(x)}$ and $i_l(x)$ annihilate any given $v\in
\CM$. Similarly, if $f$ is contained in a deep neighborhood of zero,
then $v$ is annihilated by both $f$ and $i_l^*(df)$.

\medskip

One readily checks that the above action is compatible with the Lie
algebra structure on $T(G)$. In particular, using the map
$-\Delta_{\fg}:\fg\to \CO_G\shriektimes \fg$, i.e., the embedding of
$\fg$ into $T(G)$ as left-invariant vector fields, we obtain a new
action of the Lie algebra $\fg$ on $\CM$. We will denote this action
by $\on{Lie}_r$.

We have

\begin{itemize}

\item
$[\Lie_r(x),i_l(y)]=0$, $[\Lie_r(x),i^*_l(y^*)]=0$ for $x,y\in \fg$,
$y^*\in \fg^*$.

\item
$[\Lie_r(x),i_r(y)]=i_r([x,y])$,
$[\Lie_r(x),i^*_r(y^*)]=i^*_r(\on{ad}_x(y^*))$,

\item
$[\Lie_r(x),f]=\Lie_{r(x)}(f)$ for $f\in \CO_G$.

\item
$[\Lie_r(x),a_l(y)]=0$.

\end{itemize}

\medskip

Finally, we are ready to define the action $a_r$ of $\fg_{-\can}$ on
$\CM^{DR}$.  Namely, $a_r$ is the difference of $\on{Lie}_r$ and the
canonical $\fg_{\can}$-action on $\sSpin(\fg)$.

It is easy to see that $a_r$ is indeed an action. Moreover,

\begin{itemize}

\item
$[a_r(x'),f]=\Lie_{r(x)}(f)$, for $x'\in \fg_{-\can}$ and its image
$x\in \fg$,

\item
$[a_r(x'),a_l(y)]=0$,

\item
$[a_r(x'),i_r(y)]=0$, $[a_r(x'),i_r(y^*)]=0$.

\end{itemize}

The last property implies that the $a_r$-action of $\fg_{-\can}$ on
$\CM^{DR}$ preserves the subspace $\CM$; i.e., we obtain an action of
$\fg_{-\can}$ on $\CM$ that satisfies the desired commutation
properties.

\medskip

When we view $\CM\in \fD(G)\mod$ as a $\fg_{-\can}$-module via $a_r$,
we obtain that $\CM^{DR}$ identifies with
$\fC^{\frac{\infty}{2}}(\fg,\CM)$, where $i=i_r$, $i^*=i^*_r$,
$\Lie=\Lie_r$. In particular, $\CM^{DR}$ acquires a natural
differential $d$.

{}From the above commutation properties, it follows that this
differential satisfies:

\begin{itemize}

\item
$[d,i(\xi)]=\Lie(\xi)$ for $\xi\in T(G)$,

\item
$[d,f]=i^*(df)$ for $f\in \CO_G$.

\end{itemize}

Of course, $\CM^{DR}$ depends on the choice of the Clifford module
$\sSpin(\fg)$.

\ssec{}

Note that the above construction can be inverted: we can introduce the
category $\fD(G)\mod$ to consist of $(\CO_G,\fg_{-\can})$-modules,
where the two actions are compatible in the sense of the
$\fg_{-\can}$-action on $\CO_G$ via right translations. In this case,
the vector space underlying a representation automatically acquires an
action of $\fg$, which commutes with the $\fg_{-\can}$-action and is
compatible with the action of $\CO_G$ via left translations.

\medskip

Let us also note that in the definition of $\fD(G)\mod$ we could
interchange the roles of left and right:

Let us call the category introduced above $\fD(G)_l\mod$, and let us
define the category $\fD(G)_r\mod$ to consist of
$(\CO_G,\fg)$-modules, where the two actions are compatible via the
action of $\fg$ on $\CO_G$ by right translations. We claim that the
categories $\fD(G)_l\mod$ and $\fD(G)_r\mod$ are equivalent, but this
equivalence does not respect the forgetful functor to vector spaces.

This equivalence is defined as follows. For $\CM_l\in \fD(G)_l\mod$,
the actions $i_l$, $i^*_l$ define a new action of $\sCliff(\fg,\fg^*)$
on $\CM^{DR}$. We define an object of $\CM_r\in \fD(G)_r\mod$ by
$\Hom_{\sCliff(\fg,\fg^*)}(\sSpin(\fg),\CM^{DR})$ with respect to this
new action.

Explicitly, this can be reformulated as follows. Let $G_{\can}$ be the
canonical (i.e., Tate) central extension of $G$. It can be viewed as a
line bundle $\CP_{\can}$ over $G$, whose fiber at a given point
$\bg\in G$ is the relative determinant line
$\det(\fg,\on{Ad}_\bg(\fg))$. The action of $\fg$ on $G$ be left
(resp., right) translations extends to an action of $\fg_{\can}$
(resp., $\fg_{-\can}$) on $\CP_{\can}$.

Then,
$$\CM_r\simeq \CM_l\underset{\CO_G}\otimes \CP^{-1}_{\can},$$
as $\CO_G$-modules, respecting both the $a_l$ and $a_r$ actions.

\medskip

In what follows, unless stated otherwise, we will think of
$\fD(G)\mod$ in the $\fD(G)_l\mod$ realization.

\ssec{} 

Let $H$ be a group-scheme, mapping to $G$. We claim that the category
$\fD(G)\mod$ carries a natural infinitesimally trivial action of $H$,
corresponding to the action of $H$ on $G$ by left translations.  (As
we shall see later, this is a part of a more general structure, the
latter being an action of the group ind-scheme $G\times G$ on
$\fD(G)\mod$).

\medskip

For $\CM\in \fD(G)\mod$ we set $\act_l^*(\CM)$ to be isomorphic to
$\CO_H\otimes \CM$ as an $\CO_H$-module. The action of $\CO_G$ is
given via
$$\CO_G\overset{\Delta_G}\to \CO_G\shriektimes \CO_G
\overset{\gamma\otimes \on{id}} \longrightarrow\CO_H\shriektimes
\CO_G.$$ The action $a_l$ of $\fg$ is given via the map
$$\fg\overset{\Delta_\fg}\to \CO_G\shriektimes \fg
\overset{\gamma\otimes \on{id}} \longrightarrow \CO_G\shriektimes \fg
\to \CO_H\shriektimes \fg,$$ where $\CO_H\shriektimes \fg$ acts on
$\CO_H\otimes \CM$ by $\on{id}\otimes m$.

To construct isomorphism $\act^*(\CM)|_{H^{(1)}}\simeq
p^*(\CM)|_{H^{(1)}}$ we identify both sides with $\CM\oplus
\epsilon\cdot \fh^*\otimes \CM$ as vector spaces, and the required
isomorphism is given by the action of $\fh$ on $\CM$, obtained by
restriction from $a_l$.

Note that by construction, the action of $\fg_{-\can}$ on
$\act_l^*(\CM)\simeq \CO_H\otimes \CM$ is via its action on the second
multiple.

\medskip

Let now $H'$ be another group-scheme, mapping to $G$, and let us
assume that there exists a splitting $\fh'\to \fg_{\can}$. In this
case, we claim that there exists another infinitesimally trivial
action of $H'$ on $\fD(G)\mod$, corresponding to the action of $H$ on
$G$ by right translations:

\medskip

For $\CM\in \fD(G)\mod$, we define $\act_r^*(\CM)$ to be isomorphic to
$\CM\otimes \CO_{H'}$ as an $\CO_{H'}$-module and as a
$\fg$-module. The action of $\CO_G$ is given by the co-multiplication
map $\CO_G\to \CO_G\shriektimes \CO_{H'}$. It is easy to see that the
commutation relation is satisfied. The associativity and unit
constraint are evident.

To construct the isomorphism $\act_r^*(\CM)|_{H'{}^{(1)}}\simeq
p^*(\CM)|_{H'{}^{(1)}}$, note that both sides are identified with
$\CM\oplus \epsilon\cdot \fh'{}^*\otimes \CM$ as $\fg$-modules.  The
required isomorphism is given by the action of $\fh'$ on $\CM$,
obtained by restriction from $a_r$. Again, it is easy to see that the
axioms of Harish-Chandra action hold.

Let us note that the action $a_r$ of $\fg_{-\can}$ on
$\act_r^*(\CM)\simeq \CM\otimes \CO_{H'}$ is given via the map
$$\fg_{-\can}\overset{\Delta_\fg}\longrightarrow  
\fg_{-\can}\shriektimes \CO_{G}\to \fg_{-\can}\shriektimes \CO_{H'}.$$

\medskip

Let us denote by $\fD(G)\mod^{l(H)}$ (resp., $\fD(G)\mod^{r(H')}$) the
corresponding categories of strongly equivariant objects of
$\fD(G)\mod$. Moreover, it is easy to see that the actions of $H$ and
$H'$ commute in the natural sense, i.e., we have an action of $H\times
H'$ on $\fD(G)\mod$. We will denote the resulting category by
$\fD(G)\mod^{l(H),r(H')}$.

\ssec{}    \label{D-mod on groups and quotients}

Let now $K\subset G$ be a group-subscheme such that the quotient $G/K$
exists as a strict ind-scheme of ind-finite type (in this case it is
formally smooth).  We will call such $K$ "open compact".

We will choose a particular model for the module $\sSpin(\fg)$,
denoted $\sSpin(\fg,\sk)$ by letting it be generated by a vector
$\one\in \sSpin(\fg)$, annihilated by $\sk\oplus (\fg/\sk)^*\subset\fg
\oplus \fg^*\subset \sCliff(\fg,\fg^*)$. This $\sSpin(\fg,\sk)$
carries a natural action of $K$, which gives rise to a splitting of
$G_{\can}$ over $K$.

\medskip

By the assumption on $G/K$, it makes sense to consider right D-modules
on it; we will denote this category by $\fD(G/K)\mod$.

\begin{prop}  \label{AG equivalence}
We have a canonical equivalence $\fD(G)\mod^{r(K)}\simeq
\fD(G/K)\mod$.
\end{prop}

\begin{proof}
Let $\pi$ denote the projection $G\to G/K$. For $\CF\in \fD(G/K)\mod$,
consider the $\CO_G$-module $\CM:=\Gamma(G,\pi^*(\CF))$.

For $x\in \fg$, the (negative of the) corresponding vector field
acting on $\CF$ gives rise to a map $a_l(x):\CM\to \CM$, as a vector
space, and these data satisfy the conditions for $\CM$ to be a
$\fD(G)$-module.

We claim that the action of the Lie algebra $\sk\subset \fg_{-\can}$
on $\CM$, given by $a_r$, coincides with the natural action of $\sk$
on $\pi^*(\CF)$ obtained by deriving the group action. This would
imply that $\CM$ is naturally an object of $\fD(G)\mod^{r(K)}$.

To prove the assertion, we can assume that $\CF$ is an extension of a
D-module on an affine ind-subscheme of $G/K$. Then it is sufficient to
check that the subspace $\Gamma(G/K,\CF)\subset \CM\subset \CM^{DR}$
is annihilated by the operators $\Lie_{r(x)}$ for $x\in \sk$. But this
is straightforward from the construction.

\medskip

Vice versa, let $\CM$ be an object of $\fD(G)\mod^{r(K)}$, which we
identify with the corresponding quasi-coherent sheaf on $G$. Consider
the complex of sheaves $\pi_*(\CM^{DR})$ on $G/K$; it carries an
action of the operators $i_r(x)$, $\on{Lie}_{r(x)}$, $x\in \fg$.  We
set $\CF^{DR}$ to be the subcomplex of $\pi_*(\CM^{DR})$ annihilated
by the above operators for $x\in \sk$.

Set $\CF$ to be the degree $0$ part of $\CF^{DR}$; it is easy to see
that $\CF\simeq \left(\pi_*(\CM)\right)^K$. The degree $-1$ part of
$\CF^{DR}$ identifies with $\CF\underset{\CO(G/K)}\otimes T(G/K)$, and
the differential $$d:(\CF^{DR})^{-1}\to (\CF^{DR})^0$$ defines on
$\CF$ a structure of a right $\fD$-module. Moreover, the entire
complex $\CF^{DR}$ identifies with the de Rham complex of $\CF$.

\end{proof}

Let $\delta_{K,G}$ be the object of $\fD(G)\mod^{r(K)}$ corresponding
to the delta-function $\delta_{1,G/K}$ under the equivalence of
categories of \propref{AG equivalence}. It can be constructed as
$\on{Ind}_{\sk}^{\fg}(\CO_K)$ as a module over $\fg$ and $\CO_G$. As a
module over $\fg_{-\can}$ it is also isomorphic to
$\on{Ind}_{\sk}^{\fg_{-\can}}(\CO_K)$.

\medskip

More generally, let $\CL$ be an object of $\QCoh^!_{G/K}$.  Let
$\Ind_{\CO_{G/K}}^{\fD_{G/K}}(\CL)$ be the induced D-module.
The corresponding object of $\fD(G)\mod^{r(K)}$, i.e.,,
$\Gamma\left(G,\pi^*\left(\Ind_{\CO_{G/K}}^{\fD_{G/K}}(\CL)\right)\right)$
can described as follows:

Consider the $\CO_G$-module $\Gamma(G,\pi^*(\CL))$; it is acted on
naturally by $K$. Consider the $\fg_{-\can}$-module
$\Ind_{\sk}^{\fg_{-\can}}\left(\Gamma(G,\pi^*(\CL))\right)$.  It is
naturally acted on by $\CO_G$, so that the actions of $\fg_{-\can}$
and $\CO_G$ satisfy the commutation relation with respect to the right
action of $G$ on itself. Hence,
$\Ind_{\sk}^{\fg_{-\can}}\left(\Gamma(G,\pi^*(\CL))\right)$ is an
object of $\fD(G)\mod^{r(K)}$ and we have a natural isomorphism:

\begin{equation} \label{ind D-modules}
\Gamma\left(G,\pi^*\left(\Ind_{\CO_{G/K}}^{\fD_{G/K}}(\CL)
\right)\right)\simeq
\Ind_{\sk}^{\fg_{-\can}}\left(\Gamma(G,\pi^*(\CL))\right).
\end{equation}

\ssec{The bi-equivariant situation} \label{bi-equivariant}

Let now $K_1,K_2$ be two "open compact" subgroups of $G$. Note that we
have a natural equivalence of categories
\begin{equation} \label{opposition}
\fD(G/K_1)\mod^{K_2}\to \fD(G/K_2)\mod^{K_1}:\CF\mapsto \CF^{\on{op}},
\end{equation}
defined as follows.

Assume, without loss of generality, that $\CF$ is supported on a
closed $K_2$-invariant subscheme $\CY\in G/K_1$, and let $\CY^{\on{op}}$ be
the corresponding $K_1$-invariant subscheme in $K_2\backslash G$. We
can find a normal "open compact" subgroup $K_2'\subset K_2$ such that if
we denote by $\CY'{}^{\on{op}}$ the preimage of $\CY^{\on{op}}$ in
$K'_2\backslash G$, the projection
$$\pi_1:\CY'{}^{\on{op}}\to \CY$$ is well-defined (and makes $\CY'{}^{\on{op}}$
a torsor with respect to the corresponding smooth group-scheme over
$\CY$).

Consider $\pi_1^!(\CF)$. This is a D-module on $\CY'{}^{\on{op}}$,
equivariant with respect to the action of the group $(K_2/K'_2)\times
K_1$ on this scheme. Hence, it gives rise to a $K_1$-equivariant
D-module on $K_2\backslash G$.  To obtain $\CF^{\on{op}}$ we apply the
involution $\bg\mapsto \bg^{-1}:K_2\backslash G\to G/K_2$.

\medskip

Let us now describe what this equivalence looks like in terms of the
equivalences
$$\fD(G/K_1)\mod^{K_2}\simeq \fD(G)\mod^{l(K_2),r(K_1)} \text{ and }
\fD(G/K_2)\mod^{K_1}\simeq \fD(G)\mod^{l(K_1),r(K_2)}.$$

First, the inversion on $G$ defines an equivalence
$$\fD(G)_l\mod^{l(K_2),r(K_1)}\simeq \fD(G)_r\mod^{l(K_2),r(K_1)},$$
and the sought-for equivalence is obtained from the one above
via $\fD(G)_r\mod^{l(K_2),r(K_1)}\simeq \fD(G)_l\mod^{l(K_2),r(K_1)}$.

(Note that the determinant line that played a role in the
$\fD(G)_r\mod\simeq \fD(G)_l\mod$ equivalence, appears also in
\eqref{opposition}, when we descend {\it right} D-modules from
$\CY'{}^{\on{op}}$ to $\CY^{\on{op}}$.)

\ssec{}

Consider now the DG-category $\bC(\fD(G)\mod)^{r(K)}$, and we claim
that the construction in \propref{AG equivalence} generalizes to a
functor $\bC(\fD(G)\mod)^{r(K)}\to \bC(\fD(G/K)\mod)$.

Indeed, for $\CM^\bullet\in \bC(\fD(G)\mod)^{r(K)}$ let us denote by
$(\CM^\bullet)^{DR}$ the total complex of the corresponding bicomplex.
We have several actions of $\Lambda^\bullet(\sk)$ on it, and let
$i^\flat_r$ be the sum of the one given by restricting the
$i_r$-action of $\fg$ and $i^\sharp$.  In addition,
$(\CM^\bullet)^{DR}$ carries a natural action of $K$. These
two structures combine to that of object of $\bC(\on{pt}/K)$.

Consider again the complex of sheaves $\pi_*((\CM^\bullet)^{DR})$ on
$G/K$.  It carries an action of $\Lambda^\bullet(\sk)$ coming from
$i^\flat_r$ and an action of the group-scheme $K$. Define
$$(\CF^\bullet)^{DR}:=\CHom_{\bC(\on{pt}/K)}
\left(\BC,\pi_*((\CM^\bullet)^{DR}) \right).$$ This is an
$\Omega^\bullet_{G/K}$-module on in the terminology of \cite{BD}.

Finally, we consider the functor $\bC(\fD(G)\mod)^{r(K)}\to
\bC(\fD(G/K)\mod)$ given by
$$\CM^\bullet\mapsto
\on{Ind}^{\fD(G/K)}_{\CO_{G/K}}\Bigl((\CF^\bullet)^{DR}\Bigr)\in
\bC(\fD(G/K)),$$ where $\on{Ind}^{\fD(G/K)}_{\CO_{G/K}}$ is the
induction functor from $\Omega^\bullet_{G/K}$-modules to D-modules on $G/K$,
see \cite{BD}, Sect. 7.11.12.

\begin{lem}
The resulting functor
$$\CM^\bullet\mapsto
\on{Ind}^{\fD(G/K)}_{\CO_{G/K}}\Bigl((\CF^\bullet)^{DR}\Bigr):
\bC(\fD(G)\mod)^{r(K)}\to \bC(\fD(G/K)\mod)$$ is exact.
\end{lem}

\begin{proof}
This follows from the fact that the functor
$$\CM^\bullet\mapsto (\CF^\bullet)^{DR},$$
viewed as a functor from $\bC(\fD(G)\mod)^{r(K)}$ to the DG category of
(non-quasi coherent) sheaves on $G/K$ is exact.
\end{proof}

Hence, we obtain a well-defined functor $D(\fD(G)\mod)^{r(K)}\to
D(\fD(G/K)\mod)$.

\begin{prop}   \label{equivariant and quotient}
The above functor $D(\fD(G)\mod)^{r(K)}\to D(\fD(G/K)\mod)$ is an
equivalence.  Its quasi-inverse is given by
$$D(\fD(G/K)\mod)\to D(\fD(G)\mod^{r(K)})\to D(\fD(G)\mod)^{r(K)}.$$ 
\end{prop}

As a corollary, we obtain that in this case the evident functor
$D(\fD(G)\mod^{r(K)})\to D(\fD(G/K)\mod)$ is an equivalence.

\begin{proof}

The functor
$$\CF^\bullet\mapsto \Gamma(G,\pi^*(\CF^\bullet))=:\CM^\bullet\mapsto
\on{Ind}^{\fD(G/K)}_{\CO_{G/K}}\Bigl(
\CHom_{\bC(\on{pt}/K)}\left(\BC,\pi_*((\CM^\bullet)^{DR})\right)
\Bigr)$$ is isomorphic to the composition
$$\bC(\fD(G/K)\mod) \overset{DR}\longrightarrow
\Omega^\bullet_{G/K}\mod \overset{\on{Ind}^{\fD(G/K)}_{\CO_{G/K}}}
\longrightarrow\bC(\fD(G/K)\mod),$$ and hence, on the derived level,
it induces a functor isomorphic to the identity.

\medskip

Vice versa, for $\CM^\bullet\in \bC(\fD(G)\mod^{r(K)})$ we have a
natural map
\begin{equation} \label{map from ind}
\pi^*\Bigl(\on{Ind}^{\fD(G/K)}_{\CO_{G/K}}\Bigl(\CHom_{\bC(\on{pt}/K)}
\left(\BC,\pi_*((\CM^\bullet)^{DR})\right) \Bigr)\to \CM^\bullet,
\end{equation}
and we claim that it is a quasi-isomorphism. This follows from the
fact that as a complex of vector spaces, the LHS of \eqref{map from
ind} is naturally filtered, and the associated graded is isomorphic to
$$\left(\on{Sym}(\fg/\sk)\otimes
\Lambda^\bullet(\fg/\sk)\right)\otimes \CM^\bullet,$$ where the first
multiple has the Koszul differential.

\end{proof}

We shall now establish the following:

\begin{prop}   \label{de Rham and semiinf}
For $\CM^\bullet\in D^+(\fD(G)\mod)^K$ and $\CF^\bullet\in
D^+(\fD(G/K)\mod)$, corresponding to each other under the equivalence
of \propref{equivariant and quotient}, we have a canonical
quasi-isomorphism
$$H^\bullet_{DR}(G/K,\CF^\bullet)\simeq
\fC^{\frac{\infty}{2}}(\fg;K_{red},\CM^\bullet).$$
\end{prop}

Note that by \lemref{equiv cohomology of st complex}(1),
$$\fC^{\frac{\infty}{2}}(\fg;K_{red},\CM^\bullet)\simeq
H^\bullet_{DR}\left(\on{pt}/K,(\CM^\bullet)^{DR}\right).$$

\begin{proof}

We can assume that the complex $\CM^\bullet$ is such that 
each $\CM^i$, as a $K$-equivariant $\CO_G$-module, is of the form
$\pi^*(\CL^i)$, where $\CL^i$ is a quasi-coherent
sheaf on $G/K$, which is the direct image from an affine subscheme.
Such $\CL$ is obviously loose in the sense of \cite{BD}, i.e.,
it has the property that the higher cohomologies
$H^i(G/K,\CL\otimes \CL^1)$ vanish for any quasi-coherent sheaf 
$\CL^1$ on $G/K$.

Hence, the de Rham cohomology of $\CF^\bullet$
can be computed as $\Gamma(G/K,(\CF^\bullet)^{DR})$. Note that
the latter complex can be identified by definition with
$$\fC^{\frac{\infty}{2}}(\fg;K,\CM^\bullet).$$
Hence, the assertion of the proposition follows from
\lemref{equiv cohomology of st complex}(2).

\end{proof}

\ssec{Variant: central extensions and twisting}

Let now $\fg'$ be a central extension of $\fg$ by means of $\BC$. We
will denote by $\fg'_{-\can}$ the Baer sum of $\fg_{-\can}$ and the
Baer negative of $\fg'$.

We introduce the category $\fD(G)'\mod$ to consist of (discrete)
vector spaces $\CM$, endowed with an action
$$\CO_G\startimes \CM\overset{m}\to \CM$$ as before,
and a Lie algebra action
$$a_l:\fg'\startimes \CM\to \CM,$$ (such that, of course, $1\in
\BC\subset \fg'$ acts as identity), and such that the two pieces of
data are compatible in same way as in the definition of $\fD(G)\mod$.

\medskip

We claim that in this case, the vector space, underlying an object
$\CM\in \fD(G)'\mod$ carries a canonically defined action, denoted
$a_r$, of $\fg'_{-\can}$, which commutes with $a_l$, and which
satisfies $[a_r(x'),f]=\Lie_{r(x)}(f)$ for $x'\in \fg'_{-\can}$ and
$f\in \CO_G$.

We construct $a_r$ by the same method as in the case of $\fg'=\fg$.
Namely, we tensor $\CM$ by $\sSpin(\fg)$, and show that it carries an
action of $T(G)':=\CO_G\shriektimes \fg'$, from which we produce the
desired $a_r$.

Note, however, that in this case $\CM\otimes \sSpin(\fg)$ does not
carry any differential.

\medskip

Let us assume now that $\fg'$ is a scalar multiple of an extension,
induced by some central extension of $G$ by means of $\BG_m$. Let $K$
be an "open compact" subgroup of $K$, and assume that $G'$ splits over
$K$. We can then consider the category $\fD(G/K)'\mod$ of twisted
D-modules on $G/K$.

In this case we also have a well-defined category $\fD(G)'\mod^K$
(along with its DG and triangulated versions $\bC(\fD(G)'\mod)^K$ and
$D(\fD(G)'\mod)^K$). \propref{AG equivalence} and \propref{equivariant
and quotient} generalize to the twisted context in a straightforward
way.

\ssec{D-modules with coefficients in a category}

Let $\CC$ be an abelian category, satisfying assumption (*) of
\secref{class of cat}. Then it makes sense to consider the category
$\fD(G)\mod\otimes \CC$, and all the results of the present section
carry over to this context.

In particular, for an "open compact" subgroup $K\subset G$ we can
consider the category $\fD(G/K)\mod\otimes \CC$ (see \secref{objects
over indscheme}), and we have the analogs of \propref{AG equivalence}
and \propref{equivariant and quotient}.

\section{Convolution}     \label{conv section}

\ssec{Action of group ind-schemes on categories}

We will now generalize the contents of \secref{weak action of groups}
and \secref{HCh action of groups} into the context of group
ind-schemes. Let $G$ be an affine reasonable group ind-scheme as
above.  Let $\CC$ be a category satisfying assumption (*) of
\secref{class of cat}.

\medskip

A weak action of $G$ on a $\CC$ is the data of a functor
$$\on{act}^*:\CC\to \on{QCoh}_G^*\otimes \CC,$$
and two functorial isomorphisms as in \secref{weak action of groups}.

For $X\in\CC$ and a scheme $S$ mapping to $G$ we obtain a functor
$$X\mapsto \act^*(X)|_S:\CC\to \on{QCoh}_S\otimes \CC.$$
The following assertion is proved as \lemref{act is flat}:

\begin{lem}
For any $S\to G$, the functor $X\mapsto \act^*(X)|_S$ is exact,
and its image consists of $\CO_S$-flat objects.
\end{lem}

\medskip

Let $G^{(1)}$ be the first infinitesimal neighborhood of $1\in G$. This
is a formal scheme equal to $\on{Spf}(\BC\oplus \epsilon \cdot \fg^*)$.
If $G$ acts on $\CC$ and $X\in \CC$, we obtain an object
$$X^{(1)}:=\on{act}(X)|_{G^{(1)}}\in \on{QCoh}^*_{G^{(1)}}\otimes
\CC.$$

We say that the action of $G$ on $\CC$ is of Harish-Chandra type if we
are given a functorial identification between $X^{(1)}$ and
$p^*(X)|_{G^{(1)}}$, satisfying the same compatibility conditions as
in \secref{HCh action of groups}.

\medskip

Let now $\fg'$ be a central extension of $\fg$ by means of $\BC$. Let
$G'{}^{(1)}$ be the formal scheme $\on{Spf}(\BC\oplus \epsilon \cdot
\fg'{}^*)$.  It projects onto $G^{(1)}$ and contains $\Spec(\BC\oplus
\epsilon \cdot \BC)$ as a closed subscheme.

We say that a $G$ action on $\CC$ is of twisted Harish-Chandra type
relative to $\fg'$, if for every $X\in \CC$ we have a functorial
isomorphism between $\act^*(X)|_{G'{}^{(1)}}$ and
$p^*(X)|_{G'{}^{(1)}}$ such that the induced map
$$\act^*(X)|_{G'{}^{(1)}}|_{\Spec(\BC\oplus \epsilon \BC)}\simeq
X\oplus \epsilon\cdot X \simeq p^*(X)|_{G'{}^{(1)}}|_{\Spec(\BC\oplus
\epsilon \BC)},$$ is the automorphism
$$\on{id}_X\oplus \epsilon\cdot \on{id}_X: X\oplus \epsilon\cdot X\to
X\oplus \epsilon\cdot X,$$
and which satisfies the second compatibility as in the non-twisted case.

\ssec{Example: $\fg$-modules}

Let $\bA$ be an associative topological algebra with an action of $G$
(see \secref{topological algebroids}).
Then the category $\bA\mod$ carries a weak $G$-action.

If in addition, we have a continuous map $\fg'\to \bA$ that sends
$1\in \BC\subset \fg'$ to the identity in $\bA$ such that the
commutator map $\fg\startimes \bA\to \bA$ is the dual of the map
$\bA\to \fg^*\shriektimes \bA$, obtained by deriving the $G$-action,
then the above action of $G$ on $\bA\mod$ is of twisted Harish-Chandra
type relative to $\fg'$.

We will consider some particular cases of this situation. The most
basic example is $\CC=\fg'\mod$:

\medskip

Let $M$ be a $\fg'$-module. We will denote by $a$ the action map
$\fg'\startimes M\to M$ and by $a^*:M\to \fg'{}^*\shriektimes M$ its
dual. For $S\to G$, we set $\act^*_S(M)$ to be isomorphic to
$\CO_S\otimes M$ as an $\CO_S$-module.  The $\fg'$-action on it is
given via the map
$$\fg'\overset{\Delta_\fg}\to \CO_G\shriektimes \fg'\to
\CO_S\shriektimes \fg'.$$ and the action of the latter on
$\CO_S\otimes M$ by means of $\on{id}\otimes m$.

\medskip

The restriction of $\act^*(M)$ to $G'{}^{(1)}$ identifies as a 
$(\BC\oplus \epsilon\cdot \fg'{}^*)$-module with the free module
$$M\oplus \epsilon \cdot \fg'{}^*\shriektimes M.$$ In terms of this
identification, the $\fg'$-action is given by
$$x\otimes (v_1+\epsilon \cdot v_2)\mapsto a(x\otimes v_1)+\epsilon
\cdot \left((a\otimes \on{id}_{\fg'{}^*})(\on{ad}^*(x)\otimes v_1)+
a(x\otimes v_2)\right),$$ where $\on{ad}^*$ is the map $\fg'\to
\fg'\shriektimes \fg'{}^*$, adjoint to the bracket.

We construct an isomorphism between $M^{(1)}$ and
$$p^*(M)|_{G'{}^{(1)}}\simeq M\oplus \epsilon \cdot
\fg'{}^*\shriektimes M$$ as $\fg'$-modules using the map
$$v_1+\epsilon \cdot v_2 \mapsto v_1+\epsilon \cdot (a^*(v_1)+v_2).$$

The category $\fg'\mod$ is universal in the following sense. Let $\CC$
be an abelian category as above, endowed an action of $G$ and a
functor $\sF:\CC\to \Vect$, respecting the action in the natural
sense.

Assume that the $G$ action on $\CC$ is of Harish-Chandra type relative
to $\fg'$.  Then the functor $\sF$ naturally lifts to a functor
$\CC\to \fg'\mod$.

\ssec{Example: D-modules on $G$}

Consider now the category $\fD(G)'\mod$. We claim that it carries an
action of $G$ of twisted Harish-Chandra type relative to $\fg'$,
corresponding to the action of $G$ on itself by left translations:

Let $\CM$ be an object of $\fD(G)'\mod$, and $S$ a scheme mapping to
$G$.  We define $\act^*(\CM)|_S$ to be isomorphic to $\CO_S\otimes
\CM$ as an $\CO_S$-module.  The action of $\CO_G$ is given via the
co-multiplication map $\CO_G\overset{\Delta_G}\to \CO_G\shriektimes
\CO_G\to \CO_S\shriektimes \CO_G$.  The action of $a_l$ of $\fg'$ is
also given via the map $\fg'\overset{\Delta_\fg}\to \CO_G\shriektimes
\fg\to \CO_S\shriektimes \fg'$.

Note that the action of $\fg'_{-\can}$ on $\act^*(\CM)|_S\simeq
\CO_S\otimes \CM$ is via the $a_r$-action on the second multiple.

The infinitesimal trivialization of this action is defined in the same
way as for $\fg'\mod$ via the map $a_l^*:\CM\to \fg'{}^*\shriektimes
\CM$.

\medskip

We will now define another action of $G$ on $\fD(G)'\mod$,
corresponding to the action of $G$ on itself by right translations. It
will be of twisted Harish-Chandra type relative to $\fg'_{-\can}$:

\medskip

For $\CM$ and $S$ as above, we let $\act^*(\CM)|_S$ to be again
isomorphic to $\CM\otimes \CO_S$ as a $\CO_S$-module, and the
$\CO_G$-action is given via the co-multiplication map $\CO_G\to
\CO_G\shriektimes \CO_S$.  The $a_l$-action of $\fg'$ is $a_l\otimes
\on{id}_{\CO_S}$. The resulting $a_r$-action of $\fg'_{-\can}$ is then
given by the map $\fg'\to \CO_G\shriektimes \fg'\to \CO_S\shriektimes
\fg'$.

The infinitesimal trivialization of the right action is defined in the
same way as for the category $\fg'_{-\can}\mod$ using the map
$a_r^*:\CM\to \fg'_{-\can}{}^*\shriektimes \CM$.

\medskip

It is easy to see that the two actions of $G$ on $\fD(G)'\mod$ commute
in the natural sense. Thus, we obtain an action of $G\times G$ on
$\fD(G)'\mod$, which is of twisted Harish-Chandra type relative to
$\fg\oplus \fg'_{-\can}$.

\ssec{The twisted product}

Let $\CC$ be a category equipped with an action of $G$ of twisted
Harish-Chandra type with respect to $\fg'$. Let $X$ be an object of
$\CC$ and $\CM\in \fD(G)'\mod$.  We will define an object
$\CM\tboxtimes X\in \fD(G)\mod\otimes \CC$:

As an object of $\QCoh^!_G\otimes \CC$, it is isomorphic to
$$\CM\underset{\CO_G}\otimes \act^*(X),$$
see \secref{objects over indscheme}. The action of $\fg'$ on it is
defined as follows.

Consider the ind-subscheme $G^{(1)}\times G\subset G\times G$, and let
$p_2$ denote its projection on the second multiple.  Let $\sk\subset
\fg$ be a lattice and $\Spec(\BC\oplus \epsilon \cdot \sk^*)$ the
corresponding subscheme of $G^{(1)}$.

We have to construct an isomorphism
\begin{equation}   \label{constr action}
\on{mult}^*\left(\CM\underset{\CO_G}\otimes \act^*(X)\right)
|_{\Spec(\BC\oplus \epsilon \cdot \sk^*)\times G}\simeq p_2^*
\left(\CM\underset{\CO_G}\otimes \act^*(X)\right) |_{\Spec(\BC\oplus
  \epsilon \cdot \sk^*)\times G}
\end{equation}
of objects of $\QCoh^!_{\Spec(\BC\oplus \epsilon \cdot \sk^*)\times
G}\otimes \CC$, compatible with the identification
$$\on{mult}^*\left(\CM\underset{\CO_G}\otimes \act^*(X)\right)|_{1\times G}
\simeq \CM\underset{\CO_G}\otimes \act^*(X) \simeq
p_2^*\left(\CM\underset{\CO_G}\otimes \act^*(X)\right)|_{1\times G}.$$

\medskip

Let $\sk'$ be the preimage of $\sk$ in $\fg'$, and let
$\Spec(\BC\oplus \epsilon \cdot \sk'{}^*)$ be the preimage of
$\Spec(\BC\oplus \epsilon \cdot \sk^*)$ in $G'{}^{(1)}$. We have an
isomorphism
$$\on{mult}^*(\CM)|_{\Spec(\BC\oplus \epsilon \cdot \sk'{}^*)\times
G}\simeq p_2^*(\CM)|_{\Spec(\BC\oplus \epsilon \cdot \sk'{}^*)\times
G}$$ in $\QCoh^!_{\Spec(\BC\oplus \epsilon \cdot \sk'{}^*)\times
G}\otimes \CC$, given by the $a_l$-action of $\fg'$ on $\CM$.  We also
have an isomorphism
$$\on{mult}^*(\act^*(X))|_{\Spec(\BC\oplus \epsilon \cdot
\sk'{}^*)\times G}\simeq \act^*(\act^*(X))|_{\Spec(\BC\oplus \epsilon
\cdot \sk'{}^*)\times G}\simeq p_2^*(\act^*(X))|_{\Spec(\BC\oplus
\epsilon \cdot \sk'{}^*)\times G}$$ in $\QCoh^*_{\Spec(\BC\oplus
\epsilon \cdot \sk'{}^*)\times G}\otimes \CC$ where the first arrow is
the associativity constraint for the action, and second one is the
infinitesimal trivialization.

Combining the two we obtain an isomorphism
$$\on{mult}^*\left(\CM\underset{\CO_G}\otimes \act^*(X)\right)
|_{\Spec(\BC\oplus \epsilon \cdot \sk'{}^*)\times G}\simeq p_2^*
\left(\CM\underset{\CO_G}\otimes \act^*(X)\right)
|_{\Spec(\BC\oplus \epsilon \cdot \sk'{}^*)\times G}$$ in 
$\QCoh^!_{\Spec(\BC\oplus \epsilon \cdot \sk'{}^*)\times G}\otimes \CC$, 
but it is easy to see that the two central extensions cancel out, 
and we obtain an isomorphism as in \eqref{constr action}.

\medskip

By construction, this system of isomorphisms is compatible for
different choices of $\sk$. Thus, we obtain an action of $\fg$, as a
Tate vector space, on $\CM\underset{\CO_G}\otimes \act^*(X)$,
satisfying the desired commutation relation with $\CO_G$. Moreover,
from the axioms it follows that this action of $\fg$ is compatible
with the Lie algebra structure. We will denote this action by
$\wt{a}_l$.

\medskip

Thus, $\CM\tboxtimes X$ is an object of $\fD(G)\mod\otimes \CC$; in
particular, it carries an action of $\fg'_{-\can}$, denoted
$\wt{a}_r$. Let us describe this action explicitly:

Let $\sk$ be a lattice in $\fg$ as above, and let $\sk_{-\can}$ denote
its preimage in $\fg_{-\can}$. We have to construct
$$\on{mult}^*\left(\CM\underset{\CO_G}\otimes \act^*(X)\right)
|_{G\times \Spec(\BC\oplus \epsilon \cdot \sk_{-\can}^*)}\simeq p_1^*
\left(\CM\underset{\CO_G}\otimes \act^*(X)\right) |_{G\times
  \Spec(\BC\oplus \epsilon \cdot \sk_{-\can}^*)}$$ in
$\QCoh^!_{G\times \Spec(\BC\oplus \epsilon \cdot
  \sk_{-\can}^*)}\otimes \CC$.  It is constructed as in the previous
case, using the $a_r$-action of $\fg'_{-\can}$ on $\CM$.

\medskip

Finally, let us note that we can consider an object of $\CC$ given by
$$(\CM\tboxtimes X)^{DR}\simeq
\fC^{\frac{\infty}{2}}(\fg,\CM\tboxtimes X)$$ that carries a canonical
differential. We will denote by $\wt{\Lie}_l$, $\wt{\Lie}_r$,
$\wt{i}_l$, $\wt{i}^*_l$, $\wt{i}_r$, $\wt{i}^*_r$ the corresponding
structures on it.

\ssec{Definition of convolution}    \label{convolution action}

Let now $K\subset G$ be an "open compact" group subscheme, over which
$\fg'$ (and hence also $\fg'_{-\can}$) is split. Let $X$ be an object
of $\CC^{w,K}$ (resp., $\CC^{K}$) and $\CM$ be an object of
$\fD(G)'\mod^{w,r(K)}$ (resp., $\fD(G)'\mod^{r(K)}$).  We claim that
in this case $\CM\tboxtimes X$ is naturally an object of
$\fD(G)\mod^{w,r(K)}\otimes \CC$ (resp., $\fD(G)\mod^{r(K)}\otimes
\CC$).  This follows from the description of the action $\wt{a}_r$
given above.

More generally, if $X^\bullet\in \bC(\CC)^K$ and
$\CM^\bullet\in \bC(\fD(G)'\mod)^{r(K)}$, then the complex
$\CM^\bullet\tboxtimes X^\bullet$ is naturally an object of
$\bC(\fD(G)\mod\otimes \CC)^{r(K)}$. We will denote by 
$(\CM^\bullet\tboxtimes X^\bullet)_{G/K}$ the resulting object of 
$\bC(\fD(G/K)\mod\otimes \CC)$.

We define a functor
$$\bC(\fD(G)'\mod)^{r(K)}\times \bC(\CC)^K\to \bC(\CC)$$
by
\begin{equation} \label{defn conv}
\CM^\bullet,X^\bullet\mapsto
\fC^{\frac{\infty}{2}}(\fg;K_{red},\CM^\bullet\tboxtimes X^\bullet).
\end{equation}
This functor is exact when restricted to
$\bC^+(\fD(G)'\mod)^{r(K)}\times \bC^+(\CC)^K$, and hence we obtain a
functor $D^+(\fD(G)'\mod)^{r(K)}\times D^+(\CC)^K\to D(\CC)$, denoted
$$\CM^\bullet,X^\bullet\mapsto \CM^\bullet\underset{K}\star X^\bullet.$$

By \lemref{equiv cohomology of st complex} and \propref{de Rham and
semiinf},
\begin{equation} \label{convolution as De Rham}
\CM^\bullet\underset{K}\star X^\bullet\simeq H^\bullet_{DR}(G/K,
(\CM^\bullet\tboxtimes X^\bullet)_{G/K}).
\end{equation}

\medskip

Using the equivalence $D(\fD(G)'\mod)^{r(K)}\simeq D(\fD(G/K)'\mod)$
we obtain also a functor $D^+(\fD(G/K)'\mod)\times D^+(\CC)^K\to
D(\CC)$, denoted
$$\CF^\bullet,X^\bullet\mapsto \CF^\bullet\underset{K}\star X^\bullet.$$

\medskip

Let $H\subset G$ be another group subscheme, not necessarily
"open compact", and consider the category
$\bC(\fD(G)'\mod)^{l(H),r(K)}$. The above convolution functor is
easily seen to give rise to an exact functor
$$\bC^+(\fD(G)'\mod)^{l(H),r(K)}\times \bC^+(\CC)^K\to \bC(\CC)^H,$$
and the corresponding functor
$$\bC^+(\fD(G/K)'\mod)^{H}\times \bC^+(\CC)^K\to \bC(\CC)^H.$$

Note however, that if we start with an object $\CF^\bullet\in
\bC^+(\fD(G/K)'\mod)^{H}$ that comes from an object in the naive
subcategory $\bC^+(\fD(G/K)'\mod^{H})$, the convolution
$\CF^\bullet\star X^\bullet$ is defined only as an object of
$\bC(\CC)^H$ (and not of $\bC(\CC^H)$). This is one of the reasons why
one should work with $\bC(\CC)^H$, rather than with $\bC(\CC^H)$.

\medskip

Let us denote by $\bC^{bd}(\fD(G/K)\mod)$ the subcategory of
$\bC(\fD(G/K)\mod)$ that consists of bounded from below complexes,
supported on a finite-dimensional closed subscheme of $G/K$. Let
$D^{bd}(\fD(G/K)\mod)$ be the corresponding full subcategory of
$D(\fD(G/K)\mod)$.

Let $\bC^{bd}(\fD(G)\mod)^{r(K)}$ be the subcategory of
$\bC(\fD(G)\mod)^{r(K)}$, consisting of bounded from below complexes,
supported set-theoretically on a preimage of a finite-dimensional
closed subscheme of $G/K$; let $D^{bd}(\fD(G)\mod)^{r(K)}$ be the
corresponding full subcategory of $D(\fD(G)\mod)^{r(K)}$.

One easily shows that under
the equivalence $D(\fD(G/K)\mod)\simeq D(\fD(G)\mod)^{r(K)}$ the
subcategories $D^{bd}(\fD(G/K)\mod)$ and $D^{bd}(\fD(G)\mod)^{r(K)}$
correspond to one-another.

\begin{lem}    \label{boundedness}
For $\CM^\bullet\in D^{bd}(\fD(G)\mod)^{r(K)}$ and $X^\bullet\in
D^+(\CC)^K$ (resp., $X^\bullet\in D^b(\CC)^K$), the convolution
$\CM^\bullet\underset{K}\star X^\bullet$ belongs to $D^+(\CC)$ (resp.,
$D^b(\CC)$).
\end{lem}

\begin{proof}

Under the assumptions of the lemma the $\CC$-valued complex of
D-modules $(\CM^\bullet \tboxtimes X^\bullet)_{G/K}$ is
quasi-isomorphic to one bounded from below (resp., bounded) and
supported on a finite-dimensional closed subscheme of $G/K$. Hence,
its de Rham cohomology is bounded from below (resp., bounded).

\end{proof}

\ssec{Examples} \label{semijective}

Let us consider the basic example, when $\CM$ is the $\fD(G)'$-module
$$\delta'_{K,G}\simeq \Ind_{\sk\oplus \BC}^{\fg'}(\CO_K).$$ Note that
$\delta'_{K,G}\in \bC(\fD(G)'\mod)^{l(K),r(K)}$.

\begin{prop}   \label{convolution with delta}
For $X^\bullet\in \bC(\CC)^K$, we have a canonical quasi-isomorphism
in $\bC(\CC)^K$:
$$\delta'_{K,G}\underset{K}\star X^\bullet\simeq X^\bullet$$
\end{prop}

\begin{proof}

Note that we can regard $\CC$ as a category, acted on by $K$ (rather
than $G$).  In particular, it makes sense to consider
$\delta_{K,K}\tboxtimes X^\bullet\in \bC(\fD(K)\mod\otimes
\CC)^{r(K)\times K}$.

Let us regard $\delta'_{K,G}\tboxtimes X^\bullet$ as an object of the
categories $\bC(\fg_{-\can}\mod\otimes \CC)^{K\times K}$ and
$\bC(\sk\mod\otimes \CC)^{K\times K}$.  We have a natural map
$\delta_{K,K}\tboxtimes X^\bullet\to \delta'_{K,G}\tboxtimes
X^\bullet$ in the latter category, and since as a
$\fg'_{-\can}$-module $\delta'_{K,G}\simeq
\on{Ind}^{\fg'_{-\can}}_{\sk\oplus \BC}(\CO_K)$, the latter map
induces an isomorphism
$$\on{Ind}^{\fg'_{-\can}}_{\sk\oplus \BC}\left(\delta_{K,K}\tboxtimes
X^\bullet\right)\to \delta'_{K,G}\tboxtimes X^\bullet\in
\bC(\fg_{-\can}\mod\otimes \CC)^{K\times K}.$$

Hence, as objects of $\bC(\CC)^K$,
$$\delta'_{K,G}\underset{K}\star X^\bullet:=
\fC^{\frac{\infty}{2}}(\fg;K_{red},\delta'_{K,G}\tboxtimes
X^\bullet)\overset{\on{quasi-isom}}\simeq
\fC(\sk;K_{red},\delta_{K,K}\tboxtimes X^\bullet).$$

This reduces the assertion of the proposition to the case when
$G=K$. Note that we have a natural map
$$X^\bullet\to \fC(\sk;K_{red},\delta_{K,K}\tboxtimes X^\bullet),$$
and it is easily seen to be a quasi-isomorphism, since as objects of
$\bC(\CC)$,
$$\fC(\sk;K_{red},\delta_{K,K}\tboxtimes X^\bullet)\simeq
\on{Av}_{K_u}(X^\bullet).$$

\end{proof}

More generally, let $K'\subset K$ be a group subscheme, and let
$\delta'_{K/K',G/K'}$ be the twisted D-module on $G/K'$ equal to the
direct image of $\CO_{K/K'}$ under $K/K'\to G/K'$. Arguing as above, we
obtain the following

\begin{lem}   
For $X^\bullet\in \bC(\CC)^{K'}$,
$$\delta'_{K/K',G/K'}\underset{K'}\star X^\bullet\simeq
\Av_K(X^\bullet)\in D(\CC)^K.$$
\end{lem}

\medskip

Let now $\bg$ be a point of $G$. For an object $X\in \CC$ we will
denote by $X^\bg$ (or $\delta_{\bg,G}\star \CM$) the twist of $X$ by
means of $\bg$, i.e., the restriction of $\act^*(X)$ to $\bg$.

Applying this to $\CF\in \fD(G/K)'\mod$, we obtain a $\bg$-translate
of $\CF$ with respect to the action of $G$ on $G/K$. In particular,
$(\delta_{1,G/K})^\bg\simeq (\delta_{\bg,G/K})$. The following results
from the definitions:

\begin{lem} For $\CF\in \bC(\fD(G/K)'\mod)$, $X^\bullet\in \bC(\CC)^K$,
$$\left(\CF^\bullet\underset{K}\star X^\bullet\right)^\bg\simeq
(\CF^\bullet)^\bg\underset{K}\star X^\bullet.$$
 \end{lem}

In particular, for $X^\bullet$ as above,
$$(X^\bullet)^\bg\simeq \delta_{\bg,G/K}\underset{K}\star X^\bullet.$$

\medskip

Let $G_1\subset G$ be a group subindscheme, and $K_1=K\cap G_1$, and
let $\CF_1^\bullet$ be an object of $\bC(\fD(G_1/K_1)'\mod$, and let
$\CF^\bullet\in \bC(\fD(G/K)'\mod)$ be its direct image under
$G_1/K_1\to G/K$.

The action of $G$ on $\CC$ induces an action of $G_1$; hence, for
$X^\bullet\in \bC(\CC)^K$ it makes sense to consider the object
$\CF^\bullet_1\underset{K_1}\star X^\bullet\in \bC(\CC)$.

\begin{lem}  \label{other group}
For $\CF_1^\bullet\in \bC(\fD(G_1/K_1)'\mod$ and $X^\bullet\in
\bC(\CC)^K$, the objects $\CF_1^\bullet\underset{K_1}\star X^\bullet$
and $\CF^\bullet\underset{K}\star X^\bullet$ in $\bC(\CC)$ are
canonically quasi-isomorphic.
\end{lem}

\begin{proof}

Let $\CM^\bullet$ (resp., $\CM_1^\bullet$) be the object of
$\bC(\fD(G)'\mod)^K$ (resp., $\bC(\fD(G_1)'\mod)^{K_1}$) be the object
corresponding to $\CF^\bullet$ (resp., $\CF_1^\bullet$) under the
equivalence of \propref{AG equivalence}.

Let $(\CM^\bullet\tboxtimes X^\bullet)_{G/K}$ (resp.,
$(\CM_1^\bullet\tboxtimes X^\bullet)_{G_1/K_1}$) be the corresponding
objects of the categories $\bC(\fD(G/K)'\mod\otimes \CC)$ and
$\bC(\fD(G_1/K_1)'\mod\otimes \CC)$, respectively.

The assertion follows now from the fact that $\bC(\fD(G/K)'\mod\otimes
\CC)$ is the direct image $\bC(\fD(G_1/K_1)'\mod\otimes \CC)$ under
$G_1/K_1\hookrightarrow G/K$.

\end{proof}

\medskip

Finally, let us consider the example when $\CC=\fD(\CY)'$, where $\CY$
is a strict ind-scheme, acted on by $G$, and $\fD(\CY)'$ is the
category of twisted D-modules on $\CY$, compatible with a twisting on
$G$.

Recall that in this case we have a functor
\begin{equation} \label{geom conv}
D^{bd}(\fD(G/K)')\times D^b(\fD(\CY)'\mod)^K\to D^b(\fD(\CY)'\mod)
\end{equation}
defined as follows:

Consider the ind-scheme $G\underset{K}\times \CY$, which maps to $\CY$
via the action map of $G$ on $\CY$; this ind-scheme is equipped with a
twisting, which is pulled back from the one on $\CY'$ using the above
map.

For $\CF_1\in D^{bd}(\fD(G/K)')$, $\CF_2\in D^b(\fD(\CY)'\mod)^K$ one
can form their twisted external product
$$\CF_1\tboxtimes \CF_2\in D^b(\fD(G\underset{K}\times \CY)').$$
Then $\CF_1\underset{K}\star \CF_2$ is the direct image of 
$\CF_1\tboxtimes \CF_2$ under the above map
$G\underset{K}\times \CY\to \CY$.

It follows immediately from the definitions, that the functor
\eqref{geom conv} is canonically isomorphic to the one given by
\eqref{defn conv}.

\ssec{Convolution action on Harish-Chandra modules}

We shall now study a particular case of the above situation, when
$\CC=\fg'\mod$.  First, for $\CM\in \fD(G)'\mod$ and $N\in \fg'\mod$
let us describe the object $\CM\tboxtimes N\in \fD(G)\mod\otimes
\fg'\mod$ more explicitly.

As a vector space $\CM\tboxtimes N$ is isomorphic to $\CM\otimes N$.
We will denote by $m$ the action of $\CO_G$ on $\CM$, and by $a_l$,
$a_r$ the actions of $\fg'$, $\fg'_{-\can}$ on it. We will denote by
$a$ the action of $\fg'$ on $N$.

Let $\wt{m}$, $\wt{a}_l$, $\wt{a}_r$ and $\wt{a}$ be the actions of
$\CO_G$, $\fg$, $\fg_{-\can}$ and $\fg'$, respectively, on $\CM\otimes
N$ defining on $\CM\tboxtimes N$ a structure of object of
$\fD(G)\mod\otimes \fg'\mod$. We have

\begin{itemize}

\item $\wt{m}=m\otimes \on{id}$,

\item $\wt{a}=(m\otimes a)\circ (\gamma\otimes \on{id}_\fg)\circ
\Delta_{\fg'}$,

\item $\wt{a}_l=(a_l\otimes \on{id})-(m\otimes a)\circ (\gamma\otimes
\on{id}_\fg)\circ \Delta_{\fg'}$,

\item $\wt{a}_r=a_r\otimes \on{id}+ \on{id}\otimes a$.

\end{itemize}

More generally, if $\CM^\bullet$ is an object of
$\bC(\fD(G)'\mod)^{r(K)}$ and $N^\bullet$ is an object of
$\bC(\fg'\mod)^K$, the twisted product $\CM\tboxtimes N$ is naturally
an object of $\bC(\fD(G)\mod)^{r(K)}\otimes \bC(\fg'\mod)$, where the
algebraic action of $K$ on $\CM^\bullet\otimes N^\bullet$ is the
diagonal one, and so is the action of $\sk[1]$.

In this case the convolution $\CM^\bullet\underset{K}\star N^\bullet$
is computed by means of the complex
$$\fC^{\frac{\infty}{2}}\left(\fg;K_{red}, \CM^\bullet\otimes
N^\bullet\right),$$ with respect to the diagonal action of
$\fg_{-\can}$. The $\fg'$-module structure on this complex is given by
$\wt{a}$.

Note, however, that the above complex carries a different action of
$\fg'$, namely one given by $a_l$. We will denote this other functor
$\bC(\fD(G)'\mod)^{r(K)}\times \bC(\fg'\mod)^K\to \bC(\fg'\mod)$ by
$$\CM^\bullet, N^\bullet\mapsto
\CM^\bullet\overset{\natural}{\underset{K}\star} N^\bullet.$$ Note
that if $\CM\in \bC(\fD(G)'\mod)^{H,r(K)}$ for some group-scheme $H$,
then $\CM^\bullet\overset{\natural}{\underset{K}\star} N^\bullet$ is
naturally an object of $\bC(\fg'\mod)^H$.

The two actions of $\fg'$ on $\fC^{\frac{\infty}{2}}\left(\fg;K_{red},
\CM^\bullet\otimes N^\bullet\right)$ are related by the formula
$$a_l-\wt{a}=\wt{a}_l=[d,\wt{i}_l],$$ where $\wt{i}_l$ is the action
of the annihilation operators on
$$\fC^{\frac{\infty}{2}}\left(\fg;K_{red}, \CM^\bullet\otimes
N^\bullet\right)\subset \fC^{\frac{\infty}{2}}\left(\fg,
\CM^\bullet\otimes N^\bullet\right)\simeq (\CM\tboxtimes N)^{DR}.$$

Therefore, the cohomologies of $\CM^\bullet\underset{K}\star
N^\bullet$ and $\CM^\bullet\overset{\natural}{\underset{K}\star}
N^\bullet$ are isomorphic as $\fg'$-modules.

\begin{cor}  \hfill  \label{endomorphisms of semijective}

\smallskip

\noindent{\em (1)} For $N\in (\fg',K)\mod$, the complex
$$\fC^{\frac{\infty}{2}}\left(\fg;K_{red},\delta'_{K,G}\otimes
N\right)$$ is acyclic away from degree $0$.

\smallskip

\noindent{\em (2)} When regarded as a $\fg'$-module via the
$a_l$-action on $\delta'_{K,G}$, the above $0$-th cohomology is
isomorphic to $N$.

\smallskip

\noindent{\em (3)} The image of $U(\fg')$ in
$\End_{\fg'_{-\can}}(\delta'_{K,G})$ is dense.

\end{cor}

\begin{proof}

The first two points follows from \propref{convolution with delta} and
the above comparison of $\CM\overset{\natural}{\underset{K}\star}
N^\bullet$ and $\CM^\bullet\underset{K}\star N^\bullet$.

Let $U(\fg',\sk)$ be the topological algebra of endomorphisms of the
forgetful functor
$$(\fg',K)\mod\to \Vect.$$ Evidently, the image of $U(\fg')$ in
$U(\fg',\sk)$ is dense.  We claim now that $U(\fg',\sk)$ is isomorphic
to $\End_{\fg'_{-\can}}(\delta'_{K,G})$.

The map in one direction, i.e., $U(\fg',\sk)\to
\End_{\fg'_{-\can}}(\delta'_{K,G})$, is evident: given an element in
$U(\fg',\sk)$, we obtain a {\it functorial} endomorphism of every
vector space underlying an object of $\fg'\mod$; in particular
$\delta'_{K,G}$.  This endomorphism commutes with
$\fg'\mod$-endomorphisms of $\delta'_{K,G}$, in particular, with the
action of $\fg'_{-\can}$.

To construct the map in the opposite direction, note that an
endomorphism of $\delta'_{K,G}$ as a $\fg'_{-\can}$-module defines an
endomorphism of the functor
$$N\mapsto
h^0\left(\fC^{\frac{\infty}{2}}\left(\fg;K_{red},\delta'_{K,G}\otimes
N\right)\right): \fg'\mod\to \Vect,$$ and the latter is isomorphic to
the forgetful functor.

\end{proof}

We will now study the behavior of Lie algebra cohomology under
convolution.  We shall first consider a technically simpler case, when
we will consider D-modules on a group-scheme $H$, mapping to $G$, such
that $\fg'$ splits over $\fh$. Let $K'_H,K''_H\subset H$ be group
subschemes of finite codimension.

\begin{prop}    \label{semiinf and de rham}
For $N^\bullet\in D^+(\fg\mod)^{K'_H}$ and $\CM^\bullet\in
D^+(\fD(H))^{l(K'_H),r(K''_H)}$, the complex
$\fC(\fh;K'_H{}_{red,}\CM^\bullet\underset{K'_H}\star N^\bullet)$ is
quasi-isomorphic to
$$H^\bullet_{DR}(H'_K\backslash H,\CF^\bullet)\otimes
\fC(\fh;H''_K{}_{red},\CM^\bullet),$$ where $\CF^\bullet$ is the
object of $D^+(\fD(H'_K\backslash H))$, corresponding to
$\CM^\bullet$.
\end{prop}

\begin{proof}

By \lemref{homotopy lemma}, 
$$\fC(\fh;K'_H{}_{red},\CM^\bullet\underset{K''_H}\star
N^\bullet)\simeq
\fC(\fh;K'_H{}_{red},\CM^\bullet\overset{\natural}{\underset{K''_H}\star}
N^\bullet).$$ The latter, by definition can be rewritten as
$$\fC(\fh\oplus \fh;K'_H{}_{red}\times K''_H{}_{red},
\CM^\bullet\otimes N^\bullet),$$ where the action of the first copy of
$\fh$ is via $a_l$ on $\CM$, and the action of the second copy is
diagonal with respect to $a_r$ and $a$. Hence, the above expression
can be rewritten as
$$\fC(\fh;K''_H{}_{red},\fC(\fh;K'_H{}_{red},\CM^\bullet)\otimes
N^\bullet),$$ where the $\fh$-action on
$\fC(\fh;K'_H{}_{red},\CM^\bullet)\otimes N^\bullet$ is the diagonal
one with respect to the $a_r$-action on $\CM^\bullet$ and the existing
action on $N^\bullet$.

Applying again \lemref{homotopy lemma}, we can replace the
$a_r$-action on $\CM^\bullet$ by the trivial one. Hence,
$$\fC(\fh;K''_H{}_{red},\fC(\fh;K'_H{}_{red},\CM^\bullet)\otimes
N^\bullet)\simeq \fC(\fh;K'_H{}_{red},\CM^\bullet) \otimes
\fC(\fh;K''_H{}_{red},N^\bullet),$$ which is what we had to show.

\end{proof}

We will now generalize the above proposition to the case of
semi-infinite cohomology with respect to $\fg$.

Let $N_1^\bullet$ and $N_2^\bullet$ be objects of
$D^+(\fg'_{-\can}\mod)^{K'}$ and $D^+(\fg'\mod)^{K''}$, respectively,
for some "open compact" $K,K''\subset G$.  Let $\CM^\bullet$ be an object
of $D^+(\fD(G)'\mod)^{l(K'),r(K'')}$, supported over a closed
pro-finite dimensional subscheme of $G$. In this case the convolution
$\CM^\bullet\underset{K''}\star N_2^\bullet$ makes sense as an object
of $D^+(\fg'\mod)^{K'}$. Similarly, we can consider the convolution
"on the right"
$$N_1^\bullet \underset{K'}\star \CM^\bullet\in
D^+(\fg'_{-\can}\mod)^{K'}.$$

\begin{prop}   \label{semiinf adjunction}
Under the above circumstances, 
$$\fC^{\frac{\infty}{2}}\left(\fg;K'_{red}, N_1^\bullet \otimes
\left(\CM^\bullet\underset{K''}\star N_2^\bullet\right)\right)$$ and
$$\fC^{\frac{\infty}{2}}\left(\fg;K''_{red}, \left(N_1^\bullet
\underset{K'}\star \CM^\bullet\right)\otimes N_2^\bullet\right)$$ are
quasi-isomorphic.
\end{prop}

\begin{proof}

By symmetry, it would be sufficient to show that there exists a
quasi-isomorphism between
\begin{equation} \label{asym exp}
\fC^{\frac{\infty}{2}}\left(\fg;K'_{red}, N_1^\bullet \otimes 
\left(\CM^\bullet\underset{K''}\star N_2^\bullet\right)\right)
\end{equation}
and
\begin{equation} \label{sym exp}
\fC^{\frac{\infty}{2}}\left(\fg\oplus \fg;K'_{red}\times K''_{red},
N_1^\bullet \otimes \CM^\bullet\otimes N_2^\bullet\right),
\end{equation}
where the first copy of $\fg_{-\can}$ acts diagonally on $N_1^\bullet
\otimes \CM^\bullet$ (via the existing $\fg'_{-\can}$ action on
$N_1^\bullet$ and the $a_l$-action on $\CM^\bullet$) and second copy
acts diagonally on $\CM^\bullet\otimes N_2^\bullet$ (via the
$a_r$-action on $\CM^\bullet$ and the existing $\fg'$-action on
$N_2^\bullet$).

\medskip

By \lemref{homotopy lemma}, in \eqref{sym exp} we can replace the
action of the first copy of $\fg_{-\can}$, by one where the
$\fg'$-action on $\CM^\bullet\otimes N_2^\bullet$ is given by
$\wt{a}$. The resulting expression would be equal to the one in
\eqref{asym exp} modulo the following complication:

To define \eqref{asym exp} one has to replace
$\fC^{\frac{\infty}{2}}(\fg;K''_{red}, \CM^\bullet\otimes
N_2^\bullet)$ by a quasi-isomorphic complex, which is bounded from
below.  We have to show that taking
$\fC^{\frac{\infty}{2}}\left(\fg;K'_{red}, N_1^\bullet \otimes
?\right)$ survives this quasi-isomorphism.

\medskip

Let us first consider a particular case, when $\CM^\bullet$ is induced
from an $\CO_G$-module, i.e., has the form
\begin{equation}  \label{good rows} 
\Ind^{\fg'_{-\can}}_{\sk''}(\CL^\bullet)
\end{equation} 
for some complex $\CL^\bullet$ of $K'$-equivariant $\CO_G$-modules. In
this case we have a quasi-isomorphism
$$\fC(\sk'';K''_{red}, \CL^\bullet\otimes N_2^\bullet)\to
\fC^{\frac{\infty}{2}}(\fg;K''_{red}, \CM^\bullet\otimes
N_2^\bullet)$$ of complexes of $\fg'$-modules. Moreover, the PBW
filtration defines a filtration on the RHS, of which
$\fC(\sk'';K''_{red}, \CL^\bullet\otimes N_2^\bullet)$ is the first
term, by $\fg'$-stable subcomplexes, all quasi-isomorphic to
one-another.

Since the functor $\fC^{\frac{\infty}{2}}\left(\fg;K'_{red},
N_1^\bullet \otimes ?\right)$ commutes with direct limits, the
required assertion about quasi-isomorphism holds.

\medskip

The case of a general $\CM^\bullet$ follows from the one considered
above, since the assumption on $\CM^\bullet$ implies that it can be
represented by a complex, associated with a bi-complex with finitely
many rows, each of the form \eqref{good rows}.

\end{proof}

\ssec{Convolution action on $\fD(G)$-modules}     \label{convolution
  on CADO} 

Let us now consider the case when $\CC=\fD(G)'$ with the action of $G$
by left translations.

Given two objects $\CM_1,\CM_2\in \fD(G)'\mod$ let us first describe
how $\CM_1\tboxtimes \CM_2$ looks like as an object of
$\fD(G)\mod\otimes \fD(G)'\mod$.

By construction as a vector space $\CM_1\tboxtimes \CM_2\simeq
\CM_1\otimes \CM_2$.  We will denote the $a^1_l$, $a^2_l$ (resp.,
$a^1_r$, $a^2_r$, $m^1$, $m^2$) the actions of $\fg'$ (resp.,
$\fg'_{-\can}$, $\CO_G$) on $\CF_1$ and $\CF_2$, respectively.  We
will denote by $\wt{a}^1_l$, $\wt{a}^2_l$, $\wt{a}^1_r$, $\wt{a}^2_r$,
$\wt{m}^1$, $\wt{m}^2$ the actions of $\fg$, $\fg'$, $\fg_{-\can}$,
$\fg'_{-\can}$, respectively on $\CF_1\tboxtimes \CF_2$, corresponding
to the $\fD(G)\mod\otimes \fD(G)'\mod$-structure.

The action of $\CO_G$, corresponding to the $\fD(G)\mod$-structure on
$\CM_1\tboxtimes \CM_2$ is via the first multiple in $\CM_1\otimes
\CM_2$, which we will denote by $m$. The action of $\CO_G$,
corresponding to the $\fD(G)'\mod$-structure is via the
co-multiplication map $\Delta_G:\CO_G\to \CO_G\shriektimes \CO_G$.

These actions are described as follows: 

\begin{itemize}

\item
$\wt{m}^1=m^1\otimes \on{id}$,

\item
$\wt{m}^2=(m^1\otimes m^2)\circ (\gamma\otimes \on{id})\circ
\Delta_G$,

\item
$\wt{a}^1_l=a^1_l\otimes \on{id}-(m\otimes a^2_l)\circ (\gamma\otimes
\on{id})\circ \Delta_{\fg'}$,

\item
$\wt{a}^2_l=(m\otimes a^2_l)\circ (\gamma\otimes \on{id})\circ \Delta_{\fg'}$,

\item
$\wt{a}^1_r=a^1_r\otimes \on{id}+\on{id}\otimes a^2_l$,

\item
$\wt{a}^2_r=\on{id}\otimes a^2_r$.

\end{itemize}

If $\CM^\bullet_1\in \bC^+(\fD(G)')\mod^{r(K)}$ and $\CM^\bullet_2\in
\bC^+(\fD(G)')\mod^{l(K)}$, the convolution
$\CM^\bullet_1\underset{K}\star \CM^\bullet_2$ is computed by means of
$$\fC^{\frac{\infty}{2}}(\fg;K_{red},\CM_1^\bullet\otimes
\CM_2^\bullet),$$ with respect to the diagonal (i.e.,
$\wt{a}^1_r=a^1_r+a^2_l$) action of $\fg_{-\can}$, and the actions of
$\CO_G$, $\fg'$ and $\fg'_{-\can}$, specified above.

Note, however, that the above complex carries a different
$\fD(G)'$-module structure. Namely, the action of $\CO_G$ is
$(m_1\otimes m_2)\circ (\gamma\otimes \on{id})\circ \Delta_G$ as
before, and the action of $\fg'$ is $a^1_l$. In this case the action
of $\fg'_{-\can}$ equals $(a^1_r\otimes m)\circ \Delta_\fg$.

We will denote this new functor
$$\bC^+(\fD(G)')\mod^{r(K)}\times \bC^+(\fD(G)')\mod^{l(K)}\to
\bC(\fD(G)'\mod$$
by
$$\CM^\bullet_1,\CM^\bullet_2\mapsto
\CM^\bullet_1\overset{\natural}{\underset{K}\star} \CM^\bullet_2.$$

\begin{lem}
For $\CM^\bullet_1\in D^{bd}(\fD(G)'\mod)^{r(K)}$,
$\CM^\bullet_2\in D^{bd}(\fD(G)'\mod)^{l(K)}$,
the objects
$$\CM^\bullet_1\underset{K}\star \CM^\bullet_2,
\CM^\bullet_1\overset{\natural}{\underset{K}\star} \CM^\bullet_2\in
D^b(\fD(G)'\mod)$$ are isomorphic.
\end{lem}

\begin{proof}

{}From the assumption it follows that there exist "open compact" groups
$K',K''$ such that $\CM^\bullet_1\in D(\fD(G)'\mod)^{l(K'),r(K)}$ and
$\CM^\bullet_2\in D(\fD(G)'\mod)^{l(K),r(K'')}$.

As we saw above, the convolution $\CM^\bullet_1\underset{K}\star
\CM^\bullet_2$ can be interpreted as an action of $\CF_1^\bullet \in
D^{bd}(\fD(G/K)'\mod)^{K'}$, corresponding to $\CM_1^\bullet$, on
$\CF^\bullet_2\in D^b(\fD(G/K'')'\mod)^{K}$, corresponding to
$\CM_2^\bullet$. The result is an object in
$D^b(\fD(G/K'')\mod)^{K'}$.

However, this convolution can be rewritten also as an action of
$'\CF^\bullet_2\in D^{bd}(\fD(K\backslash G)'\mod)^{K''}$ on
$'\CF_1^\bullet\in D^b(\fD(K'\backslash G)'\mod)^{K}$, with the result
being in
$$D^b(\fD(K'\backslash G)'\mod)^{K''}\simeq D^b(\fD(G/K'')\mod)^{K'}.$$

The latter convolution is manifestly the same as 
$$\CM^\bullet_1\overset{\natural}{\underset{K}\star} \CM^\bullet_2\in
D^b(\fD(G)'\mod)^{l(K'),r(K'')}.$$

\end{proof}

\ssec{Associativity of convolution}

Let now $\CM_1^\bullet$ be an object of
$\bC^{bd}(\fD(G)'\mod)^{r(K)}$, $\CM_2^\bullet\in
\bC^{bd}(\fD(G)'\mod)^{l(K),r(K')}$ and $X^\bullet\in \bC^+(\CC)^{K'}$
for a category $\CC$ as above,

\begin{prop}  \label{assoc of conv}
Under the above circumstances, there exists a canonical isomorphism in
$D^+(\CC)$
$$\left(\CM_1^\bullet\underset{K}\star \CM_2^\bullet\right)\underset{K'}
\star X^\bullet\simeq 
\CM_1^\bullet\underset{K}\star \left(\CM_2^\bullet\underset{K'}
\star X^\bullet\right),$$
compatible with three-fold convolutions.
\end{prop}

The rest of this subsection is devoted to the proof of this proposition.

Consider the bi-graded object of $\CC$ given by
\begin{equation} \label{triple product}
\left(\CM_1^\bullet\boxtimes
\CM_2^\bullet\right)\underset{\CO_{G\times G}}\otimes
\on{mult}^*\left(\act^*(X^\bullet)\right).
\end{equation}

It carries two actions of the Lie algebra $\fg_{-\can}\oplus
\fg_{-\can}$, corresponding to the two isomorphisms
$$\CM_1^\bullet \tboxtimes \left(\CM_2^\bullet\tboxtimes
X^\bullet\right)\simeq \left(\CM_1^\bullet\boxtimes
\CM_2^\bullet\right)\underset{\CO_{G\times G}}\otimes
\on{mult}^*\left(\act^*(X^\bullet)\right)\simeq \left(\CM_1^\bullet
\tboxtimes \CM_2^\bullet\right)\tboxtimes X^\bullet.$$

The action of the second copy of $\fg_{-\can}$ is the same in the two
cases.  The difference of the actions of the first copy of
$\fg_{-\can}$ is given by the $\fg$-action, coming from its
$\wt{a}_l$-action on $\CM_2^\bullet\tboxtimes X^\bullet$.

Hence, by \lemref{homotopy lemma}, the two complexes
$$\fC^{\frac{\infty}{2}}\Bigl(\fg\oplus\fg;K\times K', \CM_1^\bullet
\tboxtimes \left(\CM_2^\bullet\tboxtimes X^\bullet\right)\Bigr)$$ and
$$\fC^{\frac{\infty}{2}}\Bigl(\fg\oplus\fg;K\times K',
\left(\CM_1^\bullet \tboxtimes \CM_2^\bullet\right)\tboxtimes
X^\bullet\Bigr)$$ are isomorphic.

As in the proof of \propref{semiinf adjunction}, we have to show that
the above complexes are isomorphic in the derived category to
$\CM_1^\bullet\underset{K}\star \left(\CM_2^\bullet\underset{K'}\star
X^\bullet\right)$ and $\left(\CM_1^\bullet\underset{K}\star
\CM_2^\bullet\right)\underset{K'}\star X^\bullet$, respectively. This
is done as in the proof of \propref{semiinf adjunction} by replacing
$\CM_1^\bullet$ and $\CM_2^\bullet$ by appropriately chosen complexes,
for which the above semi-infinite complexes can be represented as
direct limits of quasi-isomorphic complexes, bounded from below.

\ssec{An adjunction in the proper case}   \label{proper conv}

Let now $K_1,K_2\in G$ be two "open compact" subgroups of $G$, and
assume that $G/K_1$ is ind-proper. Let $\CF$ be a finitely generated
object of $\fD(G/K_1)'\mod^{K_2}$. As in \secref{bi-equivariant},
we have a well-defined object $\CF^{\on{op}}$ in
$\fD(G/K_2)^{''}\mod^{K_1}$, where the superscript $''$ indicates the
twisting opposite to $'$. Then the Verdier dual $\BD(\CF^{\on{op}})$ is an
object of $\fD(G/K_2)'\mod^{K_1}$.

\begin{prop}   \label{HCh adjunction}
The functor $$D(\CC)^{K_1}\to D(\CC)^{K_2}: X_1^\bullet\mapsto
\CF\underset{K_1}\star X_1^\bullet$$ is left adjoint to the functor
$$D(\CC)^{K_2}\to D(\CC)^{K_1}: X_2\mapsto
\BD(\CF^{\on{op}})\underset{K_2}\star X_2.$$
\end{prop}

\begin{proof}

We need to construct the adjunction maps
$$\CF\underset{K_1}\star (\BD(\CF^{\on{op}})\underset{K_2}\star
X^\bullet_2)\to X^\bullet_2 \text{ and } X^\bullet_1\to
\BD(\CF^{\on{op}})\underset{K_2}\star (\CF\underset{K_1}\star
X_1^\bullet),$$ such that the identities concerning the two
compositions hold.

In view of \propref{assoc of conv}, it would suffice to construct the
maps
$$\CF\underset{K_1}\star \BD(\CF^{\on{op}})\to \delta_{1,G/K_2}\in
D(\fD(G/K_2)\mod)^{K_2}$$ and
$$\delta_{1,G/K_1}\to \BD(\CF^{\on{op}})\underset{K_2}\star \CF\in
D(\fD(G/K_1)\mod)^{K_1},$$ such that the corresponding identities
hold.

By the definition of convolution, constructing these maps is
equivalent to constructing morphisms
\begin{equation} \label{adj 1}
H^\bullet(G/K_1,\Delta_{G/K_1}^*(\CF\boxtimes \BD(\CF))\to \BC\in
D(\on{pt}/K_2)
\end{equation}
and
\begin{equation} \label{adj 2}
\BC\to H^\bullet(G/K_2,\Delta_{G/K_2}^!(\CF^{\on{op}}\boxtimes
\BD(\CF^{\on{op}}))\in D(\on{pt}/K_1),
\end{equation}
respectively, where $\Delta_{G/K}$ denotes the diagonal morphism
$G/K\to G/K\times G/K$.  (Note that in each of the cases, the
pull-back of the corresponding twisted D-module on the product under
the diagonal map is a non-twisted right D-module.)

The morphism in \eqref{adj 2} follows from Verdier duality, and
likewise for \eqref{adj 1}, using the fact that
$$H^\bullet(G/K_1,\cdot)\simeq H_c^\bullet(G/K_1,\cdot).$$

The fact that the identities concerning the compositions of adjunction
maps hold, is an easy verification.

\end{proof}


\section{Categories over topological commutative algebras}
\label{cat over top alg}

\ssec{The notion of a category flat over an algebra}

Let $\CC$ be an abelian category as in \secref{class of cat},
satisfying assumption (**), and let $Z$ be a commutative algebra,
mapping to the center of $\CC$. An example of this situation is when
$\bA$ is a topological algebra, $Z$ is a (discrete) commutative
algebra mapping to the center of $\bA$ and $\CC=\bA\mod$. Then
the functor $\sF$ factors naturally through a functor $\sF_Z:\CC\to
Z\mod$.

Note that we have
a naturally defined functor $Z\mod\times \CC\to \CC$ given by
$$M,X\mapsto M\underset{Z}\otimes X.$$ This functor is right exact in
both arguments. We have
$$\sF_Z(M\underset{Z}\otimes X)\simeq M\underset{Z}\otimes \sF_Z(X).$$
This shows, in particular, that if $M$ is $Z$-flat, then the
above functor of tensor product is exact in $Z$. 
We will denote by 
$$M^\bullet,X^\bullet\mapsto M^\bullet\overset{L}{\underset{Z}\otimes}
X^\bullet: D^-(Z\mod)\times D^-(\CC)\to D^-(\CC)$$ the corresponding
derived functor. We have
$$\sF_Z(M^\bullet\overset{L}{\underset{Z}\otimes} X^\bullet)\simeq
M^\bullet\overset{L}{\underset{Z}\otimes} \sF_Z(X^\bullet).$$

It is easy to see that for a fixed
$X^\bullet\in \bC^-(\CC)$, the derived functor of
$$M^\bullet\mapsto M^\bullet\underset{Z}\otimes
X^\bullet:\bC^-(Z\mod)\to \bC^-(\CC)$$ is isomorphic to
$M^\bullet\overset{L}{\underset{Z}\otimes} X^\bullet$.  However, this
is not, in general, true for the functor
$$X^\bullet\mapsto M^\bullet\underset{Z}\otimes
X^\bullet:\bC^-(\CC)\to \bC^-(\CC)$$ for a fixed $M^\bullet$.

\medskip

We shall say that an object $X\in \CC$ is flat over $Z$ if the functor
$$M\mapsto M\underset{Z}\otimes X:Z\mod\to \CC$$
is exact. This is equivalent to $\sF_Z(X)$ being flat as a $Z$-module.

\medskip

We shall say that $\CC$ is flat over $Z$ if every object of $X$ admits
a surjection $X'\to X$ for $X$ being flat over $Z$.

Consider the example of $\CC=\bA\mod$. Suppose there exists a family
of open left ideals $\bI\subset \bA$ such that $\bA\simeq
\underset{\longleftarrow}{\lim}\, \bA/\bI$, such that each $\bA/\bI$
is flat as a $Z$-module. Then $\CC$ is flat over $Z$.

\begin{lem}   \label{two derived tensor products}
Let $\CC$ be flat over $Z$, then for a fixed $M^\bullet\in
\bC^-(Z\mod)$ the left derived functor of $$X^\bullet\mapsto
M^\bullet\underset{Z}\otimes X^\bullet:\bC^-(\CC)\to \bC^-(\CC)$$ is
isomorphic to $M^\bullet\overset{L}{\underset{Z}\otimes} X^\bullet$.
\end{lem}

\begin{proof}

By assumption, every object in $\bC^-(\CC)$ admits a quasi-isomorphism
from one consisting of objects that are $Z$-flat. Hence, it suffices
to show that if $X^\bullet\in \bC^-(\CC)$ consists of $Z$-flat
objects, and $M^\bullet\in \bC^-(Z\mod)$ is acyclic, then
$M^\bullet\underset{Z}\otimes X^\bullet$ is acyclic as well.  However,
this is evident from the definitions.

\end{proof}

If $\phi:Z\to Z'$ is a homomorphism, we will denote by $\CC_{Z'}$ the
base-changed category, i.e., one whose objects are $X\in \CC$, endowed
with an action of $Z'$, such that the two actions of $Z$ on $X$
coincide. Morphisms in this category are $\CC$-morphisms that commute
with the action of $Z'$.

By construction, $Z'$ maps to the center of $\CC_{Z'}$. The composed
functor $\CC_{Z'}\to \CC\overset{\sF_Z}\to Z\mod$ factors naturally
through $Z'\mod$.

The forgetful functor $\CC_{Z'}\to \CC$ admits a left adjoint $\phi^*$
given by $X\mapsto Z'\underset{Z}\otimes X$. Note that this functor
sends $Z$-flat objects in $\CC$ to $Z'$-flat objects in $\CC_{Z'}$. In
particular, if $\CC$ is flat over $Z$, then so is $\CC_{Z'}$ over
$Z'$.

As in \lemref{two derived tensor products} and \lemref{der adj}, we
obtain the following

\begin{lem}    \label{* res}
Assume that $\CC$ is $Z$-flat. Then the right derived functor of
$\phi^*$
$$L\phi^*:D^-(\CC)\to D^-(\CC_{Z'})$$ is well-defined and is the left
adjoint to the forgetful functor $D(\CC_{Z'})\to D(\CC)$.  Moreover,
$$\sF_{Z'}\circ L\phi^*(X^\bullet)\simeq
\sF(X^\bullet)\overset{L}{\underset{Z}\otimes}Z'.$$
\end{lem}

In particular, we obtain that if $\CC$ is flat over $Z$ and $X\in \CC$
is $Z$-flat, then for $Y^\bullet \in \bC(\CC_{Z'})$,
$$\RHom_{D(\CC_{Z'})}(Z'\underset{Z}\otimes X, Y^\bullet)\simeq 
\RHom_{D(\CC)}(X,Y^\bullet).$$

\medskip

Let now $N$ be a $Z$-module. For $Y\in \CC$ we define the object
$\uHom_Z(N,Y)$ by
$$\Hom_{\CC}(X,\uHom_Z(N,Y)):=\Hom_\CC(N\underset{Z}\otimes X,Y).$$ If
$N$ is finitely presented, we have
$$\sF_Z(\uHom_Z(N,Y))\simeq \Hom_Z(N,\sF_Z(Y)).$$

For $Z'$ as above, which is finitely presented as a $Z$-module, we
define the functor $\phi^!:\CC\to \CC_{Z'}$ to be the right adjoint of
the forgetful functor $\CC_{Z'}\to \CC$.  It is given by $X\mapsto
\uHom_Z(Z',X)$. By definition, it maps injective objects in $\CC$ to
injectives in $\CC_{Z'}$.

We will denote by $R\phi^!:D^+(\CC)\to D^+(\CC_{Z'})$ the
corresponding right derived functor. It is easily seen to be the right
adjoint of the forgetful functor $\bC(\CC_{Z'})\to \bC(\CC)$.

\begin{prop}   \label{! restrictions}
Assume that $\CC$ is flat over $Z$, and that $Z'$ is perfect as an
object of $\bC(Z\mod)$. Then
$$R\phi^!\circ \sF_Z\simeq \sF_{Z'}\circ R\phi^!:D^+(\CC)\to
D^+(Z'\mod)$$
\end{prop}

\begin{proof}

To prove the proposition it suffices to check that if $Y^\bullet\in
\bC^+(\CC)$ is a complex, consisting of injective objects of $\CC$,
and $M^\bullet\in \bC^b(Z\mod)$ is a complex of finitely presented
modules, quasi-isomorphic to a perfect one, then
$\sF_Z(\uHom(M^\bullet,Y^\bullet))$ is quasi-isomorphic to
$\RHom_{D(Z\mod)}(M^\bullet,\sF_Z(Y^\bullet))$.

If $M^\bullet$ is a bounded complex, consisting of finitely generated
projective modules, then the assertion is evident. Hence, it remains
to show that if $M^\bullet$ is an acyclic complex of finitely
presented $Z$-modules, and $Y^\bullet$ is as above, then
$\uHom(M^\bullet,Y^\bullet)$ is acyclic. By assumption on $\CC$, it
would suffice to check that for $X\in \CC$ which is $Z$-flat,
$$\Hom_{\CC}(X,\uHom(M^\bullet,Y^\bullet))\simeq
\Hom_{\CC}(M^\bullet\underset{Z}\otimes X,Y^\bullet)$$ is acyclic. By
the flatness assumption on $Y$, the complex
$M^\bullet\underset{Z}\otimes X$ is also acyclic, and hence our
assertion follows from the injectivity assumption on $Y^\bullet$.

\end{proof}

\begin{cor}
If, under the assumptions of the proposition, $X'\in D(\CC_{Z}')$ is
quasi-perfect, then it is quasi-perfect also as an object of $D(\CC)$.
\end{cor}

\begin{proof}
This follows from the fact that the functor $R\phi^!:D(Z\mod)\to
D(Z'\mod)$ commutes with direct sums, and hence, so does the functor
$R\phi^!:D(\CC)\to D(\CC_{Z'})$.
\end{proof}

\ssec{A generalization}

Let $\CC$ be as in the previous subsection, and assume in addition
that it satisfies assumption (*) of \secref{class of cat}. Let $\bZ$
be a topological commutative algebra, which acts functorially on every
object of $\CC$. In this case we will say that $\bZ$ maps to the
center $\CC$. The functor $\sF$ naturally factors through a functor
$\sF_\bZ:\CC\to \bZ\mod$.

For every discrete quotient $Z$ of $\bZ$, let $\CC_Z$ be the
subcategory of $\CC$, consisting of objects, on which $\bZ$ acts via
$Z$. If $\bZ\twoheadrightarrow Z \twoheadrightarrow Z'$, then
$\CC_{Z'}$ is obtained from $\CC_Z$ by the procedure described in the
previous subsection.

We shall say that $\CC$ is flat over $\bZ$, if each $\CC_Z$ as above
is flat over $Z$. Equivalently, we can require that this happens for a
cofinal family of discrete quotients $Z$ of $\bZ$. Henceforth, we will
assume that $\CC$ is flat over $\bZ$.

In what follows we will make the following additional assumption on
$\bZ$. Namely, that we can present $\bZ$ as
$\underset{\underset{i}{\longleftarrow}}{\lim}\, Z_i$, such that for
$j\geq i$ the ideal of $\phi_{j,i}:Z_j\to Z_i$ is perfect as an object
of $D(Z_j\mod)$.

Recall that a discrete quotient $Z$ of $\bZ$ reasonable if for some
(equivalently, any) index $i$ such that $\bZ\to Z$ factors through
$Z_i$, the algebra $Z$ is finitely presented as a $Z_i$-module. We
shall call $Z$ admissible if the finite-presentation condition is
replaced by the perfectness one.

\medskip

Let us call an object $M\in \bZ\mod$ finitely presented if $M$ belongs
to some $Z\mod$ and is finitely presented as an object of this
category, if $Z$ is reasonable. By the assumption on $\bZ$, this
condition does not depend on a particular choice if $Z$.

For a finitely presented $M\in \bZ\mod$ and $X\in \CC$ we define
$\uHom_\bZ(M,X)\in \CC$ as
$$\underset{X_i}{\underset{\longrightarrow}{\lim}}\,
\uHom_{Z_i}(M,X_i),$$ where $X_i$ runs over the set of subobjects of
$X$ that belong to $\CC_{Z_i}$ for some discrete reasonable quotient
$Z_i$ of $\bZ$. We have
$$\sF_\bZ(\uHom_\bZ(M,X))\simeq \Hom_{\bZ}(M,\sF_\bZ(X)).$$

Consider $M=Z$ for some reasonable quotient $\phi:\bZ\to Z$. Then
$X\mapsto \uHom_\bZ(Z,X)$ defines a functor $\CC\to \CC_Z$, which we
will denote by $\phi^!$.

\begin{lem}
The functor $\phi^!$ is the right adjoint to the forgetful functor
$\CC_Z\to \CC$.
\end{lem}

\begin{proof}
By assumption (*), it suffices to check that for every finitely
generated object $Y$ of $\CC_Z$,
$$\Hom_{\CC_Z}(Y,\uHom_\bZ(Z,X))\simeq \Hom_\CC(Y,X).$$

By the finite generation assumption, we reduce the assertion to the
case when $X\in \CC_{Z_i}$ for some $\bZ\twoheadrightarrow
Z_i\twoheadrightarrow Z$, considered in the previous subsection.
\end{proof}

Evidently, the functor $\phi^!$ maps injective objects in $\CC$ to
injectives in $\CC_Z$. Let $R\phi^!$ denote the right derived functor
of $\phi^!$. By the above, it is the right adjoint to the forgetful
functor $D(\CC_Z)\to D(\CC)$.

\begin{prop}  \label{! ind restrictions}  
Assume that $Z$ is admissible. Then we have an isomorphism of functors:
$$\sF_Z\circ R\phi^!\simeq R\phi^!\circ \sF_\bZ:D^+(\CC)\to
D^+(Z\mod).$$
\end{prop}

\begin{proof}

As in the proof of \propref{! restrictions}, it suffices to show that
if $X^\bullet\in \bC^+(\CC)$ is a complex, consisting of injective
objects of $\CC$, and $M^\bullet$ is a perfect object of $D(Z\mod)$,
then $\CHom_{\bC(\bZ\mod)}(M^\bullet,\sF_\bZ(X^\bullet))$ computes
$\RHom_{D(\bZ\mod)}(M^\bullet,\sF_\bZ(X^\bullet))$.

By devissage, we can assume that $X^\bullet$ consists of a single
injective object $X\in \CC$. For every $Z_i$ such that
$\bZ\overset{\phi_i}\twoheadrightarrow Z_i\twoheadrightarrow Z$, note
that $\phi_i^!(X)$ is an injective object of $\CC_{Z_i}$, and $X\simeq
\underset{i}{\underset{\longrightarrow}{\lim}}\, X_i$.

Using \propref{! restrictions}, the assertion of the present
proposition follows from the next lemma:

\begin{lem}  \label{direct limit for modules}
For $N^\bullet\in \bC^+(\bZ\mod)$ and
$N_i^\bullet:=\phi_i^!(N^\bullet)$, the map
$$\underset{i}{\underset{\longrightarrow}{\lim}}\, \Hom_{D(Z_i\mod)}
(M^\bullet,N_i^\bullet)\to \Hom_{D(\bZ\mod)}(M^\bullet,N^\bullet)$$
is a quasi-isomorphism, provided that $M^\bullet\in D(Z\mod)$ is perfect.
\end{lem}

\end{proof}

\begin{proof} (of the Lemma) The proof follows from the next observation: 

Let $P^\bullet\to M^\bullet$ be a quasi-isomorphism, where
$P^\bullet\in \bC^-(\bZ\mod)$.  Then we can find a quasi-isomorphism
$Q^\bullet\to P^\bullet$ such that $Q^\bullet\in \bC^-(\bZ\mod)$ and
for any integer $i$, the module $Q^i$ is supported on some discrete
quotient of $\bZ$.

\end{proof}

\begin{cor}   \label{restriction and direct limits}
The functor $R\phi^!:D^+(\CC)\to D^+(\CC_Z)$ commutes with uniformly
bounded from below direct sums.
\end{cor}

\begin{proof}
This follows from the corresponding fact for the functor
$R\phi^!:D^+(\bZ\mod)\to D^+(Z\mod)$.
\end{proof}

\medskip

\begin{prop}
Let $X_1^\bullet$ be a quasi-perfect object of $\bC(\CC_Z)$ and
$X_2^\bullet$ be an object of $\bC^+(\CC_Z)$ for some discrete
quotient $Z$. Then
$$\Hom_{D(\CC)}(X_1^\bullet,X_2^\bullet)\simeq
\underset{\underset{Z_i}{\longrightarrow}}{\lim}\,
\Hom_{D(\CC_{Z_i})}(X_1^\bullet,X_2^\bullet),$$ where the direct limit
is taken over the indices $i$ such that $\bZ\to Z$ factors through
$Z_i$.
\end{prop}

\begin{proof}

We can find a system of quasi-isomorphisms $X_2^\bullet\to
Y_i^\bullet$, where each $Y_i^\bullet\in \bC(\CC_{Z_i})$ consist of
injective objects of $\CC_{Z_i}$, and such that these complexes form a
direct system with respect to the index $i$, and such that all
$Y_i^\bullet$ are uniformly bounded from below.

By \propref{! ind restrictions} and \corref{restriction and direct limits},
$R\phi^!(X_2^\bullet)$ is given by the complex 
$$\underset{i}{\underset{\longrightarrow}{\lim}}\,
\phi_i^!(Y_i^\bullet).$$ Then, by the quasi-perfectness assumption,
$$\Hom_{D(\CC)}(X_1^\bullet,X_2^\bullet)\simeq \Hom_{D(\CC_Z)}
(X_1^\bullet,R\phi^!(X_2^\bullet))\simeq 
\underset{i}{\underset{\longrightarrow}{\lim}}\, 
\Hom_{D(\CC_Z)}(X_1^\bullet, \phi_i^!(Y_i^\bullet)).$$

By \propref{! restrictions}, the latter is isomorphic to
$$\underset{i}{\underset{\longrightarrow}{\lim}}\, 
\Hom_{D(\CC_Z)}(X_1^\bullet, R\phi_i^!(Y_i^\bullet))\simeq
\underset{i}{\underset{\longrightarrow}{\lim}}\, 
\Hom_{D(\CC_{Z_i})}(X_1^\bullet, Y_i^\bullet),$$
which is what we had to show.

\end{proof}

Finally, we will prove the following assertion:

\begin{prop}    \label{quasi-perf}
Let $X^\bullet$ be an object of $\bC^-(\CC_Z)$, where $Z$ is an
admissible quotient of $\bZ$. Then $X^\bullet$ is quasi-perfect as an
object of $D(\CC_Z)$ if and only if it is quasi-perfect as an object
of $D(\CC)$.
\end{prop}

\begin{proof}

Since the functor $R\phi^!:D^+(\CC)\to D^+(\CC_Z)$ commutes with
direct sums, the implication "quasi-perfectness in $D(\CC_Z)$"
$\Rightarrow$ "quasi-perfectness in $D(\CC)$" is clear.

To prove the implication in the opposite direction, we proceed by
induction.  We suppose that the functor
$$Y\mapsto \Hom_{D(\CC_Z)}(X^\bullet,Y[i']):\CC_Z\to \Vect$$
commutes with direct sums for $i'<i$. This assumption is satisfied
for some $i$, since $X^\bullet$ is bounded from above.

Let us show that in this case the functor $Y\mapsto
\Hom_{D(\CC_Z)}(X^\bullet,Y[i])$ also commutes with direct sums. For
$\underset{\alpha}\oplus\, Y_\alpha\in \CC_Z$ consider the exact
triangle in $D^+(\CC_Z)$:
$$\underset{\alpha}\oplus\, Y_\alpha\to
R\phi^!(\underset{\alpha}\oplus\, Y_\alpha)\to
\tau^{>0}\left(R\phi^!(\underset{\alpha}\oplus\, Y_\alpha)\right),$$
where $\tau$ is the cohomological truncation.

Consider the corresponding commutative diagram:
$$ \CD \Hom_{D(\CC)}(X^\bullet,\underset{\alpha}\oplus\,
Y_\alpha[i-1]) @<<< \underset{\alpha}\oplus\,
\Hom_{D(\CC)}(X^\bullet,Y_\alpha[i-1]) \\ @VVV @VVV \\ \Hom_{D(\CC_Z)}
\left(X^\bullet,\tau^{>0}\left(R\phi^!(\underset{\alpha}\oplus\,
Y_\alpha[i-1])\right)\right) @<<< \underset{\alpha}\oplus\,
\Hom_{D(\CC_Z)}
\left(X^\bullet,\tau^{>0}\left(R\phi^!(Y_\alpha[i-1])\right)\right) \\
@VVV @VVV \\ \Hom_{D(\CC_Z)}(X^\bullet,\underset{\alpha}\oplus\,
Y_\alpha[i]) @<<< \underset{\alpha}\oplus\,
\Hom_{D(\CC_Z)}(X^\bullet,Y_\alpha[i]) \\ @VVV @VVV \\
\Hom_{D(\CC)}(X^\bullet,\underset{\alpha}\oplus\, Y_\alpha[i]) @<<<
\underset{\alpha}\oplus\, \Hom_{D(\CC)}(X^\bullet,Y_\alpha[i]) \\ @VVV
@VVV \\ \Hom_{D(\CC_Z)}
\left(X^\bullet,\tau^{>0}\left(R\phi^!(\underset{\alpha}\oplus\,
Y_\alpha[i])\right)\right) @<<< \underset{\alpha}\oplus\,
\Hom_{D(\CC_Z)}
\left(X^\bullet,\tau^{>0}\left(R\phi^!(Y_\alpha[i])\right)\right)
\endCD.
$$

The horizontal arrows in rows 1 and 4 are isomorphisms since
$X^\bullet$ is quasi-perfect in $D(\CC)$. The arrows in rows 2 and 5
are isomorphisms by the induction hypothesis.  Hence, the map in row 3
is an isomorphism, which is what we had to show.

\end{proof}

\ssec{The equivariant situation}

Assume now that the category $\CC$ as in the previous subsection is
equipped with an infinitesimally trivial action of a group-scheme
$H$. Assume that this action commutes with that of $\bZ$. The latter
means that for every $X\in \CC$, the $\bZ$-action on $\act^*(X)$ by
transport of structure coincides with the action obtained by regarding
it merely as an object of $\CC$. Then for every discrete quotient $Z$
of $\bZ$, the category $\CC_Z$ carries an infinitesimally trivial
action of $H$.

We have a functor $\phi^!:\bC^+(\CC)^H\to \bC^+(\CC_Z)^H$, and let
$R\phi^!:D^+(\CC)^H\to D^+(\CC_Z)^H$ be its right derived functor.
(Below we will show that it is well-defined.) We are going to prove
the following:

\begin{prop}   \label{! ind res equiv}
$R\phi^!:D^+(\CC)^H\to D^+(\CC_Z)^H$ is the right adjoint
to the forgetful functor $D(\CC_Z)^H\to D(\CC)^H$. Moreover,
the diagram of functors
$$
\CD
D^+(\CC)^H @>{R\phi^!}>> D^+(\CC_Z)^H \\
@V{\sF_\bZ}VV  @V{\sF_Z}VV  \\
D^+(\bZ\mod) @>{R\phi^!}>> D^+(Z\mod)
\endCD
$$
is commutative.
\end{prop}

\begin{proof}

For any quasi-isomorphism $X^\bullet\to X_1^\bullet$ in $\bC^+(\CC)^H$
we can find a quasi-isomorphism from $X_1^\bullet$ to a complex,
associated with a bi-complex $X_2^{\bullet,\bullet}$, whose rows are
uniformly bounded from below and have the form $\on{Av}_H(Y^\bullet)$,
where $Y^\bullet\in \bC^+(\CC)$ consists of injective objects.

By \propref{! ind restrictions} and \corref{restriction and direct
limits}, if we assign to $X^\bullet$ the complex in $\bC(\CC_Z)^H$
associated with the bi-complex $\phi^!(X_2^{\bullet,\bullet})$, this
is the desired right derived functor of $\phi^!$. It is clear from the
construction that the diagram of functors
$$
\CD
D^+(\CC)^H @>{R\phi^!}>> D^+(\CC_Z)^H \\
@VVV  @VVV  \\
D^+(\CC) @>{R\phi^!}>> D^+(\CC_Z),
\endCD
$$
where the vertical arrows are the forgetful functors, is commutative.

Hence, it remains to show that $R\phi^!$ satisfies the desired
adjointness property. By devissage, we are reduced to showing that for
$Y^\bullet$ as above and $Y_1^\bullet\in \bC(\CC_Z)^H$,
$$\Hom_{D(\CC)^H}(Y^\bullet_1,\on{Av}_H(Y^\bullet))\simeq
\Hom_{D(\CC_Z)^H}\left(Y_1^\bullet,
\phi^!(\on{Av}_H(Y^\bullet))\right).$$ However, the LHS is isomorphic
to $\Hom_{D(\CC)}(Y^\bullet_1,Y^\bullet)$, and the RHS is isomorphic
to
$$\Hom_{D(\CC_Z)^H}\left(Y_1^\bullet,
\on{Av}_H(\phi^!(Y^\bullet))\right)\simeq \Hom_{D(\CC_Z)}(Y_1^\bullet,
\phi^!(Y^\bullet)),$$ and, as we have seen above,
$\phi^!(Y^\bullet)\to R\phi^!(Y^\bullet)$ is an isomorphism in
$D^+(\CC_Z)$.

\end{proof}

\end{document}